\documentclass[11pt]{amsart}

\usepackage{amsmath,amssymb,amsthm,tikz}
\usetikzlibrary{patterns}
\usetikzlibrary{arrows.meta}
\usetikzlibrary{decorations.pathreplacing}

\usepackage[numbers,comma,square,sort&compress]{natbib}
\usepackage{mathabx}   %% made \widecheck work
\usepackage{datetime}
\usdate
\usepackage[normalem]{ulem}  % for striking out text
\RequirePackage{color}
\RequirePackage[colorlinks,urlcolor=my-blue,linkcolor=my-red,citecolor=my-green]{hyperref}
\usepackage{soul}

\numberwithin{figure}{section}

 \usepackage{enumitem}    %% This caused an error message 2018-04-30 

\usepackage{scalerel}    %%%%  used first for scaling \bullet  
\usepackage{stmaryrd}   %%%% for integer interval notation 
\usepackage{mathabx}   %% made \widecheck work 
\usetikzlibrary{shapes.misc}     %%% for non-coalescence 
\usetikzlibrary{calc}

\usepackage{rotating,centernot,cancel}    %% to manipulate symbols 

\definecolor{my-blue}{rgb}{0.0,0.0,0.6}
\definecolor{my-red}{rgb}{0.5,0.0,0.0}
\definecolor{my-green}{rgb}{0.0,0.5,0.0}
\definecolor{nicos-red}{rgb}{0.75,0.0,0.0}

\newtheorem{theorem}{\sc Theorem}[section]
\newtheorem{lemma}[theorem]{\sc Lemma}
\newtheorem{proposition}[theorem]{\sc Proposition}

\numberwithin{equation}{section}
\theoremstyle{remark}
\newtheorem{remark}[theorem]{Remark}
\AtBeginEnvironment{remark}{%
  \pushQED{\qed}%
}
\AtEndEnvironment{remark}{\popQED\endremark}
\newtheorem{example}[theorem]{Example}

\newcommand{\be}{\begin{equation}}
\newcommand{\ee}{\end{equation}}
\newcommand{\beq}{\begin{equation}}
\newcommand{\eeq}{\end{equation}}

\providecommand{\abs}[1]{\vert#1\vert}

\makeatletter
\newcommand{\oset}[3][0ex]{%
	\mathrel{\mathop{#3}\limits^{
			\vbox to#1{\kern-1\ex@
				\hbox{$\scriptstyle#2$}\vss}}}}
\makeatother

\newcommand{\bbP}{\mathbb{P}}

\def\cC{\mathcal{C}}  
\def\cF{\mathcal{F}}
\def\cG{\mathcal{G}}
\def\cH{\mathcal{H}}

\def\cN{\mathcal{N}}
\def\cP{\mathcal{P}}
\def\cQ{\mathcal{Q}}
\def\cR{\mathcal{R}}

\def\cS{\mathcal{S}}
\def\cI{\mathcal{I}}
\def\cK{\mathcal{K}}
\def\cB{\mathcal{B}}

\def\cA{\mathcal{A}}
 \def\cV{\mathcal{V}}

\font \mymathbb = bbold10 at 11pt
\newcommand{\one}{\mbox{\mymathbb{1}}}

\def\bE{\mathbb{E}} 
\def\bN{\mathbb{N}}
\def\bP{\mathbb{P}}

\def\bR{\mathbb{R}}
\def\bZ{\mathbb{Z}}

 \def\Z{\bZ}  \def\R{\bR}\def\N{\bN}

\def\w{\omega}

\def\e{\varepsilon}

\font \mymathbb = bbold10 at 11pt

\def\m1{\mathbf{1}}

 \def\cJ{\mathcal{J}}

\def\E{\bE}
\def\P{\bP} %% environment measure 
\def\funct lp{L} %a function needed
\def\funct lpbar{\bar L} %a function needed

\def\note#1{{\color{red}{(#1)}}}

\def\Uset{\mathcal U}

\def\cA{\mathcal A}

\def\b{{b}}

\def\cE{{\mathcal E}}

\DeclareMathOperator{\ri}{ri}    %rel interior  
\def\cH{\mathcal{H}}

\definecolor{darkgreen}{rgb}{0.0,0.5,0.0}
\definecolor{darkblue}{rgb}{0.0,0.0,0.3}
\definecolor{nicosred}{rgb}{0.65,0.1,0.1}
\definecolor{light-gray}{gray}{0.7}
\allowdisplaybreaks[1]

\usepackage{filemod}

\def\cif1{v}   %projection of CIF 
\newcommand\abullet{{\raisebox{2pt}{\scaleobj{0.5}{\bullet}}}}  %%%% smaller bullet to use inside functions 
\newcommand\bbullet{{\raisebox{0.5pt}{\scaleobj{0.6}{\bullet}}}} %% looks good for paths x_\bbullet 
\def\tsp{\hspace{0.5pt}}  %this is sometimes a good spacing in formulas when \, is too large 
\def\tspa{\hspace{0.7pt}}
\def\tspb{\hspace{0.9pt}}

\newcommand{\sig}{{\scaleobj{0.8}{\boxempty}}} 
\newcommand{\sigg}{{\scaleobj{0.9}{\boxempty}}} 

\def\timonote#1{{\color{darkgreen}{(T: #1)}}}
\soulregister\timonote{1}
\def\selfnote#1{{\color{brown}{(Note: #1)}}}
\def\firasnote#1{{\color{orange}{(F: #1)}}}
\soulregister\firasnote{1}
\def\chrisnote#1{{\color{purple}{(C: #1)}}}
\soulregister\chrisnote{1}

\def\level{\bbL}     %%%% levels on lattice 

\def\pawgt{W}  %path weight 
\def\path{x_\bbullet}   %generic path
\def\cif{\varphi}

\newcommand{\Lpp}{L}

\def\ead{e_{\mathrm{ad}}}
\def\ed{e_{\mathrm{d}}}
\def\ghatz{g(\Dir)}%{b_0}
\def\ghatprm{a_-}
\def\ghatprp{a_+}
\def\ghatpr#1{a_{#1}}   %% a_\sig, a_\pm; change together with \ghatprm, \ghatprp
\def\hcross{{\mathrm{HC}}}
\def\conn{{\mathrm{Conn}}}
\def\pbad{\bar p}               %% bound on the probability of a bad staircase 
\def\Mtop{M_+}                  %% top of the support of \mu
\def\gtop{g^{\rm top}}          %% shape function of the top-weight LPP
\def\Ltop{\Lpp^{\rm top}}       %% top-weight last-passage time
\def\Uperc{\Uset_{\rm perc}}    %% flat segment of g
\def\opc{\vec p_c}              %% critical density of oriented site percolation
\def\Bad{{\mathrm{Bad}}}          %% a staircase is bad
\def\Good{{\mathrm{Good}}}        %% good event of one grid interval and bin
\def\Stairs{{\mathfrak S}}        %% family of staircases that may be needed
\def\Lip{{\mathrm{Lip}}}

\allowdisplaybreaks[1]

\newif\ifshowgists
\showgistsfalse  % toggle between true/false to show/hide
\ifshowgists
  \long\def\gist#1{{\color{brown}{(Gist: #1)}}}
\else
  \long\def\gist#1{\ignorespaces}
\fi

\newif\ifshownotes
\shownotesfalse   % toggle between true/false to show/hide
\ifshownotes
  \long\def\note#1{{\color{brown}{(Note: #1)}}}
\else
  \long\def\note#1{\ignorespaces}
\fi

\newif\ifshowselfnotes
\showselfnotestrue  % toggle between true/false to show/hide
\ifshowselfnotes
  \long\def\selfnote#1{{\color{brown}{(#1)}}}
\else
  \long\def\selfnote#1{\ignorespaces}
\fi

\def\Elev{E}               %% crossing events E(s,n)
\def\level{\ell}           %% level map \ell(z)=z_1+z_2; level values are written I, I_1, I_2, ...
\def\rcoord{r}             %% transverse coordinate
\def\Dir{\xi}
\def\cnst{C}
\def\halfgap{\varsigma}    %% the half-gap (\hat g'(0-)-\hat g'(0+))/2 
\def\vbad{{\cV}}      %% the set of very bad vertices (Appendix app:general), and the vertex set of a path event (Section sub:bk), formerly V
\newcommand\arrowabove[3][0pt]{\vphantom{\raisebox{1.55ex}{$\scriptscriptstyle#3$}}\mathord{\ooalign{$#2$\cr\hidewidth\kern#1\smash{\raisebox{1.55ex}{$\scriptscriptstyle#3$}}\hidewidth\cr}}}
\def\Vfw{\arrowabove{\cV}{\rightarrow}}   %% trapped vertices, formerly \bar V
\def\Vbw{\arrowabove{\cV}{\leftarrow}}    %% backward trapped vertices, formerly \bar V^{\leftarrow}
\def\Afw{\arrowabove[0.2em]{A}{\rightarrow}}     %% forward reachable set, formerly A(z)
\def\Abw{\arrowabove[0.2em]{A}{\leftarrow}}      %% backward reachable set, formerly A^{\leftarrow}(z)
\def\exc{\Theta}           %% the excess (\w_z-M_0)^+ of a vertex (Remark in Section sec:witnesses and Appendix sub:onestep), formerly \varsigma
\def\Cm{C_-}               %% E\w_0^-, formerly C_1, which coincided with \cnst_1 of Lemma lem:grid
\def\clen{\Lambda}         %% characteristic length (L is taken by \Lpp)
\def\Lip{\mathrm{Lip}}        %% Lipschitz constants of g
\def\jP{\overline{\bbP}}       %% joint law of (\w,\tilde\w,\chi) (\bP is taken by the macros) 

\def\mult{\mathfrak m}        %% number of disjoint witnesses in a stage

\author[C.~Janjigian]{Christopher Janjigian}
\address{Christopher Janjigian\\ Purdue University\\  Department of Mathematics \\ 150 N University St\\ West Lafayette, IN 47901\\ USA.}
\email{cjanjigi@purdue.edu}
\urladdr{http://www.math.purdue.edu/~cjanjigi}
\thanks{C.~Janjigian was partially supported by National Science Foundation grant DMS-2125961 and Simons Foundation grant MPS-TSM-00012155. This work was also supported in part by a grant of access to OpenAI models through the ChatGPT for Academic Researchers program and by Anthropic through the Claude for Scientists program.}

\author[F.~Rassoul-Agha]{Firas Rassoul-Agha}
\address{Firas Rassoul-Agha\\ University of Utah\\  Department of Mathematics\\ 155S 1400E\\   Salt Lake City, UT 84112\\ USA.}
\email{firas@math.utah.edu}
\urladdr{http://www.math.utah.edu/~firas}
\thanks{F.\ Rassoul-Agha was partially supported by National Science Foundation grants DMS-2054630 and DMS-2450951 and Simons Foundation grant MPS-TSM-00013661. This work was also supported in part by a grant of access to Anthropic and OpenAI models through the Claude for Scientists and the ChatGPT EDU programs.}

\author[T.~Sepp\"al\"ainen]{Timo Sepp\"al\"ainen}
\address{Timo Sepp\"al\"ainen\\ University of Wisconsin--Madison\\  Department of Mathematics\\ 480 Lincoln Drive \\  Madison, WI 53706\\ USA.}
\email{seppalai@math.wisc.edu}
\urladdr{http://www.math.wisc.edu/~seppalai}
\thanks{T.~Sepp\"al\"ainen was partially supported by National Science Foundation grants  DMS-1854619, DMS-2152362, and DMS-2448375, by Simons Foundation grant 1019133, by the Wisconsin Alumni Research Foundation, and by a grant from the ChatGPT for Academic Researchers program.}


\begin{document}

\title[Differentiability of shape in planar CGM]{Differentiability of the shape function\\ in the planar corner growth model}

%\date{Last modified on \today\ at\ \filemodprinttime{\jobname}}
\makeatletter
\renewcommand{\datename}{}
\makeatother

%\date{Last updated: \today} 
\date{Last updated: November 24, 2026} 

\begin{abstract}
We prove that the shape function of planar directed last-passage percolation with general i.i.d.\ real-valued weights is differentiable on the open first quadrant, under a moment assumption of order $p''>4$ on the positive part of the weight and a nearly optimal   order $p'>1$ on the negative part. Due to known results, only directions outside a percolation cone need be considered. The proof is a multiscale path-counting argument that uses the shape theorem, oriented percolation on a coarse-grained lattice, a simple change of measure in a finite region, and the van den Berg--Kesten inequality. Our Lean-formalized proof is an AI-assisted adaptation and improvement of the OpenAI proof \cite{Ope-26-convdiff} of differentiability in planar first-passage percolation.
\end{abstract}

\maketitle

\tableofcontents

\section{Introduction}

\subsection{First-passage percolation and its limit shape}
First-passage percolation was introduced by Hammersley and Welsh \cite{Ham-Wel-65} in 1965 as a model for the flow of a fluid through a random porous medium. It has since become one of the cornerstone models of probability theory, with connections to random metric geometry, growth processes, and statistical physics; see \cite{Auf-Dam-Han-17} for a survey.

In its standard form, independent and identically distributed nonnegative random weights, the passage times, are placed on the edges (or the vertices) of the lattice $\Z^d$. The passage time of a path is the sum of the weights on it, and the first-passage time between two points is the minimal passage time over all paths joining them. This defines a random (pseudo) metric on $\Z^d$. By subadditivity, the first-passage time from the origin to $nx$ grows linearly in $n$, and its rate defines a deterministic norm-like function, the \emph{time constant}. Its unit ball is the \emph{limit shape}: the set of points reached by time $t$, scaled by $1/t$, converges to it \cite{Cox-Dur-81,Kes-86-stflour}.

\subsection{Regularity of the limit shape}
Soft properties of the limit shape, such as boundedness, convexity, and continuity of the time constant, follow quickly from subadditivity and the symmetries of the lattice. Deeper properties have been harder.  Their importance goes beyond the shape itself.

Differentiability of the boundary of the limit shape is closely tied to the existence and uniqueness of Busemann functions, and through them to the existence, uniqueness, and directedness of semi-infinite geodesics and to the structure of the geodesic tree \cite{New-95,Hof-08,Dam-Han-14,Geo-Ras-Sep-17-ptrf-1,Geo-Ras-Sep-17-ptrf-2}. Curvature of the limit shape is tied to fluctuations. First-passage percolation is believed to belong to the Kardar--Parisi--Zhang (KPZ) universality class \cite{Kar-Par-Zha-86,Cor-12}, in which, in two dimensions, the passage time fluctuates on the scale $n^{1/3}$ and geodesics wander on the scale $n^{2/3}$. Under suitable assumptions, including curvature of the limit shape, these exponents are known to satisfy the KPZ scaling relation $\chi=2\xi-1$ \cite{Cha-13,Auf-Dam-14}.

Yet, apart from a few special cases, even differentiability of the limit shape, the most basic of these properties, has been open since the introduction of the model. Prior known results concern flat pieces of the boundary created by an atom at the minimum of the passage-time distribution, the \emph{percolation cone} \cite{Dur-Lig-81,Mar-02}, and differentiability at the edges of this cone \cite{Auf-Dam-13}.

Recently, a manuscript produced by OpenAI \cite{Ope-26-convdiff} claimed a proof of the differentiability of the planar limit shape for gamma-distributed edge weights. In this paper we prove the analogous result for a closely related model, directed last-passage percolation, for general weights under moment  assumptions. The changes to the argument that remove the distributional assumptions are not specific to last-passage percolation.  We would expect the same extensions  in first-passage percolation.

\subsection{Directed last-passage percolation}\label{intro:LPP}
In this paper we work with a close cousin of first-passage percolation: planar directed last-passage percolation (LPP), also known as the corner growth model (CGM). Here the weights sit on the vertices of $\Z^2$, paths are restricted to up-right steps, and the last-passage time between two ordered points is the \emph{maximal} weight of a path joining them. The limit of the last-passage time in direction $\xi$, the \emph{shape function} $g(\xi)$, is concave and homogeneous \cite{Mar-04}.  The regularity questions above have exact analogues for it.

Which of the two models one studies is largely a matter of taste, since each has its own advantages. In first-passage percolation every vertex can be reached from every other, which removes boundary effects and gives the model a rich symmetry, at the price of a more complicated geometry of paths. In directed last-passage percolation the geometry is simpler, since paths are monotone and two of them can cross only by sharing a vertex. The model is also connected to many other systems: series of queues \cite{Gly-Whi-91,Sep-97-aap-1,OCo-99}, interacting particle systems such as the totally asymmetric simple exclusion process \cite{Ros-81,Sep-98-mprf-2}, and random matrices \cite{Joh-00}.

In contrast with FPP, the CGM 
has known exactly solvable cases. For exponential or geometric weights, fortuitous coincidences allow explicit computations, including an explicit formula for the shape function \cite{Ros-81,Joc-Pro-Sho-26-,Sep-98-mprf-1,Sep-98-mprf-2}. The conjectured regularity properties of the shape function  can be verified by calculus.

\subsection{Main result of this paper}\label{sec:i-main}
Theorem~\ref{thm:main} states that, assuming i.i.d.\ non-degenerate real-valued vertex weights, $\E[(\w_0^-)^{p'}]<\infty$ for some $p'>1$, and $\E[(\w_0^+)^{p''}]<\infty$ for some $p''>4$, the shape function $g$ is differentiable at every point of the open first quadrant $\R_{>0}^2$. Other than these moment assumptions, no further regularity requirements on the marginal  weight distribution are imposed.  The proof of Theorem~\ref{thm:main} adapts the first-passage percolation argument of \cite{Ope-26-convdiff} to directed last-passage percolation and sharpens it, as described in Section~\ref{sec:written} below.  The lower-tail assumption is close to optimal, while the upper-tail assumption comes from the method and is presumably not needed beyond $p''>2$.  Sections \ref{sub:lb_mom_chat} and \ref{sub:ub_mom-chat} expound more on this. 

When the maximum of the support is an atom, $g$ has a flat segment in the directions in which the top weights percolate. Differentiability on this flat segment, including at its endpoints and in the critical case, is treated separately in Section~\ref{sec:flat}, using \cite[Theorem~2.10]{Geo-Ras-Sep-15}.  The results on the flat segment are adaptations of the first-passage percolation work in \cite{Dur-Lig-81,Mar-02,Auf-Dam-13}. 

\subsection{Motivation of this work}
\label{sec:motive}
The last 30 years have seen enormous progress in understanding the properties of exactly solvable planar LPP models, while results for  models with general weight distributions lag far behind. 
  Three open  problems have remained  persistent roadblocks to progress on the general models: differentiability of the shape function, strict concavity of the shape function away from the percolation cone, and the nonexistence of bi-infinite geodesics.  Once OpenAI posted solutions to these questions for planar first-passage percolation, the way forward in directed  LPP became clear: transport the OpenAI FPP results to the CGM, verify them, and write them up in  a way that is accessible to the probability community. The present paper achieves this for the differentiability question, based on the FPP solution   in  \cite{Ope-26-convdiff}. Strict concavity of the shape function and the nonexistence of bi-infinite geodesics for the general CGM are treated in companion papers in preparation.

Once these three open problems have been resolved for the general CGM, we expect that the program of our articles  
\cite{Jan-Ras-Sep-23, Jan-Ras-Sep-26}
on the Busemann process, the geometry of semi-infinite geodesics, competition interfaces,  instability graphs, and stationary distributions can be carried much further.

\subsection{How this paper was written}\label{sec:written}
The starting point for this work was the OpenAI manuscript \cite{Ope-26-convdiff}. 
We had Claude Fable 5.1 (High)  port the argument of  \cite{Ope-26-convdiff} to the CGM.   In the process it observed that the assumption of gamma weights in \cite{Ope-26-convdiff} can be relaxed to continuous weights with more than six moments, by changing the way the weights are modified in the change of measure argument. We first restricted the proof to bounded continuous weights   to minimize technicalities that obscure the main flow of the proof. 

The next mathematical step involved generalization to arbitrary i.i.d.\ weights under reduced  moment assumptions, again with the help of Claude. 
 The continuity of the weight distribution originally appeared   in   the change of measure in Lemmas~\ref{lem:kernel} and~\ref{lem:tilt}. Removing this assumption required a modification of that  argument, together with a separate treatment of the flat segment created by an atom at the top of the support.  Boundedness of the weights was used in exactly three places: the attachments bound \eqref{attachments}, the staircase bound \eqref{eq:badstair}, and the one-step bound \eqref{eq:one-step}. We have kept the main text under the boundedness assumption  for clarity of exposition. Appendix~\ref{app:general} shows how the proofs  are extended to obtain Theorem~\ref{thm:main} under the moment assumptions stated above in Section \ref{sec:i-main}.  

The mathematics produced with Claude served as a draft. We then went through the arguments line by line, reorganized the proofs, filled in details, added explanations of steps and constructions, corrected mistakes, and chose the notation and terminology used here. We rewrote a number of sections ourselves, and Claude rewrote others along the same lines under our direction, after which we checked them again. The figures were drafted by Claude and edited by us. The paper in its present form is one that we have verified in full and whose exposition reflects our choices. As an independent check, the results have also been formalized in Lean, as described next.

\subsection{Lean formalization}
The results of this paper have been formalized in Lean~4. The formalization is available at
%\begin{center}
\url{https://github.com/cjanjigi/LPP}
%\end{center}
as part of a larger repository that formalizes basic results about directed last-passage percolation. It was produced with a slightly modified version of Scott Armstrong's Lean AutoFormalization skills, available at
\begin{center}
\url{https://github.com/scottnarmstrong/LeanAutoformalizationSkills}
\end{center}
with ChatGPT 6.1 Sol (xHigh) and Claude Opus 5.5 (High) as orchestrators and most sub-tasks delegated to cheaper models within those two families.

\subsection{Organization of the paper}
Section~\ref{sec:setting} describes the setting and states the main result. Section \ref{sec:strategy} summarizes the main points of the proof.   Section~\ref{sec:crossings} proves estimates for crossings of thin tubes, and Section~\ref{sec:scales} introduces the characteristic lengths through a change of measure. Section~\ref{sec:endpoints} extends the crossing estimates to arbitrary endpoints, and Section~\ref{sec:witnesses} concludes the proof through the van den Berg--Kesten inequality. Appendices~\ref{app:oriented} and~\ref{app:general} contain auxiliary  results on oriented percolation and the extensions of the main proofs  to unbounded weights. Appendix~\ref{app:dictionary} is a glossary of the symbols used in the proofs, with pointers to their definitions.  

\subsection{Notation and conventions} We collect here some standard and bespoke notation and terminology  for quick reference. New symbols and terms are also defined when they first appear in the text, and the glossary in Appendix~\ref{app:dictionary} collects them with pointers to their definitions. 

$\Z_{\ge0}=\{0,1,2,\dotsc\}$ and $\N=\Z_{>0}=\{1,2,3,\dotsc\}$.  The open and closed first quadrants of the plane are $\R^2_{>0}=(0,\infty)^2$ and $\R^2_{\ge0}=[0,\infty)^2$.  
The line segment between points  $u,v\in\R^2$ is denoted by $[u,v]=\{tu+(1-t)v: t\in[0,1]\}$. 
On $\R^2$ we use  the   norms $\abs z_1=\abs{z_1}+\abs{z_2}$ and $\abs z_\infty=\abs{z_1}\vee\abs{z_2}$. 
The coordinatewise order on $\R^2$ is $u\le v$ if $u_1\le v_1$ and $u_2\le v_2$. For real $a$ and $b$, $a\vee b=\max\{a,b\}$ and $a\wedge b=\min\{a,b\}$, and for real $t$, $t^+=t\vee0$ and $t^-=(-t)\vee0$. The indicator of an event or a set $A$ is $\one_A$, also written $\one\{\dotsm\}$. The symbol $\sigg$ stands for the sign $+$ or $-$.

Throughout the proof we work with a single fixed direction $\Dir$ in the open first quadrant, normalized by $\abs\Dir_1=1$. Constants written as $\cnst$ with parameters in a subscript or in parentheses, such as $\cnst_A$, $\cnst_{\!A,R}$, and $\cnst(D,\kappa,\theta_0)$, depend on the indicated parameters and on the fixed data $\mu$ and $\Dir$, but not on the scale parameters $s$, $j$, and $k$. Their values may change from one occurrence to the next. The two constants $\cnst_*$ and $\cnst_1$ are exceptions: they are fixed in Section~\ref{sec:endpoints} and do not vary.

For an event of the form $E=\{\text{there exists }x\text{ with property }A\}$, an $x$ with property $A$ is said to \emph{witness} $E$. The ending of a numbered remark is marked by a $\triangle$.

Some letters are used for more than one purpose, always locally. The letter $\alpha$ denotes an attachment path in the proof of Lemma~\ref{lem:intervals} and in Appendix~\ref{sub:attachments}, a parameter in Lemma~\ref{lem:kernel}, and a generic staircase in Section~\ref{sec:endpoints}. The letter $\eta$ is a parameter in Lemmas~\ref{lem:travel} and~\ref{lem:intervals}, a constant in the proof of Lemma~\ref{lem:tilt}, and the parameter of Section~\ref{sec:witnesses}, fixed in \eqref{eq:etac}. The letter $\delta$ is a distance in Lemmas~\ref{lem:travel} and~\ref{lem:intervals} and a probability in Section~\ref{sec:scales}. The letter $k$ is the scale index in Sections~\ref{sec:scales} and~\ref{sec:endpoints} and a generic level in Section~\ref{sec:witnesses}. The letter $h$ is the real argument of $\hat g$ in Section~\ref{sub:extreme}, a bin index in Section~\ref{sec:endpoints}, and the second coordinate of a block index $(p,h)$ in Section~\ref{sec:witnesses}. Lattice regions that are used inside one proof are denoted by $\cR$ with or without a subscript.

\section{Setting and main result}
\label{sec:setting}

\subsection{Corner growth model}\label{sub:cgm}
Let $\w=(\w_z)_{z\tsp\in\tsp\Z^2}$ be i.i.d.\ real-valued  random weights defined  on a probability space $(\Omega,\cF,\P)$, with common marginal law $\mu$. Expectation with respect to $\P$ is denoted by $\E$.  Throughout assume
\begin{enumerate} [label={\rm(A\arabic*)}] \itemsep=2pt 
\item\label{A1}\label{ass:mu} $\mu$ is not a point mass,
\item\label{A2} $\E\bigl[(\w_0^-)^{p'}\bigr]<\infty$ for some $p'>1$, and 
\item\label{A3} $\E\bigl[(\w_0^+)^{p''}\bigr]<\infty$ for some $p''>4$.
\end{enumerate}
No other regularity assumptions on the weight distribution are needed. The moment assumptions are discussed in   Sections \ref{sub:lb_mom_chat} and  \ref{sub:ub_mom-chat}.

Special vectors are denoted by $e_1=(1,0)$, $e_2=(0,1)$, $\ed=(1,1)$, and $\ead=(1,-1)$. 
The \emph{anti-diagonal level} of $z=(z_1,z_2)\in\R^2$ is \begin{equation}\label{eq:level}
\level(z)=z\cdot\ed=z_1+z_2.
\end{equation}
Note that $\abs z_1=\level(z)$ for $z$ in the closed first quadrant $\R^2_{\ge0}$.

\medskip
\emph{Paths.} An \emph{admissible path} is a sequence $x_{m:n}=(x_i)_{m\le i\le n}$ of lattice points, indexed by an integer interval, with
$x_{i+1}-x_i\in\{e_1,e_2\}$
and $\level(x_i)=i$ for all $i$.
Every vertex is thus indexed by its level. The path $x_{m:n}$ runs from $x_m$ to $x_n$, or from level $m$ to level $n$. Its \emph{length} is $n-m$. For $m\le k\le l\le n$ the \emph{segment} $x_{k:l}=(x_i)_{k\le i\le l}$ is again an admissible path. When $u\le v$ are vertices on the path we also write $x_{u:v}=x_{\level(u):\level(v)}$ for the segment from $u$ to $v$, and we mix the two conventions freely, as in $x_{u:n}$. In particular, a path from $u$ to $v$ is a path $x_{u:v}$ with $x_{\level(u)}=u$ and $x_{\level(v)}=v$. Such a path exists if and only if $u\le v$. If $x_{k:m}$ and $y_{m:n}$ satisfy $x_m=y_m$, their \emph{concatenation} is the path $x_{k:m}\,y_{m:n}$ whose segments between levels $k,m$ and $m,n$ are $x_{k:m}$ and $y_{m:n}$. 

\medskip
\emph{Weights.} The \emph{weight} of an admissible path $x_{m:n}$ is
\begin{equation}\label{eq:weight}
\pawgt(x_{m:n})=\sum_{i=m}^{n-1}\w_{x_i},
\end{equation}
the terminal vertex $x_n$  excluded. To highlight this point, the vertices $x_m,\dots,x_{n-1}$ are called the \emph{collected} vertices of the path $x_{m:n}$. Thus $\pawgt(x_{m:m})=0$ for a path of length zero, and weight is  additive over concatenation:  $\pawgt(x_{m:n})=\pawgt(x_{m:k})+\pawgt(x_{k:n})$ for $m\le k\le n$. For $u\le v$ in $\Z^2$ 
(coordinatewise order: $u\le v$ if $u_1\le v_1$ and $u_2\le v_2$) the \textit{last-passage time} from $u$ to $v$ is
\be\label{def:LPP}
\Lpp_{u,v}=\max\{\pawgt(x_{u:v}):\ x_{u:v}\text{ an admissible path from $u$ to $v$}\}, \ \text{ with } \  \Lpp_{u,u}=0.
\ee
A maximizing path is a \emph{geodesic}.

\subsection{Differentiability of the shape function}\label{sub:shape}
The shape theorem \cite[Theorem~5.1 and Remark~5.1(ii)]{Mar-04}, whose hypotheses hold under \ref{A2} and \ref{A3}, provides a deterministic, finite, concave, positively $1$-homogeneous function $g$ on $\R^2_{\ge0}$, continuous up to the boundary, such that almost surely
\begin{equation}\label{eq:shape}
\sup_{z\tsp\in\tsp\Z^2_{\ge0},\ \abs z_1\ge R}\frac{\abs{\Lpp_{0,z}-g(z)}}{\abs z_1}\longrightarrow0\qquad\text{as }R\to\infty.
\end{equation}
We record the consequences of \eqref{eq:shape} that are used below.
\begin{enumerate}[label={\rm(S\arabic*)}]
\item\label{S1} By translation invariance, \eqref{eq:shape} holds almost surely with $\Lpp_{0,z}$ replaced by $\Lpp_{u,u+z}$, for every fixed $u\in\Z^2$, and hence simultaneously for all $u$ in any fixed finite set.
\item\label{S2} For every $\e>0$,
\[
\sup_{u\le v,\ \abs{v-u}_1\ge R}\P\bigl\{\abs{\Lpp_{u,v}-g(v-u)}>\e\abs{v-u}_1\bigr\}\longrightarrow0\qquad\text{as }R\to\infty.
\]
%Indeed, by translation invariance the probability equals $\P\{\abs{\Lpp_{0,z}-g(z)}>\e\abs z_1\}$ with $z=v-u$, which is at most $\P\{\sup_{\abs{z'}_1\ge R}\abs{\Lpp_{0,z'}-g(z')}/\abs{z'}_1>\e\}$, and the latter tends to zero by \eqref{eq:shape}.
\item\label{S3} For every $z\in\R^2_{\ge0}$, $g(z)\ge\level(z)\,\E\w_0$. Indeed, for   $z\in\Z^2$ and $n\ge1$, $\Lpp_{0,nz}$ dominates the weight of any fixed admissible path from $0$ to $nz$, which is a sum of $n\level(z)$ i.i.d.\ weights. By the law of large numbers and \eqref{eq:shape}, $g(z)=\lim_n\Lpp_{0,nz}/n\ge\level(z)\E\w_0$. Rational points follow by homogeneity and all points by continuity.
\item\label{S4} Since $g$ is finite and concave on the open first quadrant, it is Lipschitz on every compact subset of the open first quadrant, hence (by homogeneity) on every cone $\{z\in\R_{\ge0}^2:\level(z)>0,\ \abs{\rcoord(z)}\le(\theta_0-\theta)\level(z)\}$ inside it, in the notation of Section~\ref{sub:coords} below. 
Lipschitz continuity  up to the axes is not expected.\footnote{By the boundary asymptotics of \cite[Theorem~2.4]{Mar-04}, when $\w_0$ has finite variance $\sigma^2>0$, $g(1,\alpha)=\E\w_0+2\sigma\sqrt\alpha+o(\sqrt\alpha)$ as $\alpha\downarrow0$, so $g$ is not Lipschitz up to the axes. Our paper does not need these boundary asymptotics and hence our assumption  \ref{A2} can allow an infinite  variance.}  
\end{enumerate}

Let $\Uset=\{(t,1-t):0\le t\le1\}$ and $\ri\Uset$ its relative interior. Here  is the main result of the paper.

\begin{theorem}\label{thm:main}
Under \ref{A1}--\ref{A3}, $g$ is differentiable at every point of $\ri\Uset$. Equivalently, $g$ is differentiable on the open quadrant $\R^2_{>0}$.
\end{theorem}
The equivalence is a consequence of homogeneity: $g$ is differentiable at $z\in\R_{>0}^2$ if and only if it is differentiable at $z/\abs z_1$.

\subsection{Flat segment}\label{sec:flat} 
When the maximum of the support of $\mu$ is an atom with a large enough mass, $g$ is linear on a closed cone of directions around the diagonal, called the \emph{percolation cone}. In these directions the differentiability of $g$ has already been established.  This section takes care of the percolation cone case by quoting  the relevant results from   \cite[Section 2.10]{Geo-Ras-Sep-15}.  These results  were adapted to the CGM from the oriented percolation review \cite{Dur-84} and from the  FPP work \cite{Dur-Lig-81,Mar-02,Auf-Dam-13}. 

Define 
\[\Mtop=\sup(\operatorname{supp}\mu)\in(-\infty,\infty]
\quad\text{and}\quad
p_1=\mu\{\Mtop\}.\] 
Let $\Ltop_{u,v}$ be the last-passage time under  the Bernoulli  weights $\one\{\w_z=\Mtop\}$, that is, the maximal number of collected vertices with weight $\Mtop$ along an admissible path from $u$ to $v$.    Its shape function $\gtop$ satisfies  $0\le\gtop(z)\le\level(z)$. We call
\[\Uperc=\{z\in\Uset:\ \gtop(z)=1\}\]
the \emph{flat segment}. If $\Mtop=\infty$, which \ref{A3} allows, then no vertex has weight $\Mtop$, so $p_1=0$, $\Ltop\equiv0$, $\gtop\equiv0$, and $\Uperc=\emptyset$, and everything in the rest of this subsection is trivially true or vacuous. 

Since $\gtop$ is concave and continuous on $\Uset$ and at most $1$ there, $\Uperc$ is a closed subsegment of $\Uset$, possibly empty, and it is symmetric about $(\tfrac12,\tfrac12)$ because the law of $\w$ is invariant under reflection $(i,j)\mapsto(j,i)$ across the diagonal.  The flat segment is described by the following theorem, which is a short step from \cite[Theorem~2.10]{Geo-Ras-Sep-15}. Let $\opc\in(0,1)$ be the critical density of oriented site percolation on $\Z^2$ with steps $\{e_1,e_2\}$. For $p_1>\opc$,   the edge speed of the cluster of the origin exists and equals a constant $\beta_{p_1}\in(\tfrac12,1)$ (\cite{Dur-84}, Section 2.10 in \cite{Geo-Ras-Sep-15}).

\begin{theorem}\label{thm:flat-top}
Assume $\Mtop<\infty$ and set $p_1=\mu\{\Mtop\}$. Then 
\begin{equation}\label{eq:flat-top}
\Uperc=\{z\in\Uset:\ g(z)=\Mtop\}
\end{equation}
and $g(z)<\Mtop$ for $z\in\Uset\setminus\Uperc$. 
Furthermore, 
\begin{enumerate}[label={\rm(\alph*)}]\itemsep=3pt 
\item\label{flat:a} If $p_1<\opc$, then $\Uperc=\emptyset$.
\item\label{flat:b} If $p_1=\opc$, then $\Uperc=\{(\tfrac12,\tfrac12)\}$ and $g$ is differentiable at $(\tfrac12,\tfrac12)$.
\item\label{flat:c} If $p_1>\opc$, then $\Uperc=\{z\in\Uset:\ 1-\beta_{p_1}\le z\cdot e_1\le\beta_{p_1}\}$ and $g$ is differentiable at the two endpoints of this segment.
\end{enumerate}
\end{theorem}

\begin{proof}
The weights $\bar\w_z=\w_z-\Mtop+1$ and $\one\{\w_z=\Mtop\}$ are i.i.d., $\le1$, equal to $1$ with probability $p_1<1$ by \ref{A1}, and have integrable negative parts by \ref{A2}. Their shape functions are $\bar g(z)=g(z)-(\Mtop-1)\level(z)$ and $\gtop$, respectively. These  equal $1$ at $z\in\Uset$ exactly when $g(z)=\Mtop$ and when $z\in\Uperc$, respectively. In each of the three cases, \cite[Theorem~2.10]{Geo-Ras-Sep-15} gives the same set for both, which proves \eqref{eq:flat-top}. Since $\w_z\le\Mtop$, every admissible path from $0$ to $z\in\Z^2_{\ge0}$ has weight at most $\Mtop\level(z)$. Hence $g(z)\le\Mtop\level(z)$ for all $z\in\R^2_{\ge0}$ by homogeneity and continuity, so $g\le\Mtop$ on $\Uset$, and \eqref{eq:flat-top} gives $g(z)<\Mtop$ for $z\in\Uset\setminus\Uperc$. Parts \ref{flat:a}--\ref{flat:c} are the statements of \cite[Theorem~2.10]{Geo-Ras-Sep-15} for $\bar\w$.  Since $g$ and $\bar g$ differ by a linear function, they are differentiable at the same points.
\end{proof}

By \eqref{eq:flat-top} and homogeneity, $g(z)=\Mtop\level(z)$ is linear on the cone spanned by $\Uperc$, so $g$ is differentiable at the points of the relative interior of $\Uperc$. Together with Theorem~\ref{thm:flat-top}\ref{flat:b}--\ref{flat:c}, this shows that $g$ is differentiable at every point of $\Uperc\cap\ri\Uset$. The remaining  task  is   to prove differentiability at the points $\Dir\in(\ri\Uset)\setminus\Uperc$, that is, at the points with $\gtop(\Dir)<1$. 
%If $\mu$ has no atom at $\Mtop$, then $\gtop\equiv0$ and $\Uperc=\emptyset$. Conversely, if $p_1>0$ then $\gtop(z)\ge p_1\level(z)$, since a fixed path from $0$ to $z$ collects on average $p_1\level(z)$ vertices of weight $\Mtop$.

\subsection{Coordinates adapted to a direction}\label{sub:coords}
By Section~\ref{sec:flat}, it suffices to prove Theorem~\ref{thm:main} at the points of $\ri\Uset$ outside the closed  flat segment. Fix $\Dir=(\Dir_1,1-\Dir_1)\in(\ri\Uset)\setminus\Uperc$, $0<\Dir_1<1$, for the rest of the paper, so that $\gtop(\Dir)<1$. 

\begin{figure}[htb]
\centering
\begin{tikzpicture}[scale=2.6,>={Stealth[scale=1.3]}]
  % cone of theta-interior vectors, theta_0-theta=0.15
  \fill[gray!18] (0,0) -- (0.88,1.32) -- (0.22,1.98) -- cycle;
  \draw[gray!60] (0,0) -- (0.88,1.32);
  \draw[gray!60] (0,0) -- (0.22,1.98);
  \node[gray!90!black,anchor=west] at (0.62,2.12) {$\{\abs\rcoord\le(\theta_0-\theta)\level\}$};
  % axes
  \draw[->] (0,0) -- (2.35,0) node[right] {$z_1$};
  \draw[->] (0,0) -- (0,2.35) node[above] {$z_2$};
  % level lines
  \draw[very thick] (1,0) -- (0,1) node[pos=0.2,above right=-4pt] {$\Uset=\{\level=1\}$};
  \draw[dashed] (2.2,-0.2) -- (-0.2,2.2) node[pos=1,left=1pt] {$\{\level=\level(z)\}$};
  % ray through xi
  \draw[dotted,thick] (0,0) -- (0.5,1.5);
  \draw[->,thick] (0,0) -- (0.25,0.75);
  \draw[->,thick] (0,0) -- (0.5,-0.5) node[right] {$\ead$};
  % decomposition of z
  \draw[->,thick] (0.5,1.5) -- (1.2,0.8) node[midway,above right=-2pt] {$\rcoord(z)\,\ead$};
  % points
  \fill (0.25,0.75) circle (0.033) node[left=2pt] {$\Dir$};
  \fill (0.5,1.5) circle (0.033) node[left=2pt] {$\level(z)\,\Dir$};
  \fill (1.2,0.8) circle (0.033) node[right=1pt] {$z=\level(z)\,\Dir+\rcoord(z)\,\ead$};
  \fill (0,0) circle (0.026) node[below left] {$0$};
\end{tikzpicture}
\caption{\small The coordinates $\level(z)$ and $\rcoord(z)$ of a point $z$, here with $\Dir_1=0.25$, $\level(z)=2$, and $\rcoord(z)=0.7$. The level $\level$ is constant on antidiagonal lines. The ray from $0$ through the fixed direction $\Dir$ meets the level line of $z$ at $\level(z)\Dir$, and $\rcoord(z)$ is the signed distance, in units of $\ead$, from $\level(z)\Dir$ to $z$. Shaded: the cone of $\theta$-interior vectors $\{\level>0,\ \abs\rcoord\le(\theta_0-\theta)\level\}$ for $\theta_0-\theta=0.15$. It is symmetric about the ray through $\Dir$ and, as $\theta\downarrow0$, fills the open quadrant only up to the ray through $\Dir-\theta_0\ead$, here the $z_2$-axis. The point $z$ drawn is not $\theta$-interior.}
\label{fig:coords}
\end{figure}
Besides the level \eqref{eq:level}, for $z=(z_1,z_2)\in\R^2$ define the \emph{transverse coordinate} by
\begin{equation}\label{eq:coords}
\rcoord(z)=(1-\Dir_1)z_1-\Dir_1z_2,
\end{equation}
so that $\level(\Dir)=1$, $\rcoord(\Dir)=0$, $\level(\ead)=0$, $\rcoord(\ead)=1$.  Every $z\in\R^2$ is expressed uniquely as
\begin{equation}\label{eq:decomp}
z=\level(z)\tspa\Dir+\rcoord(z)\tspa\ead ,\quad\text{that is,}\quad z_1=\Dir_1\level(z)+\rcoord(z)\quad\text{and}\quad
z_2=(1-\Dir_1)\level(z)-\rcoord(z). 
\end{equation}
 Figure~\ref{fig:coords} illustrates.  

Both $\level$ and $\rcoord$ are linear functionals. An admissible step increases $\level$ by $1$ and changes $\rcoord$ by $\rcoord(e_1)=1-\Dir_1$ or $\rcoord(e_2)=-\Dir_1$. In particular, $\abs{\rcoord(x_{i+1})-\rcoord(x_i)}\le1$ along admissible paths, and
\begin{equation}\label{eq:lip-r}
\abs{\rcoord(z)-\rcoord(z')}\le\abs{z-z'}_1\qquad\text{for all }z,z' .
\end{equation}
By \eqref{eq:decomp}, a vector $z$ with $\level(z)=n>0$ lies in the open quadrant if and only if $-\Dir_1n<\rcoord(z)<(1-\Dir_1)n$, and on the axes if and only if $\rcoord(z)\in\{-\Dir_1n,(1-\Dir_1)n\}$. On a fixed level $I\in\Z$ the lattice points are $z=(j,I-j)$, $j\in\Z$, with $\rcoord(z)=j-\Dir_1I$: their transverse coordinates form the set $\Z-\Dir_1I$, so every half-open unit interval contains exactly one of them.

\begin{figure}[htb]
\centering
\begin{tikzpicture}[scale=0.62,>=Stealth]
  %% lattice
  \foreach \x in {-1,...,5} \foreach \y in {-3,...,8} \fill[gray!50] (\x,\y) circle (0.05);
  %% axes
  \draw[->,gray] (-1.5,0) -- (5.6,0) node[right,black] {$z_1$};
  \draw[->,gray] (0,-3.5) -- (0,8.6) node[above,black] {$z_2$};
  %% blocks: lines r=rho and r=rho-1, direction xi=(0.25,0.75)
  \draw[dotted,thick] (-1.25,0.25) -- (1.25,7.75) node[above right=-3pt] {$\rcoord=-1$};
  \draw[densely dotted,thick] (2.25,6.75) -- (-0.25,-0.75) node[below left=-4pt] {$\rcoord=0$};
  \draw[dotted,thick] (4.75,8.25) -- (1,-3) node[below] {$\rcoord=1.5$};
  \draw[densely dotted,thick] (2.25,-3.25) -- (6,8) node[above right=-2pt] {$\rcoord=2.5$};
  %% staircase rho=0
  \draw[very thick] (0,0) -- (0,3) -- (1,3) -- (1,6) -- (2,6);
  \foreach \pt in {(0,0),(0,1),(0,2),(0,3),(1,3),(1,4),(1,5),(1,6),(2,6)} \fill \pt circle (0.12);
  %\node[below left=-1pt] at (0,0) {$\hat z(0,0)$};
  %\node[right=3pt] at (2,6) {$\hat z(8,0)$};
  %% staircase rho=2.5
  \draw[very thick] (2,-2) -- (2,-1) -- (3,-1) -- (3,2) -- (4,2) -- (4,4);
  \foreach \pt in {(2,-2),(2,-1),(3,-1),(3,0),(3,1),(3,2),(4,2),(4,3),(4,4)} \fill \pt circle (0.12);
  %\node[right] at (2.1,-2.2) {$\hat z(0,2.5)$};
  %\node[right=4pt] at (4,4.45) {$\hat z(8,2.5)$};
  %% direction xi
  \draw[->,thick] (0,0) -- (0.25,0.75);
  \node[above right=-1pt] at (0.1,0.55) {$\Dir$};
\end{tikzpicture}
\caption{\small The staircases $I\mapsto \hat z(I,\rho)$ for $\rho=0$ (left) and $\rho=2.5$ (right), with $\Dir_1=0.25$, between levels $0$ and $8$. Each is the unique admissible path with all vertices in the half-open strip $\{\rho-1<\rcoord\le\rho\}$ bounded by the dotted lines (the line $\rcoord=\rho$ included), which is the block of transverse width $1$ along the line through $\rho\ead$ in direction $\Dir$. The vertex of level $I$ is $(j,I-j)$ with $j=\lfloor\rho+\Dir_1I\rfloor$, so the $e_1$ steps occur at the levels where $\rho+\Dir_1I$ crosses an integer. All staircases look alike up to translation and the placement of the $e_1$ steps.}
\label{fig:staircase}
\end{figure}
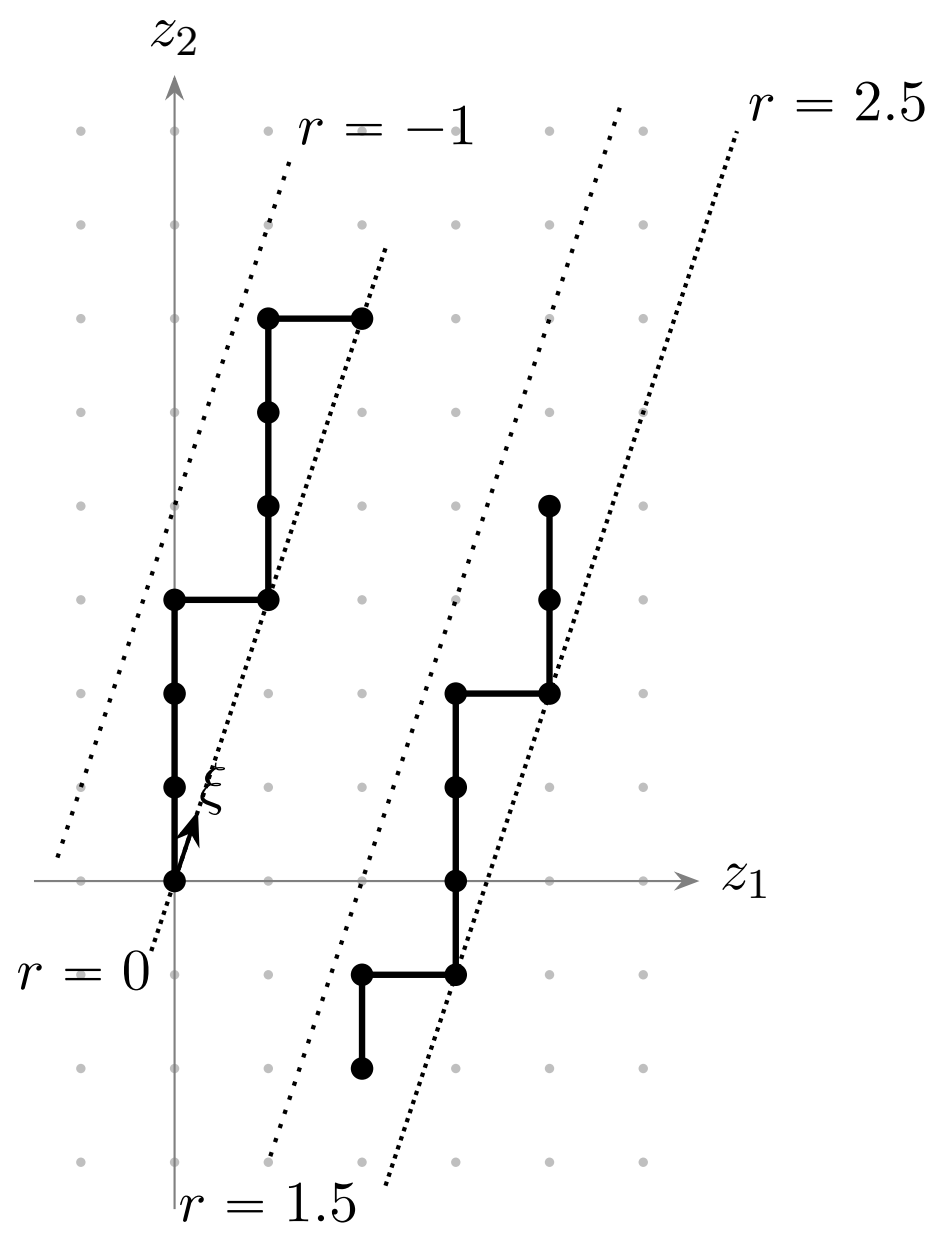
For $\rho\in\R$ and level index $I\in\Z$,
\begin{align}\label{stairdef}
\text{$\hat z(I,\rho)$ denotes the unique lattice point $z$ with $\level(z)=I$ and $\rcoord(z)\in(\rho-1,\rho]$.}
\end{align}
We observe that consecutive points $\hat z(I,\rho),\hat z(I+1,\rho)$ differ by $e_1$ or $e_2$.  By definition, their difference has level $1$ and transverse coordinate in $(-1,1)$.  Hence for some $j\in\Z$, $\hat z(I+1,\rho)-\hat z(I,\rho)=(j,1-j)$ with $\rcoord=j-\Dir_1$.
This transverse coordinate lies in $(-1,1)$ only for $j\in\{0,1\}$. 

Thus $I\mapsto \hat z(I,\rho)$ is an admissible path we call  the \emph{staircase at $\rho$}.  %All its vertices have transverse coordinate  $\rcoord\in(\rho-1,\rho]$.  
Figure~\ref{fig:staircase} illustrates.

For the remainder of the paper, fix the parameter $\theta_0=\min\{\Dir_1,1-\Dir_1\}$. For $\theta\in(0,\theta_0)$, we say
\[
\text{$z\in\R^2$ with $\level(z)>0$ is \emph{$\theta$-interior} if $\abs{\rcoord(z)}\le(\theta_0-\theta)\level(z)$.}
\]
Then, by \eqref{eq:decomp}, $z_1=\Dir_1\level(z)+\rcoord(z)\ge(\Dir_1-\theta_0+\theta)\level(z)\ge\theta\level(z)$, and likewise $z_2\ge\theta\level(z)$.
Consequently, $z$ lies in the open first quadrant and $z/\abs z_1=\Dir+(\rcoord(z)/\level(z))\ead$ lies in the compact subset $\{\Dir+h\ead:\abs h\le\theta_0-\theta\}$ of $\ri\Uset$. Throughout, $\Lip_\theta$ denotes a Lipschitz constant of $g$ (with respect to $\abs\cdot_1$) on the cone of $\theta$-interior vectors, which exists by \ref{S4}.

\subsection{Extreme supporting functionals}\label{sub:extreme}

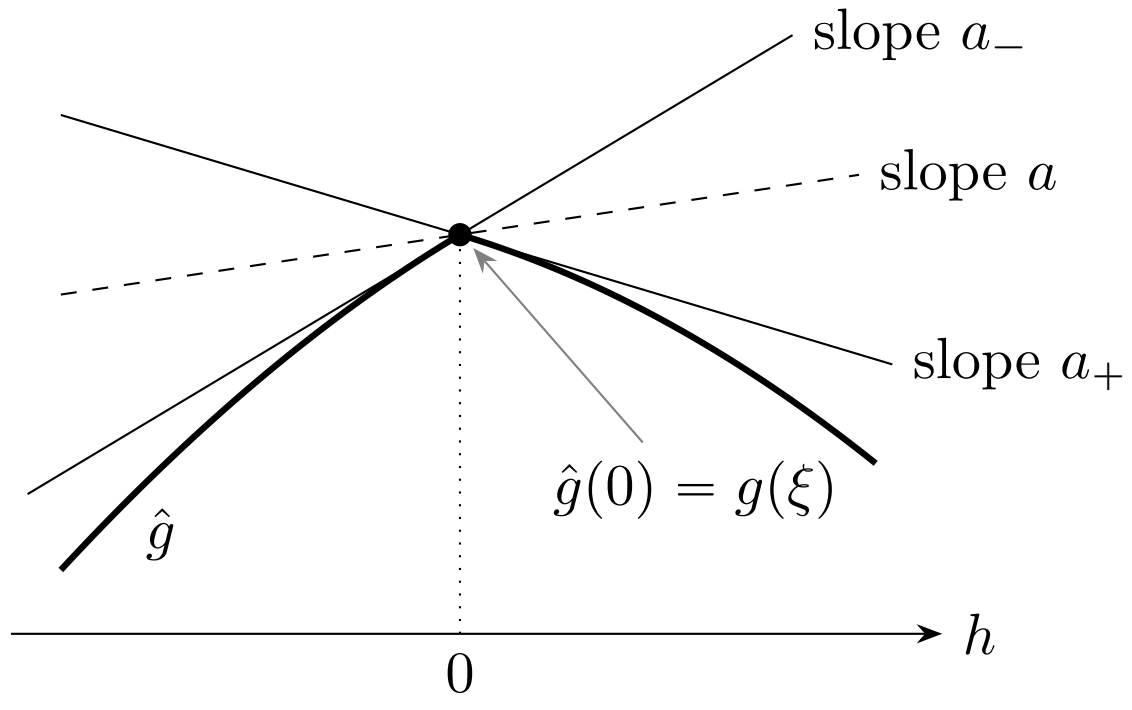
\begin{figure}[htb]
\centering
\begin{tikzpicture}[scale=1.1,>=Stealth]
  \draw[->] (-2.7,-0.4) -- (2.9,-0.4) node[right] {$h$};
  \draw[dotted] (0,-0.4) -- (0,2);
  \node[below] at (0,-0.4) {$0$};
  \draw[thin] (-2.6,0.44) -- (2.0,3.2) node[right] {slope $\ghatprm$};
  \draw[thin] (-2.4,2.72) -- (2.6,1.22) node[right] {slope $\ghatprp$};
  \draw[thin,dashed] (-2.4,1.64) -- (2.4,2.36) node[right] {slope $a$};
  \draw[very thick,domain=-2.4:0,samples=40] plot (\x,{2+0.6*\x-0.1*\x*\x});
  \draw[very thick,domain=0:2.5,samples=40] plot (\x,{2-0.3*\x-0.1*\x*\x});
  \fill (0,2) circle (0.07);
  \node at (-1.8,0.2) {$\hat g$};
  \draw[->,gray] (1.1,0.75) -- (0.08,1.92);
  \node[below right=-2pt] at (0.5,0.7) {$\hat g(0)=g(\Dir)$};
\end{tikzpicture}
\caption{\small The figure depicts the situation that our proof will rule out, namely,  a  corner of $\hat g$ at $0$. The two extreme supporting lines, with slopes $\ghatprm$ and $\ghatprp$, touch $\hat g$ to first order on their own sides. The dashed line has the average slope $a$ and touches $\hat g$ only at $0$.}
\label{fig:corner}
\end{figure}
Let $\hat g(h)=g(\Dir+h\ead)$ for $-\Dir_1\le h\le1-\Dir_1$. In particular, $\hat g(0)=g(\Dir)$.
By \eqref{eq:decomp} this parametrizes the segment $\Uset$, and $\hat g$ is a finite concave continuous function. Write $\ghatpr{\pm}$ for its one-sided derivatives at $0$:
\[
\ghatprm=\hat g'(0-)=\lim_{h\uparrow0}\frac{\hat g(h)-\hat g(0)}{h}\ \ge\ \ghatprp=\hat g'(0+)=\lim_{h\downarrow0}\frac{\hat g(h)-\hat g(0)}{h}.
\]
Both limits exist and are finite by concavity.   Figure~\ref{fig:corner} illustrates. A finite concave function on an open subset of $\R^2$ is differentiable at a point if and only if its two partial derivatives (in any two independent directions) exist there \cite[Theorem~25.2]{Roc-70}. Since $g(\Dir+h\Dir)=(1+h)g(\Dir)$ is differentiable in $h$, differentiability of $g$ at $\Dir$ is equivalent to $\ghatprm=\ghatprp$. 
Define the linear functionals
\begin{equation}\label{eq:supports}
\begin{split}
&f_\sig(z)=\ghatz\,\level(z)+\ghatpr{\sig}\,\rcoord(z)\quad\text{for signs } \ \sigg\in\{-,+\},\\[3pt] 
&f=\tfrac12(f_-+f_+),\quad a=\tfrac{\ghatprm+\ghatprp}2,
\quad\text{and}\quad \halfgap =\tfrac{\ghatprm-\ghatprp}{2}\ge0. 
\end{split}
\end{equation}
Then  
\be\label{t.f(z)} f(z)=\ghatz\level(z)+a\rcoord(z) \ee
and
\begin{equation}\label{eq:fpm}
f_\pm(z)=f(z)\mp \halfgap \,\rcoord(z).
\end{equation}
In the nondifferentiable case, $\ghatprp<a<\ghatprm$ and $\halfgap >0$.

%\eqref{eq:contact}: small error $o(h)$ requires that the sign of $h$ match sign in $f_\pm$.

\begin{lemma}\label{lem:contact}
For $\sigg\in\{-,+\}$, $f_\sig\ge g$ on $\R^2_{\ge0}$, $f_\sig(\Dir)=f(\Dir)=g(\Dir)$, and
\begin{equation}\label{eq:contact}
f_\sig(\Dir+h\ead)-g(\Dir+h\ead)=o(\abs h)\qquad\text{as }h\to0\text{ with }\sigg h>0.
\end{equation}
\end{lemma}

\begin{proof}
Since $\ghatpr{\sig}$ lies between the two one-sided derivatives of the concave function $\hat g$ at $0$, the line $h\mapsto\hat g(0)+\ghatpr{\sig} h$ is a supporting line: $\hat g(h)\le\hat g(0)+\ghatpr{\sig} h$ for all $h\in[-\Dir_1,1-\Dir_1]$. For $z\in\R^2_{\ge0}$ with $\level(z)=n>0$, homogeneity and \eqref{eq:decomp} give $g(z)=n\hat g(\rcoord(z)/n)$ with $\rcoord(z)/n\in[-\Dir_1,1-\Dir_1]$, hence $g(z)\le n \ghatz+\ghatpr{\sig}\rcoord(z)=f_\sig(z)$, while for $z=0$ both sides vanish. The equalities at $\Dir$ are immediate, and \eqref{eq:contact} is the definition of the one-sided derivative $\ghatpr{\sig}=\hat g'(0\sigg)$:
\begin{align*}  \text{for $\sigg h>0$, } \quad   f_\sig(\Dir+h\ead)-g(\Dir+h\ead) &= \hat g(0) + \hat g'(0\sigg)\tsp h - \hat g(h)
\\
& %=\ghatpr{\sig}h-\bigl(\hat g(h)-\hat g(0)\bigr)
=h\bigl[\hat g'(0\sigg)-\tfrac{\hat g(h)-\hat g(0)}h\bigr]=o(\abs h).
\qedhere
\end{align*}
\end{proof}

\smallskip

\subsection{Path weights   and supporting functionals}
Throughout, the weight $\pawgt(x_{m:n})$ of a path is compared with $f(x_n-x_m)$ and $f_\sig(x_n-x_m)$, $\sigg\in\{-,+\}$. A bound of the form $\pawgt\le f_\sig+b$ or $\pawgt\ge f_\sig-b$, stated for a path $x_{m:n}$, means $\pawgt(x_{m:n})\le f_\sig(x_n-x_m)+b$ or $\pawgt(x_{m:n})\ge f_\sig(x_n-x_m)-b$, and likewise with $f$.
Since $f_\sig$ and $f$ are linear, $f_\sig(x_n-x_m)=f_\sig(x_k-x_m)+f_\sig(x_n-x_k)$ for $m\le k\le n$, likewise for $f$, and \eqref{eq:weight} gives the same for $\pawgt$. Hence bounds of the form $\pawgt\le f_\sig+b$ or $\pawgt\ge f_\sig-b$ for pieces of a path add up to a bound of the same form for the whole path, with the sum of the constants $b$. 

By \eqref{eq:fpm}, $f_\pm(x_n-x_m)=f(x_n-x_m)\mp \halfgap \,\rcoord(x_n-x_m)$, so that
\begin{equation}\label{eq:defect-min}
\min_{\sig\tsp\in\tsp\{-,+\}}f_\sig(x_n-x_m)=f(x_n-x_m)-\halfgap \,\abs{\rcoord(x_n-x_m)}\ \le\ f(x_n-x_m) .
\end{equation}
The gap between one of $f_\pm$ and $f$ created by nondifferentiability $(\halfgap >0)$ at $\Dir$ is one crucial lever of the final stage of the proof of Theorem \ref{thm:main}.  

By concavity  (Lemma~\ref{lem:contact}) $g(x_n-x_m)\le f_\sig(x_n-x_m)$, and by the definition \eqref{def:LPP} of last-passage times $\pawgt(x_{m:n})\le\Lpp_{x_m,x_n}$. So each admissible path $x_{m:n}$ started at  $x_m=u$ satisfies 
\begin{equation}\label{eq:defect-lower}
\pawgt(x_{m:n})\ \le\ f_\sig(x_n-x_m)+(\Lpp_{u,x_n}-g(x_n-u)) .
\end{equation}
Combined with \ref{S1}, this says: almost surely, for every fixed $u$ and every $\e>0$ there exists $R(\e,u,\w)<\infty$ such that $\pawgt(x_{m:n})\le f_\sig(x_n-x_m)+\e(n-m)$ for every admissible path from $u$ of length $n-m\ge R$.

On the other hand, along the direction $\Dir$ the geodesic collects weight $f+o(n)$, relative to the \emph{average} supporting functional $f$, because $f(\Dir)=g(\Dir)$. To see this let $z_n$ be the lattice point on level $n$ with $\rcoord(z_n)\in(-1,0]$ and $x^{(n)}_{0:n}$ a geodesic from $0$ to $z_n$. Then $\abs{z_n-n\Dir}_1=\abs{\rcoord(z_n)}\,\abs\ead_1\le2$ and hence $z_n$, for $n\ge2/\theta_0$, is $\theta_0/2$-interior. Consequently, $\abs{g(z_n)-n\ghatz}\le2\tsp\Lip_{\theta_0/2}$, while $f(z_n)=n\ghatz+a\rcoord(z_n)$ differs from $n\ghatz$ by at most $\abs a$. Hence, almost surely, by \eqref{eq:shape},
\begin{equation}\label{eq:geodesic-defect}
\pawgt(x^{(n)}_{0:n})=\Lpp_{0,z_n}=g(z_n)+o(n)=f(z_n)+o(n).
\end{equation}

 The proof of Theorem~\ref{thm:main} shows that if 
 $\halfgap >0$, then every long admissible path $x_{0:n}$ has $\pawgt(x_{0:n})\le f(x_n-x_0)-Cn$ for a constant $C>0$, contradicting \eqref{eq:geodesic-defect}. The nondifferentiability  hypothesis $\halfgap >0$ enters exactly once, in Lemma~\ref{lem:fulldefect} at the end of the proof, through \eqref{eq:defect-min}: a path that moves transversally distance $s$ pays $\halfgap s$ relative to $f$, and (by the results of Sections~\ref{sec:crossings}--\ref{sec:endpoints}) it cannot recover this by exceeding the favorable extreme supporting functional.

The argument does not require $\Dir$ to be an isolated corner or $g$ to be strictly concave near $\Dir$.  Only the contact property \eqref{eq:contact} is used, and it holds whether or not $\Dir$ is the endpoint of a linear piece of $g$.  If $g$ is linear on one side of $\Dir$, the contact on that side is exact. 

\section{Strategy of the proof}\label{sec:strategy}

This section summarizes the main  points of the argument.  Further discussion appears in Section \ref{sec:remarks} after the main proof. 
A direction $\Dir\in\ri\Uset$   is fixed and lattice  paths are viewed through the $\Dir$-adapted $(\level, \rcoord)$  coordinates of \eqref{eq:decomp}.  Assume that $g$ is not differentiable at $\Dir$, equivalently, that $\halfgap >0$ in \eqref{eq:supports}. 

\begin{figure}[htb]
\centering
\begin{tikzpicture}[scale=0.26,>=Stealth]
  \filldraw[fill=gray,fill opacity=0.22,draw=gray!70,thin] (-1.00,1.00) -- (10.20,5.80) -- (12.20,3.80) -- (1.00,-1.00) -- cycle;
  \filldraw[fill=gray,fill opacity=0.22,draw=gray!70,thin] (10.00,6.00) -- (21.20,10.80) -- (23.20,8.80) -- (12.00,4.00) -- cycle;
  \filldraw[fill=gray,fill opacity=0.22,draw=gray!70,thin] (14.00,6.00) -- (25.20,10.80) -- (27.20,8.80) -- (16.00,4.00) -- cycle;
  \filldraw[fill=gray,fill opacity=0.22,draw=gray!70,thin] (25.00,11.00) -- (36.20,15.80) -- (38.20,13.80) -- (27.00,9.00) -- cycle;
  \draw[<->,gray!80!black,thin] (-2.40,0.40) -- (-0.40,-1.60);
  \node[gray!80!black,anchor=north,font=\scriptsize] at (-2,-0.50) {$2s_*$};
  \draw[<->,gray!80!black,thin] (-2.40,2.40) -- (8.80,7.20);
  \node[gray!80!black,anchor=south,font=\scriptsize] at (3.20,5.00) {$2n_*$};
  \draw[black,line width=1.7pt,line join=round] (0.00,0.00) -- (1.00,0.00) -- (2.00,0.00) -- (3.00,0.00) -- (3.00,1.00) -- (4.00,1.00) -- (5.00,1.00) -- (5.00,2.00) -- (6.00,2.00) -- (7.00,2.00) -- (7.00,3.00) -- (8.00,3.00) -- (9.00,3.00) -- (9.00,4.00) -- (10.00,4.00) -- (11.00,4.00) -- (11.00,5.00);
  \draw[black,line width=0.7pt,line join=round] (11.00,5.00) -- (12.00,5.00) -- (13.00,5.00) -- (14.00,5.00) -- (15.00,5.00);
  \draw[black,line width=1.7pt,line join=round] (15.00,5.00) -- (15.00,6.00) -- (16.00,6.00) -- (17.00,6.00) -- (17.00,7.00) -- (18.00,7.00) -- (19.00,7.00) -- (19.00,8.00) -- (20.00,8.00) -- (21.00,8.00) -- (22.00,8.00) -- (22.00,9.00) -- (23.00,9.00) -- (24.00,9.00) -- (25.00,9.00) -- (25.00,10.00) -- (26.00,10.00);
  \draw[black,line width=0.7pt,line join=round] (26.00,10.00) -- (26.00,11.00) -- (27.00,11.00) -- (28.00,11.00) -- (29.00,11.00) -- (29.00,12.00) -- (30.00,12.00);
  \fill (0.00,0.00) circle (0.35);
  \fill (11.00,5.00) circle (0.35);
  \fill (15.00,5.00) circle (0.35);
  \fill (26.00,10.00) circle (0.35);
  \fill (30.00,12.00) circle (0.35);
  \node[anchor=east,font=\scriptsize] at (0,0) {$0$};
  \node[anchor=west,font=\scriptsize] at (30,12.00) {$x_n$};
\end{tikzpicture}
\caption{\small A path $x_{0:n}$ with $x_0=0$ decomposed into stages. Each stage is the piece of the path inside a parallelogram of antidiagonal width $2s_*$ and diagonal length $2n_*$ centered at the stage's first vertex (shaded), and it ends either when the path leaves the parallelogram through its top side, a longitudinal stage (thick), or through one of its long sides, a transverse stage (thin). The last piece is the remainder.}
\label{fig:strategy}
\end{figure}
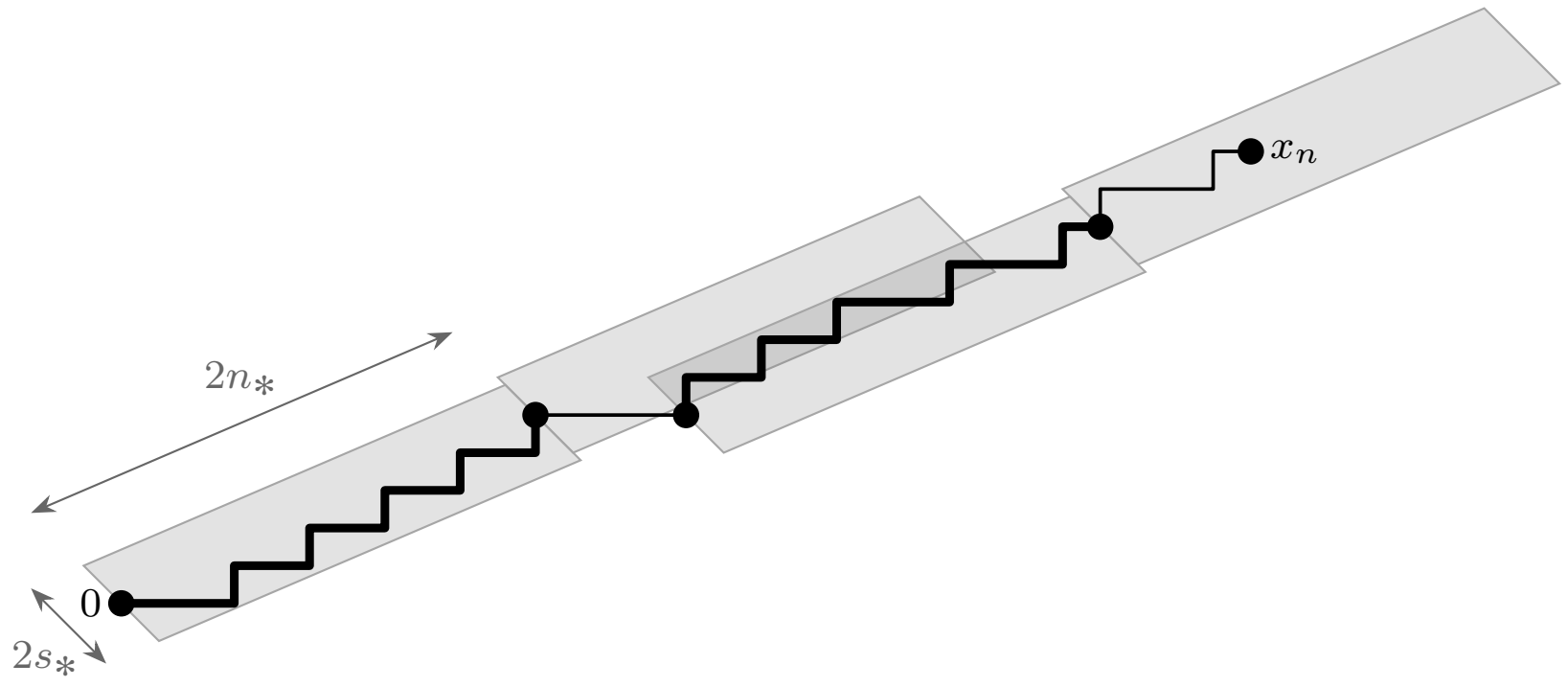

\medskip 

\emph{Characteristic scale.}  Integers   $n_*\gg s_*$ are found  so that statements \ref{strat.a}--\ref{strat.b} below hold for parallelograms of antidiagonal width of order $s_*$ and diagonal length of order $n_*$.
\begin{enumerate}[label={\rm(\alph*)}, ref={\rm(\alph*)}] \itemsep=2pt
\item\label{strat.a}  With high probability, every admissible path   inside such a parallelogram satisfies $\pawgt\le f_\pm+\tfrac{\halfgap}2\tsp s_*$   (Proposition~\ref{prop:endpoints}).
\item\label{strat.b} The probability that some admissible path crosses the parallelogram along its full length with weight $\pawgt>f-Ks_*$, where $K$ is a fixed large constant, is  small  (Section~\ref{sec:scales}).
\end{enumerate}

Property \ref{strat.a} comes from the shape theorem through a coarse-grained oriented percolation argument (Section~\ref{sec:crossings}), which first gives it for paths that cross the parallelogram from bottom to top.  A multiscale construction (Section~\ref{sec:endpoints}) extends it to paths with arbitrary endpoints.

Property \ref{strat.b}  defines $n_*$. A change of measure (Section~\ref{sec:scales}) shows that crossings with $\pawgt>f-Ks_*$ are rare once the length is of order $s_*^3$, and $n_*$ is defined as the minimal length at which such crossings are rare. The multiscale construction of \ref{strat.a} uses the minimality in this definition on a whole sequence of smaller scales. Crossings shorter than the characteristic length are not rare, and they serve as its building blocks. 
\medskip

\emph{Stages.} Let $x_{0:n}$ be an admissible path from the origin and cut it into stages, as in Figure~\ref{fig:strategy}.  Each stage starts where the previous one ended, and it ends when the path has advanced $2n_*$ levels or has moved by $\pm s_*$ in the antidiagonal direction. A stage therefore stays in a parallelogram of width about $2s_*$ and length $2n_*$, and it is one of two kinds. 
\begin{itemize} \itemsep=2pt 
    \item 
A \emph{longitudinal} stage runs the full length $2n_*$. It contains a crossing of a parallelogram of length $n_*$, so by \ref{strat.b}, with high probability, $\pawgt\le f-Ks_*$ along that crossing, and \ref{strat.a} controls the rest; altogether $\pawgt\le f-\tfrac{\halfgap}2\tsp s_*$ for the stage. 
\item  A \emph{transverse} stage leaves through one of the long sides, so it has antidiagonal displacement of absolute value   $\ge s_*$. By \ref{strat.a}, $\pawgt\le f_\pm+\tfrac{\halfgap}2\tsp s_*$, and since $f_\pm=f\mp \halfgap \,\rcoord$ we choose the sign for which $f_\pm=f-\halfgap \abs\rcoord\le f-\halfgap s_*$. Again $\pawgt\le f-\tfrac{\halfgap}2\tsp s_*$. This is   where the assumption $\halfgap >0$ enters.  Without it, $\halfgap =0$ and a transverse move costs nothing.
\end{itemize}

For a stage at which \ref{strat.a} or \ref{strat.b} fails  the bound $\pawgt\le f-\tfrac{\halfgap}2\tsp s_*$ may fail, and $\pawgt$ may even exceed $f$. But each excess of order $s_*$ over $f$ requires another rare event inside the stage and the sum of excesses  over the entire path is small (Section \ref{sec:witnesses}).  

\medskip

\emph{Completion of the proof.} A path of length $n$ has about $n/(2n_*)$ stages. Section~\ref{sec:witnesses} shows  that the stages at which \ref{strat.a} or \ref{strat.b} fails are few, uniformly over all paths with a given number of stages. The rare events responsible for the failures in distinct stages involve disjoint sets of vertices, so their joint occurrence can be controlled by   the van den Berg--Kesten inequality.    Choice of $(s_*,n_*)$  fixes the  failure probabilities  in \ref{strat.a}--\ref{strat.b} small enough to overcome the exponential number of  arrangements of   stages along a path. Bad stages cost at most $p\,\frac{n}{2n_*}\cdot\frac{\halfgap}2\tsp s_*$ in total, with $p$ small. Adding up,
\[
\pawgt(x_{0:n})\ \le\ f(x_n)-(1-p)\frac{n}{2n_*}\cdot\frac{\halfgap}2\tsp s_*+p\,\frac{n}{2n_*}\cdot\frac{\halfgap}2\tsp s_*\ =\ f(x_n)-Cn
\]
with a constant $C>0$ that depends on $\halfgap $, $s_*$, and $n_*$ but not on $n$. Apply this to the geodesic in \eqref{eq:geodesic-defect} to produce a  contradiction. The source of the contradiction is the assumption $\halfgap >0$. Thus $\halfgap =0$ and Theorem \ref{thm:main} is proved. 

\medskip

\emph{Flat segment.}  The nondifferentiability assumption implies that $\Dir$ lies off the closed flat segment $\Uperc$. The proof uses this assumption in one place only.   Namely, the modification argument of 
Lemma~\ref{lem:tilt} requires that every path in a $\Dir$-directed  tube has a positive fraction of vertices whose weight is below the maximum.

\medskip 

\emph{Bounded weights versus moment assumptions.} The main part of the paper proves Theorem \ref{thm:main} for bounded weights $\abs{\w_x}\le M$. The proof under moment assumptions  is   in Appendix~\ref{app:general}. Boundedness of the weights is used in two ways: sums of weights along deterministic paths obey Hoeffding's inequality \eqref{eq:hoeffding}, and a single step changes the weight of a path by at most $M$.

\section{Localized crossings}
\label{sec:crossings}

The development of the proof of Theorem \ref{thm:main} begins.  
This section needs no assumption about the percolation cone and so  $\Dir\in\ri\Uset$ is an arbitrary  fixed direction.   The condition $\Dir\notin\Uperc$ enters first in Lemma~\ref{lem:tilt} in Section \ref{sec:cubic}.

\begin{figure}[htb]
\centering
\begin{tikzpicture}[scale=0.34,>=Stealth]
  \clip (-9.6,-4.6) rectangle (13.6,16.2);
  %% lattice grid
  \foreach \x in {-10,...,14} \draw[gray!30,very thin] (\x,-5) -- (\x,17);
  \foreach \y in {-5,...,17} \draw[gray!30,very thin] (-10,\y) -- (14,\y);
  %% level lines I_1 and I_2 (dotted)
  \draw[gray!60,densely dotted] (-4.60,4.60) -- (4.60,-4.60);
  \node[gray,anchor=east,font=\scriptsize] at (-4.70,4.60) {$\level=I_1$};
  \draw[gray!60,densely dotted] (2.60,15.40) -- (11.80,6.20);
  \node[gray,anchor=east,font=\scriptsize] at (2.50,15.40) {$\level=I_2$};
  %% the block: parallelogram with sides along xi and e_ad
  \filldraw[fill=gray!18,draw=black,thick] (-2.00,2.00) -- (5.20,12.80) -- (9.20,8.80) -- (2.00,-2.00) -- cycle;
  %% lattice points of the block
  \foreach \p in {(-2,2),(-1,1),(-1,2),(-1,3),(0,0),(0,1),(0,2),(0,3),(0,4),(0,5),(1,-1),(1,0),(1,1),(1,2),(1,3),(1,4),(1,5),(1,6),(2,-2),(2,-1),(2,0),(2,1),(2,2),(2,3),(2,4),(2,5),(2,6),(2,7),(2,8),(3,0),(3,1),(3,2),(3,3),(3,4),(3,5),(3,6),(3,7),(3,8),(3,9),(4,1),(4,2),(4,3),(4,4),(4,5),(4,6),(4,7),(4,8),(4,9),(4,10),(4,11),(5,3),(5,4),(5,5),(5,6),(5,7),(5,8),(5,9),(5,10),(5,11),(5,12),(6,4),(6,5),(6,6),(6,7),(6,8),(6,9),(6,10),(6,11),(6,12),(7,6),(7,7),(7,8),(7,9),(7,10),(7,11),(8,7),(8,8),(8,9),(8,10),(9,9)} \fill \p circle (0.09);
  %% a crossing: admissible path from level I_1 to level I_2 inside the block
  \draw[line width=1.5pt,line join=round] (1,-1) -- (2,-1) -- (2,0) -- (2,1) -- (2,2) -- (2,3) -- (2,4) -- (2,5) -- (3,5) -- (4,5) -- (4,6) -- (5,6) -- (6,6) -- (6,7) -- (6,8) -- (7,8) -- (7,9) -- (7,10) -- (8,10);
  \fill (1,-1) circle (0.14); \fill (8,10) circle (0.14);
  \node[anchor=west,font=\scriptsize] at (-0.8,-1) {$x_{I_1}$};
  \node[anchor=west,font=\scriptsize] at (8,10) {$x_{I_2}$};
  %% dimensions
  \draw[<->,gray!80!black,thin] (2.90,-2.90) -- (10.10,7.90);
  \node[gray!80!black,anchor=west,font=\scriptsize] at (6.80,2.50) {length $I_2-I_1$};
  \draw[<->,gray!80!black,thin] (-2.78,0.98) -- (1.22,-3.02);
  \node[gray!80!black,anchor=north east,font=\scriptsize] at (-0.58,-0.92) {width $2Rs$};
  \draw[->,gray] (-7,9) -- (-5.40,11.40) node[above right=-2pt,font=\scriptsize] {$\Dir$};
  \draw[->,gray] (-7,9) -- (-5.4,7.4) node[right,font=\scriptsize] {$\ead$};
\end{tikzpicture}
\caption{\small A block $\{I_1\le\level\le I_2,\ \abs{\rcoord-\rho}\le Rs\}$ in the original coordinates, drawn to scale with $\Dir=(0.4,0.6)$, $I_2-I_1=18$, and $Rs=2$, over the lattice $\Z^2$ (light grid). The dots are the lattice points of the block. The sides of the parallelogram are along $\Dir$ and $\ead$, its length $I_2-I_1$ is measured in levels along $\Dir$, and its width $2Rs$ in units of $\ead$. The thick path is a crossing: an admissible path from level $I_1$ to level $I_2$ all of whose vertices lie in the block.}
\label{fig:block}
\end{figure}
A \emph{block} is a set of the form $\{z\in\Z^2:\ I_1\le\level(z)\le I_2,\ \abs{\rcoord(z)-\rho}\le Rs\}$, a parallelogram in the original coordinates with two sides along $\Dir$ and two along $\ead$. The  \emph{length} of the block  is $I_2-I_1$ and its \emph{width}  $2Rs$. A \emph{crossing} of it is an admissible path $x_{I_1:I_2}$ all of whose vertices lie in it. See Figure~\ref{fig:block}.
 The main result of this section, Proposition~\ref{prop:rarefast}, says that a crossing of a block of width of order $s$ and of any length $n=I_2-I_1$  cannot exceed either extreme support functional  by $\e s$, that is, that such a path $x_{m:m+n}$  
cannot satisfy $\pawgt(x_{m:m+n}) > f_\sig(x_{m+n}-x_m)+\e s$. 
Its proof turns a positive probability of such a crossing into paths that exceed $f_\sig$ by a positive multiple of their length, contradicting \eqref{eq:defect-lower}.  

We  first construct  near-geodesics that  can be attached to prescribed endpoints and stay close to a straight-line segment.

\begin{lemma}\label{lem:travel}
Fix $D\ge1$, $\theta\in(0,\theta_0)$, and $\eta,\delta>0$. There exist  $s_1<\infty$ and a function $p$ with $p(s)\to1$ as $s\to\infty$, depending only on $D,\theta,\eta,\delta,\mu,\Dir$, with the following property. For $s\ge s_1$ and every deterministic pair $u\le v$ in $\Z^2$ with $s\le\level(v-u)\le Ds$ and $v-u$ $\theta$-interior, with probability at least $p(s)$ there exists an admissible path $x_{u:v}$ with
\[
\pawgt(x_{u:v})\ge g(v-u)-\eta s
\]
all of whose vertices are within $\abs{\cdot}_1$-distance $\delta s$ of the line segment $[u,v]$.
\end{lemma}

\gist{This proof uses the shape theorem and Lipschitz $g$. The distance from $u$ to $v$ is decomposed into blocks that are microscopically large to take advantage of the shape theorem, but macroscopically small enough to keep the concatenation of the block geodesics within $\delta s$ of $[u,v]$. Lip $g$ ensures that $g$ is close to summable over the segments.}

\begin{figure}[htb]
\centering
\begin{tikzpicture}[scale=0.62,>=Stealth]
  %% tube of radius delta s around [u,v]
  \fill[gray!14] (-1.45,1.1) -- (4.55,9.1) -- (7.45,6.9) -- (1.45,-1.1) -- cycle;
  \draw[gray!60] (-1.45,1.1) -- (4.55,9.1);
  \draw[gray!60] (1.45,-1.1) -- (7.45,6.9);
  \draw[<->,gray] (3,4) -- (4.46,2.91) node[midway,below=5pt, right=-6pt,sloped] {\footnotesize $\delta s$};
  %% segment [u,v]
  \draw[dashed] (0,0) -- (6,8);
  %% lattice
  \foreach \x in {0,...,6} \foreach \y in {0,...,8} \fill[gray!50] (\x,\y) circle (0.05);
  %% rectangles {u_j <= z <= u_{j+1}}
  \foreach \a/\b/\c/\d in {0/0/1/2, 1/2/3/4, 3/4/4/6, 4/6/6/8} \fill[gray!40,opacity=0.6] (\a,\b) rectangle (\c,\d);
  %% concatenated near-geodesic
  \draw[very thick] (0,0) -- (0,1) -- (1,1) -- (1,2) -- (2,2) -- (2,3) -- (3,3) -- (3,4) -- (4,4) -- (4,5) -- (4,6) -- (5,6) -- (5,7) -- (6,7) -- (6,8);
  \foreach \pt in {(0,0),(1,2),(3,4),(4,6),(6,8)} \fill \pt circle (0.13);
  \node[below right] at (-1,0) {$u=u_0$};
  \node[left] at (1,2) {$u_1$};
  \node[left] at (3,4) {$u_2$};
  \node[left] at (4,6) {$u_3$};
  \node[above left] at (7,8) {$u_4=v$};
\end{tikzpicture}
\caption{\small Lemma~\ref{lem:travel} with $B=4$. The points $u_j=u+\lfloor jz/B\rfloor$ lie within $\abs\cdot_1$-distance $2$ of the segment $[u,v]$ (dashed). The darker boxes are the rectangles $\{u_j\le z\le u_{j+1}\}$: every admissible path from $u_j$ to $u_{j+1}$, in particular a geodesic, stays inside the box, whose $\abs\cdot_1$-diameter is $\abs{u_{j+1}-u_j}_1\le Ds/B+2\le\delta s/4+2$. Hence the concatenation of the four geodesics (thick path) stays in the tube of radius $\delta s$ around $[u,v]$ (lighter shading). Each piece has weight at least $g(u_{j+1}-u_j)-\eta s/(2B)$ by \ref{S2}, and $\sum_jg(u_{j+1}-u_j)\ge g(v-u)-2\tsp\Lip_{\theta/2}B$ by the Lipschitz property of $g$.}
\label{fig:travel}
\end{figure}
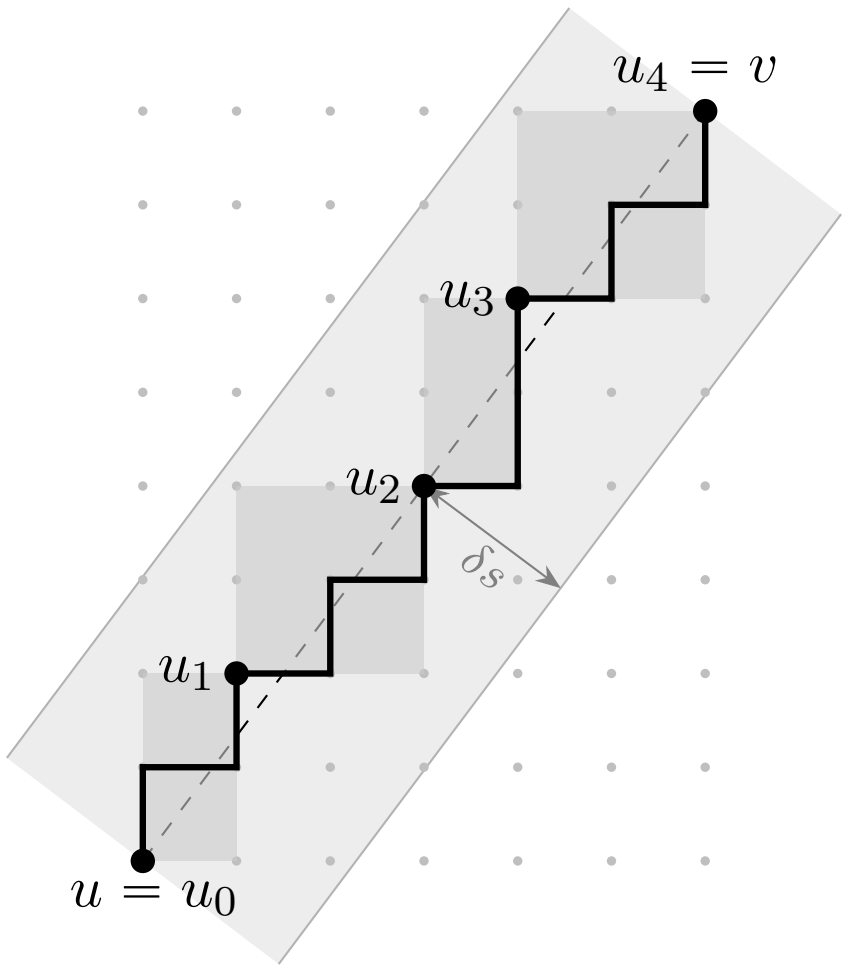
\begin{proof}
Write $z=v-u$ and let $B$ be a fixed integer with $B\ge4D/\delta$ (see Figure~\ref{fig:travel}). Let $u_j=u+\lfloor jz/B\rfloor$ for $0\le j\le B$, the floor taken coordinatewise. Then $u_0=u$, $u_B=v$, and $u_j\le u_{j+1}$ (as $jz/B$ is coordinatewise nondecreasing in $j$ and the floor is monotone). Since each coordinate of $\lfloor jz/B\rfloor-jz/B$ lies in $(-1,0]$, $u_j$ is within $\abs\cdot_1$-distance $2$ of the point $u+jz/B$ of $[u,v]$, and consequently $\abs{u_{j+1}-u_j-z/B}_1\le2$, as two vectors with coordinates in $(-1,0]$ differ by at most $1$ in each coordinate.

Since $z$ is $\theta$-interior, $\abs{\rcoord(z/B)}\le(\theta_0-\theta)\level(z)/B$, and by \eqref{eq:lip-r}
\begin{align*}
\abs{\rcoord(u_{j+1}-u_j)}&\le\abs{\rcoord(z/B)}+\abs{\rcoord(u_{j+1}-u_j-z/B)}\\
&\le(\theta_0-\theta)\level(z)/B+\abs{u_{j+1}-u_j-z/B}_1\le(\theta_0-\theta)\level(z)/B+2
\end{align*}
while $\level(u_{j+1}-u_j)\ge\level(z)/B-2$. Hence $u_{j+1}-u_j$ is $\theta/2$-interior as soon as $(\theta/2)\level(z)/B\ge2+2(\theta_0-\theta/2)$, which holds for $s\ge s_1'=s_1'(B,\theta,\theta_0)$ because $\level(z)\ge s$. The vector $z/B$ itself is $\theta$-interior by linearity of $\level$ and $\rcoord$.

%$\theta$-interior assumption is used in this step to control the Lip constant of $g$. No other use in this proof.

With $\Lip_{\theta/2}$ the Lipschitz constant from \ref{S4}, homogeneity gives $\sum_{j=0}^{B-1}g(z/B)=g(z)$, so
\be\label{t.aux13}
\sum_{j=0}^{B-1}g(u_{j+1}-u_j)\ \ge\ \sum_{j=0}^{B-1}\bigl[g(z/B)-\Lip_{\theta/2}\abs{u_{j+1}-u_j-z/B}_1\bigr]\ \ge\ g(z)-2\tsp\Lip_{\theta/2}B .
\ee

Each pair $(u_j,u_{j+1})$ is deterministic with $\abs{u_{j+1}-u_j}_1=\level(u_{j+1}-u_j)\in[s/B-2,Ds/B+2]$. Apply \ref{S2} with $\eta'=\eta/(4D)$: with probability at least $1-\e(s)$, where $\e(s)\to0$ as $s\to\infty$ and $\e(s)$ depends only on $D,B,\eta,\mu,\Dir$, a geodesic $x^{(j)}_{u_j:u_{j+1}}$ has
\be\label{t.aux15}
\Lpp_{u_j,u_{j+1}}=\pawgt(x^{(j)}_{u_j:u_{j+1}})\ \ge\ g(u_{j+1}-u_j)-\frac{\eta}{4D}\Bigl(\frac{Ds}{B}+2\Bigr)\ \ge\ g(u_{j+1}-u_j)-\frac{\eta s}{2B}
\ee
for $s\ge2B/D$. The last inequality follows from $\frac{\eta}{4D}(\frac{Ds}{B}+2)=\frac{\eta s}{4B}+\frac{\eta}{2D}$ and $\frac{\eta}{2D}\le\frac{\eta s}{4B}$ for $s\ge2B/D$.  Every admissible path from $u_j$ to $u_{j+1}$ stays in the rectangle $\{u_j\le z'\le u_{j+1}\}$, whose $\abs\cdot_1$-diameter is $\abs{u_{j+1}-u_j}_1\le Ds/B+2\le\delta s/4+2$. So with $s\ge16/(3\delta)$, every vertex of $x^{(j)}_{u_j:u_{j+1}}$ is within $\delta s/4+2+2\le\delta s$ of the segment $[u,v]$.

On the intersection of the $B$ events above, which has probability at least $p(s)=1-B\e(s)$, let $x_{u:v}$ be the path with $x_{u_j:u_{j+1}}=x^{(j)}_{u_j:u_{j+1}}$ for every $j$. Then
\[\pawgt(x_{u:v})=\sum_{j=0}^{B-1}\pawgt(x^{(j)}_{u_j:u_{j+1}})
\overset{\eqref{t.aux13},\eqref{t.aux15}}\ge g(z)-2\tsp\Lip_{\theta/2}B-\frac{\eta s}2\ge g(z)-\eta s\quad\text{for }s\ge4\tsp\Lip_{\theta/2}B/\eta .\]
Altogether, the lemma holds with $s_1$ the largest of $s_1'(B,\theta,\theta_0)$, $2B/D$, $16/(3\delta)$, and $4\tsp\Lip_{\theta/2}B/\eta$, which depends only on $D,\theta,\eta,\delta,\mu,\Dir$ (recall that $B$ is determined by $D$ and $\delta$).
\end{proof}

 Call a   set of the form
\[\cJ=\{z\in\Z^2:\level(z)=I,\ \rcoord(z)\in[\rho_1,\rho_2]\}\]
a \emph{level interval} with  $\rcoord$-width $\rho_2-\rho_1$. We extend  Lemma \ref{lem:travel} from a fixed pair of points to all pairs of points from one level interval to another, positioned relative to each other in a $\theta$-interior manner. Their mutual  distance and $\rcoord$-widths grow. The combinatorial explosion of pairs is tamed by connecting the levels to macroscopically nearby \textit{nets}  of bounded numbers of points. Lemma \ref{lem:travel} is then applied  to the constant number of connections between the nets. 

It is helpful to observe that 
the lattice points reachable from $u\in\Z^2$ in exactly $m$ steps are $u+(j,m-j)$ with transverse coordinates $\rcoord(u)-\Dir_1m+j$, for  $ 0\le j\le m$. 
These are exactly the lattice points on level $\level(u)+m$  with $\rcoord\in[\rcoord(u)-\Dir_1m,\rcoord(u)+(1-\Dir_1)m]$, an interval of length $m$. 

\begin{lemma}\label{lem:intervals}
Fix $D\ge1$, $\theta\in(0,\theta_0)$, and $\eta,\delta>0$. There exist $s_1<\infty$ and $p(s)\to1$ as $s\to\infty$, depending only on $D,\theta,\eta,\delta,\mu,\Dir$, with the following property. Let $\cJ_1,\cJ_2$ be deterministic level intervals on levels $I_1<I_2$ with $s\le I_2-I_1\le Ds$ and $\rcoord$-widths at most $Ds$, such that $v-u$ is $\theta$-interior for all $u\in \cJ_1$ and $v\in \cJ_2$. For $s\ge s_1$, with probability at least $p(s)$, every pair $u\in \cJ_1$ and $v\in \cJ_2$ is joined by an admissible path $x_{u:v}$ of weight at least $g(v-u)-\eta s$, all of whose vertices lie within $\abs{\cdot}_1$-distance $\delta s$ of $[u,v]$.
\end{lemma}

We prove   Lemma \ref{lem:intervals} here under the assumption that the weights are bounded:
\begin{align}\label{ass:bounded}
\text{$\mu([-M,M])=1$ for some $M<\infty$.}
\end{align}
Assumption \eqref{ass:bounded} enters the proof below in one place, explained after the proof  in  Remark~\ref{rk:bounded1}. The proof under assumption \ref{A2} is   in Appendix~\ref{sub:attachments}.

Under \eqref{ass:bounded}, sums of weights along deterministic paths obey Hoeffding's inequality \cite[Theorem~2]{Hoe-63}:   if $S_m$ is a sum of $m$ i.i.d.\ copies of $\w_0$, then for every $t>0$
\begin{equation}\label{eq:hoeffding}
\P\bigl\{S_m<m\E\w_0-t\bigr\}\le\exp\Bigl\{-\frac{t^2}{2mM^2}\Bigr\}. 
\end{equation}

\gist{Put two nets $\cN_1$ and $\cN_2$ of constant  numbers of vertices at distance $m\approx\kappa s$ of $\cJ_1$ and $\cJ_2$. Control passage times from $\cJ_i$ to $\cN_i$ with Hoeffding: polynomially many exp bounds goes to zero. Apply Lemma \ref{lem:travel} to the constant number of connections from $\cN_1$ to $\cN_2$.}

\begin{figure}[htb]
\centering
\begin{tikzpicture}[scale=0.9,>=Stealth]
  %% axes: horizontal = r, vertical = level (not to scale)
  \draw[->,gray] (-2.2,-0.5) -- (-2.2,8.4) node[above] {\footnotesize $\level$};
  \draw[->,gray] (-2.2,-0.5) -- (9.2,-0.5) node[right] {\footnotesize $\rcoord$};
  %% level intervals J_1, J_2
  \draw[very thick] (0,0) -- (5,0) node[above=7pt,right=-70pt] {\footnotesize $\cJ_1=\{\level=I_1,\ \rcoord\in[\rho_1,\rho_2]\}$};
  \draw[very thick] (2,8) -- (7,8) node[midway, left=10pt, above=-1pt] {\footnotesize $\cJ_2$};
  %% u, its reachable cone, u-hat
  \fill (2,0) circle (0.08) node[below=1pt] {\footnotesize $u$};
  \fill[gray!25] (2,0) -- (1.7,1.3) -- (2.9,1.3) -- cycle;
  \draw[gray] (2,0) -- (1.7,1.3);  \draw[gray] (2,0) -- (2.9,1.3);
  \draw[thick] (2,0) -- (2.0,1.3);
  \node[right] (alab) at (0,0.8) {\footnotesize $\alpha_{u:\hat u}$};
   \draw[->,gray,shorten >=1pt] (alab) -- (2.0,0.65);
  %% v, its backward cone, v-hat
  \fill (5,8) circle (0.08) node[above=1pt] {\footnotesize $v$};
  \fill[gray!25] (5,8) -- (4.1,6.7) -- (5.3,6.7) -- cycle;
  \draw[gray] (5,8) -- (4.1,6.7);  \draw[gray] (5,8) -- (5.3,6.7);
  \draw[thick] (4.9,6.7) -- (5,8);
  \node[left] at (4.1,7.4) {\footnotesize $\alpha_{\hat v:v}$};
  \node[left] (blab) at (4.1,7.4) {\footnotesize $\alpha_{\hat v:v}$};
   \draw[->,gray,shorten >=1pt] (blab) -- (5,7.3);
  %% connection: tube and path
  \fill[gray!12] (1.5,1.5) -- (4.4,6.9) -- (5.4,6.5) -- (2.5,1.1) -- cycle;
  \draw[dashed] (2.0,1.3) -- (4.9,6.7);
   \draw[thick] (2.0,1.3) .. controls (2.42,2.82) and (3.76,3.84) .. (3.6,4.27) .. controls (3.72,5.25) and (4.73,5.75) .. (4.9,6.7);
  %% level lines
  \foreach \y/\lab in {0/$I_1$, 1.3/$I_1+m$, 6.7/$I_2-m$, 8/$I_2$} {\draw[gray!50] (-1.8,\y) -- (9,\y); \node[right,gray] at (9,\y) {\footnotesize \lab};}
  %% nets
  \foreach \x in {-0.8,-0.1,0.6,1.3,2.0,2.7,3.4,4.1,4.8,5.5,6.2} \fill (\x,1.3) circle (0.06);
  \node[right] at (6.4,1.3) {\footnotesize $\cN_1$};
  \foreach \x in {0.7,1.4,2.1,2.8,3.5,4.2,4.9,5.6,6.3,7.0,7.7} \fill (\x,6.7) circle (0.06);
  \node[right] at (7.9,6.7) {\footnotesize $\cN_2$};
  \fill (2.0,1.3) circle (0.08) node[above left=-2pt] {\footnotesize $\hat u$};
  \fill (4.9,6.7) circle (0.08) node[below right=-2pt] {\footnotesize $\hat v$};
  \node[right] at (4.3,3.6) {\footnotesize Lemma~\ref{lem:travel}: $\pawgt\ge g(\hat v-\hat u)-\eta s/4$};
\end{tikzpicture}
\caption{\small Lemma~\ref{lem:intervals}, drawn in the coordinates $(\rcoord,\level)$ and not to scale. The nets $\cN_1$ and $\cN_2$ are deterministic sets of lattice points on the levels $I_1+m$ and $I_2-m$ with transverse gaps between $\theta_0m/4$ and $\theta_0m/2$. From a point $u\in\cJ_1$ the lattice points reachable in $m$ steps form a level interval of $\rcoord$-width $m$ (the shaded cone), which contains a net point $\hat u$. Symmetrically, $v\in\cJ_2$ is reachable in $m$ steps from a net point $\hat v$. The attachments $\alpha_{u:\hat u}$ and $\alpha_{\hat v:v}$ are deterministic paths of length $m$, whose weights are controlled by Hoeffding's inequality, and the boundedly many pairs $(\hat u,\hat v)$ are joined by the near-geodesics of Lemma~\ref{lem:travel}, which stay in the tube of radius $\delta s/2$ around $[\hat u,\hat v]$ (shaded).}
\label{fig:intervals}
\end{figure}
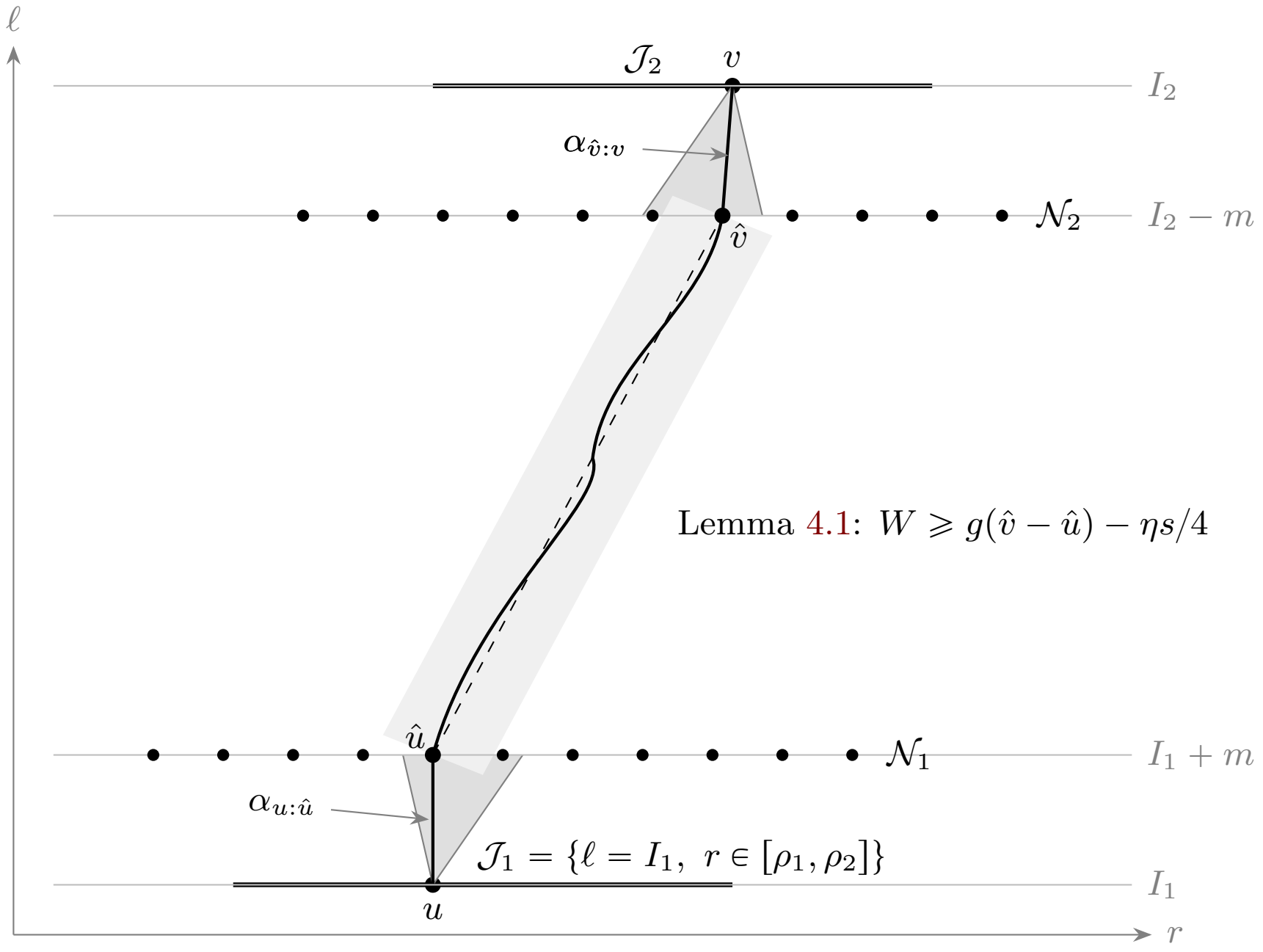
\begin{proof}[Proof of Lemma \ref{lem:intervals} under Assumption \eqref{ass:bounded}]
Let $\Lip_{\theta/2}$ be the Lipschitz constant from \ref{S4}  and fix %Note: the $+1$ in the denomenator is only to avoid dividing by zero when $\Lip_{\theta/2}=\E\w_0=0$. The estimates use just $\kappa(\Lip_{\theta/2}+\abs{\E\w_0})\le\eta/16$.
\[
\kappa=\min\Bigl\{\frac{\theta}{16},\ \frac{\delta}{4},\ \frac{\eta}{16(\Lip_{\theta/2}+\abs{\E\w_0}+1)}\Bigr\}\quad\text{and}\quad
m=\lceil\kappa s\rceil .
\]
The level intervals are (see Figure~\ref{fig:intervals})   \[ \cJ_1=\{y\in\Z^2: \level(y)=I_1,\ \rcoord(y)\in[\rho_1,\rho_2]\}  \ \text{ and } \ \cJ_2=\{y\in\Z^2: \level(y)=I_2,\ \rcoord(y)\in[\rho_1',\rho_2']\}. \]

\medskip
\emph{Requirements on $s_1$.} The proof imposes finitely many lower bounds on $s$, each depending only on $D,\theta,\eta,\delta,\mu,\Dir$, and $s_1$ is the largest of them. The first is $\theta_0\kappa s_1\ge4$, so that $\theta_0m/4\ge1$. The second is $s_1\ge16$, so that $2m\le2\kappa s+2\le s/8+2\le s/4$ (recall $\kappa\le\theta/16<1/16$). In particular $I_2-m>I_1+m$, since $I_2-I_1\ge s>2m$. The third requirement is $s_1\ge16/\theta$, so that together with $4m\le4\kappa s+4\le\theta s/4+4$ we have $\theta s/2\ge4m$ and hence $(\theta/2)s\ge2m(1+\theta_0-\theta/2)$ (as $\theta_0\le\tfrac12$). The fourth requirement is $s_1\ge16(\Lip_{\theta/2}+\abs{\E\w_0})/\eta$, and the fifth is $s_1\ge4/\delta$.

\medskip
On level $I_1+m$ choose a set $\cN_1$ of lattice points whose transverse coordinates increase with consecutive gaps in $[\theta_0m/4,\theta_0m/2]$, exactly one of them at most $\rho_1-\Dir_1m+\theta_0m/2$, and exactly one at least $\rho_2+(1-\Dir_1)m-\theta_0m/2$. This is possible because the lattice points on a level have transverse spacing $1\le\theta_0m/4$ (by the first requirement on $s_1$). Set
\[\cR_1=[\rho_1-\Dir_1m+\theta_0m/2,\ \rho_2+(1-\Dir_1)m-\theta_0m/2].\] By construction $\cR_1$ is contained in the interval spanned by the transverse coordinates of $\cN_1$. We claim that every interval $[\gamma,\gamma+m]$ contained in $[\rho_1-\Dir_1m,\rho_2+(1-\Dir_1)m]$ contains the transverse coordinate of at least one point of $\cN_1$. Indeed, $\cR_1$ is obtained from $[\rho_1-\Dir_1m,\rho_2+(1-\Dir_1)m]$ by removing a piece of length $\theta_0m/2$ at each end, so $[\gamma,\gamma+m]\cap\cR_1$ is an interval of length at least $m-2\cdot\theta_0m/2=(1-\theta_0)m\ge\theta_0m/2$, the maximal net  gap. 

The lattice points reachable from $u\in\cJ_1$ in $m$ steps are those on level $I_1+m$ with transverse coordinate in $[\rcoord(u)-\Dir_1m,\rcoord(u)+(1-\Dir_1)m]$, an interval of length $m$ contained in $[\rho_1-\Dir_1m,\rho_2+(1-\Dir_1)m]$. By the claim proved in the preceding paragraph, it contains a point $\hat u\in \cN_1$. Thus $\hat u$ is reachable from $u$ in $m$ steps, $\abs{\hat u-u}_1=m$, and, by \eqref{eq:lip-r},
\begin{align}\label{aux397}
\abs{\rcoord(\hat u-u)}\le m.
\end{align}

Similarly, on level $I_2-m$ choose a set $\cN_2$ of lattice points with consecutive transverse gaps in $[\theta_0m/4,\theta_0m/2]$ spanning at least the range $[\rho_1'-(1-\Dir_1)m+\theta_0m/2,\ \rho_2'+\Dir_1m-\theta_0m/2]$. The points from which $v\in\cJ_2$ is reachable in $m$ steps are those on level $I_2-m$ with transverse coordinate in $[\rcoord(v)-(1-\Dir_1)m,\rcoord(v)+\Dir_1m]$ (reverse the steps), and the same argument gives $\hat v\in \cN_2$ from which $v$ is reachable in $m$ steps, with $\abs{v-\hat v}_1=m$ and
\begin{align}\label{aux403}
\abs{\rcoord(v-\hat v)}\le m.
\end{align}

 Both nets have cardinality bounded by a constant  $\cnst(D,\kappa,\theta_0)$. To see this, consider that all but the two extreme points of $\cN_1$ have coordinates in $\cR_1$, an interval of length $\rho_2-\rho_1+(1-\theta_0)m\le Ds+m-\theta_0m$, and consecutive gaps are at least $\theta_0m/4$. From these,   
\[
\abs{\cN_1}\le3+\tfrac{Ds+m-\theta_0m}{\theta_0m/4}<4\tfrac{Ds+m}{\theta_0m} \le 4\tfrac{Ds+\kappa s+1}{\theta_0\kappa s} \le\cnst(D,\kappa,\theta_0). \]
Fix for each $u$ one admissible path $\alpha_{u:\hat u}$ from $u$ to $\hat u$ (say the one taking its $e_1$ steps first) and for each $v$ one path $\alpha_{\hat v:v}$ from $\hat v$ to $v$. These are deterministic.

For $u\in\cJ_1$, $v\in\cJ_2$ write $z=v-u$. By assumption,  $z$ is $\theta$-interior. Let $\tilde z=\hat v-\hat u=z-(\hat u-u)-(v-\hat v)$. Then, by $\level(z)\in[s,Ds]$ and the second requirement on $s_1$,
\[\level(\tilde z)=\level(z)-2m\ge s-s/4\ge s/2>0\quad\text{and}\quad\level(\tilde z)\le\level(z)\le Ds.\]
Furthermore, using \eqref{aux397} and \eqref{aux403}, we have
\[\abs{\rcoord(\tilde z)}\le(\theta_0-\theta)\level(z)+2m=(\theta_0-\theta)(\level(\tilde z)+2m)+2m.\]
Since $\level(z)\ge s$, the third requirement on $s_1$ implies $(\theta/2)\level(z)\ge2m(1+\theta_0-\theta/2)$, which rearranges to $(\theta_0-\theta)\level(z)+2m\le(\theta_0-\theta/2)(\level(z)-2m)$. This gives $\abs{\rcoord(\tilde z)}\le(\theta_0-\theta/2)\level(\tilde z)$ for $s\ge s_1$, proving that $\tilde z$ is $\theta/2$-interior.

The weight of $\alpha_{u:\hat u}$ is a sum of $m$ i.i.d.\ weights. By Hoeffding's inequality  \eqref{eq:hoeffding} with $t=\eta s/4$ and $m\le s$ (by the second requirement on $s_1$),
%$m\le\kappa s+1\le s/16+1\le s$ since $s\ge16$ and $\kappa<1/16$
\begin{align}\label{attachments}
\P\bigl\{\pawgt(\alpha_{u:\hat u})<m\E\w_0-\eta s/4\bigr\}\le\exp\Bigl\{-\frac{\eta^2s^2}{32mM^2}\Bigr\}\le\exp\Bigl\{-\frac{\eta^2s}{32M^2}\Bigr\}.
\end{align}
The same bound holds for $\alpha_{\hat v:v}$.
There are at most $\abs{\cJ_1}+\abs{\cJ_2}\le2Ds+2$ attachments, so with probability at least $1-(2Ds+2)e^{-\eta^2s/(32M^2)}$, which tends to one as $s\to\infty$, every attachment has weight at least $m\E\w_0-\eta s/4$.

Apply Lemma~\ref{lem:travel} at scale $s/2$ with $D'=2D$, $\theta'=\theta/2$, $\eta'=\eta/2$, and $\delta'=\delta$ to every pair $(\hat u,\hat v)\in \cN_1\times \cN_2$ (the displacement $\hat v-\hat u$ is $\theta/2$-interior with level in $[s/2,2D(s/2)]$). With probability at least $1-\cnst(D,\kappa,\theta_0)^2(1-p(s/2))\to1$, every such pair is joined by a path of weight at least $g(\hat v-\hat u)-\eta s/4$ within $\delta s/2$ of $[\hat u,\hat v]$.

On the intersection of these events, given $u,v$ as in the claim of the lemma, let $x_{u:v}$ be the path with $x_{u:\hat u}=\alpha_{u:\hat u}$, $x_{\hat u:\hat v}$ the net connection from $\hat u$ to $\hat v$, and $x_{\hat v:v}=\alpha_{\hat v:v}$. Since $z$ and $\tilde z$ are both $\theta/2$-interior and $\abs{z-\tilde z}_1\le2m$, the Lipschitz property gives $g(\tilde z)\ge g(z)-2\tsp\Lip_{\theta/2}\tsp m$. Hence
\[
\pawgt(x_{u:v})\ \ge\ 2\bigl(m\E\w_0-\tfrac{\eta s}4\bigr)+g(\tilde z)-\tfrac{\eta s}4\ \ge\ g(z)-2m\bigl(\Lip_{\theta/2}+\abs{\E\w_0}\bigr)-\tfrac{3\eta s}{4}\ \ge\ g(z)-\eta s ,
\]
because
\[ 2m(\Lip_{\theta/2}+\abs{\E\w_0})\le2(\kappa s+1)(\Lip_{\theta/2}+\abs{\E\w_0})\le\eta s/8+2(\Lip_{\theta/2}+\abs{\E\w_0})\le\eta s/4 \] by the fourth requirement on $s_1$. As for the location: the vertices of $\alpha_{u:\hat u}$ and $\alpha_{\hat v:v}$ are within $m$ of $u$ and $v$. Every point $\hat u+\tau(\hat v-\hat u)$ of $[\hat u,\hat v]$ is within $(1-\tau)\abs{\hat u-u}_1+\tau\abs{\hat v-v}_1=m$ of the point $u+\tau(v-u)$ of $[u,v]$. So all vertices of $x_{u:v}$ are within $\delta s/2+m\le\delta s$ of $[u,v]$, as $m\le\delta s/4+1\le\delta s/2$ by the fifth requirement on $s_1$.
\end{proof}

\begin{remark}[Bounded vs.\ unbounded weights]\label{rk:bounded1}  
The boundedness of the weights enters Lemma~\ref{lem:intervals} only in the bound \eqref{attachments}. The attachment paths $\alpha_{u:\hat u}$, $\alpha_{\hat v:v}$ are deterministic paths of length $m\approx\kappa s$.   Hoeffding's inequality \eqref{eq:hoeffding} makes the probability that one of them has weight  $S_m<m\E\w_0-\eta s/4$ exponentially small in $s$, which survives the union bound over the $O(s)$ attachments.

For unbounded weights this line of reasoning fails without a high enough moment bound.  Assume   a finite  mean $\E\abs{\w_0}<\infty$ and left tail 
$\P(\w_0<-t)\asymp t^{-p'''}$.  
Then an inclusion-exclusion argument gives, for large $s>0$ and still with $m\approx\kappa s$,  
\be\label{t.aux.bdd78} 
\P\{ S_m<m\E\w_0-\eta s/4\} \ge \tfrac m4\,\P(\w_0<-\eta s/2)  \asymp s^{1-p'''} 
\ee
which is not $o(1/s)$  unless $p'''>2$. 

To deduce \eqref{t.aux.bdd78}, let $q=\P(\w_0<-\eta s/2)$, and for $1\le i\le m$ let $A_i=\{\w_i<-\eta s/2\}$ and $B_i=\{\sum_{j\ne i}\w_j\le(m-1)\E\w_0+\eta s/8\}$. On $A_i\cap B_i$, 
\[ S_m<m\E\w_0-\eta s/2+\eta s/8+\abs{\E\w_0}\le m\E\w_0-\eta s/4 \quad 
\text{once } \abs{\E\w_0}\le\eta s/8. \]
The events $A_i$ and $B_i$ are independent, and $\P(B_i)\ge\tfrac12$ for large $s$.
%by the law of large numbers, since $m\approx\kappa s$. 
By inclusion-exclusion, 
\begin{align*}
 \P\{ S_m<m\E\w_0-\eta s/4\} &\ge   \P\bigl({\textstyle\bigcup}_i(A_i\cap B_i)\bigr)\ge\sum_i\P(A_i)\P(B_i)-\sum_{i<j}\P(A_i\cap A_j)\\[-7pt] 
 &\ge\tfrac{mq}2-\tfrac{m^2q^2}2\ge\tfrac{mq}4
\end{align*} 
when $mq\le\tfrac12$, which holds for large $s$ because $q=o(1/s)$ when $\E\w_0^-<\infty$.

The  fix is to abandon deterministic paths.  Attachments are constructed to avoid \emph{very bad} 
vertices with weights $\w_x<-t'$.  Under the  assumption $\E[(\w_0^-)^{p'}]<\infty$ with  $p'>1$,  very bad vertices 
are isolated at the relevant scale and the costs of the detours can be controlled.   The one-sided Bernstein inequality for the truncated weights $\w_x\vee(-t')$ replaces Hoeffding. The details   are in Appendix~\ref{sub:attachments}. Figure~\ref{fig:attach-general} there shows an attachment that avoids the very bad vertices, and Figure~\ref{fig:intervals-general} shows how Figure~\ref{fig:intervals} changes. 
\end{remark}

The next proposition is the main result of this section. It is proved by an oriented percolation argument  on a coarse-grained lattice. 
%The percolation result, Lemma~\ref{lem:oriented} in Section~\ref{sub:oriented}, requires only a uniform bound on the failure probability of a site and bounded-range dependence.  
The upper tail of $\mu$ enters the proof through the shape theorem, through Lemma~\ref{lem:travel} and in step 6.

\gist{The proposition is an oriented percolation argument that ``amplifies'' a hypothetical crossing that is localized in a thin tube (width $Rs$, length $n\gg s$) and beats the supporting line by $\e s$, i.e.\ $\pawgt(x_{0:n})>f_\sig(x_n-x_0)+\e s$.  Lemma~\ref{lem:intervals} is used for gluing these crossings, along an oriented path of good sites in the coarse-grained lattice, into a single path. The gluing is cheap, at cost below $\e s/2$ per layer, because the connectors drift to the side where $f_\sig$ touches $g$: The coarse lattice is drawn so that it can only move in the direction $\sig$, the blocks are narrow enough that every connector shifts in that direction, and the connectors are long enough that the shift is a small angle, where tangency makes the gap $f_\sig-g$ negligible.}

%For the percolation argument, a site of the coarse lattice corresponds to a block in the original lattice that has $T$ disjoint blocks in it, each of which contains a crossing with $\pawgt>f_\sig+\e s$ (equivalently, up to $o(s)$, $>g+\e s$) with probability at least $p_0$, by the contradiction hypothesis. The site is good if one of its $T$ blocks has such a crossing and the connectors to the two successor sites exist. Failure has probability at most $(1-p_0)^T+2(1-p(s))$, which can be made smaller than the threshold of oriented percolation. With positive probability there is then an infinite oriented path of good sites, and concatenating the crossings and connectors along it gives an admissible semi-infinite path such that, after $t$ layers, additivity gives $\pawgt>f_\sig+t\e s-t\e s/2=f_\sig+t\e s/2\ge g+t\e s/2$. [The loss of $t\e s/2$ is from the connectors.] On the other hand, $\pawgt\le\Lpp\le g+o(t)$ by the shape theorem, with $s$ and $n$ fixed. This is a contradiction for $t$ large.

\begin{proposition}\label{prop:rarefast}
Fix $R,\e>0$ and $\sigg\in\{-,+\}$. Then
\begin{equation}\label{eq:rarefast}
\begin{aligned}
\lim_{s\to\infty}\ \sup_{n\tsp\in\tsp\N}\ &\P\bigl\{\text{there exists an admissible path $x_{0:n}$ such that  } \\[-4pt]
&\qquad\qquad
\pawgt(x_{0:n})>f_\sig(x_n-x_0)+\e s \text{ and }\abs{\rcoord(x_i)}\le Rs\text{ for }0\le i\le n\bigr\}=0 . 
\end{aligned}
\end{equation}
\end{proposition}

%The moment assumption on $\w_0^+$ is used in the proof of the proposition only through the hypotheses of the shape theorem: in Lemma~\ref{lem:travel}, hence in Lemma~\ref{lem:intervals} (step 4), and in step 6, where \eqref{eq:defect-lower} and \ref{S1} rule out an infinite path that beats $f_\sig$ by a positive multiple of its length. This covers $n=1$ as well. There the event is that a vertex in a window of about $2Rs$ lattice points has weight above $f_\sig(x_1-x_0)+\e s$, and if this had probability at least $p_0$ for arbitrarily large $s$, the percolation argument would produce an infinite path that collects such a vertex in every layer of $n+d\le D_0s+2$ levels, so that it beats $f_\sig$ by about $\e/(2D_0)$ per level, contradicting the shape theorem. Consistently, the hypotheses of the shape theorem give the case of bounded $n$ directly: its positive-tail condition $\int_0^\infty\P(\w_0>t)^{1/2}\,dt<\infty$ implies $\P(\w_0>t)=o(t^{-2})$, since $\P(\w_0>t)^{1/2}\,t/2\le\int_{t/2}^t\P(\w_0>r)^{1/2}\,dr\to0$. For $n\le n_0$ fixed, a path $x_{0:n}$ has $\pawgt(x_{0:n})\le n_0\max_i\w_{x_i}$ and $f_\sig(x_n-x_0)\ge-n_0(\abs{\ghatz}+K)$, so the event requires a vertex among the $O(s)$ vertices of the block with weight above $\e s/(2n_0)$ for large $s$, which has probability $O(s)\,\P\{\w_0>\e s/(2n_0)\}\to0$. Uniformity over all $n$ is what needs the proof.

We state some remarks before the proof. 

\begin{remark}[The role of $n$]    For  fixed $n$ the proposition is a direct consequence of the hypotheses of the shape theorem: since $\int_0^\infty\P(\w_0>t)^{1/2}dt<\infty$ implies $\P(\w_0>t)=o(t^{-2})$, a union bound shows that, with probability tending to one, no vertex of an  $n\times 2Rs$ block  has weight above $\e s/(2n)$.  Then no crossing of the block exceeds $f_\sig$ by $\e s$ once $s$ is large.

The proposition is applied in two places below  to blocks where  $n\gg s$:  namely, in the proof of Lemma~\ref{lem:tilt} with $n=\lceil s^3\rceil$, and in the proof of Proposition~\ref{prop:endpoints} where $n$ is the distance between two coarse-grained  levels. 
\end{remark}

\begin{remark}[Translation by real $\rho$] \label{rk:translate}
By translation invariance, Proposition~\ref{prop:rarefast} holds for every lattice translate of the event in \eqref{eq:rarefast}, that is, for the crossings $x_{I:I+n}$ of the block $\{z:\ I\le\level(z)\le I+n,\ \abs{\rcoord(z)-\rcoord(z_0)}\le Rs\}$, where $I$ is any level and $z_0$ any lattice point on level $I$. A block $\{z:\ I\le\level(z)\le I+n,\ \abs{\rcoord(z)-\rho}\le Rs\}$ centered at an arbitrary real $\rho$ need not be of this form, but it is contained in one with a larger width. Indeed, every unit interval contains a point of $\Z-\Dir_1I$, so there is a lattice point $z_0$ on level $I$ with $\abs{\rcoord(z_0)-\rho}\le1$, and for $s\ge1$ the block centered at $\rho$ is then contained in $\{z:\ I\le\level(z)\le I+n,\ \abs{\rcoord(z)-\rcoord(z_0)}\le(R+1)s\}$, which is the translate by $z_0$ of the block in \eqref{eq:rarefast} with $R+1$ in place of $R$. 
So the proposition also holds for crossings $x_{I:I+n}$ of blocks centered at an arbitrary real $\rho$, with $\pawgt(x_{I:I+n})>f_\sig(x_{I+n}-x_I)+\e s$ and with $R+1$ in place of $R$, uniformly in $I$ and $\rho$.
\end{remark}

\begin{figure}[htb]
\centering
\begin{tikzpicture}[scale=0.55,>=Stealth]
  \draw[gray!70] (-1.60,1.60) -- (2.60,-2.60);
  \node[left,gray] at (-1.60,1.60) {\footnotesize $I_{0}$};
  \draw[gray!40,dashed] (-1.60,5.10) -- (5.80,-2.30);
  \node[left,gray] at (-1.60,5.10) {\footnotesize $I_{0}+n$};
  \draw[gray!70] (-1.60,6.40) -- (6.32,-1.52);
  \node[left,gray] at (-1.60,6.40) {\footnotesize $I_{1}$};
  \draw[gray!40,dashed] (-1.60,9.90) -- (7.72,0.58);
  \node[left,gray] at (-1.60,9.90) {\footnotesize $I_{1}+n$};
  \draw[gray!70] (-1.60,11.20) -- (8.24,1.36);
  \node[left,gray] at (-1.60,11.20) {\footnotesize $I_{2}$};
  \draw[gray!40,dashed] (-1.60,14.70) -- (9.64,3.46);
  \node[left,gray] at (-1.60,14.70) {\footnotesize $I_{2}+n$};
  \draw[gray!70] (-1.60,16.00) -- (10.16,4.24);
  \node[left,gray] at (-1.60,16.00) {\footnotesize $I_{3}$};
  \draw[dotted] (0,0) -- (6.00,9.00) node[above] {\footnotesize $\{\rcoord=0\}$};
  \fill[gray!25] (-1.10,1.10) -- (0.30,3.20) -- (0.90,2.60) -- (-0.50,0.50) -- cycle;
  \draw[gray!60,thin] (-1.10,1.10) -- (0.30,3.20) -- (0.90,2.60) -- (-0.50,0.50) -- cycle;
  \fill[pattern=north east lines,pattern color=gray!70] (-0.30,0.30) -- (1.10,2.40) -- (1.70,1.80) -- (0.30,-0.30) -- cycle;
  \draw[gray!60,thin] (-0.30,0.30) -- (1.10,2.40) -- (1.70,1.80) -- (0.30,-0.30) -- cycle;
  \fill[gray!25] (0.50,-0.50) -- (1.90,1.60) -- (2.50,1.00) -- (1.10,-1.10) -- cycle;
  \draw[gray!60,thin] (0.50,-0.50) -- (1.90,1.60) -- (2.50,1.00) -- (1.10,-1.10) -- cycle;
  \fill[pattern=north east lines,pattern color=gray!70] (1.82,2.98) -- (3.22,5.08) -- (3.82,4.48) -- (2.42,2.38) -- cycle;
  \draw[gray!60,thin] (1.82,2.98) -- (3.22,5.08) -- (3.82,4.48) -- (2.42,2.38) -- cycle;
  \fill[gray!25] (2.62,2.18) -- (4.02,4.28) -- (4.62,3.68) -- (3.22,1.58) -- cycle;
  \draw[gray!60,thin] (2.62,2.18) -- (4.02,4.28) -- (4.62,3.68) -- (3.22,1.58) -- cycle;
  \fill[gray!25] (3.42,1.38) -- (4.82,3.48) -- (5.42,2.88) -- (4.02,0.78) -- cycle;
  \draw[gray!60,thin] (3.42,1.38) -- (4.82,3.48) -- (5.42,2.88) -- (4.02,0.78) -- cycle;
  \fill[gray!25] (4.74,4.86) -- (6.14,6.96) -- (6.74,6.36) -- (5.34,4.26) -- cycle;
  \draw[gray!60,thin] (4.74,4.86) -- (6.14,6.96) -- (6.74,6.36) -- (5.34,4.26) -- cycle;
  \fill[pattern=north east lines,pattern color=gray!70] (5.54,4.06) -- (6.94,6.16) -- (7.54,5.56) -- (6.14,3.46) -- cycle;
  \draw[gray!60,thin] (5.54,4.06) -- (6.94,6.16) -- (7.54,5.56) -- (6.14,3.46) -- cycle;
  \fill[gray!25] (6.34,3.26) -- (7.74,5.36) -- (8.34,4.76) -- (6.94,2.66) -- cycle;
  \draw[gray!60,thin] (6.34,3.26) -- (7.74,5.36) -- (8.34,4.76) -- (6.94,2.66) -- cycle;
  \draw[very thick] (-0.10,0.10) -- (0.53,0.37) -- (0.55,1.25) -- (1.21,1.49) -- (1.35,2.15);
  \draw[very thick,densely dashed] (1.35,2.15) -- (2.02,2.78);
  \draw[very thick] (2.02,2.78) -- (2.65,3.05) -- (2.67,3.93) -- (3.33,4.17) -- (3.47,4.83);
  \draw[very thick,densely dashed] (3.47,4.83) -- (5.74,3.86);
  \draw[very thick] (5.74,3.86) -- (6.37,4.13) -- (6.39,5.01) -- (7.05,5.25) -- (7.19,5.91);
  \fill (0,0) circle (0.08) node[below left=-1pt] {\footnotesize $0$};
  \node[right] (bl) at (5.0,-1.6) {\footnotesize block of $\hcross_{v_0,2}$};
  \draw[->,gray,shorten >=1pt] (bl) -- (1.30,1.20);
  \node[right] (cn) at (6.6,0.9) {\footnotesize connector};
  \draw[->,gray,shorten >=1pt] (cn) -- (2.06,2.09);
  \node[right] (b2) at (8.4,3.6) {\footnotesize blocks at $v_2$};
  \draw[->,gray,shorten >=1pt] (b2) -- (7.18,3.02);
  \draw[->,gray] (7.00,9.00) -- (10.00,6.00) node[midway,above right=-2pt] {\footnotesize drift of $Y_{v_t}$ ($\protect\sigg=+$)};
    \draw[->,gray] (-2.2,0) -- (10.5,0) node[right] {$z_1$};
  \draw[->,gray] (0,-2.4) -- (0,14.4) node[above] {$z_2$};
\end{tikzpicture}
\caption{\small The coarse-graining in the proof of Proposition~\ref{prop:rarefast}, in the original Cartesian $\Z^2$ coordinates but not to scale, for $\protect\sigg=+$ and $T=3$. The levels are cut into layers $[I_t,I_{t+1})$ of $n+d$ levels each. In the first $n$ levels of layer $t$, the site $v$ of the coarse-grained lattice with $t(v)=t$ owns $T$ disjoint blocks (parallelograms with two sides along $\Dir$ and two along $\ead$), centered at transverse coordinates $z_{v,\protect\bbullet}$.  The hatched block is the one whose heavy-crossing event is used and the thick path is the heavy crossing.  The dashed segment is the  connector, built in the last $d$ levels of the layer, that joins the exit vertex of the heavy crossing to the entry vertex of the heavy crossing  in the next layer. Along the oriented path $v_0,v_1, v_2,\dots$ the centers $Y_{v_t}$ drift in the direction $\protect\sigg\ead$.}
\label{fig:coarse}
\end{figure}
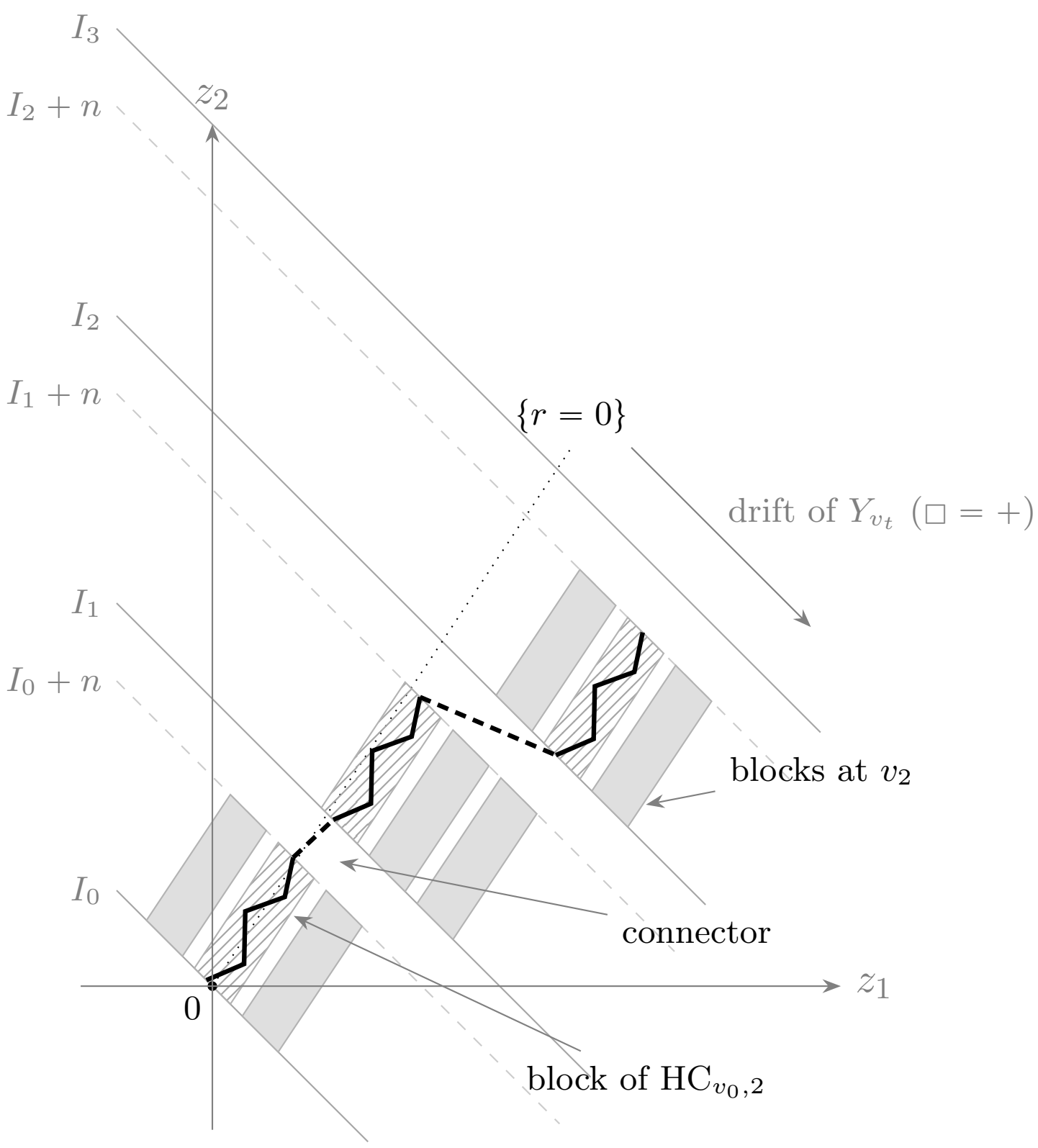

\begin{remark}[$g$ instead of $f_\sig$]  \label{rk:gversion}
If the assumption of Proposition \ref{prop:rarefast} is strengthened to $s\to\infty$ and $n/s\to\infty$,  
then \eqref{eq:rarefast} holds also with $g(x_n-x_0)$ in place of $f_\sig(x_n-x_0)$, for the same $R$ and $\e$. Since $g\le f_\sig$ by Lemma~\ref{lem:contact}, the version with $g$ implies \eqref{eq:rarefast} for both signs along such sequences. 
We prove here that the $g$ version  follows from  having \eqref{eq:rarefast}  for both  signs $\pm$ and  with $\e/2$ in place of $\e$. 

For $z\in\R^2_{\ge0}$ with $\level(z)=n>0$, homogeneity and \eqref{eq:decomp} give 
$g(z)=g(n\Dir+\rcoord(z)\tspa\ead)=n\hat g(\rcoord(z)/n)$, 
while $f_\sig(z)=n\hat g(0)+\ghatpr\sig\rcoord(z)$ by \eqref{eq:supports}. Hence, 
\[
\text{with $h=\rcoord(z)/n$,} \quad f_\sig(z)-g(z)=n\bigl[\hat g(0)+\hat g'(0\sigg)\tsp h-\hat g(h)\bigr] .\]
By \eqref{eq:contact} the bracket is $o(\abs h)$ as $h\to0$ with $\sigg h>0$. So given $\e>0$ there exists  $h_0>0$ such that the bracket is at most $\frac{\e}{4R}\abs h$ whenever $\sigg h>0$ and $\abs h\le h_0$. A path $x_{0:n}$ with $\abs{\rcoord(x_i)}\le Rs$ for all $i$ has $z=x_n-x_0$ with $\abs{\rcoord(z)}\le2Rs$, hence $\abs h=\abs{\rcoord(z)/n}\le2Rs/n\le h_0$ once $n/s\ge2R/h_0$. Choosing $\sigg$ with $\sigg\rcoord(z)\ge0$ (either sign if $\rcoord(z)=0$) then gives $f_\sig(z)-g(z)\le\frac{\e}{4R}\abs{\rcoord(z)}\le\e s/2$. Therefore, far enough along a sequence with $n/s\to\infty$, every such path with $\pawgt(x_{0:n})>g(x_n-x_0)+\e s$ satisfies $\pawgt(x_{0:n})>f_\sig(x_n-x_0)+\e s/2$ for at least one of the two signs.  Thereby the event with $g$ is contained in the union of the two signed events of \eqref{eq:rarefast}  with $\e/2$. The hypothesis $n/s\to\infty$ makes the slope $h$ of the crossings small, so that the first-order contact \eqref{eq:contact} applies. 

The convergence \eqref{eq:rarefast}  does not hold in general with $ng(\Dir)$ in place of $f_\sig(x_n-x_0)$ or $g(x_n-x_0)$.  Since $f_\sig(z)-ng(\Dir)=\ghatpr\sig\rcoord(z)$, a crossing with $\ghatpr\sig\rcoord(z)$ positive and of order $s$, and with weight close to $f_\sig(z)$, exceeds $ng(\Dir)$ by an amount of order $s$. The proof below works with $f_\sig$ rather than $g$ because $f_\sig$ is linear, so that bounds for consecutive pieces of a path add up, while $g$ is not additive.
\end{remark}
\begin{proof}[Proof of Proposition \ref{prop:rarefast}]
Assume  the assertion fails. Then  there exist $p_0>0$ and a sequence of pairs $(s,n)=(s,n_s)$, with $s\to\infty$ and integers $n_s\ge1$, 
along which the probability in \eqref{eq:rarefast} is at least $p_0$. We fix such a sequence and derive a contradiction. The argument is the same for $\sigg=-$ and $\sigg=+$. Figure~\ref{fig:coarse} shows the coarse-graining in the original coordinates, and Figure~\ref{fig:arrangement} the construction of steps 1--3 in the coordinates $(\rcoord,\level)$.

\medskip
\emph{1.\ Set-up.} Set $\varkappa=10$ and let $\zeta=\zeta(\varkappa)$ be given by Lemma~\ref{lem:oriented}. Fix an integer $T$ with $(1-p_0)^T<\zeta/2$, and set
\[
S_0=100T(R+1)\lceil s\rceil,\quad C_0=200T(R+1),\quad\text{and}\quad d=\lceil D_0s\rceil,
\]
where $D_0=D_0(T,R, \e,g)$ is chosen in step 4, with no dependence on $s$ or $n$. Note that $S_0\le C_0s$. For $t\in\Z_{\ge0}$ set $I_t=t(n+d)$ and call the interval $[I_t,I_{t+1})$  of levels  the $t$-th \emph{layer}. It consists of $n$ levels, in which the heavy crossings below take place, followed by $d$ levels reserved for the connectors.   The layers, together with transverse positions spaced by $S_0$, are indexed by the sites of a \emph{coarse-grained lattice}  $\Z^2$. The last step of the proof  is  an oriented percolation argument on the coarse-grained lattice. This produces a contradiction   for one fixed pair $(s,n)$ with large enough $s$.  

%No relation between $n$ and $s$ is needed. The length $n$ enters the construction only as the number of levels that a layer reserves for the heavy crossings. The connectors have length $d$ of order $s$ whatever $n$ is, and the contradiction in step 6 is reached for one fixed pair $(s,n)$ with large enough $s$.  

Levels and transverse coordinates on the coarse-grained lattice are not denoted by $\level$ and $\rcoord$ but by the functions $t$ and $Y$ below. The terms  \emph{site} and \emph{oriented path} refer to the coarse-grained lattice, while \emph{vertex} and \emph{admissible path} refer to the original one. To a site $v=(v_1,v_2)$ of the coarse-grained lattice assign its layer index and transverse center
\[
t(v)=v_1+v_2\qquad\text{and}\qquad Y_v=\sigg\bigl(t(v)+v_2\bigr)S_0=\sigg(v_1+2v_2)S_0 .
\]
Then $I_{t(v)}$ is the first level of the layer of $v$. An oriented step $v\to v+e_1$ or $v\to v+e_2$ moves to layer $t(v)+1$  and changes $Y_v$ by $\sigg S_0$ or $\sigg2S_0$. The transversal  drift is towards  the side on which $f_\sig$ has first-order contact with $g$, as determined by the sign $\sigg$ in \eqref{eq:contact}.

\medskip

%The heavy crossings are produced from the contradiction assumption $p_0>0$.

\begin{figure}[htb]
\centering
\begin{tikzpicture}[scale=0.62,>=Stealth]
  %% coordinates: horizontal r, vertical level; layer t of site v, drawn for sigma=+, T=3 (not to scale)
  \draw[->,gray] (-9.2,-0.4) -- (-9.2,7.4) node[above] {\footnotesize $\level$};
  \draw[->,gray] (-9.2,-0.4) -- (10.4,-0.4) node[right] {\footnotesize $\rcoord$};
  \foreach \y/\lab in {0/$I_{t(v)}$, 4/$I_{t(v)}+n$, 6/$I_{t(v)+1}$} {\draw[gray!50] (-8.8,\y) -- (10,\y); \node[right,gray] at (10,\y) {\footnotesize \lab};}
  %% region R_v
  \draw[dashed,gray] (-8.5,0) rectangle (8.5,6);
  \node[gray,below right] at (-8.5,6) {\footnotesize $\cR_v$};
  %% trial blocks
  \foreach \x in {-2.6,0,2.6} \fill[gray!20] (\x-0.7,0) rectangle (\x+0.7,4);
  \draw[<->,gray] (-3.3,0.4) -- (-1.9,0.4) node[midway,above] {\footnotesize $2Rs$};
  \node[below] (astrip) at (0,-0.75) {\footnotesize blocks of $\hcross_{v,\bbullet}$, centered at $\rcoord(z_{v,\bbullet})$};
  \draw[->,gray,shorten >=1pt] (astrip) -- (0.2,0);
  \draw[->,gray,shorten >=1pt] (astrip) -- (-2.5,0);
  \draw[->,gray,shorten >=1pt] (astrip) -- (2.5,0);
  \draw[<->,gray] (-8.5,1.3) -- (8.5,1.3) node[pos=0.8,above] {\footnotesize $8S_0$};
  %% a successful crossing in the middle block
  \draw[thick] (0.3,0) .. controls (-0.4,1.4) and (0.5,2.6) .. (-0.2,4);
  \fill (-0.2,4) circle (0.07);
  \node[above left=-2pt] at (-0.2,4) {\footnotesize $z_t$};
  %% J_v at level I_t+n
  \draw[very thick] (-3.3,4) -- (3.3,4);
  \node[below right] at (3.3,4) {\footnotesize $\cJ^{\rm out}_v$};
  %% J' for the two forward neighbors at level I_{t+1}
  \draw[very thick] (2.0,6) -- (4.6,6); \node[above] at (3.3,6) {\footnotesize $\cJ^{\rm in}_{v+e_1}$};
  \draw[very thick] (6.0,6) -- (8.6,6); \node[above] at (7.3,6) {\footnotesize $\cJ^{\rm in}_{v+e_2}$};
  %% connectors
  \draw[thick] (-0.2,4) .. controls (0.8,4.6) and (2.4,5.2) .. (3.0,6);
  \draw[thick] (-0.2,4) .. controls (2.5,4.4) and (5.5,5.3) .. (7.2,6);
  \node[left] (acon) at (-0.6,5.2) {\footnotesize connectors from $z_t$};
  \draw[->,gray,shorten >=1pt] (acon) -- (1.9,5.3);
  \draw[->,gray,shorten >=1pt] (acon) -- (3.0,4.7);
  %% center marks
  \draw[dotted] (0,0) -- (0,4);
  \draw[<->,gray] (-3.3,2.2) -- (3.3,2.2) node[pos=0.75,above] {\footnotesize $S_0/2$};
\end{tikzpicture}
\caption{\small The arrangement at a site $v$ of the coarse-grained lattice, for $\protect\sigg=+$ and $T=3$, drawn in the coordinates $(\rcoord,\level)$ and not to scale. In the first $n$ levels of the layer, $T$ disjoint blocks of width $2Rs$ centered at $\rcoord(z_{v,\protect\bbullet})$ host the heavy-crossing events $\hcross_{v,\protect\bbullet}$, which are independent translates of the event in \eqref{eq:rarefast}. All of them lie in $\{\abs{\rcoord-Y_v}\le S_0/4\}$, so a heavy crossing exits on the level interval $\cJ^{\rm out}_v$. The connection events ask for paths (connectors) from every point of $\cJ^{\rm out}_v$ to every point of $\cJ^{\rm in}_{v'}$, for the two forward neighbors $v'=v+e_1,v+e_2$ (centered at $Y_v+S_0$ and $Y_v+2S_0$), with weight at least $f_\protect\sig-\e s/2$ and all vertices in the region $\cR_v$ (dashed), which is centered at $Y_v$. The next layer begins at level $I_{t(v)+1}$; the blocks of $v+e_i$ are centered at $\rcoord(z_{v+e_i,\protect\bbullet})$, within $\tfrac12$ of $Y_{v+e_i}+h_\protect\bbullet$.}
\label{fig:arrangement}
\end{figure}
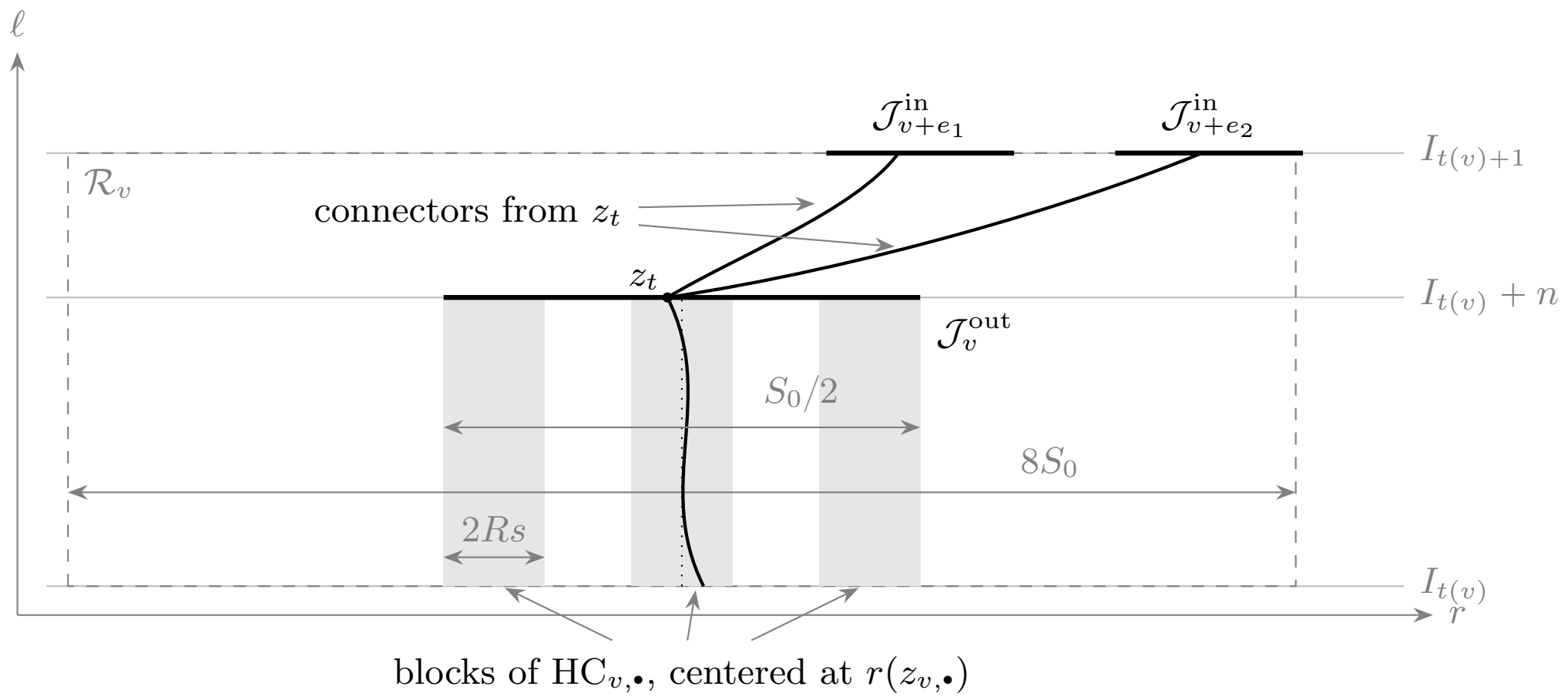
\emph{2.\ Heavy crossings.} At each site $v$ define $T$ events as follows. For $1\le j\le T$ let $h_j=(2j-T-1)\cdot3(R+1)\lceil s\rceil$ and let $z_{v,j}$ be the vertex of the original lattice on level $I_{t(v)}$ that satisfies  \[ \abs{\rcoord(z_{v,j})-Y_v-h_j}\le\tfrac12.\]
Such a vertex $z_{v,j}$ exists, since, by the paragraph above \eqref{stairdef},
the transverse coordinates of the lattice points on level $I_{t(v)}$ make up the set $\Z-\Dir_1I_{t(v)}$. The $j$th \textit{heavy-crossing event $\hcross_{v,j}$ at $v$}  is the translate by $z_{v,j}$ of the event in \eqref{eq:rarefast}: namely, that there exists a crossing $x_{I_{t(v)}:I_{t(v)}+n}$ of the block
\be\label{vj-block}
\{z:I_{t(v)}\le\level(z)\le I_{t(v)}+n,\ \abs{\rcoord(z)-\rcoord(z_{v,j})}\le Rs\}\ee
with weight
\be\label{t.aux89} \pawgt(x_{I_{t(v)}:I_{t(v)}+n})>f_\sig(x_{I_{t(v)}+n}-x_{I_{t(v)}})+\e s.
\ee
Consecutive centers $\rcoord(z_{v,j})$ differ by at least $6(R+1)\lceil s\rceil-1>2Rs$ for every $s\ge1$, 
%as $6(R+1)\lceil s\rceil\ge6Rs+6$
so the $T$ blocks \eqref{vj-block} are pairwise disjoint and thereby the events $\{\hcross_{v,j}\}_{1\le j\le T}$ are independent. By translation invariance $\P(\hcross_{v,j})\ge p_0$
and hence $\P\bigl(\textstyle\bigcap_{j=1}^T \hcross_{v,j}^c\bigr)\le (1-p_0)^T$.
All $T$ blocks lie in $\{z: \abs{\rcoord(z)-Y_v}\le S_0/4\}$, since \[ \abs{h_j}+Rs+\tfrac12\le3T(R+1)\lceil s\rceil+(R+1)\lceil s\rceil\le S_0/4. \]

\medskip

%Connections are produced by Lemma~\ref{lem:intervals}.

\emph{3.\ Connections.}
The crossings of the blocks at $v$ end on the level interval
\begin{align*}
\cJ^{\rm out}_v=\{z:\ \level(z)=I_{t(v)}+n,\ \abs{\rcoord(z)-Y_v}\le S_0/4\} .
\end{align*}
For each of the two forward neighbors $v'\in\{v+e_1,v+e_2\}$ of $v$, consider the level intervals
\begin{align*}
\cJ^{\rm in}_{v'}=\{z:\ \level(z)=I_{t(v)+1},\ \abs{\rcoord(z)-Y_{v'}}\le S_0/4\}.
\end{align*}
The \emph{connection event} $\conn_{v,v'}$ is the event that each pair  $u\in\cJ^{\rm out}_v$ and $u'\in\cJ^{\rm in}_{v'}$ is joined by an admissible path $x_{u:u'}$ with $\pawgt(x_{u:u'})\ge f_\sig(u'-u)-\e s/2$, all of whose vertices $x_i$ satisfy $\abs{\rcoord(x_i)-Y_v}\le4S_0$. Its levels lie in $[I_{t(v)}+n,I_{t(v)+1}]$.

\medskip
\emph{4.\ Lower bound on the connection probability.}
A displacement $z=u'-u$ from step 3 with $u\in\cJ^{\rm out}_v$ and $u'\in\cJ^{\rm in}_{v'}$ has $\level(z)=d$ and, since $Y_{v'}-Y_v\in\{\sigg S_0,\sigg2S_0\}$, the transverse displacement satisfies $\sigg\rcoord(z)\in[S_0/2,5S_0/2]$. By \eqref{eq:contact} there exists $h_0=h_0(\e,T,R)>0$ with
\[
0\le f_\sig(\Dir+h\ead)-g(\Dir+h\ead)\le\frac{\e}{12C_0}\abs h\qquad\text{for }0<\sigg h\le h_0 .
\]
Choose $D_0$ with $3C_0/D_0\le\min\{h_0,\theta_0/2,6\}$; in particular $D_0\ge C_0/2\ge1$. Then $h=\rcoord(z)/d$ satisfies
\[
0<\sigg h\le\frac{5S_0}{2D_0s}\le\frac{5C_0}{2D_0}<h_0.
%\qquad\text{and}\qquad
%\abs{\rcoord(z)}\le\frac{5C_0}{2D_0}\,d\le\frac{\theta_0}{2}\level(z),
\]
From this,
\be\label{t.aux27}
\abs{\rcoord(z)}=\abs{h}d
\le\frac{5C_0}{2D_0}\,d\le\frac{\theta_0}{2}\level(z)  \qquad\text{and}\qquad
d\abs{h}=\abs{\rcoord(z)}\le\frac{5S_0}2\le \frac{5C_0s}2.
\ee
The first inequality of \eqref{t.aux27} says that
 $z=u'-u$ is $\theta_0/2$-interior. Since this holds for all $u\in\cJ^{\rm out}_v$ and $u'\in\cJ^{\rm in}_{v'}$, part of the hypotheses of Lemma~\ref{lem:intervals} is verified. 
 
 By the first equality in \eqref{eq:decomp}, $z=d(\Dir+h\ead)$, so by homogeneity 
\be\label{t.hom}
f_\sig(z)-g(z)=d\bigl[f_\sig(\Dir+h\ead)-g(\Dir+h\ead)\bigr]\le\frac{\e}{12C_0}\,d\abs h
%=\frac{\e}{12C_0}\abs{\rcoord(z)}
\overset{\eqref{t.aux27}}\le\frac{\e}{12C_0}\cdot\frac{5C_0s}{2}<\frac{\e s}4.
\ee

Set $D=D_0+1$. The levels of $\cJ^{\rm out}_v,\cJ^{\rm in}_{v'}$ differ by $d=\lceil D_0s\rceil\in[D_0s,D_0s+1]\subset[s,Ds]$ for $s\ge1$, and their $\rcoord$-widths are at most $S_0/2\le C_0s/2\le D_0s\le Ds$.
Now apply Lemma~\ref{lem:intervals} to $\cJ^{\rm out}_v,\cJ^{\rm in}_{v'}$ with this $D$, $\theta=\theta_0/2$, $\eta=\e/4$, and $\delta=1/2$.  With probability at least $p(s)\to1$ as $s\to\infty$, uniformly in $v$, every pair $u\in\cJ^{\rm out}_v$ and  $u'\in\cJ^{\rm in}_{v'}$ is joined by a path $\pi^{u,u'}$ all of whose vertices are within $s/2$ of $[u,u']$ and whose weight satisfies  $\pawgt(\pi^{u,u'})\ge g(z)-\e s/4$. Use \eqref{t.hom} to compare this with $f_\sig$:
\be\label{t.aux91}
\pawgt(\pi^{u,u'})\ge g(z)-\e s/4> f_\sig(z)-\e s/2=f_\sig(u'-u)-\e s/2.  \ee
Since $\rcoord$ is linear, every point of $[u,u']$ has
\begin{align*}
    \abs{\rcoord-Y_v}&\le\max\{\abs{\rcoord(u)-Y_v},\abs{\rcoord(u')-Y_v}\}\\
    &\le\max\{\abs{\rcoord(u)-Y_v},\abs{\rcoord(u')-Y_{v'}}\} + \abs{Y_{v}-Y_{v'}} \le S_0/4+2S_0=\tfrac94S_0.
\end{align*}
By \eqref{eq:lip-r},  $\abs{\rcoord(\pi^{u,u'}_i)-Y_v}\le9S_0/4+s/2\le4S_0$.  The conditions of  the definition of $\conn_{v,v'}$ have been verified and thereby  $\P(\conn_{v,v'})\ge p(s)$.

\medskip
\emph{5.\ Good and bad sites on the coarse-grained lattice.}
Call $v$ \emph{good} if some heavy-crossing event $\hcross_{v,j}$ and both connection events $\conn_{v,v+e_1}$ and $\conn_{v,v+e_2}$ occur, and \emph{bad} otherwise. Then
\[\P(v\text{ bad})\le(1-p_0)^T+2(1-p(s))\le\zeta,\]
for all $v\in\Z^2_{\ge0}$ and all members $(s,n)$ of the sequence with $s\ge s_2$, where $s_2\ge s_1$ 
%$s_1$ from Lemma~\ref{lem:intervals} with the above parameters 
is such that $2(1-p(s))\le\zeta/2$ for $s\ge s_2$ (recall that $T$ was already chosen to satisfy $(1-p_0)^T<\zeta/2$).  The events $\hcross_{v,\bbullet}$ and $\conn_{v,v+e_i}$, $i\in\{1,2\}$, are determined by the weights in the region
\[
\cR_v=\{z:\ I_{t(v)}\le\level(z)\le I_{t(v)+1},\ \abs{\rcoord(z)-Y_v}\le4S_0\},
\]
because the blocks defining $\hcross_{v,\bbullet}$ lie in it and the events $\conn_{v,v+e_i}$ were defined in step 3 through connectors with all vertices in $\cR_v$. So the event $\{v\text{ good}\}$ is measurable with respect to $(\w_z)_{z\tsp\in\tsp\cR_v}$.

If $\cR_v\cap\cR_{v'}\ne\emptyset$ then $\abs{t(v)-t(v')}\le1$ and $\abs{Y_v-Y_{v'}}\le8S_0$, equivalently, $\abs{(t(v)+v_2)-(t(v')+v_2')}\le8$. Hence $\abs{v_2-v_2'}\le9$ and $\abs{v_1-v_1'}\le10$.
Consequently, if $\abs{v-v'}_\infty>\varkappa=10$, then $\cR_v\cap\cR_{v'}=\varnothing$. Hence for any finite family of sites $v$ at pairwise sup-distance greater than $\varkappa$ the regions $\cR_v$ are pairwise disjoint and the events $\{v\text{ good}\}$ are independent, as required in Lemma~\ref{lem:oriented}. Note that the choice of $\varkappa=10$ does not depend on $T$, $D_0$, or the scale $s$.

\medskip
\emph{6.\ Oriented percolation on the coarse-grained lattice and completion of the proof.}
Fix one member $(s,n)$ of the sequence with $s\ge s_2$, so that $\P(v\text{ bad})\le\zeta$ for all $v\in\Z^2_{\ge0}$. By Lemma~\ref{lem:oriented}, the event $\Gamma_\infty$ that there exists an infinite oriented path $v_{0:\infty}$ with $v_0=0$ of good sites has positive probability. On $\Gamma_\infty$ we build an infinite admissible path $x_{0:\infty}$ by concatenating heavy crossings and connectors, layer by layer, as follows. For each $t\ge0$, since $v_t$ is good and $t(v_t)=t$, some heavy-crossing event $\hcross_{v_t,j}$ occurs; let $x_{I_t:I_t+n}$ be a crossing of its block with $\pawgt(x_{I_t:I_t+n})>f_\sig(x_{I_t+n}-x_{I_t})+\e s$,
as in \eqref{t.aux89}. Its exit vertex $x_{I_t+n}$ lies in $\cJ^{\rm out}_{v_t}$ and the entry vertex $x_{I_{t+1}}$ of the crossing chosen in the next layer lies in $\cJ^{\rm in}_{v_{t+1}}$, so the connection event $\conn_{v_t,v_{t+1}}$ provides a connector $x_{I_t+n:I_{t+1}}$ from $x_{I_t+n}$ to $x_{I_{t+1}}$ with
$\pawgt(x_{I_t+n:I_{t+1}})\ge f_\sig(x_{I_{t+1}}-x_{I_t+n})-\e s/2$,
as in \eqref{t.aux91}. By additivity, for every $t\ge0$,
\be\label{t.3.4.h}
\pawgt(x_{0:I_t+n})>f_\sig(x_{I_t+n}-x_0)+(t+1)\e s-t\,\e s/2=f_\sig(x_{I_t+n}-x_0)+\frac{t\,\e s}2+\e s .
\ee
The vertex $x_0$ lies on level $0$ with $\abs{\rcoord(x_0)}\le S_0/4$, so it is one of at most $S_0/2+1$ lattice points. Call this finite deterministic set $\cI_0$. Moreover $x_0\le x_{I_t+n}$, so $\abs{x_{I_t+n}-x_0}_1=\level(x_{I_t+n}-x_0)=t(n+d)+n\ge t(n+d)$.

By \eqref{eq:defect-lower} and \ref{S1} applied to the finitely many origins $u\in\cI_0$, there exists an almost sure event $\Omega_0$ on which, for every $u\in\cI_0$, every admissible path $x_{0:l}$ from $u$ of length $l\ge R(\w)$ satisfies $\pawgt(x_{0:l})\le f_\sig(x_l-x_0)+\frac{\e s}{4(n+d)}\,l$ (here $s$ and $n+d$ are fixed numbers). On $\Gamma_\infty\cap\Omega_0$, which has positive probability, we get for all $t$ with $t(n+d)\ge R(\w)$:
\[
\frac{t\,\e s}2+\e s\ \overset{\eqref{t.3.4.h}}<\ \pawgt(x_{0:I_t+n})-f_\sig(x_{I_t+n}-x_0)\ \le\ \frac{\e s}{4(n+d)}\bigl(t(n+d)+n\bigr)\ \le\ \frac{\e s}{4}\,(t+1),
\]
which is false for every $t\ge0$. This contradiction proves \eqref{eq:rarefast}.
%With $\cI_0$, one sentence locates $x_0$ ($\abs{\rcoord(x_0)}\le S_0/4$ on level $0$), so the set of possible origins is finite and S1 applies exactly as stated. As a bonus, $R(\w)$ can be taken as a maximum over that finite set, so it doesn’t depend on which origin occurs. Without $\cI_0$, you skip locating $x_0$ but use S1 for a countable set, which also holds since a countable intersection of almost sure events is almost sure. Then the threshold is $R_u(\w)$ for each origin $u$, and you apply it with $u=x_0$, which is fixed once $\w$ is.
\end{proof}

\section{Characteristic scale}
\label{sec:scales}

Proposition~\ref{prop:rarefast} concerns crossings whose weight exceeds $f_\sig$ by $\e s$. We now find a length at which crossings are rare even when we allow the weight to fall below $g$ by the order $s$ of the block width. The tool is a change of measure in a region of size $s^3\times s$ that perturbs each weight with probability of order $s^{-2}$.  This keeps the  $L^2$ cost bounded.

The title of the section is borrowed from finite-size scaling in percolation. 
The correlation length is the side length at which the probability of a box crossing passes a fixed threshold.  See Kesten \cite{Kes-87} for the subcritical phase and Chayes, Chayes, Grimmett, Kesten, and Schonmann \cite{Cha-etal-89} for the supercritical phase. 

In our setting  the role of the box crossing   is played by a path  that crosses a block of width $16s$ and length $n$ and has  weight above $\ghatz n-Ks$. This event $\Elev(s,n)$ is defined below in Section \ref{sec:cubic}.   The characteristic length $\clen$ associated with the transverse width $s$ is the length $n$ at which the  probability $\P(\Elev(s,n))$ first drops below a fixed threshold $q_0$ (see \eqref{eq:clen-def}). 

This section does not directly use  the boundedness assumption \eqref{ass:bounded}. The  inputs are the shape theorem, Proposition~\ref{prop:rarefast}, and the assumption $\Dir\notin\Uperc$.    
This is the only section that explicitly uses   the assumption $\Dir\notin\Uperc$. It enters through the event $V$ in the proof of Lemma~\ref{lem:tilt}: since $\gtop(\Dir)<1$, every path in a tube has a positive fraction of vertices with weight below $\Mtop$, which the modification of the environment can increase.  After this section the remaining proof depends on the assumption $\Dir\notin\Uperc$ only through the properties of the scales $(s_k, n_k, \clen_k)$.  

\subsection{Local modification of the environment}\label{sub:localmod}
Let $F_\mu $ be the distribution function of $\mu$ and $F_\mu^{-1}(y)=\inf\{t:F_\mu(t)\ge y\}$, $0<y<1$, its quantile function. For  $t\in\R$ and $0<y<1$,  $F_\mu(t)\ge y$ if and only if $t\ge F_\mu^{-1}(y)$. Write $F_\mu(t-)=\mu((-\infty,t))$. For a finite set of sites $\cR$ and $\delta\in(0,1)$ define a random modification $\w'$ of $\w$ as follows. Let $\tilde\w=(\tilde\w_v)_{v\tsp\in\tsp\Z^2}$ be an independent copy of $\w$ and $\chi=(\chi_v)_{v\tsp\in\tsp\Z^2}$ i.i.d.\ Bernoulli($\delta$) variables, independent of $(\w,\tilde\w)$, and set
\[
\w'_v=\begin{cases}\max(\w_v,\tilde\w_v)&\text{if }v\in\cR\text{ and }\chi_v=1,\\ \w_v&\text{otherwise.}\end{cases}
\]
In words: independently at each site of $\cR$, with probability $\delta$ the weight is replaced by the maximum of itself and an independent resampled weight. Write $\jP$ for the joint law of $(\w,\tilde\w,\chi)$, so that $\w$ has law $\P$ under $\jP$, and $\P'$ for the law of $\w'$ under $\jP$.

\begin{lemma}\label{lem:kernel} The following hold.
\begin{enumerate}[label={\rm(\roman*)}]
\item $\w'\ge\w$ coordinatewise, and $\P'$ is a product measure, absolutely continuous with respect to $\P$, with
\[
\E\Bigl[\Bigl(\frac{d\P'}{d\P}\Bigr)^{\!2}\,\Bigr]\le e^{3\delta^2\abs\cR} .
\]
\item For $\alpha\in(0,1)$ and $t\in\R$ with $F_\mu(t)\le1-\alpha$ let $\psi_\alpha(t)=F_\mu^{-1}(F_\mu(t)+\alpha/2)-t$. Then $\psi_\alpha(t)>0$ for all such $\alpha$ and $t$. Furthermore, for every $\beta>0$ there exist $\alpha\in(0,\tfrac12)$ and $c'>0$ such that the set
\[
\cI=\{t<\Mtop:\ F_\mu(t)\le1-\alpha,\ \psi_\alpha(t)\ge c'\}
\]
satisfies $\mu\bigl(\R\setminus(\cI\cup\{\Mtop\})\bigr)\le\beta$. If $v\in\cR$ and $\w_v\in\cI$, then, conditionally on $\w$, with $\jP$-probability at least $\delta \alpha/2$ we have $\w'_v\ge\w_v+c'$, and these events are conditionally independent across $v$ given $\w$.
\end{enumerate}
\end{lemma}
\begin{proof}
(i) Clearly $\w'\ge\w$, and the coordinates of $\w'$ are independent because the triples $(\w_v,\tilde\w_v,\chi_v)$ are. For $v\notin\cR$ the marginal of $\P'$ is $\mu$. For $v\in\cR$ it is $(1-\delta)\mu+\delta\mu_{\max}$, where $\mu_{\max}$ is the law of the maximum of two independent $\mu$-distributed variables $X,Y$. For a Borel set $B$,
\[
\mu_{\max}(B)=\P\{X\in B,\ Y\le X\}+\P\{Y\in B,\ X<Y\}=\int_B\bigl(F_\mu(t)+F_\mu(t-)\bigr)\,\mu(dt),
\]
so $\mu_{\max}$ has density $h=F_\mu +F_\mu(\,\bbullet\,-)$ with respect to $\mu$, with $0\le h\le2$ and $\int h\,d\mu=\mu_{\max}(\R)=1$. Hence $d\P'/d\P=\prod_{v\tsp\in\tsp\cR}\bigl(1-\delta+\delta h(\w_v)\bigr)$, and
\[
\int\bigl(1-\delta+\delta h\bigr)^2d\mu=(1-\delta)^2+2\delta(1-\delta)+\delta^2\int h^2d\mu\le1-\delta^2+4\delta^2\le e^{3\delta^2}.
\]
By independence the second moment of $d\P'/d\P$ is the product of the $\abs\cR$ factors.
\medskip

(ii) Let $F_\mu(t)\le1-\alpha$ and set $y=F_\mu(t)+\alpha/2$, so that $0<y<1$. Since $F_\mu(t)<y$ and $F_\mu $ is right-continuous, there exists $t''>t$ with $F_\mu(t'')<y$. As $F_\mu $ is nondecreasing, every $t'$ with $F_\mu(t')\ge y$ satisfies $t'>t''$, hence $F_\mu^{-1}(y)\ge t''>t$, i.e.\ $\psi_\alpha(t)>0$. Next let $\beta>0$. The complement of $\cI\cup\{\Mtop\}$ is, up to the $\mu$-null set $(\Mtop,\infty)$, the union of $\cA_1=\{t<\Mtop:\ F_\mu(t)>1-\alpha\}$ and $\cA_2=\{t<\Mtop:\ F_\mu(t)\le1-\alpha,\ \psi_\alpha(t)<c'\}$. As $\alpha\downarrow0$ the sets $\cA_1$ decrease to $\{t<\Mtop:\ F_\mu(t)=1\}$, which is empty because $t<\Mtop=\sup\operatorname{supp}\mu$ implies $\mu((t,\infty))>0$. So there exists $\alpha\in(0,\tfrac12)$ with $\mu(\cA_1)\le\beta/2$. For this $\alpha$, as $c'\downarrow0$ the sets $\cA_2$ decrease to $\{t<\Mtop:\ F_\mu(t)\le1-\alpha,\ \psi_\alpha(t)\le0\}$, which is empty by the first assertion. So there exists $c'>0$ with $\mu(\cA_2)\le\beta/2$.  

Finally let $v\in\cR$ with $\w_v=t\in\cI$. On the event $\{\chi_v=1,\ \tilde\w_v\ge F_\mu^{-1}(F_\mu(t)+\alpha/2)\}$, whose conditional probability given $\w$ is
\[
\delta\tsp\mu\{t':F_\mu(t')\ge F_\mu(t)+\alpha/2\}\ge\delta(1-F_\mu(t)-\alpha/2)\ge\delta\alpha/2,
\]
%for $0<y<1$, $\mu\{t':F(t')\ge y\}=\mu\bigl([F^{-1}(y),\infty)\bigr)=1-F(F^{-1}(y)-)\ge1-y$, since $F(t')<y$ for $t'<F^{-1}(y)$ 
we have $\tilde\w_v\ge F_\mu^{-1}(F_\mu(t)+\alpha/2)>t=\w_v$, hence $\w'_v=\max(\w_v,\tilde\w_v)=\tilde\w_v\ge t+\psi_\alpha(t)\ge t+c'$.  The conditional independence across $v$ holds because the pairs $(\tilde\w_v,\chi_v)$ are independent of each other and of $\w$.
\end{proof}

\subsection{Cubic length scale} \label{sec:cubic}
Fix the constant 
\be\label{t.Koo}   K=100(1+\abs{\ghatprm}+\abs{\ghatprp})\ee 
for the remainder of the paper. 
For real $s\ge1$, integer $n\ge1$, and $z\in\Z^2$ define the event 
\be\label{def:Ens}\begin{aligned}
\Elev_z(s,n)= &\bigl\{ \exists \text{ admissible path } x_{\level(z):\level(z)+n} \text{ such that } \abs{\rcoord(x_i)-\rcoord(z)}\le8s  \\
&\quad \text{ for } \level(z)\le i\le\level(z)+n  \,\text{  and }  \, \pawgt(x_{\level(z):\level(z)+n})>\ghatz n-Ks 
\bigr\} .  
\end{aligned}\ee
Write $\Elev(s,n)=\Elev_0(s,n)$. Thus $\Elev_z(s,n)$ is the translate of $\Elev(s,n)$ by $z$, and $\P(\Elev_z(s,n))=\P(\Elev(s,n))$ for every $z$.

Here and below, for an event of the form $E=\{\text{there exists }x\text{ with property }A\}$, an $x$ with property $A$ is said to \emph{witness} $E$, or to be a \emph{witness} of $E$. 
This  terminology is from logic.

\gist{The idea below is that raising the weights turns the event $\pawgt(x_{0:n})>ng(\Dir)-Ks$ into $\pawgt'(x_{0:n})>g(x_n-x_0)+\e s$ (in fact with $f_+$ in place of $g$, but this is immaterial). Each site of the box is resampled independently with probability $\delta$ of order $s/n$, and when it is, its weight goes up by an amount of order one. The boost does not target any path: it is independent of $\w$, so whichever path witnesses the event gains, with high probability, about $\delta n$ times a constant, i.e.\ a constant times $s$ (this uses the events $U$ and $V$, which guarantees many improvable sites on every path: $V$ bounds the number of top-weight vertices on every path through a bounded-differences inequality, and $U$ bounds the number of vertices with weight in $\cA$ through a binomial tail and a union bound over the $2^n$ paths). This achieves the job since $(g(x_n-x_0)+\e s)-(ng(\Dir)-Ks)$ is at most of order $s$. The price is that the raised environment has a different law $\P'$, but its $L^2$ density with respect to $\P$ is bounded, since the cost is of order $\abs{\cR}\delta^2\approx ns(s/n)^2=s^3/n$. So the event of Proposition~\ref{prop:rarefast}, which is rare under $\P$, stays rare under $\P'$. [VERY high level hand waving: the new measure increases the weights at each site of $\cR$ by roughly $\delta$. The total variance is then about $\delta^2|\cR|$. So for the two measures to be close, one needs this $O(1)$.] This is why the length is $n\approx s^3$.}

\begin{lemma}\label{lem:tilt}
$\P(\Elev(s,\lceil s^3\rceil))\to0$ as $s\to\infty$.
\end{lemma}
\begin{proof}
Let $n=\lceil s^3\rceil$, so that $s^3\le n\le2s^3$, and \[ \cR=\{z\in\Z^2:\ 0\le\level(z)\le n,\ \abs{\rcoord(z)}\le8s\}, \]
a region with $n+1$ levels and at most $16s+1$ lattice points per level, so $\abs\cR\le(n+1)(16s+1)\le34ns$ for $s\ge1$. Any  path witnessing $\Elev(s,n)$ lies in $\cR$. Since $\Dir\notin\Uperc$, $\gtop(\Dir)<1$, and we set
\[
\varphi=\tfrac14\bigl(1-\gtop(\Dir)\bigr)\in(0,\tfrac14] \quad\text{and}\quad \beta=16^{-1/\varphi}. 
\]
Let $\alpha,c',\cI$ be as in Lemma~\ref{lem:kernel}(ii) with this $\beta$. Now  $\cA=\R\setminus(\cI\cup\{\Mtop\})$ has $\mu(\cA)\le\beta$. Set 
\be\label{t.aux.eta}
\eta=\frac12(K+16\abs{\ghatprp}+1)\quad\text{and}\quad
\delta=\frac{8\eta s}{\alpha\varphi c'n}.\ee 
Note that $\delta\le8\eta/(\alpha\varphi c's^2)$, which is less than $1$ for large $s$. 

Let $V$ be the event that every admissible path $x_{0:n}$ with all vertices in $\cR$ has at most $(1-2\varphi)n$ collected vertices with weight $\Mtop$. For $u$ on level $0$ and $v$ on level $n$, both in $\cR$ and with $u\le v$, let $\Ltop_{u,v}(\cR)$ be the maximal number of collected vertices with weight $\Mtop$ along an admissible path from $u$ to $v$ with all vertices in $\cR$. Then $\Ltop_{u,v}(\cR)\le\Ltop_{u,v}$, and $\E\Ltop_{u,v}\le\gtop(v-u)$ by superadditivity. 
%$\E\Ltop_{0,mz}\ge m\,\E\Ltop_{0,z}$ for $z=v-u$, and $\E\Ltop_{0,mz}/m\to\gtop(z)$ as $m\to\infty$ by the shape theorem and bounded convergence.
Then $z=v-u$ satisfies $z=n\bigl(\Dir+(\rcoord(z)/n)\ead\bigr)$ with $\abs{\rcoord(z)}\le16s$ and $n\ge s^3$, so $\abs{\rcoord(z)}/n\le16s^{-2}\to0$. By homogeneity and continuity of $\gtop$ at $\Dir$, $\gtop(z)\le n\bigl(\gtop(\Dir)+\varphi\bigr)=(1-3\varphi)n$ for all such $z$ once $s$ is large. The variable $\Ltop_{u,v}(\cR)$ is a function of the $\abs\cR$ independent indicators $\one\{\w_z=\Mtop\}$ across  $z\in\cR$. Changing one of them changes $\Ltop_{u,v}(\cR)$  by at most $1$. By the bounded differences inequality \cite[Theorem~6.2]{Bou-Lug-Mas-13}, 
\begin{align*}
\P\bigl\{\Ltop_{u,v}(\cR)>(1-2\varphi)n\bigr\}&\le\P\bigl\{\Ltop_{u,v}(\cR)-\E\Ltop_{u,v}(\cR)>\varphi n\bigr\}\\
&\le\exp\Bigl\{-\frac{2\varphi^2n^2}{\abs\cR}\Bigr\}\le\exp\Bigl\{-\frac{\varphi^2n}{17s}\Bigr\}\le e^{-\varphi^2s^2/17}.
\end{align*}
There are at most $(16s+1)^2\le289s^2$ pairs $(u,v)$, so $\P(V^c)\le289s^2e^{-\varphi^2s^2/17}\to0$. 

%The event $V$ is where $\Dir\notin\Uperc$ is used, through $\varphi=\tfrac14(1-\gtop(\Dir))>0$ in the definition of $V$.  In a direction inside the flat segment, a path could run almost entirely through sites of weight $\Mtop$, which can never be raised, and the modification below would have nothing to act on. The restriction to paths inside $\cR$ is what makes the bounded differences inequality useful, since $\Ltop_{u,v}$ itself depends on about $n^2$ sites. The event $U$ below uses only the choice of $\beta$.

Let $U$ be the event that every admissible path $x_{0:n}$ with all vertices in $\cR$ has at most $\varphi n$ collected vertices with weight in $\cA$. For a fixed such path the indicators $\one\{\w_{x_i}\in\cA\}$, $0\le i<n$, are i.i.d.\ Bernoulli with parameter $\mu(\cA)\le\beta$, and
\[
\P\{\mathrm{Binomial}(n,\beta)>\varphi n\}\le\binom{n}{\lceil\varphi n\rceil}\beta^{\lceil\varphi n\rceil}\le2^n\beta^{\varphi n}=2^n16^{-n}=8^{-n}.
\]
There are at most $(16s+1)2^n$ such paths (a starting point on level $0$ in $\cR$ and $n$ steps), so $\P(U^c)\le17s\,4^{-n}\to0$. On $U\cap V$, every admissible path $x_{0:n}$ with all vertices in $\cR$ has at least $n-(1-2\varphi)n-\varphi n=\varphi n$ collected vertices with weight in $\cI$, because, $\mu$-almost surely, every weight lies in exactly one of $\cI$, $\cA$, and $\{\Mtop\}$. 

On $\Elev(s,n)\cap U\cap V$ let $x_{0:n}$ be a witnessing path, chosen as a measurable function of $\w$ (for example, the lexicographically first one in the finite set of candidates). On $U\cap V$ it has at least $\varphi n$ collected vertices with weight in $\cI$. Let $B$ be a $\mathrm{Binomial}(\lceil\varphi n\rceil,\delta\alpha/2)$ variable. Its mean $\lambda$ satisfies $\lambda\ge\delta\alpha\varphi n/2=4\eta s/c'\to\infty$ as $s\to\infty$, and its variance is at most $\lambda$. By Lemma~\ref{lem:kernel}(ii), conditionally on $\w$, the number $N$ of collected vertices of $x_{0:n}$ with weight in $\cI$ at which $\w'\ge\w+c'$ stochastically dominates $B$. By this domination and Chebyshev's inequality,
\[ \jP\{N<\lambda/2\mid\w\}\le\mathrm P\{B<\lambda/2\}\le\mathrm P\{\abs{B-\lambda}>\lambda/2\}\le\frac{\lambda}{(\lambda/2)^2}=\frac4\lambda\to0 \quad\text{uniformly in $\w$.} \]
On $\{N\ge\lambda/2\}$ the weight of the path increases by at least $c'\lambda/2\ge c'\delta\alpha\varphi n/4=2\eta s$ at these vertices, and no weight decreases. So 
\[ \pawgt'(x_{0:n})\ge\pawgt(x_{0:n})+2\eta s>\ghatz n-Ks+2\eta s, 
\]
where $\pawgt'$ denotes the weight in the environment $\w'$. Since   $f_+(x_n-x_0)=\ghatz n+\ghatprp\rcoord(x_n-x_0)$ and $\abs{\rcoord(x_n-x_0)}\le16s$, in the environment $\w'$ the path satisfies  
\begin{align*}
\pawgt'(x_{0:n})\ &>\ \ghatz n-Ks+2\eta s\ \ge\ f_+(x_n-x_0)-16\abs{\ghatprp}s-Ks+2\eta s\\
&=\ f_+(x_n-x_0)+(2\eta-K-16\abs{\ghatprp})s\
\overset{\eqref{t.aux.eta}}=\ f_+(x_n-x_0)+s. 
\end{align*}  
Let $\cE(s,n)$ be the event that some admissible path $x_{0:n}$ with $\abs{\rcoord(x_i)}\le8s$ for all $i$ has $\pawgt(x_{0:n})>f_+(x_n-x_0)+s$. We have shown that \[ \Elev(s,n)\cap U\cap V\cap\{N\ge\lambda/2\}\subset\{\w'\in\cE(s,n)\}.\] 
Hence, writing $\Elev^\circ=\Elev(s,n)\cap U\cap V$, 
\begin{align*}
\P(\Elev^\circ)
&\le\jP\bigl(\Elev^\circ\cap\{N\ge\lambda/2\}\bigr)
+\jP\bigl(\Elev^\circ\cap\{N<\lambda/2\}\bigr)\\
&\le\jP\bigl(\Elev^\circ\cap\{\w'\in\cE(s,n)\}\bigr)
+\E\bigl[\one_{\Elev^\circ}\,\jP\{N<\lambda/2\mid\w\}\bigr]\\
&\le\jP\bigl(\Elev^\circ\cap\{\w'\in \cE(s,n)\}\bigr)+\sup_{\w}\jP\{N<\lambda/2\mid\w\}\le\P'(\cE(s,n))+\frac4\lambda.
\end{align*}

By Lemma~\ref{lem:kernel}(i) and the Cauchy--Schwarz inequality,
\begin{align*}
\P'(\cE(s,n))=\E\Bigl[\one_{\cE(s,n)}\frac{d\P'}{d\P}\Bigr]&\le\Bigl(\E\Bigl[\Bigl(\frac{d\P'}{d\P}\Bigr)^2\Bigr]\Bigr)^{1/2}\P(\cE(s,n))^{1/2}\\
&\le\exp\Bigl\{\frac32\cdot\frac{64\eta^2s^2}{\alpha^2\varphi^2(c')^2n^2}\cdot34ns\Bigr\}\P(\cE(s,n))^{1/2}\\
&=\exp\Bigl\{\frac{3264\,\eta^2}{\alpha^2\varphi^2(c')^2}\cdot\frac{s^3}{n}\Bigr\}\P(\cE(s,n))^{1/2}.
\end{align*}
Since $n\ge s^3$, the exponential factor is bounded. The event $\cE(s,n)$ is the event in Proposition~\ref{prop:rarefast} with $\sigg=+$, $R=8$, and $\e=1$, so $\P(\cE(s,n))\to0$ as $s\to\infty$.   Finally,  $\P(\Elev(s,n))\le\P(U^c)+\P(V^c)+\P'(\cE(s,n))+4/\lambda\to0$ as $s\to\infty$.
\end{proof}

\subsection{Characteristic length}\label{sub:charlength}
Consider a parameter $q_0\in(0,1)$, whose value is specified in Section~\ref{sec:witnesses}. For the rest of the paper, write $s_k=2^k$.  The \emph{characteristic length} $\clen_k$ defined in the next lemma is the longitudinal length at which crossings of blocks of width of order $s_k$ with weight allowance $Ks_k$ stop being likely. Lemma~\ref{lem:tilt} bounds  $\clen_k\le s_k^3$, and \eqref{eq:short} shows that it is much larger than $s_k$. Only these two facts, and the threshold property \eqref{eq:Lbounds}, are used later.

\gist{The point now is, given $q_0$ [which will be small], to find the right length of the box that makes the probability of a crossing $\approx q_0$. We know if we take it $s^3$, then the probability goes to $0$ and if we take it $s$, the probability goes to $1$. So the right scaling is somewhere between $s$ and $s^3$. Dyadics $s_k=2^k$ are used so that one can run a multi-scale argument, in the next section.}

\begin{proposition}\label{prop:charlength}
For every $C>0$,
\begin{equation}\label{eq:short}
\inf\bigl\{\P(\Elev(s,n)):\ n\in\N,\ 1\le n\le Cs\bigr\}\longrightarrow1\quad\text{ as }s\to\infty.
\end{equation}
Consequently there exists a positive integer $k_0=k_0(q_0,\mu,\Dir)$ such that for every $k\ge k_0$ the first dyadic integer $\clen_k\in\{1,2,4,\dots\}$ for which
\begin{equation}\label{eq:clen-def}
\min\bigl\{\P(\Elev(s_k,n)):\ n\in\N,\ \clen_k/2<n\le\clen_k\bigr\}\le q_0
\end{equation}
exists and satisfies
\begin{equation}\label{eq:Lbounds}
\clen_k\le2^{3k},\quad 2^{-k}\clen_k\to\infty,\quad\text{and}\quad \P(\Elev(s_k,n))>q_0\ \text{ for every integer }n<\clen_k/2 .
\end{equation}
\end{proposition}

\begin{proof} Let $z_n$ be the lattice point on level $n$ with $\rcoord(z_n)\in(-1,0]$. Then $z_n$ is $\theta_0/2$-interior for $n\ge2/\theta_0$, and, using the first equation in \eqref{eq:decomp}, $\abs{z_n-n\Dir}_1=\abs{\rcoord(z_n)}\abs\ead_1\le2$, so $\abs{g(z_n)-\ghatz n}\le2\tsp\Lip_{\theta_0/2}$. 

Apply Lemma~\ref{lem:travel} at scale $n$ (in place of $s$) with $D=1$, $\theta=\theta_0/2$, $\eta=K/(2C)$, and $\delta=\min\{1/4,4/C\}$. Then  for $n\ge n_1$, with probability at least $p(n)\to1$ as $n\to\infty$, there exists a path $x_{0:n}$ from $0$ to $z_n$ such that 
\begin{itemize} \itemsep=3pt
 \item    
  $\pawgt(x_{0:n})\ge g(z_n)-\eta n\ge \ghatz n-2\tsp\Lip_{\theta_0/2}-Ks/2$ (using $n\le Cs$), and 
  \item   vertices of $x_{0:n}$  are within $\delta n\le4s$ of the segment $[0,z_n]$, on which $\abs\rcoord\le1$. 
  \end{itemize}
  By \eqref{eq:lip-r} these vertices have $\abs\rcoord\le4s+1\le8s$, as $s\ge1$. If moreover $s>4\tsp\Lip_{\theta_0/2}/K$, then $2\tsp\Lip_{\theta_0/2}<Ks/2$ and this path $x_{0:n}$ witnesses $\Elev(s,n)$.

Given $\epsilon>0$ choose $n_0\ge n_1$ with $p(n)\ge1-\epsilon$ for $n\ge n_0$. Then $\P(\Elev(s,n))\ge1-\epsilon$ for all $n_0\le n\le Cs$ once $s>4\tsp\Lip_{\theta_0/2}/K$. For $1\le n<n_0$, the staircase $I\mapsto \hat z(I,0)$ (see \eqref{stairdef}) is an admissible path from $0=\hat z(0,0)$ to $z_n=\hat z(n,0)$ with $\abs\rcoord<1$ along it and with weight $S_n$, a sum of $n$ i.i.d.\ weights, satisfying $\min_{1\le n<n_0}\P\{S_n>\ghatz n-Ks\}\to1$ as $s\to\infty$. This proves \eqref{eq:short}. 
%With bounded weights, the probability is $1$ once $s$ is large. But saying $\to1$ survives the generalization to unbounded weights. 

By Lemma~\ref{lem:tilt}, for all large $k$ the integer $n=2^{3k}=s_k^3$ has $\P(\Elev(s_k,n))\le q_0$. It lies in $(l/2,l]$ for the dyadic $l=2^{3k}$, so \eqref{eq:clen-def} is satisfied with $l$ in place of $\clen_k$. Thus for large $k$,  $\clen_k$ exists and $\clen_k\le2^{3k}$. 

If $2^{-k}\clen_k$ did not tend to infinity, then along a subsequence $\clen_k\le C'2^k$ for some $C'$, and some  $n_k\le\clen_k\le C's_k$ would satisfy $\P(\Elev(s_k,n_k))\le q_0<1$, contradicting \eqref{eq:short}.  Finally, if an integer $n<\clen_k/2$ had $\P(\Elev(s_k,n))\le q_0$, the dyadic $l$ with $n\in(l/2,l]$ would satisfy $l<2n<\clen_k$ and \eqref{eq:clen-def}, contradicting the minimality of $\clen_k$.
\end{proof}

\begin{lemma}\label{lem:records}
There exists an infinite set $\cK\subset\{k\ge k_0\}$ such that for $k\in\cK$ and $0\le j\le k-k_0$,
\begin{equation}\label{eq:record}
\clen_{k-j}\ge2^{-4j}\clen_k .
\end{equation}
\end{lemma}
\begin{proof}
The positive sequence $\vartheta_k=\clen_k2^{-4k}\le2^{-k}$, $k\ge k_0$, tends to zero, so it has infinitely many strict record minima, i.e.\ indices $k$ with $\vartheta_k<\vartheta_{k'}$ for all $k_0\le k'<k$. Let $\cK$ be the set of these indices. For $k\in\cK$ and $0\le j\le k-k_0$ we have $\clen_k2^{-4k}\le\clen_{k-j}2^{-4(k-j)}$, which is \eqref{eq:record}.
\end{proof}

In the sequel  $n_k$ denotes some choice of integer that attains the minimum in \eqref{eq:clen-def}.  So
\be\label{t.n_k}
\P(\Elev(s_k,n_k))\le q_0 \quad\text{and}\quad  \clen_k/2<n_k\le\clen_k. \ee

\section{Control over a parallelogram}
\label{sec:endpoints}

 %The goal  Proposition \ref{prop:endpoints}  improves Proposition \ref{prop:rarefast} to  ALL paths in a block, not only ones that traverse the whole block. The scale $\clen_k$ comes in and  multi-scale surgery is needed.

Proposition~\ref{prop:rarefast} controls crossings between two prescribed levels. Section~\ref{sec:witnesses}, which completes the proof of Theorem~\ref{thm:main}, needs more: with high probability, every admissible path contained in a region of longitudinal size $O(\clen_k)$ and transverse size $O(s_k)$ has $\pawgt\le f_\pm+\e s_k$ for both signs. 
A union bound over the endpoints is not available, since there are $O(\clen_k^2s_k^2)$ pairs of them and Proposition~\ref{prop:rarefast} comes with no rate. 

To get around this 
we connect every vertex $u$ of the region to two \emph{coarse levels}, one below and one above $u$, by paths whose weight falls short of $f_\sig$ by at most a small multiple of $s_k$. Only boundedly many pairs of coarse levels then remain, and Proposition~\ref{prop:rarefast} applies to each pair. 

The connecting paths are assembled along a nested family of grids of levels. At each scale of the family, the grid intervals are either 
\begin{itemize}\itemsep=2pt
    \item 
long enough to host the crossing events $\Elev$ of \eqref{def:Ens},  whose probabilities are bounded below by the record scale property \eqref{eq:record}, or 

\item short enough that a deterministic staircase will do. 
\end{itemize}
Figures~\ref{fig:grid} and~\ref{fig:entrance} illustrate the construction.

Directedness forces the connecting paths to advance in level, which is why the short grid intervals are bridged by deterministic staircases. This is the point at which the construction differs from the one for the undirected first-passage model in \cite{Ope-26-convdiff}.   

Throughout this section, $s_k=2^k$,  $\clen_k=2^{m_k}$ is the characteristic scale chosen in Proposition~\ref{prop:charlength}, and $n_k$ satisfies \eqref{t.n_k}. The set 
$\cK$ is from  Lemma~\ref{lem:records} and  the scale index  $k\to\infty$ along $\cK$.   By \eqref{eq:Lbounds}, $m_k\le3k$ and $m_k-k\to\infty$.

\begin{remark}[Bounded vs.\ unbounded weights]\label{rk:overview-general}
The main result of this section is Proposition~\ref{prop:endpoints}, and nothing else from this section is used in the sequel. Its statement holds under \ref{A1}--\ref{A3}, but unbounded weights introduce additional technical difficulties, beginning with the bound \eqref{eq:badstair} on the weight of a deterministic staircase. We therefore carry out the exposition under the boundedness assumption \eqref{ass:bounded}, from the point where the staircases are introduced, after \eqref{eq:bincount}, to the end of Section~\ref{sec:endpoints}.  Whenever a step needs a change under \ref{A2}--\ref{A3}, a remark says what changes, so that the reader can continue with the modified definition or bound, and points to the place in Appendix~\ref{sub:staircases} where the details are given.

Here is an overview of the three  changes, each forced by the one before it.
\begin{enumerate}[label={\rm(\roman*)}] \itemsep=2pt
\item \emph{The staircases are replaced.} Under \ref{A2} a deterministic staircase may collect a vertex with a very negative weight, with a probability  too large for the union bounds. Each staircase is therefore replaced by a path that avoids these vertices, and follows the same line up to a constant transverse distance (see Figure~\ref{fig:stair-geometry}).  The weights of these paths are bounded below by a Bernstein bound, \eqref{eq:weights} in Lemma~\ref{lem:weights}, which takes the place of the Hoeffding bound \eqref{eq:badstair}. 
See Remarks~\ref{rk:stair-general} and~\ref{rk:fine-general}.
\item \emph{Paths are trimmed.} The construction attaches the replacement paths to other paths: to the crossings provided by the events $\Elev$ of \eqref{def:Ens} on the long grid intervals, and to the path $x_{u:u'}$ in the proof of Proposition~\ref{prop:endpoints}. This is not always possible at an endpoint of such a path, because there may be vertices with very negative weights at that endpoint or close to it. So a bounded number of vertices are removed from each end of these paths. See Remarks~\ref{rk:coarse-general} and~\ref{rk:endpoints-general}.
\item \emph{The removed weights are bounded above.} Removing vertices from a path must not lower its weight by much, so the weights of the removed vertices have to be bounded above. This is where the upper tail assumption \ref{A3} enters this section. See Remark~\ref{rk:heavy-general}.
\end{enumerate}
Apart from this, the construction and the estimates are the same, up to the values of the constants. In particular, Lemma~\ref{lem:grid} holds with a different constant $\cnst_1$ and with the grid points $\hat z(I,\rho)$ replaced by vertices of the new paths, see Remark~\ref{rk:grid-general}, and the proof of Proposition~\ref{prop:endpoints} uses that lemma but not the value of $\cnst_1$. A reader  content with this overview can skip the other bounded vs.~unbounded remarks of this section.
\end{remark}

\subsection{Grid scales}\label{sub:gridscales}
Fix $A,R>0$. Fix also an integer $J\ge1$, which selects the coarsest grid. In the probability estimates $J$ is kept fixed while $k\to\infty$ in $\cK$. At the very end $J\to\infty$. The \emph{scale indices} are the integers $j$ with \begin{equation}\label{eq:jstar}
J\le j\le j_*(k)=\lceil(m_k-1)/4\rceil . 
\end{equation}
The scale with index $j$ has \emph{longitudinal mesh} and \emph{transverse width}
\begin{equation}\label{eq:meshes}
l_{j,k}=2^{\max\{0,\,m_k-4j-1\}}\qquad\text{and}\qquad w_{j,k}=2^{k-j}=s_{k-j}, 
\end{equation}
both decreasing in $j\in[J, j_*(k)]$.  
For $j<j_*(k)$ we have $4j<m_k-1$, so $l_{j,k}=2^{m_k-4j-1}>1$ and $l_{j,k}/l_{j+1,k}\in\{2,4,8,16\}$, while $l_{j_*(k),k}=1$. Thus the meshes are nested powers of two, each a multiple of the next, with successive ratios at most $16$, and the finest mesh is $l_{j_*(k),k}=1$. A \emph{grid interval} at scale $j$ is a level interval $[I,I+l_{j,k}]$ with $I$ in
\begin{equation}\label{eq:gridset}
\cG_{j,k}=\{I\in l_{j,k}\Z:\ \abs I\le(A+4)\clen_k\}.
\end{equation}
When $j<j_*(k)$, $\clen_k/l_{j,k}=2^{4j+1}$. Since $\clen_k=2^{m_k}\le2^{4j_*(k)+1}$ and $l_{j_*(k),k}=1$, the bound $\clen_k/l_{j,k}\le2^{4j+1}$ holds also when $j=j_*(k)$. Thus the cardinality of the set of grid intervals satisfies 
\begin{equation}\label{eq:gridcount}
\abs{\cG_{j,k}}\le2(A+4)\clen_k/l_{j,k}+1\le\cnst_A2^{4j}\quad\text{for all }J\le j\le j_*(k).
\end{equation}

\medskip
\emph{Requirement on $k$.} Throughout this section we take $k\in\cK$ with
\begin{equation}\label{eq:kreq}
k\ge4k_0+3,
\end{equation}
where $k_0$ is the positive integer of Proposition~\ref{prop:charlength}. Then, since $m_k\le3k$, the index of the finest scale satisfies
\begin{equation}\label{t.aux32}
j_*(k)\le\tfrac14(m_k-1)+1\le\tfrac34(k+1)\overset{\eqref{eq:kreq}}\le k ,
\end{equation}
and every scale index $j\le j_*(k)$ satisfies
\begin{equation}\label{t.aux34}
k-j\overset{\eqref{t.aux32}}\ge k-\tfrac34(k+1)=\tfrac14(k-3) \overset{\eqref{eq:kreq}}\ge k_0 .
\end{equation}
Hence \eqref{eq:record} applies to the pair $(k,j)$ and gives $l_{j,k}=2^{-4j-1}\clen_k\le\tfrac12\clen_{k-j}$ when $l_{j,k}>1$. Since $s_{k-j}=2^{k-j}=w_{j,k}$, the last assertion of \eqref{eq:Lbounds} at index $k-j$ yields
\begin{equation}\label{eq:trialprob}
\P(\Elev(w_{j,k},l))>q_0\quad\text{for every pair of integers $(l,j)$ with $1\le l<l_{j,k}$ and $J\le j\le j_*(k)$.}
\end{equation}
In words: a crossing $x_{I:I+l}$ of a block of width $16w_{j,k}$ and of any length  $l<l_{j,k}$, with weight allowance $\pawgt(x_{I:I+l}) >  lg(\Dir)-Kw_{j,k}$, occurs with probability at least $q_0$. This is the only place where the record scales $\{\clen_i: i\in\cK\}$ are used.

Set $\cnst_*=50/\theta_0$. A scale of index $j$ is \emph{coarse} if $l_{j,k}\ge8\cnst_*j^2w_{j,k}$ and \emph{fine} otherwise. Note that
\begin{equation}\label{eq:j2w}
j^2w_{j,k}=j^22^{k-j}\le s_k\max_{i\ge1}i^22^{-i}\le2s_k\qquad\text{for all }j\ge1 . 
\end{equation}
Looking ahead to the utility of this bound, 
the transverse errors of order $j^2w_{j,k}$ allowed below in \eqref{eq:gridcost} are at most $2s_k$ uniformly over scales $j$, and their sum over  $j$ is of order $s_k=2^k$. Up to a constant factor, this sum is the amount by which $\pawgt$ can fall short of $f_\sig$ in \eqref{eq:extcost} and \eqref{eq:extcost-exit} in the proof of Proposition~\ref{prop:endpoints}.

\subsection{Staircases and crossings at a single  scale}\label{sub:onescale}
This section covers preparatory work for the proof of Lemma~\ref{lem:grid} below.  It constructs, at a single scale, the events on which the grid points can be connected. Throughout, $k\in\cK$ satisfying \eqref{eq:kreq} is fixed, and $k\to\infty$ only once we reach  Lemma~\ref{lem:grid}\ref{lem:grid:prob}. Fix a scale index $j$ with $J\le j\le j_*(k)$. Fix a grid interval $[I,I+l_{j,k}]$, that is, $I\in\cG_{j,k}$ from \eqref{eq:gridset}. The transverse coordinates $\rho$ with $\abs\rho\le Rs_k$ are grouped into \emph{bins}: the bin with index $h\in\Z$ is the interval $[hw_{j,k},(h+1)w_{j,k})$.  Only the bins meeting $[-Rs_k,Rs_k]$ matter.  So the bin indices form the set
\begin{align}\label{eq:bincount}
\begin{split}
\cH_{j,k}&=\bigl\{h\in\Z:\ [hw_{j,k},(h+1)w_{j,k})\cap[-Rs_k,Rs_k]\ne\emptyset\bigr\}\quad\text{that satisfies}\\
\abs{\cH_{j,k}}&\le2Rs_k/w_{j,k}+2\le\cnst_R2^j.
\end{split}
\end{align}
Fix a bin index $h\in\cH_{j,k}$. The task at $(j,I,h)$ is to connect $\hat z(I,\rho)$ to $\hat z(I+l_{j,k},\rho)$ for every $\rho$ in the bin $h$, by a path that stays within transverse distance $O(j^2w_{j,k})$ of $\rho$ and has $\pawgt\ge f_\sig-O(j^2w_{j,k})$.  

\medskip

\emph{Staircases.} Let $\e_0>0$ be a parameter, fixed in Lemma~\ref{lem:grid} below, and set
\begin{equation}\label{eq:tk}
t_k=\e_0s_k/k .
\end{equation}
In the rest of this section, a \emph{staircase} is an admissible path fixed deterministically in terms of its endpoints. This encompasses  the staircases $I\mapsto\hat z(I,\rho)$ of \eqref{stairdef}, and, for $u\le u'$ in $\Z^2$, the path from $u$ to $u'$ that takes all its $e_1$ steps first. The weight of a staircase $\alpha$ of length $n$ is a sum of $n$ i.i.d.\ weights. We call $\alpha$ \emph{bad} if $\pawgt(\alpha)<n\E\w_0-t_k$, and we write $\Bad(\alpha)$ for this event.

Under the boundedness assumption $\abs{\w_z}\le M$ in \eqref{ass:bounded}, Hoeffding's inequality \eqref{eq:hoeffding} implies that for every staircase $\alpha$ of length $n\le8\cnst_*k^2s_k$,
\begin{equation}\label{eq:badstair}
\P(\Bad(\alpha))\le\pbad_k=\exp\Bigl\{-\frac{t_k^2}{16\cnst_*k^2s_kM^2}\Bigr\}=\exp\Bigl\{-\frac{\e_0^2\,s_k}{16\cnst_*M^2k^4}\Bigr\} .
\end{equation}

All staircases used below have length at most $8\cnst_*j^2w_{j,k}\le8\cnst_*k^2s_k$, by \eqref{t.aux32} and $w_{j,k}\le s_k$. 

\begin{remark}\label{rk:length}
By \eqref{eq:j2w}, the staircases in fact have length at most $16\cnst_*s_k$. We keep the weaker bound $8\cnst_*k^2s_k$ because it is the one used in \eqref{eq:badstair}, in Lemma~\ref{lem:stair-general} and in item~\ref{ch:family} of Appendix~\ref{sub:staircases}, and the factor $k^2$ is harmless there: it only turns $s_k$ into $s_k/k^4$ in the exponent of $\pbad_k$.
\end{remark}

\begin{remark}[Bounded vs.\ unbounded weights]\label{rk:stair-general}
The bound \eqref{eq:badstair} uses the bounded weight assumption \eqref{ass:bounded}. With only \ref{A2} to control $\w_0^-$, a staircase may collect a vertex of weight below $-t_k/\cnst$, for a constant $\cnst$, with probability too large for the union bounds below. Call such vertices \emph{very bad}. The staircases are then replaced by paths that avoid very bad vertices. In particular, for a vertex $u$, the staircase $I\mapsto\hat z(I,\rcoord(u))$ is replaced by a path $\gamma^u$ through $u$ that avoids very bad vertices and stays within a constant transverse distance of $\rcoord(u)$, and the grid points $\hat z(I,\rcoord(u))$ are replaced by the vertices $\gamma^u(I)$ of this path on the grid levels. Figure~\ref{fig:stair-geometry} shows how such a path gets around the very bad vertices. Such a path $\alpha$ is called bad if $\pawgt(\alpha)<-2\abs\alpha\,\E\w_0^--t_k/2$, a lower threshold than the one in the definition of $\Bad(\alpha)$ above. In place of \eqref{eq:badstair} one has
\begin{equation}\label{eq:badstair-general}
\P\bigl\{\Bad(\alpha)\text{ for some path }\alpha\text{ used in the construction}\bigr\}\longrightarrow0\qquad\text{as }k\to\infty .
\end{equation}
The paths are constructed in Appendix~\ref{sub:staircases}, items~\ref{ch:gamma}, \ref{ch:coarse} and~\ref{ch:endpoints}. The family of paths used in the construction is listed in item~\ref{ch:family}, and the precise form of \eqref{eq:badstair-general} is Lemma~\ref{lem:stair-general}, applied to this family. Like the staircases, the paths of item~\ref{ch:family} have length at most $8\cnst_*k^2s_k$ for large $k$.
\end{remark}

\emph{Fine scales.} Suppose scale $j$ is fine, that is, $l_{j,k}<8\cnst_*j^2w_{j,k}$. For $\rho$ in the bin $h$ let $y^\rho_{I:I+l_{j,k}}$ be the staircase $\hat z(I,\rho),\hat z(I+1,\rho),\dots,\hat z(I+l_{j,k},\rho)$ at $\rho$. It connects $\hat z(I,\rho)$ to $\hat z(I+l_{j,k},\rho)$, it has length $l_{j,k}$, its displacement has coordinates $\level=l_{j,k}$ and $\abs\rcoord<1$ (both endpoints have $\rcoord\in(\rho-1,\rho]$), and all its vertices satisfy $\abs{\rcoord-\rho}<1$, by \eqref{stairdef}. In particular $f_\sig(y^\rho_{I+l_{j,k}}-y^\rho_I)\le\ghatz l_{j,k}+\abs{\ghatpr{\sig}}$. On the event $\{\Bad(y^\rho_{I:I+l_{j,k}})\}^c$,
\be\label{eq:finecost}\begin{aligned}
\pawgt(y^\rho_{I:I+l_{j,k}})&\ge l_{j,k}\E\w_0-t_k\ \ge\ f_\sig(y^\rho_{I+l_{j,k}}-y^\rho_I)-(\ghatz-\E\w_0)l_{j,k}-\abs{\ghatpr{\sig}}-t_k
\\
&\ge f_\sig(y^\rho_{I+l_{j,k}}-y^\rho_I)-8\cnst_*(\ghatz-\E\w_0)j^2w_{j,k}-K-t_k ,
\end{aligned}\ee
using $\ghatz\ge\E\w_0$  (from  \ref{S3}) and $\abs{\ghatpr{\sig}}\le K$. 

The staircase $y^\rho_{I:I+l_{j,k}}$ is determined by the points $\hat z(I',\rho)$, $I\le I'\le I+l_{j,k}$. For a fixed $I'$, $\hat z(I',\rho)$ changes only when $\rho$ crosses a point of $\Z-\Dir_1I'$, of which the bin $[hw_{j,k},(h+1)w_{j,k})$ contains exactly $w_{j,k}$. So as $\rho$ ranges over the bin, the staircase changes at most $(l_{j,k}+1)w_{j,k}$ times. Hence the staircases needed at $(j,I,h)$ form a family $\Stairs^{(j)}_{I,h}$ with cardinality 
\begin{equation}\label{eq:finecount}
\abs{\Stairs^{(j)}_{I,h}}\le (l_{j,k}+1)w_{j,k}+1\le9\cnst_*j^2w_{j,k}^2 , 
\end{equation}
 using $l_{j,k}+1\le8\cnst_*j^2w_{j,k}+1$ and $\cnst_*j^2w_{j,k}^2\ge w_{j,k}+1$. 
Define the event 
\begin{equation}\label{eq:goodfine}
\Good^{(j)}_{I,h}=\bigcap_{\alpha\tsp\in\tsp\Stairs^{(j)}_{I,h}}\{\Bad(\alpha)\}^c\qquad\text{when scale $j$ is fine}.
\end{equation}

\begin{remark}[Bounded vs.\ unbounded weights]\label{rk:fine-general}
Under \ref{A2}, for $\rho=\rcoord(u)$ the staircase $y^\rho_{I:I+l_{j,k}}$ is replaced by the segment of $\gamma^u$ between $\gamma^u(I)$ and $\gamma^u(I+l_{j,k})$. Then \eqref{eq:finecost} holds with $\ghatz-\E\w_0$ replaced by $\ghatz+2\E\w_0^-$, with $t_k$ replaced by $t_k/2$, and with a larger constant in place of $K$. There is no family $\Stairs^{(j)}_{I,h}$: the event of \eqref{eq:badstair-general} covers all these paths at once and is added to the good event of Lemma~\ref{lem:grid}, and $\Good^{(j)}_{I,h}$ is the whole space when scale $j$ is fine. See items~\ref{ch:gamma}, \ref{ch:family}, \ref{ch:cost} and~\ref{ch:prob} of Appendix~\ref{sub:staircases}.
\end{remark}

\begin{figure}[htb]
\centering
\begin{tikzpicture}[scale=0.8,>=Stealth]
  %% coordinates: horizontal r, vertical level (not to scale); coarse scale, j^2=3 blocks
  \draw[->,gray] (-2.4,-0.5) -- (-2.4,8.6) node[above] {\footnotesize $\level$};
  \draw[->,gray] (-2.4,-0.5) -- (10.4,-0.5) node[right] {\footnotesize $\rcoord$};
  \foreach \y/\lab in {0/{$I$}, 1.6/{$I+m'_{j,k}$}, 6.4/{$I+l_{j,k}-m'_{j,k}$}, 8/{$I+l_{j,k}$}} {\draw[gray!50] (-2.2,\y) -- (10,\y); \node[right,gray] at (10,\y) {\footnotesize \lab};}
  %% bin on level I
  \draw[very thick] (-1.0,0) -- (0.2,0);
  \node[below] at (1.2,-0.6) {\footnotesize bin $h$: $\rcoord\in[hw_{j,k},(h+1)w_{j,k})$};
  \draw[->,gray,shorten >=1pt] (-0.5,-0.7) -- (-0.5,-0.1);
  \draw[very thick] (-1.0,8) -- (0.2,8);
  %% blocks
  \foreach \x/\n in {2.0/1, 5.0/2, 8.0/3} {
    \fill[gray!20] (\x-0.8,1.6) rectangle (\x+0.8,6.4);
    \draw[gray!60,thin] (\x-0.8,1.6) rectangle (\x+0.8,6.4);
    \fill (\x,1.6) circle (0.06);
    \node[below=1pt] at (\x,1.6) {\footnotesize $z^{(j)}_{I,h,\n}$};
  }
  \draw[<->,gray] (4.2,7.0) -- (5.8,7.0) node[midway,above] {\footnotesize $16w_{j,k}$};
  \draw[<->,gray] (2.0,0.6) -- (5.0,0.6) node[midway,below] {\footnotesize $30w_{j,k}$};
  %% chosen crossing in block 1
  \draw[very thick] (1.7,1.6) .. controls (2.6,3.0) and (1.4,5.0) .. (2.3,6.4);
  \fill (1.7,1.6) circle (0.07); \node[above left=-2pt] at (1.7,1.6) {\footnotesize $u$};
  \fill (2.3,6.4) circle (0.07); \node[above right=-2pt] at (2.3,6.4) {\footnotesize $u'$};
  %% staircases
  \fill (-0.3,0) circle (0.07); \node[above left=-1pt] at (-0.3,0) {\footnotesize $\hat z(I,\rho)$};
  \fill (-0.3,8) circle (0.07); \node[above right=-1pt] at (-0.5,8) {\footnotesize $\hat z(I+l_{j,k},\rho)$};
  \draw[very thick,densely dashed] (-0.3,0) -- (1.7,1.6);
  \draw[very thick,densely dashed] (2.3,6.4) -- (-0.3,8);
  \node[above right=-3pt] at (0.9,7.3) {\footnotesize staircase};
  \node[above left=0pt] at (0.9,0.75) {\footnotesize staircase};
  \node[right] (cr) at (6.0,4.0) {\footnotesize crossing of $\Elev_{z^{(j)}_{I,h,1}}(w_{j,k},l'_{j,k})$};
  \draw[->,gray,shorten >=1pt] (cr) -- (2.1,4.1);
\end{tikzpicture}
\caption{\small The construction at a coarse scale index $j$, for one grid interval $[I,I+l_{j,k}]$ and one bin $h$, drawn in the coordinates $(\rcoord,\level)$ and not to scale, with $j^2=3$. The $j^2$ blocks are centered at the vertices $z^{(j)}_{I,h,\nu}$ on level $I+m'_{j,k}$, whose transverse coordinates are $(hw_{j,k}+30\nu w_{j,k})$ up to $\tfrac12$, and they occupy the levels from $I+m'_{j,k}$ to $I+l_{j,k}-m'_{j,k}$. Each hosts the crossing event $\Elev_{z^{(j)}_{I,h,\nu}}\!\!(w_{j,k},l'_{j,k})$, and the events are independent because the blocks are disjoint. If some $\Elev_{z^{(j)}_{I,h,\nu}}\!\!(w_{j,k},l'_{j,k})$ occurs (here $\nu=1$), its crossing from $u$ to $u'$ (thick) is attached to $\hat z(I,\rho)$ and to $\hat z(I+l_{j,k},\rho)$ by two deterministic staircases of length $m'_{j,k}=\lceil\cnst_*j^2w_{j,k}\rceil$ (dashed), for every $\rho$ in the bin. The blocks lie to the $e_1$ side of the bin, at transverse distance at most $39j^2w_{j,k}\le\theta_0m'_{j,k}$, which is what makes the staircases admissible.}
\label{fig:grid}
\end{figure}
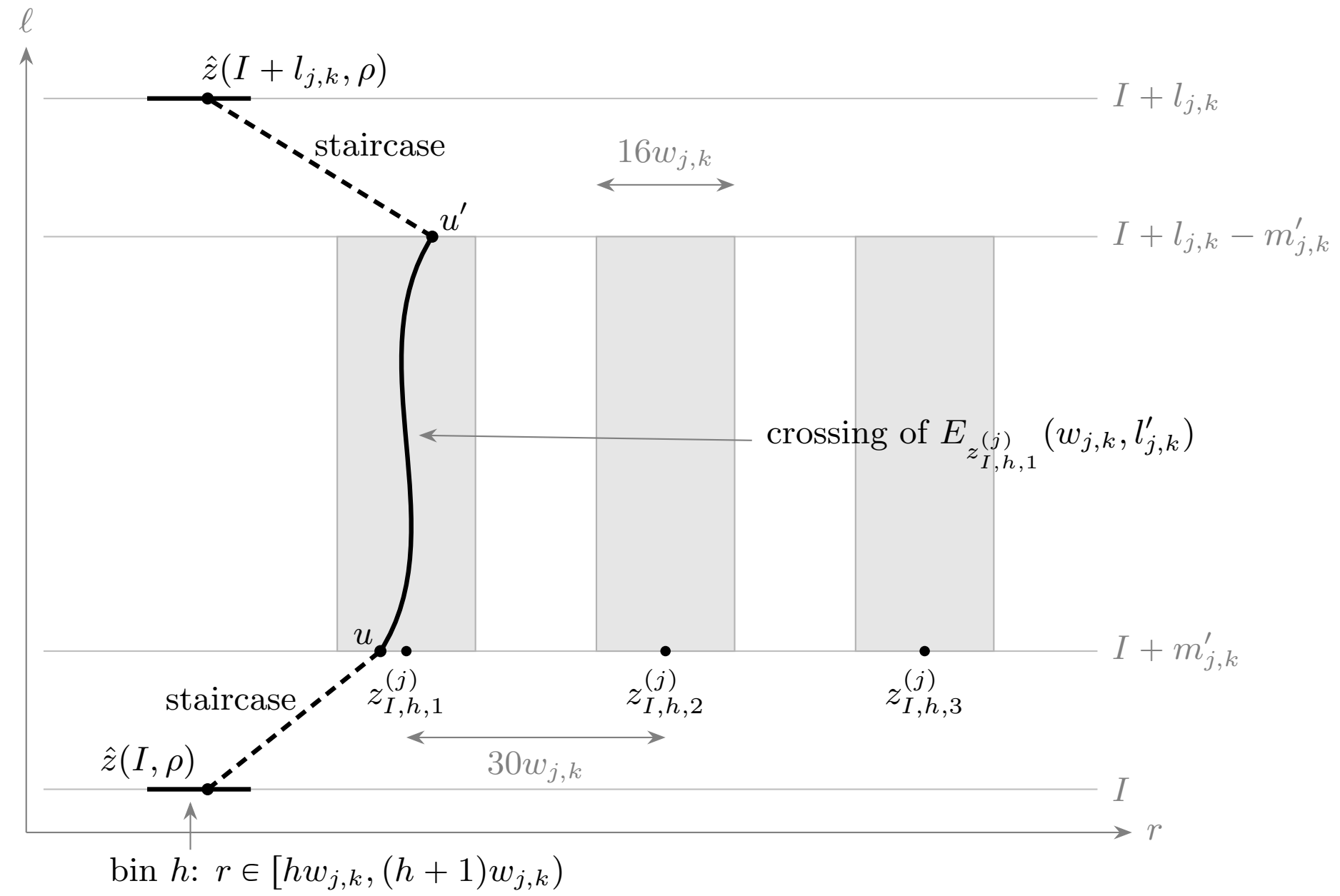
\emph{Coarse scales: crossing events.}  Suppose scale  $j$ is coarse, that is, $l_{j,k}\ge8\cnst_*j^2w_{j,k}$. See Figure~\ref{fig:grid}. Set
\begin{equation}\label{eq:mprime}
m'_{j,k}=\lceil\cnst_*j^2w_{j,k}\rceil\qquad\text{and}\qquad l'_{j,k}=l_{j,k}-2m'_{j,k} .
\end{equation}
Since $\cnst_*j^2w_{j,k}\ge1$ and $l_{j,k}\ge8\cnst_*j^2w_{j,k}\ge8$, we have
\[m'_{j,k}\le2\cnst_*j^2w_{j,k}\le l_{j,k}/4
\quad\text{and}\quad l_{j,k}/2\le l'_{j,k}<l_{j,k}.\]
For an index $\nu\in\{1,\dots,j^2\}$ let $z^{(j)}_{I,h,\nu}$ be a vertex, chosen deterministically (independently of $\w$), with
\begin{equation}\label{eq:rhonu}
\level(z^{(j)}_{I,h,\nu})=I+m'_{j,k}\qquad\text{and}\qquad\abs{\tsp\rcoord(z^{(j)}_{I,h,\nu})-(hw_{j,k}+30\nu w_{j,k})\tsp}\le\tfrac12 ,
\end{equation}
one of which exists by \eqref{eq:coords}. Note the transverse offset $30\nu w_{j,k}$ of these vertices from the bin $[hw_{j,k},(h+1)w_{j,k})$: the blocks \eqref{eq:crossblock} defined below lie to the $e_1$ side of the bin, which is what makes the staircases attaching a crossing to $\hat z(I,\rho)$ and to $\hat z(I+l_{j,k},\rho)$ admissible, as shown in the paragraph on attachments below.  

Recall that $\Elev_{z^{(j)}_{I,h,\nu}}\!\!(w_{j,k},l'_{j,k})$ is the event that there exists a crossing $x_{I+m'_{j,k}:I+l_{j,k}-m'_{j,k}}$ of the block
\begin{equation}\label{eq:crossblock}
\{z:\ I+m'_{j,k}\le\level(z)\le I+l_{j,k}-m'_{j,k},\ \abs{\rcoord(z)-\rcoord(z^{(j)}_{I,h,\nu})}\le8w_{j,k}\}
\end{equation}
with
\begin{equation}\label{success55}
\pawgt(x_{I+m'_{j,k}:I+l_{j,k}-m'_{j,k}})>\ghatz l'_{j,k}-Kw_{j,k} .
\end{equation}
Consecutive centers $\rcoord(z^{(j)}_{I,h,\nu})$ are at least $30w_{j,k}-1$ apart and each block has half-width $8w_{j,k}$, so the $j^2$ blocks are pairwise disjoint and the events $\Elev_{z^{(j)}_{I,h,\nu}}\!\!(w_{j,k},l'_{j,k})$, $1\le\nu\le j^2$, are independent. Since $1\le l'_{j,k}<l_{j,k}$, each has probability greater than $q_0$ by \eqref{eq:trialprob}. Hence
\begin{equation}\label{eq:allfail}
\P\Bigl(\,\bigcap_{\nu=1}^{j^2}\bigl\{\Elev_{z^{(j)}_{I,h,\nu}}\!\!(w_{j,k},l'_{j,k})\bigr\}^c\Bigr)\le(1-q_0)^{j^2}\le e^{-q_0j^2}.
\end{equation}

\begin{remark}[Bounded vs.\ unbounded weights]\label{rk:heavy-general}
Although \eqref{ass:bounded} is not used in the above paragraph,  under \ref{A2}--\ref{A3} the crossing needs to be trimmed by a constant number of vertices at each end, as explained in Remark~\ref{rk:coarse-general} below. The weights removed then have to be bounded above. For this, let $H^{(j)}_{I,h,\nu}$ be the event that no vertex of the block \eqref{eq:crossblock} has weight above $w_{j,k}$, and replace $\Elev_{z^{(j)}_{I,h,\nu}}\!\!(w_{j,k},l'_{j,k})$ by $\Elev_{z^{(j)}_{I,h,\nu}}\!\!(w_{j,k},l'_{j,k})\cap H^{(j)}_{I,h,\nu}$. By \ref{A3} and a union bound, these events still have probability at least $q_0/2$ for large $k$, and they are still independent over $\nu$, so in place of \eqref{eq:allfail} one has
\begin{equation}\label{eq:allfail-general}
\P\Bigl(\,\bigcap_{\nu=1}^{j^2}\bigl\{\Elev_{z^{(j)}_{I,h,\nu}}\!\!(w_{j,k},l'_{j,k})\cap H^{(j)}_{I,h,\nu}\bigr\}^c\Bigr)\le e^{-q_0j^2/2}\qquad\text{for large }k .
\end{equation}
See Lemma~\ref{lem:heavy} and item~\ref{ch:heavy} of Appendix~\ref{sub:staircases}.
\end{remark}

\medskip
\emph{Coarse scales: attachments.} Suppose $\Elev_{z^{(j)}_{I,h,\nu}}\!\!(w_{j,k},l'_{j,k})$ occurs with crossing $x_{I+m'_{j,k}:I+l_{j,k}-m'_{j,k}}$. Label its endpoints as    $u=x_{I+m'_{j,k}}$ and $u'=x_{I+l_{j,k}-m'_{j,k}}$. Both satisfy $\abs{\rcoord-(hw_{j,k}+30\nu w_{j,k})}\le8w_{j,k}+\tfrac12$. For $\rho$ in the bin $h$, the point $\hat z(I,\rho)$ has $\rcoord\in(\rho-1,\rho]\subset(hw_{j,k}-1,(h+1)w_{j,k})$, so the displacement $u-\hat z(I,\rho)$ has $\level=m'_{j,k}$ and
\[
\rcoord\bigl(u-\hat z(I,\rho)\bigr)\in\bigl[(30\nu-9)w_{j,k}-\tfrac12,\ (30\nu+8)w_{j,k}+\tfrac32\bigr]\subset\bigl(0,\ 39j^2w_{j,k}\bigr]\subset(0,\theta_0m'_{j,k}],
\]
because $w_{j,k}\ge2$, %k-j\ge k_0\ge1
$\nu\le j^2$ and $\theta_0m'_{j,k}\ge\theta_0\cnst_*j^2w_{j,k}=50j^2w_{j,k}$. By Section~\ref{sub:coords} (below \eqref{eq:lip-r}), a vector with level $m'_{j,k}>0$ and $0<\rcoord\le\theta_0m'_{j,k}\le(1-\Dir_1)m'_{j,k}$ lies in the closed first quadrant, so $\hat z(I,\rho)\le u$ and the staircase from $\hat z(I,\rho)$ to $u$ exists. Likewise $\hat z(I+l_{j,k},\rho)-u'$ has level $m'_{j,k}$ and $\rcoord\in[-39j^2w_{j,k},0)$, so the staircase from $u'$ to $\hat z(I+l_{j,k},\rho)$ exists.

Let $y^\rho_{I:I+l_{j,k}}$ be the concatenation of the staircase from $\hat z(I,\rho)$ to $u$, the crossing $x_{I+m'_{j,k}:I+l_{j,k}-m'_{j,k}}$, and the staircase from $u'$ to $\hat z(I+l_{j,k},\rho)$. Its displacement has $\level=l_{j,k}$ and $\abs\rcoord<1$. If neither staircase is bad, then by \eqref{success55}
\[
\pawgt(y^\rho_{I:I+l_{j,k}})\ge(m'_{j,k}\E\w_0-t_k)+(\ghatz l'_{j,k}-Kw_{j,k})+(m'_{j,k}\E\w_0-t_k),
\]
and hence, using $f_\sig(y^\rho_{I+l_{j,k}}-y^\rho_I)\le\ghatz l_{j,k}+\abs{\ghatpr{\sig}}$ and $l_{j,k}-l'_{j,k}=2m'_{j,k}\le4\cnst_*j^2w_{j,k}$,
\be\label{eq:coarsecost}
\begin{aligned}
\pawgt(y^\rho_{I:I+l_{j,k}})&\ge f_\sig(y^\rho_{I+l_{j,k}}-y^\rho_I)-\ghatz l_{j,k}-\abs{\ghatpr{\sig}}+\ghatz l'_{j,k}-Kw_{j,k}+2m'_{j,k}\E\w_0-2t_k\\[3pt]
&=f_\sig(y^\rho_{I+l_{j,k}}-y^\rho_I)-2m'_{j,k}(\ghatz-\E\w_0)-Kw_{j,k}-\abs{\ghatpr{\sig}}-2t_k\\[3pt]
&\ge f_\sig(y^\rho_{I+l_{j,k}}-y^\rho_I)-4\cnst_*(\ghatz-\E\w_0)j^2w_{j,k}-2Kj^2w_{j,k}-2t_k. 
\end{aligned}\ee
In the last step $\abs{\ghatpr{\sig}}\le K\le Kj^2w_{j,k}$, as $j\ge1$ and $w_{j,k}\ge2$.

As for the transverse location: the vertices of the crossing have $\abs{\rcoord-(hw_{j,k}+30\nu w_{j,k})}\le8w_{j,k}+\tfrac12$ and, since $\rho\in[hw_{j,k},(h+1)w_{j,k})$ and $\nu\le j^2$, $\abs{(hw_{j,k}+30\nu w_{j,k})-\rho}\le30j^2w_{j,k}$.   Together these give that  vertices of the crossing satisfy \[ \abs{\rcoord-\rho}\le(30j^2+9)w_{j,k}. \]
The vertices of a staircase of length $m'_{j,k}$ starting at $\hat z(I,\rho)$, or ending at $\hat z(I+l_{j,k},\rho)$, are within $m'_{j,k}$ steps of a point with $\rcoord\in(\rho-1,\rho]$, so they have $\abs{\rcoord-\rho}\le m'_{j,k}+1\le2\cnst_*j^2w_{j,k}+1$. Since $(30j^2+9)w_{j,k}\le39j^2w_{j,k}\le\cnst_*j^2w_{j,k}$, all vertices of $y^\rho_{I:I+l_{j,k}}$ satisfy the transverse bound 
\begin{equation}\label{eq:coarseloc}
\abs{\rcoord-\rho}\le2\cnst_*j^2w_{j,k}+1\le3\cnst_*j^2w_{j,k} .
\end{equation}

\begin{remark}[Bounded vs.\ unbounded weights]\label{rk:coarse-general}
As explained in Remark~\ref{rk:stair-general}, under \ref{A2} the staircases are replaced by paths that avoid very bad vertices. Such a path cannot always be attached to an endpoint of the crossing, because there may be very bad vertices around  that endpoint. The crossing is therefore trimmed: its start is moved forward along it, and its end backward, by at most a constant number of vertices, to vertices at which these paths can be attached. For $\rho=\rcoord(u)$ the two staircases are then replaced by paths that avoid the very bad vertices, from $\gamma^u(I)$ to the trimmed crossing and from the trimmed crossing back to $\gamma^u(I+l_{j,k})$. The first of these is an arrival path, as in Figure~\ref{fig:stair-geometry}. On the event $H^{(j)}_{I,h,\nu}$ of Remark~\ref{rk:heavy-general} the weights of the removed vertices are bounded above. Then \eqref{eq:coarsecost} and \eqref{eq:coarseloc} hold with $\ghatz-\E\w_0$ replaced by $\ghatz+2\E\w_0^-$, with $2t_k$ replaced by $3t_k/2$, and with larger constants. See items~\ref{ch:trim}, \ref{ch:coarse}, and~\ref{ch:cost} of Appendix~\ref{sub:staircases}.
\end{remark}

The staircases that may be needed at $(j,I,h)$ are determined by a starting point $\hat z(I,\rho)$ with $\rho$ in the bin (a lattice point on level $I$ with $\rcoord\in(hw_{j,k}-1,(h+1)w_{j,k})$: at most $w_{j,k}+1$ choices) and a vertex $u$ on level $I+m'_{j,k}$ in one of the $j^2$ blocks (at most $16w_{j,k}+2$ choices per block), and symmetrically on the exit side. They form a family $\Stairs^{(j)}_{I,h}$ with cardinality  (recall $w_{j,k}\ge2$)  
\begin{equation}\label{eq:staircount}
\abs{\Stairs^{(j)}_{I,h}}\le 
2(w_{j,k}+1)j^2(16w_{j,k}+2)\le60j^2w_{j,k}^2. 
\end{equation}
%16w^2+18w+2\le16w^2+9w^2+w^2=26w^2\le30w^2
 The good event at $(j,I,h)$ is
\begin{equation}\label{eq:goodcoarse}
\Good^{(j)}_{I,h}=\Bigl(\,\bigcup_{\nu=1}^{j^2}\Elev_{z^{(j)}_{I,h,\nu}}\!\!(w_{j,k},l'_{j,k})\Bigr)\cap\bigcap_{\alpha\tsp\in\tsp\Stairs^{(j)}_{I,h}}\Bad(\alpha)^c\qquad\text{when scale $j$ is coarse}.
\end{equation}
In summary, by \eqref{eq:badstair}, \eqref{eq:finecount}, \eqref{eq:goodfine}, \eqref{eq:allfail}, and \eqref{eq:staircount}, for $k\in\cK$ satisfying \eqref{eq:kreq}, $J\le j\le j_*(k)$, $I\in\cG_{j,k}$, and $h\in\cH_{j,k}$,  
\begin{equation}\label{eq:goodprob}
\P\bigl(\bigl\{\Good^{(j)}_{I,h}\bigr\}^c\tspb\bigr)\le
\begin{cases}
e^{-q_0j^2}+60j^2w_{j,k}^2\,\pbad_k&\text{if scale $j$ is coarse,}\\[3pt]
9\cnst_*j^2w_{j,k}^2\,\pbad_k&\text{if scale $j$ is fine.}
\end{cases}
\end{equation}

\begin{remark}[Bounded vs.\ unbounded weights]\label{rk:good-general}
Under \ref{A2}--\ref{A3} there are no families $\Stairs^{(j)}_{I,h}$, at the coarse scales as at the fine ones: the weights of all the paths that replace the staircases are controlled at once by the event in \eqref{eq:badstair-general}, which is added to the good event of Lemma~\ref{lem:grid}. Accordingly, with the events $H^{(j)}_{I,h,\nu}$ of Remark~\ref{rk:heavy-general}, \eqref{eq:goodcoarse} is replaced by 
\[
\Good^{(j)}_{I,h}=\bigcup_{\nu=1}^{j^2}\Bigl(\Elev_{z^{(j)}_{I,h,\nu}}\!\!(w_{j,k},l'_{j,k})\cap H^{(j)}_{I,h,\nu}\Bigr)\qquad\text{when scale $j$ is coarse},
\]
while $\Good^{(j)}_{I,h}$ is the whole space when scale $j$ is fine. By \eqref{eq:allfail-general}, in place of \eqref{eq:goodprob}, we have, for large $k$,
\[
\P\bigl(\bigl\{\Good^{(j)}_{I,h}\bigr\}^c\tspb\bigr)\le
\begin{cases}
e^{-q_0j^2/2}&\text{if scale $j$ is coarse,}\\[3pt]
0&\text{if scale $j$ is fine.}
\end{cases}
\]
See item~\ref{ch:prob} of Appendix~\ref{sub:staircases}. \end{remark}

\medskip

\subsection{All scales simultaneously}\label{sub:allscales}
The event on which all the connections of Section~\ref{sub:onescale} are simultaneously  available  is
\begin{equation}\label{eq:Gdef}
G_{k,J}=\bigcap_{j=J}^{j_*(k)}\ \bigcap_{I\tsp\in\tsp\cG_{j,k}}\ \bigcap_{h\tsp\in\tsp\cH_{j,k}}\Good^{(j)}_{I,h}.
\end{equation}

\begin{lemma}\label{lem:grid}
Fix $A,R>0$, $\e_0>0$, and $J\ge1$. There exists a constant $\cnst_1$, depending on $\mu$, $\Dir$, and $K$ only, such that the following hold.
\begin{enumerate}[label={\rm(\roman*)}]
\item\label{lem:grid:path} On $G_{k,J}$, for every scale index $J\le j\le j_*(k)$, every $I\in\cG_{j,k}$ and every real
$\rho\in[-Rs_k, Rs_k]$, there exists an admissible path $y^\rho_{I:I+l_{j,k}}$ from $\hat z(I,\rho)$ to $\hat z(I+l_{j,k},\rho)$ that satisfies 
\begin{equation}\label{eq:gridcost} 
\begin{aligned}
\pawgt(y^\rho_{I:I+l_{j,k}}) &\ge f_\sig(y^\rho_{I+l_{j,k}}-y^\rho_I)-\cnst_1j^2w_{j,k}-2t_k \quad\text{for } \ \sigg=\pm  \\[2pt] \text{and}\qquad
\abs{\rcoord(y^\rho_i)-\rho}&\le\cnst_1j^2w_{j,k}\quad\text{for } \ I\le i\le I+l_{j,k}.
\end{aligned}
\end{equation}
\item\label{lem:grid:prob} There exists a sequence $\beta_J\to0$ as $J\to\infty$, depending on $A,R,q_0$ only, and for each $J$ there exists $k_1=k_1(J)$, depending also on $A$, $R$, $\e_0$, $\mu$, and $\Dir$, such that $\P(G_{k,J}^c)\le2\beta_J$ for all $k\in\cK$ with $k\ge k_1$.
%$k_1(J)$ depends on $J$ through $\beta_J$: the staircase contribution \eqref{eq:stairterm} does not depend on $J$, and $k_1(J)$ is the index from which on it is at most $\beta_J$. Since $\beta_J\downarrow0$, $k_1(J)$ grows with $J$, which is why $k\to\infty$ is taken before $J\to\infty$ in Proposition~\ref{prop:endpoints}.
\end{enumerate}
\end{lemma}

\begin{remark}[Bounded vs.\ unbounded weights]\label{rk:grid-general}
Under \ref{A2}--\ref{A3}, Lemma~\ref{lem:grid} holds with a different constant $\cnst_1$, a larger $k_1(J)$, and one change in part~\ref{lem:grid:path}: for $\rho=\rcoord(u)$, the path connects $\gamma^u(I)$ to $\gamma^u(I+l_{j,k})$ instead of $\hat z(I,\rho)$ to $\hat z(I+l_{j,k},\rho)$, as in Remarks~\ref{rk:fine-general} and~\ref{rk:coarse-general}. This form of part~\ref{lem:grid:path} is asserted only for the vertices $u$ and the grid intervals that occur in the proof of Proposition~\ref{prop:endpoints}, which are specified in Remark~\ref{rk:endpoints-general}. Nothing after this lemma depends on the value of $\cnst_1$. 
\end{remark}

\begin{proof}[Proof of Lemma~\ref{lem:grid}]
\ref{lem:grid:path} Fix $j$, $I$, and $\rho$, and let $h$ be the index of the bin containing $\rho$. On $\Good^{(j)}_{I,h}$ the path $y^\rho_{I:I+l_{j,k}}$ constructed in Section~\ref{sub:onescale} exists, with weight bounded below by \eqref{eq:finecost} at a fine scale $j$ and by \eqref{eq:coarsecost} at a coarse scale $j$. Using $K\le2Kj^2w_{j,k}$ and $t_k\le2t_k$ in \eqref{eq:finecost}, and $4\cnst_*\le8\cnst_*$ in \eqref{eq:coarsecost}, in both cases
\begin{align*}
     \pawgt(y^\rho_{I:I+l_{j,k}}) &\ge f_\sig(y^\rho_{I+l_{j,k}}-y^\rho_I)-8\cnst_*(\ghatz-\E\w_0)j^2w_{j,k}-2Kj^2w_{j,k}-2t_k \\
     &\ge f_\sig(y^\rho_{I+l_{j,k}}-y^\rho_I)-\cnst_1j^2w_{j,k}-2t_k 
     \end{align*}
with
\begin{equation}\label{eq:Ccone}
\cnst_1=8\cnst_*(\ghatz-\E\w_0)+2K+3\cnst_* .
\end{equation}
The location bound is $\abs{\rcoord-\rho}<1\le\cnst_1j^2w_{j,k}$ at a fine scale and \eqref{eq:coarseloc} at a coarse scale, both at most $\cnst_1j^2w_{j,k}$ because of the extra term $3\cnst_*$ in \eqref{eq:Ccone}.  

\smallskip

\ref{lem:grid:prob} By \eqref{eq:Gdef} and a union bound,
\begin{equation}\label{eq:Gunion}
\P(G_{k,J}^c)\le\sum_{j=J}^{j_*(k)}\abs{\cG_{j,k}}\,\abs{\cH_{j,k}}\max_{I\tsp\in\tsp\cG_{j,k},\,h\tsp\in\tsp\cH_{j,k}}\P\bigl((\Good^{(j)}_{I,h})^c\bigr).
\end{equation}
By \eqref{eq:gridcount} and \eqref{eq:bincount}, $\abs{\cG_{j,k}}\abs{\cH_{j,k}}\le\cnst_{\!A,R}2^{5j}$. Insert \eqref{eq:goodprob} and separate the contribution of the failed crossing events, which comes from the coarse scales only, from the contribution of the bad staircases:
\be\label{eq:Gsplit}\begin{aligned}
\P(G_{k,J}^c)&\le\cnst_{\!A,R}\sum_{\substack{J\le j\le j_*(k)\\ j\text{ coarse}}}2^{5j}e^{-q_0j^2}\\
&\qquad +\cnst_{\!A,R}\,\pbad_k\biggl\{\,\sum_{\substack{J\le j\le j_*(k)\\ j\text{ coarse}}}60j^22^{5j}w_{j,k}^2
+\sum_{\substack{J\le j\le j_*(k)\\ j\text{ fine}}}9\cnst_*j^22^{5j}w_{j,k}^2\biggr\}.
\end{aligned}\ee

The first sum is at most
\begin{equation}\label{eq:betaJ}
\beta_J=\cnst_{\!A,R}\sum_{j\ge J}2^{5j}e^{-q_0j^2},
\end{equation}
which is finite, tends to $0$ as $J\to\infty$, and only depends on $A$, $R$, and $q_0$.

For the sums inside the braces in \eqref{eq:Gsplit}, we bound the powers of two in terms of $k$. At a coarse scale, $2^{m_k-4j-1}=l_{j,k}\ge8\cnst_*j^2w_{j,k}\ge2^{k-j+3}$ gives $3j\le m_k-k-4$. Furthermore,  $j\le k$ by \eqref{t.aux32} and $m_k\le 3k$ by \eqref{eq:Lbounds}. Terms in the first sum in braces satisfy 
\begin{equation}\label{eq:coarseterm}
60j^22^{5j}w_{j,k}^2=60j^2\,2^{3j+2k}\le60k^2\,2^{m_k+k}\le60k^2\,2^{4k}=60k^2s_k^4 .
\end{equation}
At every scale, $4j\le4j_*(k)\le m_k+3\le3k+3$ by \eqref{eq:jstar}, so $3j\le\tfrac34(3k+3)\le3k$ as $k\ge3$, and terms in the second sum in braces satisfy
\begin{equation}\label{eq:fineterm}
9\cnst_*j^22^{5j}w_{j,k}^2=9\cnst_*j^2\,2^{3j+2k}\le9\cnst_*k^2\,2^{5k}=9\cnst_*k^2s_k^5 .
\end{equation}
There are at most $j_*(k)\le k$ scales, so the two sums inside the braces in \eqref{eq:Gsplit} add up to at most $(60+9\cnst_*)k^3s_k^5$, and by \eqref{eq:badstair} and $s_k=2^k$,
\begin{equation}\label{eq:stairterm}
\cnst_{\!A,R}\,\pbad_k\,(60+9\cnst_*)k^3s_k^5
=(60+9\cnst_*)\,\cnst_{\!A,R}\,k^3s_k^5\exp\Bigl\{-\frac{\e_0^2}{16\cnst_*M^2}\cdot\frac{s_k}{k^4}\Bigr\}\longrightarrow0\qquad\text{as }k\to\infty .
\end{equation}

The left-hand side of \eqref{eq:stairterm} does not depend on $J$. Let $k_1=k_1(J)$ be such that it is at most $\beta_J$ for all $k\ge k_1$. Then \eqref{eq:Gsplit} gives $\P(G_{k,J}^c)\le2\beta_J$ for all $k\in\cK$ with $k\ge k_1$.
\end{proof}

\begin{remark}[Bounded vs.\ unbounded weights]\label{rk:gridprob-general}
Under \ref{A2}--\ref{A3}, in the proof of Lemma~\ref{lem:grid}\ref{lem:grid:prob} the good event $G_{k,J}$ is intersected with the event in \eqref{eq:badstair-general} and with two further events of probability tending to one, from Lemmas~\ref{lem:stair-geometry} and~\ref{lem:heavy}. By Remark~\ref{rk:good-general}, the terms with $\pbad_k$ in \eqref{eq:Gsplit} disappear and $\beta_J$ in \eqref{eq:betaJ} has $e^{-q_0j^2/2}$ in place of $e^{-q_0j^2}$, while the three added events cost an additional $o(1)$, so part~\ref{lem:grid:prob} of the lemma holds as stated. See items~\ref{ch:cost} and~\ref{ch:prob} of Appendix~\ref{sub:staircases}.
\end{remark}

\subsection{Arbitrary endpoints}
Define the region
\[\cR_k=\{z:\abs{\level(z)}\le A\clen_k,\ \abs{\rcoord(z)}\le Rs_k\}.\]

\begin{proposition}\label{prop:endpoints}
For every fixed $A,R,\e>0$, as $k\to\infty$ in $\cK$,
\be\label{eq:endpoints}\begin{aligned}
&\P\bigl\{\text{there exists an admissible path $x_{m:n}$ with all vertices in }\cR_k\\
&\qquad 
\text{and such that }\pawgt(x_{m:n})>f_\sig(x_n-x_m)+\e s_k\text{ for some sign }\sigg\in\{-,+\}\bigr\}\longrightarrow0. 
\end{aligned} \ee
\end{proposition}

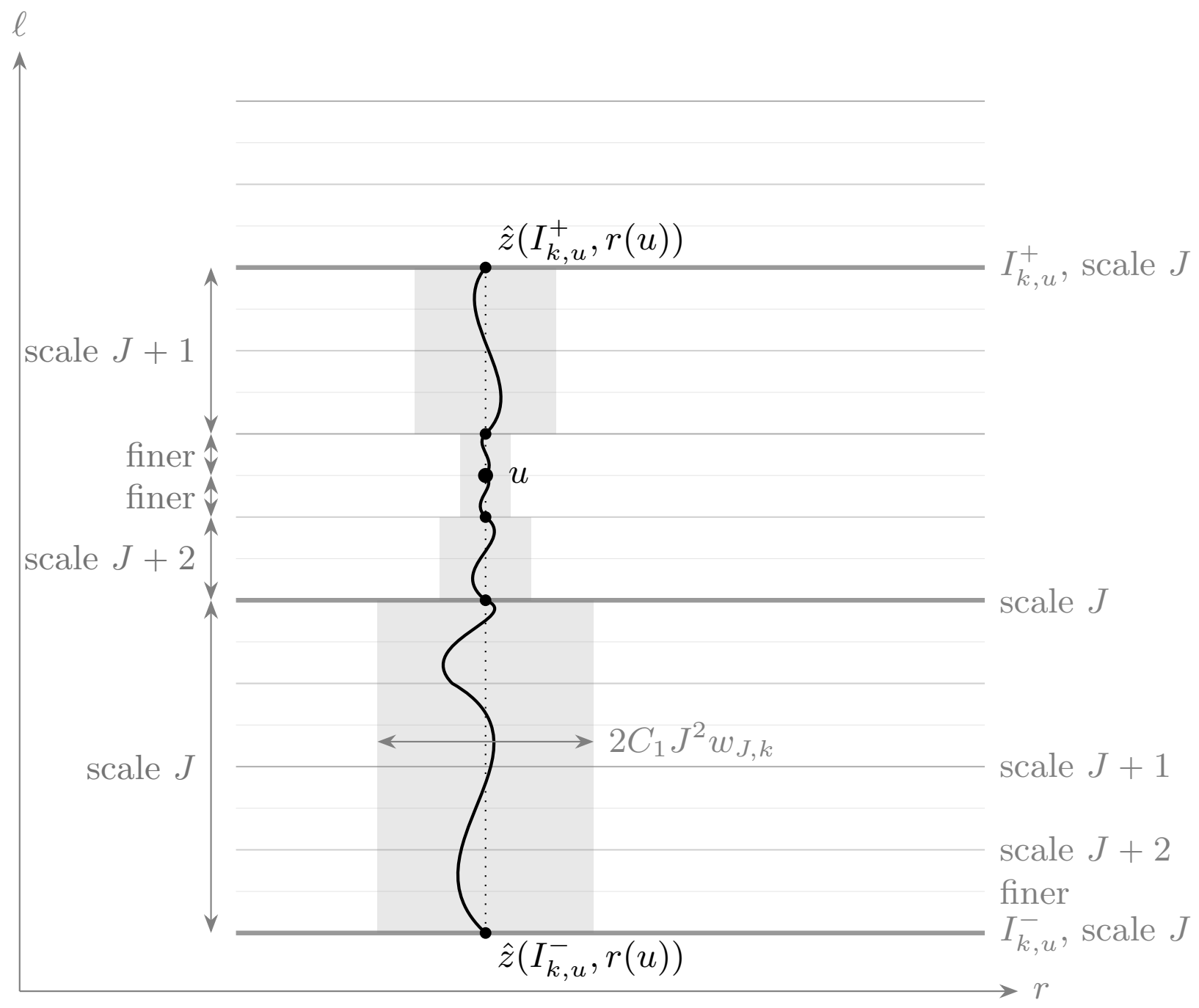
\begin{figure}[htb]
\centering
\begin{tikzpicture}[scale=0.75,>=Stealth]
  %% coordinates: horizontal r, vertical level; nested grids (not to scale)
  \draw[->,gray] (-5.6,-0.7) -- (-5.6,10.6) node[above] {\footnotesize $\level$};
  \draw[->,gray] (-5.6,-0.7) -- (6.4,-0.7) node[right] {\footnotesize $\rcoord$};
  %% scale J lines (thick), J+1 (medium), J+2 (thin)
  \foreach \y in {0,4,8} \draw[gray!80,very thick] (-3,\y) -- (6,\y);
  \foreach \y in {2,6,10} \draw[gray!60] (-3,\y) -- (6,\y);
  \foreach \y in {1,3,5,7,9} \draw[gray!35,thin] (-3,\y) -- (6,\y);
  \foreach \y in {0.5,1.5,2.5,3.5,4.5,5.5,6.5,7.5,8.5,9.5} \draw[gray!20,very thin] (-3,\y) -- (6,\y);
  \node[right,gray] at (6,0) {\footnotesize $I_{k,u}^-$, scale $J$};
  \node[right,gray] at (6,8) {\footnotesize $I_{k,u}^+$, scale $J$};
  \node[right,gray] at (6,4) {\footnotesize scale $J$};
  \node[right,gray] at (6,2) {\footnotesize scale $J+1$};
  \node[right,gray] at (6,1) {\footnotesize scale $J+2$};
  \node[right,gray] at (6,0.5) {\footnotesize finer};
  %% shaded boxes: the piece over a grid interval at scale j stays within transverse distance C_1 j^2 w_{j,k} of the dotted line
  \fill[gray,fill opacity=0.18] (-1.3,0) rectangle (1.3,4);
  \fill[gray,fill opacity=0.18] (-0.55,4) rectangle (0.55,5);
  \fill[gray,fill opacity=0.18] (-0.3,5) rectangle (0.3,5.5);
  \fill[gray,fill opacity=0.18] (-0.3,5.5) rectangle (0.3,6);
  \fill[gray,fill opacity=0.18] (-0.85,6) rectangle (0.85,8);
  %% vertex u at level 5.5
  \fill (0,5.5) circle (0.09); \node[right=2pt] at (0,5.5) {\footnotesize $u$};
  %% entrance path: [0,4] at scale J, [4,5] at scale J+2 (mesh 1), [5,5.5] finer
  \draw[thick] (0,0) .. controls (-1.1,1.0) and (1.0,2.2) .. (-0.4,3.0) .. controls (-0.9,3.5) and (0.5,3.8) .. (0,4);
  \draw[thick] (0,4) .. controls (-0.5,4.4) and (0.4,4.7) .. (0,5);
  \draw[thick] (0,5) .. controls (-0.2,5.2) and (0.15,5.35) .. (0,5.5);
  %% exit path: [5.5,6] finer, [6,8] at scale J+1
  \draw[thick] (0,5.5) .. controls (0.15,5.7) and (-0.15,5.85) .. (0,6);
  \draw[thick] (0,6) .. controls (0.6,6.6) and (-0.5,7.4) .. (0,8);
  \fill (0,0) circle (0.07); \fill (0,4) circle (0.07); \fill (0,5) circle (0.07); \fill (0,6) circle (0.07); \fill (0,8) circle (0.07);
  \node[below right=-3pt] at (0.1,0) {\footnotesize $\hat z(I_{k,u}^-,\rcoord(u))$};
  \node[above right=-3pt] at (0.1,8) {\footnotesize $\hat z(I_{k,u}^+,\rcoord(u))$};
  \draw[dotted] (0,0) -- (0,8);
  %% width annotations
  \draw[<->,gray] (-3.3,0) -- (-3.3,4); \node[left,gray!90!black] at (-3.3,2) {\footnotesize scale $J$};
  \draw[<->,gray] (-3.3,4) -- (-3.3,5); \node[left,gray!90!black] at (-3.3,4.5) {\footnotesize scale $J+2$};
  \draw[<->,gray] (-3.3,5) -- (-3.3,5.5); \node[left,gray!90!black] at (-3.3,5.25) {\footnotesize finer};
  \draw[<->,gray] (-3.3,6) -- (-3.3,8); \node[left,gray!90!black] at (-3.3,7) {\footnotesize scale $J+1$};
  \draw[<->,gray] (-3.3,5.5) -- (-3.3,6); \node[left,gray!90!black] at (-3.3,5.75) {\footnotesize finer};
  \draw[<->,gray] (-1.3,2.3) -- (1.3,2.3) node[pos=1,right] {\footnotesize $2\cnst_1J^2w_{J,k}$};
\end{tikzpicture}
\caption{\small The entrance and exit paths of a vertex $u$, drawn in the coordinates $(\rcoord,\level)$ and not to scale. The horizontal lines are the grid levels $l_{j,k}\Z$ at the scales $J$ (thick), $J+1$, $J+2$, and finer. The level interval $[I_{k,u}^-,\level(u)]$ is decomposed into one grid interval at scale $J$ followed by at most $15$ grid intervals at each finer scale (the brackets on the left show the scales of the pieces), and $[\level(u),I_{k,u}^+]$ symmetrically from the top. Along each grid interval the path of Lemma~\ref{lem:grid}(i) at $\rho=\rcoord(u)$ is used, so that consecutive pieces share the endpoints $\hat z(I,\rcoord(u))$ on the dotted line and the last piece ends exactly at $u=\hat z(\level(u),\rcoord(u))$. The pieces at scale $j$ stay within transverse distance $\cnst_1j^2w_{j,k}$ of the dotted line, that is, inside the shaded box of width $2\cnst_1j^2w_{j,k}$ over their grid interval, and have weight at least $f_{\protect\sig}-\cnst_1j^2w_{j,k}-2t_k$ for both signs; both quantities decrease geometrically in $j$.}
\label{fig:entrance}
\end{figure}

\gist{A bad path in the sense of $\pawgt>f_\sig+\e s_k$ cannot connect two coarse grid levels because (i) $J$ is fixed so there is a bounded number of them and (ii) such bad paths between fixed levels were ruled out by Proposition~\ref{prop:rarefast} (applies because coarse blocks have length $\clen_k 2^{-4j}\gg s_k$). Given an arbitrarily positioned bad path $x_{m:n}$, augment it with an initial segment from a $J$-grid level to $x_m$: Go down from $x_m$ to the first coarse level and then one more. Start from there. This guarantees the coarse piece [hence also the total patched up path] is of length $\ge\clen_k 2^{-4j}\gg s_k$. Use Lemma \ref{lem:grid} paths at different scales down to the finest scale $1$. Similarly, augment the path with a final segment from $x_n$ to a $J$-grid level. The Lemma \ref{lem:grid} upper bound \eqref{eq:gridcost} on $f_\sig-\pawgt$ is small enough that the augmented path still has $W>f_\sig+\e s_k/2$. But now it connects $J$-grid levels so Proposition~\ref{prop:rarefast} kills it.} 

\gist{The way Lemma~\ref{lem:grid} works: the fine scales are made of deterministic staircases and, even though these paths have a very low passage time, they are short and have a small $f_\sig-\pawgt$ value [basically because $f_\sig$ charges $g(\Dir)$ per level, while a fixed path collects about $\E\w_0$ per step, so a staircase causes $f_\sig-\pawgt$ to go down by about $\ghatz-\E\w_0\ge0$ per step]. Their cost is at most of order $j^22^{k-j}$ [the thickness of a fine strip]. The piece in the coarse scale is a stitching of two short attachment staircases, which again cause no trouble, and a crossing that comes from \eqref{eq:trialprob}. A tricky step here: \eqref{eq:trialprob} gives a crossing with $\pawgt>n\,g(\Dir)-Ks$. Here, $s=2^{k-j}$ and $n=2^{-4j}\clen_k-j^22^{k-j}$. To get $\pawgt>f_\sig-O(j^22^{k-j})$, the $f_\sig$ is NOT computed over the central crossing alone, but over the whole attachment-crossing-attachment path. This path is arranged to have transverse displacement $O(1)$, so $f_\sig$ of its displacement is $g(\Dir)2^{-4j}\clen_k+O(1)$. The two attachment staircases collect at least $j^22^{k-j}\E\w_0$ each.
Summed over all three pieces,
\begin{align*}
    \pawgt&\ge(2^{-4j}\clen_k-2j^22^{k-j})g(\Dir)-K2^{k-j}+2j^22^{k-j}\E\w_0\\
&=g(\Dir)2^{-4j}\clen_k-2j^22^{k-j}\bigl(g(\Dir)-\E\w_0\bigr)-K2^{k-j}
\ge f_\sig-O(j^22^{k-j}).
\end{align*}
Summing over all the pieces (coarse and fine) one has a loss of order $\sum_{j\ge J}j^22^{k-j}=s_k\sum_{j\ge J}j^22^{-j}$ [plus a fluctuation term controlled by the choice of $t_k\approx s_k/k$], which can be made smaller than $\e s_k/10000$.
Why coarse crossings exist (reason for multi-scales): A single event with probability $>q_0$ is not enough, so there are $j^2$ disjoint copies of the blocks at the coarse scale, and not having a crossing has probability at most $(1-q_0)^{j^2}$. This is summable against the $C^j$ grid intervals and bins at scale $j$, so with high probability we ALWAYS have a successful crossing, and it becomes the middle piece of the coarse path.}

\begin{proof}
Fix $A,R,\e$. Let $\cnst_1$ be the constant of Lemma~\ref{lem:grid} and choose $J$ so large that
\begin{equation}\label{eq:Jchoice}
16\cnst_1\sum_{j\ge J}j^22^{-j}\le\frac{\e}{8}.
\end{equation}
Apply Lemma~\ref{lem:grid} with this $J$ and with $\e_0=\e/512$, so that $t_k=\e s_k/(512k)$, and work on the event $G_{k,J}$.

See Figure~\ref{fig:entrance} for an illustration of the following construction. For a vertex $u$ of the region $\cR_k$ set
\begin{equation}\label{eq:Ipm}
I_{k,u}^-=\bigl(\lfloor\level(u)/l_{J,k}\rfloor-1\bigr)l_{J,k}\qquad\text{and}\qquad I_{k,u}^+=\bigl(\lceil\level(u)/l_{J,k}\rceil+1\bigr)l_{J,k} ,
\end{equation}
so that $I_{k,u}^\pm\in l_{J,k}\Z$ and both $\level(u)-I_{k,u}^-$ and $I_{k,u}^+-\level(u)$ lie in $[l_{J,k},2l_{J,k})$. In particular, 
%$l_{J,k}\le 2^{m_k}=\clen_k$.  So below, $A+2$ would do. The 4 is to match the definition in \eqref{eq:gridset}, which was chosen with slack!
\begin{equation}\label{f.1170}
\abs{I_{k,u}^\pm}\le |\level(u)|+2l_{J,k}\le A\clen_k+2l_{J,k}\le(A+4)\clen_k.
\end{equation}

Decompose the level interval $[I_{k,u}^-,\level(u)]$ into grid intervals as follows. The first is $[I_{k,u}^-,I_{k,u}^-+l_{J,k}]$, a grid interval at scale $J$. The remaining interval $[I_{k,u}^-+l_{J,k},\level(u)]$ has length less than $l_{J,k}$ and left endpoint in $l_{J,k}\Z\subset l_{J+1,k}\Z$.  It is covered, from the left, by at most $l_{J,k}/l_{J+1,k}-1\le15$ grid intervals at scale $J+1$, followed by a remainder of length less than $l_{J+1,k}$ with left endpoint in $l_{J+1,k}\Z$. Continuing, at most $15$ grid intervals are used at each scale $j>J$. Since the finest scale at $k\in\cK$ is $l_{j_*(k),k}=1$, the decomposition ends exactly at $\level(u)$. All these grid intervals have left endpoints in the sets $\cG_{j,k}$, since $\abs I\le(A+4)\clen_k$ along the way. 

In each of these grid intervals construct the path of Lemma~\ref{lem:grid}(i) at $\rho=\rcoord(u)$.  Consecutive paths share the endpoints $\hat z(I,\rcoord(u))$, and the last one ends at $\hat z(\level(u),\rcoord(u))=u$, because $u$ is the lattice point on its level with $\rcoord\in(\rcoord(u)-1,\rcoord(u)]$. Do the same on $[\level(u),I_{k,u}^+]$, with the decomposition carried out from top to bottom.  The scale-$J$ interval $[I_{k,u}^+-l_{J,k},I_{k,u}^+]$ is the topmost, and the remainder $[\level(u),I_{k,u}^+-l_{J,k}]$ is covered from its upper endpoint, which lies in $l_{J,k}\Z$, by at most $15$ grid intervals of each finer scale. 

At the conclusion of these steps we have an admissible path $y^u_{I_{k,u}^-:I_{k,u}^+}$ through $u$, whose pieces are paths $y^{\rcoord(u)}$ of Lemma~\ref{lem:grid}.
%The pieces are the paths $y^{\rcoord(u)}$ of Lemma~\ref{lem:grid}. For the concatenation we use $y^u$, because it is not determined by $\rcoord(u)$ and its index range: two vertices $u\ne u''$ with $\rcoord(u)=\rcoord(u'')$ and levels in the same interval of $l_{J,k}\Z$ have the same $I^\pm_{k,u}$ but different decompositions, hence different paths $y^u_{I^-_{k,u}:I^+_{k,u}}$ and $y^{u''}_{I^-_{k,u}:I^+_{k,u}}$. The superscript $u$ carries both $\rcoord(u)$ and $\level(u)$.
We call $y^u_{I_{k,u}^-:u}$ the \emph{entrance path} of $u$ and $y^u_{u:I_{k,u}^+}$ its \emph{exit path}.   By additivity and \eqref{eq:gridcost}, with at most $16$ grid intervals at each of at most $j_*(k)\le k$ scales, and then by $w_{j,k}=2^{k-j}=s_k2^{-j}$, the choice of $t_k$, and the bound \eqref{eq:Jchoice}, the entrance path satisfies, for both signs, \begin{equation}\label{eq:extcost} \begin{aligned}
\pawgt(y^u_{I_{k,u}^-:u})\ &\ge\ f_\sig(u-y^u_{I_{k,u}^-})-16\sum_{j\ge J}\cnst_1j^2w_{j,k}-16k\cdot2t_k\\
&=f_\sig(u-y^u_{I_{k,u}^-})-16\cnst_1s_k\sum_{j\ge J}j^22^{-j}-\frac{\e s_k}{16} \; \ge \; f_\sig(u-y^u_{I_{k,u}^-})-\frac{\e s_k}{8}-\frac{\e s_k}{16}.
\end{aligned}
\end{equation}
Similarly, the exit path satisfies
\begin{equation}\label{eq:extcost-exit}
\pawgt(y^u_{u:I_{k,u}^+})\ \ge\ f_\sig(y^u_{I_{k,u}^+}-u)-\frac{\e s_k}{8}-\frac{\e s_k}{16}.
\end{equation}
By \eqref{f.1170}, all vertices on these paths satisfy $\abs\level\le(A+4)\clen_k$ and, by \eqref{eq:j2w} and \eqref{eq:gridcost}, $\abs{\rcoord-\rcoord(u)}\le\cnst_1\max_jj^2w_{j,k}\le2\cnst_1s_k$, hence $\abs\rcoord\le(R+2\cnst_1)s_k$.

Now let $x_{u:u'}$ be an admissible path in the region $\cR_k$ with $\pawgt(x_{u:u'})>f_\sig(u'-u)+\e s_k$ for some $\sigg\in\{-,+\}$. Let $\tilde x_{I_{k,u}^-:I_{k,u'}^+}$ be the path with $\tilde x_{I_{k,u}^-:u}=y^u_{I_{k,u}^-:u}$, $\tilde x_{u:u'}=x_{u:u'}$, and $\tilde x_{u':I_{k,u'}^+}=y^{u'}_{u':I_{k,u'}^+}$. By additivity, \eqref{eq:extcost}, and \eqref{eq:extcost-exit},
\[
\pawgt(\tilde x_{I_{k,u}^-:I_{k,u'}^+})>f_\sig(\tilde x_{I_{k,u'}^+}-\tilde x_{I_{k,u}^-})+\e s_k-2\Bigl(\frac{\e s_k}8+\frac{\e s_k}{16}\Bigr)>f_\sig(\tilde x_{I_{k,u'}^+}-\tilde x_{I_{k,u}^-})+\frac{\e s_k}2 .
\]
Thus, on $G_{k,J}$, the event in \eqref{eq:endpoints} implies that for some pair of \emph{coarse levels} $I^-<I^+$ in $l_{J,k}\Z\cap[-(A+4)\clen_k,(A+4)\clen_k]$ and some $\sigg\in\{-,+\}$ there exists a crossing $x_{I^-:I^+}$ of the block
\[\{z:\ I^-\le\level(z)\le I^+,\ \abs{\rcoord(z)}\le(R+2\cnst_1)s_k\}\]
with $\pawgt(x_{I^-:I^+})>f_\sig(x_{I^+}-x_{I^-})+\e s_k/2$. The number of such pairs is at most
\[
\cnst_{\!A,J}=\bigl(2(A+4)\clen_k/l_{J,k}+1\bigr)^2=\bigl(2(A+4)2^{4J+1}+1\bigr)^2 ,
\]
which does not depend on $k$.  For each $k\in\cK$ choose one of these pairs and a sign for which the probability of such a crossing is largest, and let $n'_k=I^+-I^-$ for this pair. By Proposition~\ref{prop:rarefast} and Remark~\ref{rk:translate} (with $s=s_k$, $n=n'_k$, the lower level $I^-$ in place of $0$, $R+2\cnst_1$ in place of $R$, and $\e/2$ in place of $\e$) this largest probability tends to $0$ as $k\to\infty$ in $\cK$.  Since there are at most $\cnst_{\!A,J}$ pairs and two signs,
\[\varlimsup_{k\to\infty,\ k\tsp\in\tsp\cK}\P(\text{event in \eqref{eq:endpoints}})\le\varlimsup_{k\to\infty,\ k\tsp\in\tsp\cK}\P(G_{k,J}^c)
\le2\beta_J \underset{J\to\infty}\longrightarrow0\quad\text{by Lemma \ref{lem:grid}\ref{lem:grid:prob}}. \qedhere \]
\end{proof}

\smallskip

\begin{remark}[Bounded vs.\ unbounded weights]\label{rk:endpoints-general}
Under \ref{A2}--\ref{A3}, Proposition~\ref{prop:endpoints} holds as stated, with the following changes in its proof. The endpoints of a path $x_{u:u'}$ in $\cR_k$ with $\pawgt(x_{u:u'})>f_\sig(u'-u)+\e s_k$ are first moved inward along the path by at most a constant number of steps, to vertices $\tilde u$ and $\tilde u'$ through which the paths $\gamma^{\tilde u}$ and $\gamma^{\tilde u'}$ can be run. This lowers $\pawgt(x_{u:u'})-f_\sig(u'-u)$ by at most $\e s_k/64$ for large $k$, which is absorbed into the margin $\e s_k/2$ in the reduction to crossings between coarse levels. The entrance path of $\tilde u$ is then built as in the proof above from the paths of Lemma~\ref{lem:grid} at $\rho=\rcoord(\tilde u)$, which connect consecutive vertices $\gamma^{\tilde u}(I)$ and end exactly at $\tilde u$ because $\gamma^{\tilde u}$ passes through $\tilde u$, and likewise for the exit path of $\tilde u'$. These are the only paths of Lemma~\ref{lem:grid}\ref{lem:grid:path} that are needed under \ref{A2}--\ref{A3}: at $\rho=\rcoord(\tilde u)$, on the grid intervals contained in $[I^-_{k,\tilde u},\level(\tilde u)]$, where the part of $\gamma^{\tilde u}$ below $\tilde u$ exists, and at $\rho=\rcoord(\tilde u')$, on the grid intervals contained in $[\level(\tilde u'),I^+_{k,\tilde u'}]$, where the part of $\gamma^{\tilde u'}$ above $\tilde u'$ exists, with $I^\pm$ as in \eqref{eq:Ipm}. See items~\ref{ch:gamma}, \ref{ch:trim}, \ref{ch:coarse} and~\ref{ch:endpoints} of Appendix~\ref{sub:staircases}. 
\end{remark}

\section{Disjoint witnesses and proof of the nonexistence of corners}
\label{sec:witnesses}

Here is an informal overview of this final proof section. Fix a parameter $\eta\in(0,1]$, the value of the parameter $q_0\in(0,1)$ of Section~\ref{sec:scales}, which is fixed in \eqref{eq:q0} below, and one scale $(s_*,\clen_*)$ at which the local events of the previous sections are rare. The plane is covered by \emph{blocks} $\{\cB_{p,h}:(p,h)\in\Z^2\}$, of length $3n_*$ and width $8s_*$, where $n_*\in(\clen_*/2,\clen_*]$ is the integer $n_k$ of \eqref{t.n_k}. 

Two rare events are attached to each block: the \emph{heavy-path event}, that some path in the block has weight above $f_\sig+\eta s_*$ for one of the two extreme supporting functionals, which is rare by Proposition~\ref{prop:endpoints}, and the \emph{crossing event} $\Elev_{v}(s_*,n_*)$ at the entrance of the block, whose probability is at most $q_0$ by the choice of $n_*$. Every admissible path from the origin is cut into \emph{stages}, each contained in one block. Lemma~\ref{lem:chunk-witness} describes the full stages that contain no witness of a rare event in  their block, and Lemma~\ref{lem:multiplicity} shows that the number of witnesses with pairwise disjoint sets of collected vertices in a stage controls how far its weight can exceed $f$. The BK inequality for path events \eqref{eq:bk} and a count of the possible paths of blocks 
then show that long paths contain only few disjoint witnesses (Proposition~\ref{prop:longpaths}).

Up to this point nothing depends on whether $g$ is differentiable at $\Dir$. The corner hypothesis $\halfgap >0$ from \eqref{eq:supports} enters only in Section~\ref{sub:proofmain}.  Through the identity \eqref{eq:defect-min}, a path stage that moves transversally by $s_*$ pays $\halfgap s_*$ relative to $f$, so with $\eta\le \halfgap /2$ every witness-free full stage has weight at most $f-\eta s_*$.  Thereby long paths fall short of $f$ by a linear amount, which contradicts \eqref{eq:geodesic-defect} for geodesics.

\subsection{Path events}\label{sub:bk}
Call an event $\cA$ a \emph{path event} if there exist a finite set $\vbad \subset\Z^2$, a finite family $\Pi$ of admissible paths with all collected vertices in $\vbad$, and a function $\lambda:\Pi\to\R$, such that $\cA$ is the event that $\pawgt(x_{m:n})>\lambda(x_{m:n})$ for some $x_{m:n}\in\Pi$. We allow a family to contain several copies of a particular path with different $\lambda$-values. A \emph{witness} of $\cA$ is a path $x_{m:n}$ that satisfies $\pawgt(x_{m:n})>\lambda(x_{m:n})$. Whether a path is a witness depends only on the weights at its collected vertices, and a path event is increasing in $\w$. For path events $\cA_1,\dots,\cA_q$ (not necessarily distinct) let $D(\cA_1,\dots,\cA_q)$ be the event that there exist witnesses $x^{(1)},\dots,x^{(q)}$ of $\cA_1,\dots,\cA_q$ with pairwise disjoint sets of collected vertices.
The van den Berg--Kesten inequality for path events, Lemma~\ref{lem:bk} in Appendix~\ref{sub:bk-appendix}, states that
\begin{equation}\label{eq:bk}
\P\bigl(D(\cA_1,\dots,\cA_q)\bigr)\le\prod_{j=1}^q\P(\cA_j) .
\end{equation}

\subsection{Scale, blocks, and rare events}\label{sub:scale}
Fix $\eta\in(0,1]$. Choose $q_0\in(0,1)$ small enough so that 
\begin{equation}\label{eq:q0}
60\tspa(2q_0)^{1/32}\le\tfrac12.  
\end{equation}
This $q_0$ value is used    in Proposition~\ref{prop:charlength} and Lemma~\ref{lem:records}.   A scale index $k\in\cK$ is fixed in \eqref{eq:rare-prob}--\eqref{eq:rare-K} below. For this $k$ abbreviate the scales  $s_*=s_k$, $\clen_*=\clen_k$, and $n_*=n_k$.  By  \eqref{t.n_k},  $n_*\in(\clen_*/2,\clen_*]$ and  $\P(\Elev(s_*,n_*))\le q_0$.

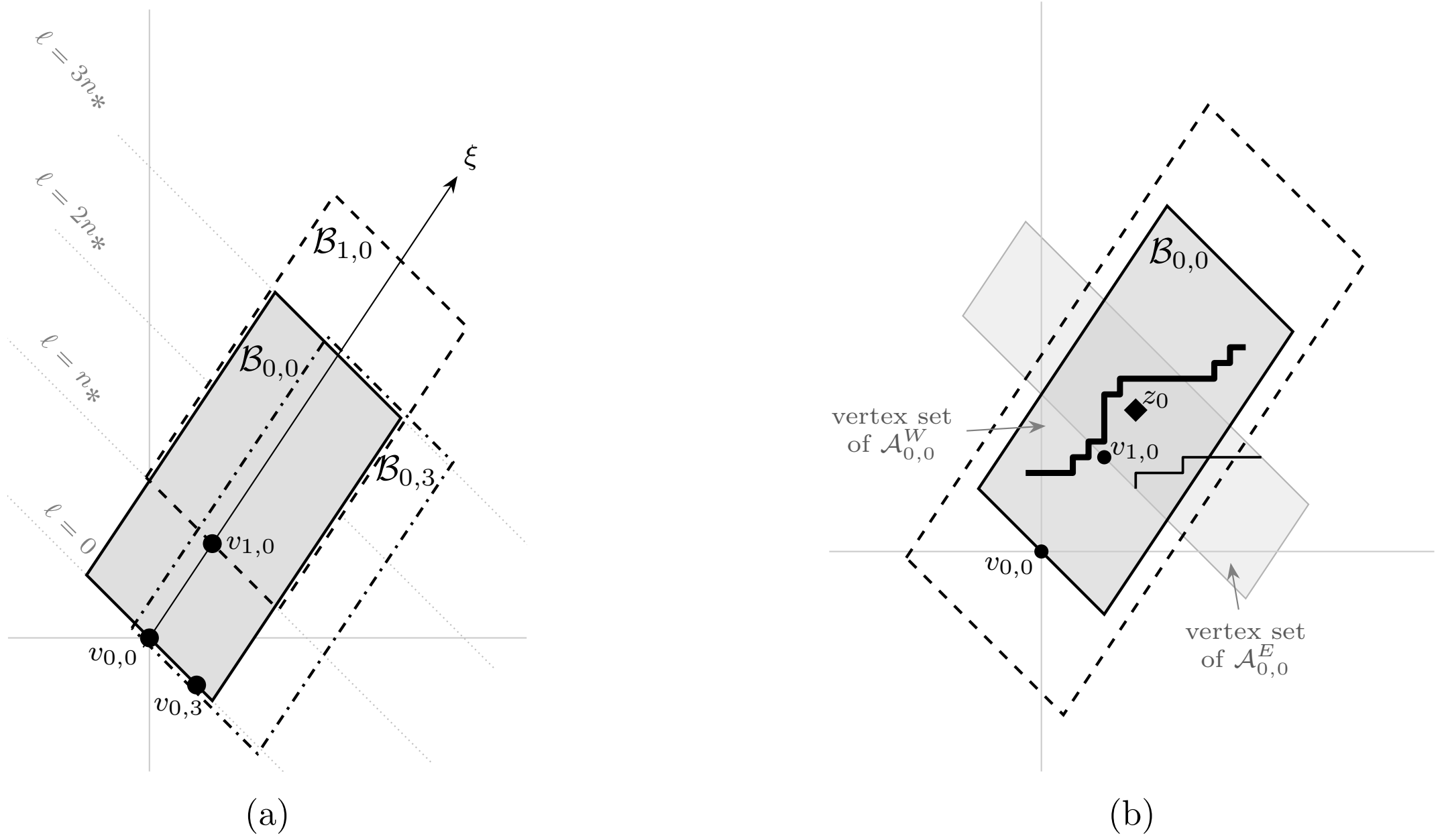
\begin{figure}[htb]
\centering
\begin{minipage}[b]{0.49\textwidth}\centering
\begin{tikzpicture}[scale=0.155,>=Stealth]
  \clip (-9,-8.5) rectangle (24,40);
  \draw[gray!40] (-9,0) -- (24,0); \draw[gray!40] (0,-8.5) -- (0,40);
  %\node[gray,below left] at (24,0) {\scriptsize $z_1$}; \node[gray,below left] at (0,40) {\scriptsize $z_2$};
  \draw[gray!50,densely dotted] (-14.00,14.00) -- (14.00,-14.00);
  \node[gray,rotate=-45,anchor=south] at (-6.00,6.00) {\tiny $\ell=0$};
  \draw[gray!50,densely dotted] (-10.00,20.00) -- (18.00,-8.00);
  \node[gray,rotate=-45,anchor=south] at (-6.00,16.00) {\tiny $\ell=n_*$};
  \draw[gray!50,densely dotted] (-6.00,26.00) -- (22.00,-2.00);
  \node[gray,rotate=-45,anchor=south] at (-6.00,26.00) {\tiny $\ell=2n_*$};
  \draw[gray!50,densely dotted] (-2.00,32.00) -- (26.00,4.00);
  \node[gray,rotate=-45,anchor=south] at (-6.00,35.00) {\tiny $\ell=3n_*$};
  \filldraw[fill=gray!25,draw=black,thick] (-4.00,4.00) -- (8.00,22.00) -- (16.00,14.00) -- (4.00,-4.00) -- cycle;
  \draw[black,thick,dashed] (-0.22,10.22) -- (11.78,28.22) -- (20.22,19.78) -- (8.22,1.78) -- cycle;
  \draw[black,thick,dash dot] (-1.12,0.62) -- (11.32,19.18) -- (19.32,11.18) -- (6.88,-7.38) -- cycle;
  \draw[->] (-0.40,-0.60) -- (19.60,29.40) node[above right=-2pt] {\scriptsize $\Dir$};
  \fill (0.00,0.00) circle (0.6);
  \node[anchor=west] at (-4.77,-1.30) {\scriptsize $v_{0,0}$};
  \fill (4.00,6.00) circle (0.6);
  \node[anchor=west] at (4.00,5.80) {\scriptsize $v_{1,0}$};
  \fill (3.00,-3.00) circle (0.6);
  \node[anchor=west] at (-0.6,-4.50) {\scriptsize $v_{0,3}$};
  \node[anchor=south east] at (10.5,15.50) {\small $\cB_{0,0}$};
  \node[anchor=south west] at (9.50,23.00) {\small $\cB_{1,0}$};
  \node[anchor=west] at (13.5,10.50) {\small $\cB_{0,3}$};
\end{tikzpicture}\\[2pt] {\small (a)}
\end{minipage}\hfill
\begin{minipage}[b]{0.49\textwidth}\centering
\begin{tikzpicture}[scale=0.155,>=Stealth]
  \clip (-13.5,-14) rectangle (25,40);
  \filldraw[fill=gray!12,draw=gray!60] (-5.00,15.00) -- (-1.00,21.00) -- (17.00,3.00) -- (13.00,-3.00) -- cycle;
  \draw[gray!40] (-13.5,0) -- (25,0); \draw[gray!40] (0,-14) -- (0,35);
  %\node[gray,below left] at (25,0) {\scriptsize $z_1$}; \node[gray,below left] at (0,40) {\scriptsize $z_2$};
  \filldraw[fill=gray!40,fill opacity=0.55,draw=black,thick] (-4.00,4.00) -- (8.00,22.00) -- (16.00,14.00) -- (4.00,-4.00) -- cycle;
  \draw[thick,dashed] (-8.60,-0.40) -- (10.60,28.40) -- (20.60,18.40) -- (1.40,-10.40) -- cycle;
  \draw[black,line width=1.7pt,line join=round] (-1.00,5.00) -- (0.00,5.00) -- (1.00,5.00) -- (2.00,5.00) -- (2.00,6.00) -- (3.00,6.00) -- (3.00,7.00) -- (4.00,7.00) -- (4.00,8.00) -- (4.00,9.00) -- (4.00,10.00) -- (5.00,10.00) -- (5.00,11.00) -- (6.00,11.00) -- (7.00,11.00) -- (8.00,11.00) -- (9.00,11.00) -- (10.00,11.00) -- (11.00,11.00) -- (11.00,12.00) -- (12.00,12.00) -- (12.00,13.00) -- (13.00,13.00);
  \draw[black,line width=0.7pt,line join=round] (6.00,4.00) -- (6.00,5.00) -- (7.00,5.00) -- (8.00,5.00) -- (9.00,5.00) -- (9.00,6.00) -- (10.00,6.00) -- (11.00,6.00) -- (12.00,6.00) -- (13.00,6.00) -- (14.00,6.00);
  \fill (0.00,0.00) circle (0.45);
  \fill (4.00,6.00) circle (0.45);
  \node[anchor=north west] at (-4.5,0.5) {\scriptsize $v_{0,0}$};
  \node[anchor=west] at (3.50,6.30) {\scriptsize $v_{1,0}$};
  \fill (6.00,9.00) +(0,0.75) -- +(0.75,0) -- +(0,-0.75) -- +(-0.75,0) -- cycle; \node[anchor=south west] at (5.50,8.30) {\scriptsize $z_0$};
  \node[anchor=north] at (8.80,20.80) {\small $\cB_{0,0}$};
  \node[gray!70!black,anchor=north,align=center,font=\tiny] (vsH) at (-9.5,10) {vertex set\\ of $\cA^{\pawgt}_{0,0}$};
  \draw[->,gray,shorten >=1pt] (vsH) -- (0.5,8);
  \node[gray!70!black,anchor=north,align=center,font=\tiny] (vsE) at (13.00,-3.70) {vertex set\\ of $\cA^{\Elev}_{0,0}$};
  \draw[->,gray,shorten >=1pt] (vsE) -- (12,-0.5);
\end{tikzpicture}\\[2pt] {\small (b)}
\end{minipage}
\caption{\small The blocks of Section~\ref{sub:scale} in the original coordinates $(z_1,z_2)$, drawn with $\Dir=(0.4,0.6)$, $s_*=1$, $n_*=10$, and $\clen_*=12$ for legibility. In the proof $\clen_*\ge n_*$ is of order up to $s_*^3$, so the blocks are much longer than wide. (a) The block $\cB_{0,0}$ (shaded) and two of its neighbors: $\cB_{1,0}$ (dashed), shifted by $n_*$ levels along $\Dir$, and $\cB_{0,3}$ (dash-dotted), shifted by $3s_*$ in the transverse direction $\protect\ead$. The two neighbors are drawn slightly enlarged so that their sides do not fall on the sides of $\cB_{0,0}$. The bullets are the lattice points $v_{p,h}$, and the light dotted lines are the levels $\level=pn_*$. Each block shares $2n_*$ of its $3n_*$ levels with the next block in $p$ and $7s_*$ of its width $8s_*$ with the next block in $h$, so every vertex lies in boundedly many blocks. (b) The block $\cB_{0,0}$ (shaded), the point $z_0$ (diamond), and the region $z_0+\{\abs\level\le2\clen_*,\ \abs\rcoord\le5s_*\}$ (dashed) to which Proposition~\ref{prop:endpoints} is applied. A witness of the heavy-path event $\cA^{\pawgt}_{0,0}$ (thick path) is any admissible path in $\cB_{0,0}$ with $\pawgt>f_{\protect\sig}+\eta s_*$ for some $\protect\sigg$, with arbitrary endpoints. The vertex set of this event is $\cB_{0,0}$. A witness of the crossing event $\cA^{\Elev}_{0,0}$ (thin path) runs from level $n_*$ to level $2n_*$ within transverse distance $8s_*$ of $v_{1,0}$ and has weight above $\ghatz n_*-Ks_*$. The vertex set of this event is the lightly shaded region $\{n_*\le\level<2n_*,\ \abs\rcoord\le8s_*+1\}$, which is wider than the block, so a witness can leave $\cB_{0,0}$.}
\label{fig:blocks}
\end{figure}
\medskip
\emph{Blocks.} For each coarse-grained index  $(p,h)\in\Z^2$ let $v_{p,h}$ be the lattice point with $\level(v_{p,h})=pn_*$ and $\rcoord(v_{p,h})\in(hs_*-1,hs_*]$. The set 
\begin{equation}\label{eq:Bblock}
\cB_{p,h}=\{z\in\Z^2:\ pn_*\le\level(z)\le pn_*+3n_*,\ \abs{\rcoord(z)-hs_*}\le4s_*\}
\end{equation}
is the \emph{block at $(p,h)$}. See Figure~\ref{fig:blocks}.

\medskip
\emph{Rare events.} Two path events, in the sense of Section~\ref{sub:bk}, are attached to the block at $(p,h)$. See Figure \ref{fig:blocks}.
\begin{enumerate}[label={\rm(\alph*)}] \itemsep=3pt
\item\label{Aa} The \emph{heavy-path event} $\cA^{\pawgt}_{p,h}$ is the event that some admissible path $x_{m:n}$ with all vertices in $\cB_{p,h}$ has $\pawgt(x_{m:n})>f_\sig(x_n-x_m)+\eta s_*$ for some $\sigg\in\{-,+\}$. Its vertex set is $\cB_{p,h}$ and its family consists of all admissible paths with all vertices in $\cB_{p,h}$, each listed twice, with the weight lower bounds $f_-(x_n-x_m)+\eta s_*$ and $f_+(x_n-x_m)+\eta s_*$. 

\item\label{Ab} The \emph{crossing event} $\cA^{\Elev}_{p,h}$ is the event $\Elev_{v_{p+1,h}}(s_*,n_*)$ of Section~\ref{sec:scales}: some admissible path $x_{(p+1)n_*:(p+2)n_*}$ with $\abs{\rcoord(x_i)-\rcoord(v_{p+1,h})}\le8s_*$ for all $i$ has $\pawgt(x_{(p+1)n_*:(p+2)n_*})>\ghatz n_*-Ks_*$. Its vertex set is $\{z:\ (p+1)n_*\le\level(z)<(p+2)n_*,\ \abs{\rcoord(z)-hs_*}\le8s_*+1\}$ and its family consists of the paths just described, each with the lower weight  bound $\ghatz n_*-Ks_*$.
\end{enumerate}
The union of two path events with vertex sets $\vbad,\vbad'$ is the path event with vertex set $\vbad \cup \vbad'$ and the union of the families. Hence
\begin{equation}\label{eq:Adef}
\cA_{p,h}=\cA^{\pawgt}_{p,h}\cup\cA^{\Elev}_{p,h}
\end{equation}
is a path event, and a witness of $\cA_{p,h}$ is a witness of $\cA^{\pawgt}_{p,h}$ or of $\cA^{\Elev}_{p,h}$.

\medskip
\emph{Probabilities.} See Figure \ref{fig:blocks}. Let $z_0$ be the lattice point with $\level(z_0)=pn_*+\lceil3n_*/2\rceil$ and $\rcoord(z_0)\in(hs_*-1,hs_*]$. Then $\cB_{p,h}\subset z_0+\{\abs\level\le2\clen_*,\ \abs\rcoord\le5s_*\}$ (as $3n_*/2\le2\clen_*$ and $4s_*+1\le5s_*$), so by Proposition~\ref{prop:endpoints} with $A=2$, $R=5$, $\e=\eta$, $\P(\cA^{\pawgt}_{p,h})\to0$ as $k\to\infty$ in $\cK$, uniformly in $(p,h)$. By translation invariance and the choice of $n_*$, $\P(\cA^{\Elev}_{p,h})=\P(\Elev(s_*,n_*))\le q_0$ for all $(p,h)$ and all $k\in\cK$.
\medskip

\emph{Scale.} Fix one scale index $k\in\cK$ so large that
\begin{align}
&\P(\cA_{p,h})\le2q_0\quad\text{for all }(p,h)\in\Z^2,\label{eq:rare-prob}\\
&\abs{\ghatz}+\abs a+M\le\eta s_*/2,\quad\text{and}\label{eq:rare-step}\\
&3\eta s_*+\abs a(2s_*+2)<Ks_* .\label{eq:rare-K}
\end{align}
The constant $\eta$ and the scale depend on each other only through these conditions, and $\eta$ is arbitrary in $(0,1]$. The third condition holds for all $s_*\ge1$ by the definition \eqref{t.Koo} of $K$: 
\[ 3\eta s_*+\abs a(2s_*+2)\le(3+4\abs a)s_*<100(1+\abs a)s_*\le Ks_*. \]

\begin{remark}[Bounded vs.\ unbounded weights]
This is the first point in this section at which the bounded and unbounded cases differ. The boundedness assumption \eqref{ass:bounded} is used in this section only through the requirement \eqref{eq:rare-step} above, which gives the one-step bound \eqref{eq:one-step} in the proof of Lemma~\ref{lem:multiplicity}. We therefore assume \eqref{ass:bounded} from \eqref{eq:rare-step} to the end of Section~\ref{sec:witnesses}, and, as in Section~\ref{sec:endpoints}, a remark at each place where it is used says what changes under \ref{A3} and points to the details in Appendix~\ref{sub:onestep}. Under \ref{A3}, the scale is fixed with $M$ replaced by $M_0=\eta s_*/4$ in \eqref{eq:rare-step}, and with the additional requirement \eqref{eq:rare-excess}, which holds for all large $k$ because $p''>4$. Condition \eqref{eq:rare-prob} is available under \ref{A2}--\ref{A3} because Proposition~\ref{prop:endpoints} is.
\end{remark}

The scale $(s_*,\clen_*,n_*)=(s_k,\clen_k,n_k)$  is now fixed, and only the lengths of paths tend to infinity. The letter $k$ is not needed again in this section and is used below for generic levels, as in Section~\ref{sec:setting}.

\subsection{Stages and witnesses}\label{sub:chunks}
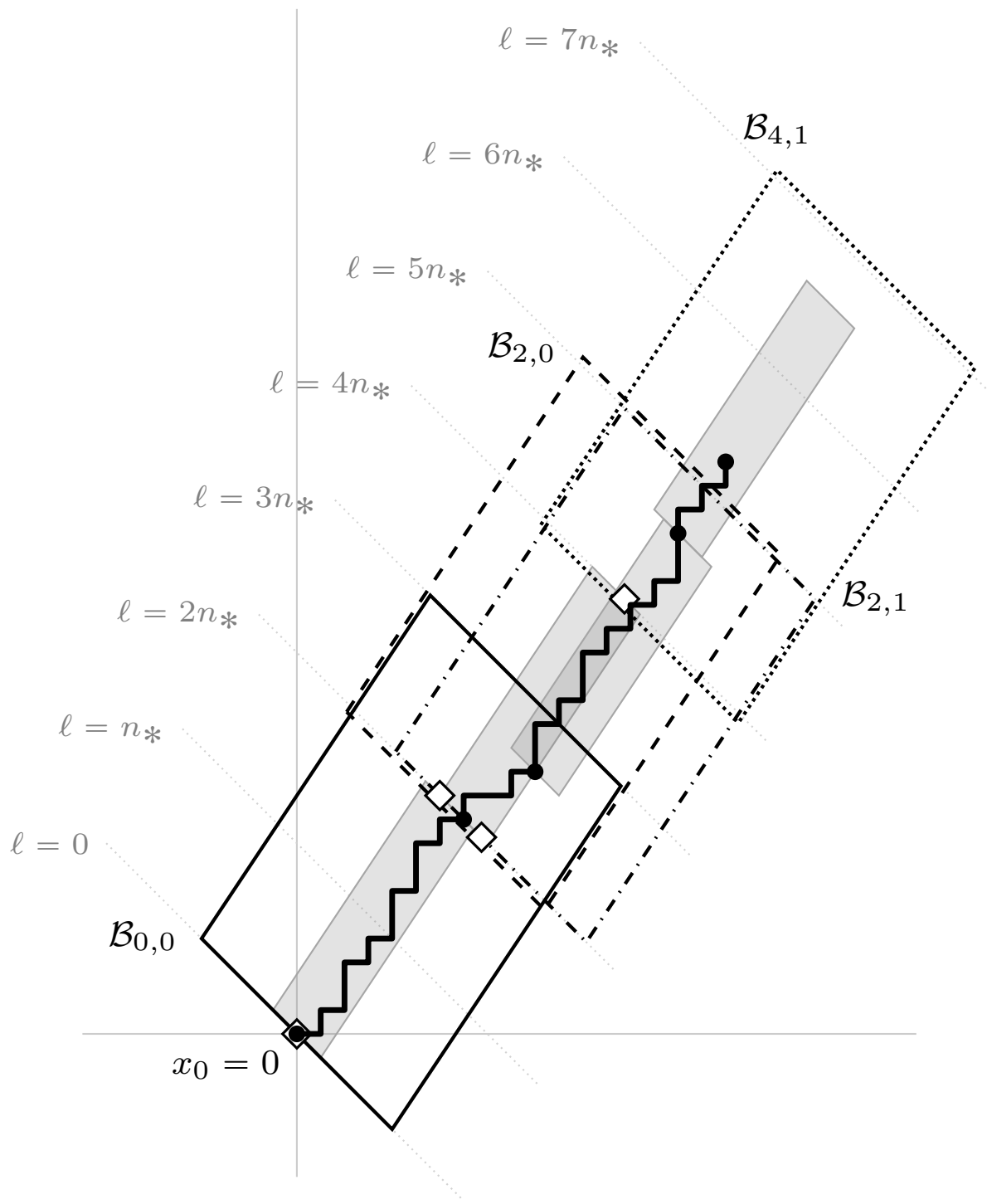
\begin{figure}[htb]
\centering
\begin{tikzpicture}[scale=0.19,>=Stealth]
  \draw[gray!40] (-9,0) -- (26,0); \draw[gray!40] (0,-6) -- (0,43);
  %\node[gray,below left] at (26,0) {\scriptsize $z_1$}; \node[gray,below left] at (0,46) {\scriptsize $z_2$};
  \draw[gray!35,densely dotted] (-8.00,8.00) -- (7.00,-7.00);
  \node[gray,anchor=east,font=\tiny] at (-8.00,8.00) {$\ell=0$};
  \draw[gray!35,densely dotted] (-4.80,12.80) -- (10.20,-2.20);
  \node[gray,anchor=east,font=\tiny] at (-4.80,12.80) {$\ell=n_*$};
  \draw[gray!35,densely dotted] (-1.60,17.60) -- (13.40,2.60);
  \node[gray,anchor=east,font=\tiny] at (-1.60,17.60) {$\ell=2n_*$};
  \draw[gray!35,densely dotted] (1.60,22.40) -- (16.60,7.40);
  \node[gray,anchor=east,font=\tiny] at (1.60,22.40) {$\ell=3n_*$};
  \draw[gray!35,densely dotted] (4.80,27.20) -- (19.80,12.20);
  \node[gray,anchor=east,font=\tiny] at (4.80,27.20) {$\ell=4n_*$};
  \draw[gray!35,densely dotted] (8.00,32.00) -- (23.00,17.00);
  \node[gray,anchor=east,font=\tiny] at (8.00,32.00) {$\ell=5n_*$};
  \draw[gray!35,densely dotted] (11.20,36.80) -- (26.20,21.80);
  \node[gray,anchor=east,font=\tiny] at (11.20,36.80) {$\ell=6n_*$};
  \draw[gray!35,densely dotted] (14.40,41.60) -- (29.40,26.60);
  \node[gray,anchor=east,font=\tiny] at (14.40,41.60) {$\ell=7n_*$};
  \filldraw[fill=gray,fill opacity=0.22,draw=gray!70,thin] (-1.00,1.00) -- (5.40,10.60) -- (7.40,8.60) -- (1.00,-1.00) -- cycle;
  \filldraw[fill=gray,fill opacity=0.22,draw=gray!70,thin] (6.00,10.00) -- (12.40,19.60) -- (14.40,17.60) -- (8.00,8.00) -- cycle;
  \filldraw[fill=gray,fill opacity=0.22,draw=gray!70,thin] (9.00,12.00) -- (15.40,21.60) -- (17.40,19.60) -- (11.00,10.00) -- cycle;
  \filldraw[fill=gray,fill opacity=0.22,draw=gray!70,thin] (15.00,22.00) -- (21.40,31.60) -- (23.40,29.60) -- (17.00,20.00) -- cycle;
  \draw[black,solid,thick] (-4.00,4.00) -- (5.60,18.40) -- (13.60,10.40) -- (4.00,-4.00) -- cycle;
  \draw[black,dashed,thick] (2.09,13.51) -- (12.01,28.39) -- (20.31,20.09) -- (10.39,5.21) -- cycle;
  \draw[black,dash dot,thick] (4.15,11.85) -- (13.75,26.25) -- (21.75,18.25) -- (12.15,3.85) -- cycle;
  \draw[black,densely dotted,thick] (10.24,21.36) -- (20.16,36.24) -- (28.46,27.94) -- (18.54,13.06) -- cycle;
  \fill[white] (0.00,0.00) +(0,0.6) -- +(0.6,0) -- +(0,-0.6) -- +(-0.6,0) -- cycle; \draw[semithick] (0.00,0.00) +(0,0.6) -- +(0.6,0) -- +(0,-0.6) -- +(-0.6,0) -- cycle;
  \fill[white] (6.00,10.00) +(0,0.6) -- +(0.6,0) -- +(0,-0.6) -- +(-0.6,0) -- cycle; \draw[semithick] (6.00,10.00) +(0,0.6) -- +(0.6,0) -- +(0,-0.6) -- +(-0.6,0) -- cycle;
  \fill[white] (7.75,8.25) +(0,0.6) -- +(0.6,0) -- +(0,-0.6) -- +(-0.6,0) -- cycle; \draw[semithick] (7.75,8.25) +(0,0.6) -- +(0.6,0) -- +(0,-0.6) -- +(-0.6,0) -- cycle;
  \fill[white] (13.75,18.25) +(0,0.6) -- +(0.6,0) -- +(0,-0.6) -- +(-0.6,0) -- cycle; \draw[semithick] (13.75,18.25) +(0,0.6) -- +(0.6,0) -- +(0,-0.6) -- +(-0.6,0) -- cycle;
  \draw[black,line width=1.3pt,line join=round] (0.00,0.00) -- (1.00,0.00) -- (1.00,1.00) -- (2.00,1.00) -- (2.00,2.00) -- (2.00,3.00) -- (3.00,3.00) -- (3.00,4.00) -- (4.00,4.00) -- (4.00,5.00) -- (4.00,6.00) -- (5.00,6.00) -- (5.00,7.00) -- (5.00,8.00) -- (6.00,8.00) -- (6.00,9.00) -- (7.00,9.00);
  \draw[black,line width=1.3pt,line join=round] (7.00,9.00) -- (7.00,10.00) -- (8.00,10.00) -- (9.00,10.00) -- (9.00,11.00) -- (10.00,11.00);
  \draw[black,line width=1.3pt,line join=round] (10.00,11.00) -- (10.00,12.00) -- (10.00,13.00) -- (11.00,13.00) -- (11.00,14.00) -- (12.00,14.00) -- (12.00,15.00) -- (12.00,16.00) -- (13.00,16.00) -- (13.00,17.00) -- (14.00,17.00) -- (14.00,18.00) -- (15.00,18.00) -- (15.00,19.00) -- (16.00,19.00) -- (16.00,20.00) -- (16.00,21.00);
  \draw[black,line width=1.3pt,line join=round] (16.00,21.00) -- (16.00,22.00) -- (17.00,22.00) -- (17.00,23.00) -- (18.00,23.00) -- (18.00,24.00);
  \fill (0.00,0.00) circle (0.35);
  \fill (7.00,9.00) circle (0.35);
  \fill (10.00,11.00) circle (0.35);
  \fill (16.00,21.00) circle (0.35);
  \fill (18.00,24.00) circle (0.35);
  \node[anchor=east,font=\scriptsize] at (-4.30,4.00) {$\cB_{0,0}$};
  \node[anchor=east,font=\scriptsize] at (11.61,28.69) {$\cB_{2,0}$};
  \node[anchor=west,font=\scriptsize] at (22.15,18.25) {$\cB_{2,1}$};
  \node[anchor=south,font=\scriptsize] at (20.16,36.44) {$\cB_{4,1}$};
  \node[anchor=north east,font=\scriptsize] at (0,0) {$x_0=0$};
\end{tikzpicture}
\caption{\small A path $x_{0:42}$ with $x_0=0$ cut into stages, in the original coordinates, drawn with $\Dir=(0.4,0.6)$, $s_*=1$, and $n_*=8$ for legibility. In the proof $n_*$ is much larger than $s_*$, so the blocks are much longer than wide. The bullets mark the start and end of each stage, and the open diamonds are the lattice points $v_{p,h}$ of the four blocks, with $v_{0,0}=x_0$. The shaded parallelogram attached to each piece is the region $\{k\le\level\le k+2n_*,\ \abs{\rcoord-\rcoord(x_k)}\le s_*\}$ determined by its first vertex $x_k$: a piece with a longitudinal cut leaves it through its far short side, a piece with a transverse cut leaves it through one of its long sides, and the remainder ends inside it. The first piece runs from level $0$ to level $16=2n_*$ without moving transversally by $s_*$, so its cut is longitudinal. The second piece starts at level $16$ and moves transversally by $s_*$ within $5$ levels, so its cut is transverse. The third piece, from level $21$ to level $37$, is again longitudinal, and the last $5$ levels form the remainder. By \eqref{eq:chunkindex} the indices of the four pieces are $(0,0)$, $(2,0)$, $(2,1)$, and $(4,1)$, and each piece lies in the block with its index, drawn solid, dashed, dash-dotted, and dotted, respectively. For legibility, $\cB_{2,0}$ and $\cB_{4,1}$ are drawn slightly enlarged, so that no two block sides fall on one another. Consecutive stages may have the same $p$ (the second and third pieces here), and their indices differ by at most $2$ in each coordinate.}
\label{fig:stages}
\end{figure}
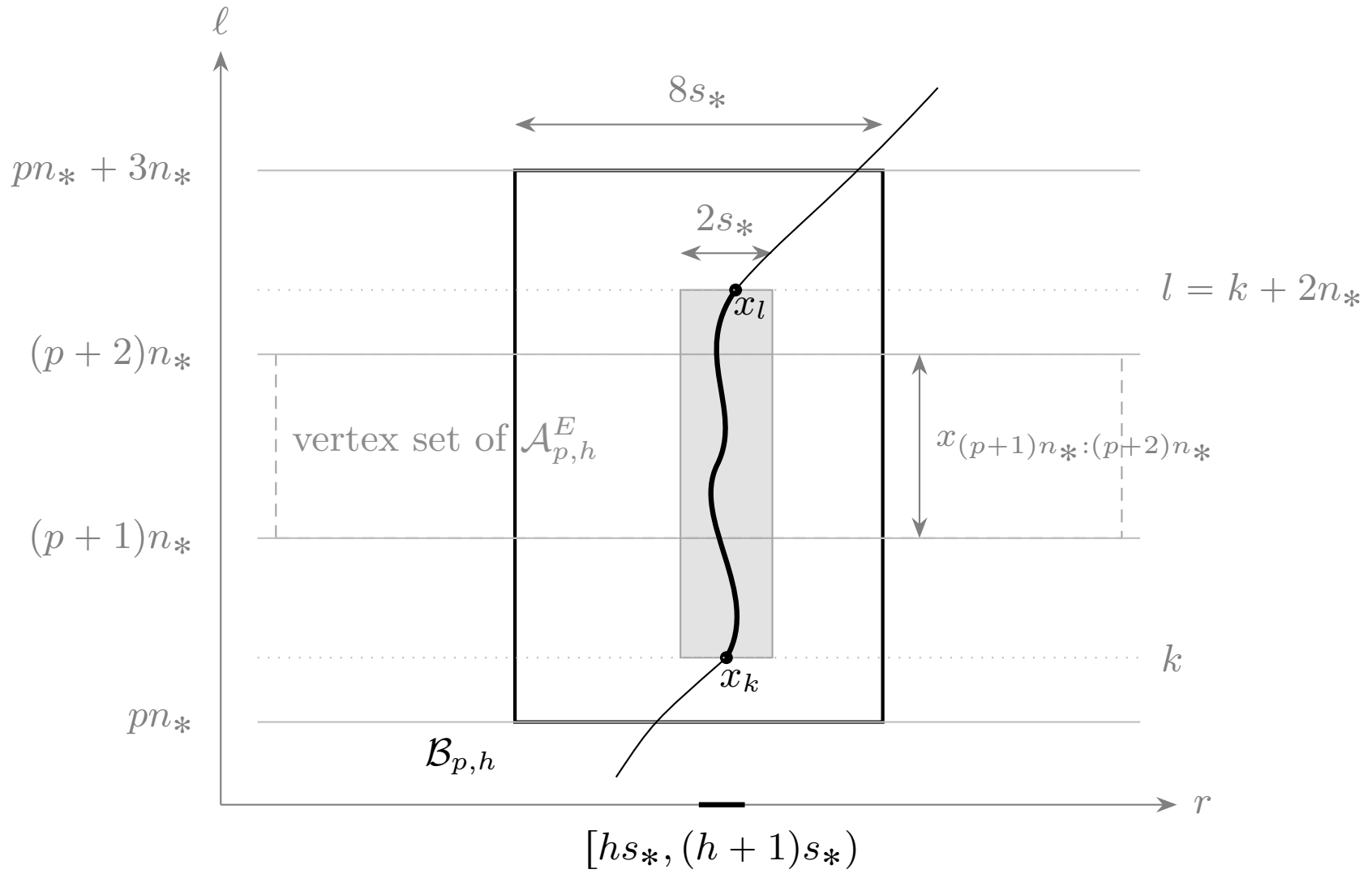
\begin{figure}[htb]
\centering
\begin{tikzpicture}[scale=0.75,>=Stealth]
  %% coordinates: horizontal r, vertical level (not to scale); s = 0.5 units, n_* = 2 units
  %% block B_{p,h}: levels p n_* .. p n_* + 3 n_* (0..6), |r - hs| <= 4s (|r| <= 2), with hs at r=0
  \draw[thick] (-2,0) rectangle (2,6);
  %% vertex set of the crossing event: levels (p+1)n_* .. (p+2)n_*, |r - hs| <= 8s+1
  \draw[gray!70,densely dashed] (-4.6,2) rectangle (4.6,4);
  \node[gray!90,below left] at (-0.9,3.5) {\footnotesize vertex set of $\cA^{\Elev}_{p,h}$};
  %% axes and level lines
  \draw[->,gray] (-5.2,-0.9) -- (-5.2,7.3) node[above] {\footnotesize $\level$};
  \draw[->,gray] (-5.2,-0.9) -- (5.2,-0.9) node[right] {\footnotesize $\rcoord$};
  \foreach \y/\lab in {0/{$pn_*$}, 2/{$(p+1)n_*$}, 4/{$(p+2)n_*$}, 6/{$pn_*+3n_*$}} {\draw[gray!50] (-4.8,\y) -- (4.8,\y); \node[left,gray] at (-5.3,\y) {\footnotesize \lab};}
  %% bin of r(x_k) on the r axis: [hs,(h+1)s) = [0,0.5)
  \draw[very thick] (0,-0.9) -- (0.5,-0.9);
  \node[below] at (0.25,-1.0) {\footnotesize $[hs_*,(h+1)s_*)$};
  %% the stage: starts at level k = 0.7 at r(x_k)=0.3, ends at level l = k + 2n_* = 4.7 (longitudinal cut)
  \filldraw[fill=gray,fill opacity=0.22,draw=gray!70,thin] (-0.2,0.7) rectangle (0.8,4.7);
  \draw[<->,gray] (-0.2,5.1) -- (0.8,5.1) node[midway,above] {\footnotesize $2s_*$};
  %% path: previous stage (thin), the stage (thick), next stage (thin)
  \draw[thin] (-0.9,-0.6) .. controls (-0.5,0.0) .. (0.3,0.7);
  \draw[very thick] (0.3,0.7) .. controls (0.7,1.4) and (-0.1,2.2) .. (0.2,2.8)
       .. controls (0.5,3.4) and (-0.1,4.0) .. (0.4,4.7);
  \draw[thin] (0.4,4.7) .. controls (0.9,5.3) and (1.4,5.6) .. (2.6,6.9);
  \fill (0.3,0.7) circle (0.07); 
  \node[right=1pt] at (0,0.45) {\footnotesize $x_k$};
  \fill (0.4,4.7) circle (0.07); \node[right=1pt] at (0.15,4.5) {\footnotesize $x_l$};
  \node[right,gray] at (4.85,0.7) {\footnotesize $k$};
  \node[right,gray] at (4.85,4.7) {\footnotesize $l=k+2n_*$};
  \draw[gray!50,dotted] (-4.8,0.7) -- (4.8,0.7);
  \draw[gray!50,dotted] (-4.8,4.7) -- (4.8,4.7);
  %% the segment used in Lemma stage-witness
  \draw[<->,gray] (2.4,2) -- (2.4,4) node[midway,right] {\scriptsize $x_{(p+1)n_*:(p+2)n_*}$};
  \node[below left] at (-2,0) {\footnotesize $\cB_{p,h}$};
  \draw[<->,gray] (-2,6.5) -- (2,6.5) node[midway,above] {\footnotesize $8s_*$};
\end{tikzpicture}
\caption{\small A stage $x_{k:l}$ with index $(p,h)$ and its block $\cB_{p,h}$ (thick outline), drawn in the coordinates $(\rcoord,\level)$ and not to scale. The stage starts at level $k\in[pn_*,(p+1)n_*)$ at a vertex $x_k$ with $\rcoord(x_k)\in[hs_*,(h+1)s_*)$, stays within transverse distance $s_*+1$ of $\rcoord(x_k)$, and ends at the first level $l$ at which it has moved transversally by $s_*$ or longitudinally by $2n_*$, here the latter. The shaded region is $\{k\le\level\le k+2n_*,\ \abs{\rcoord-\rcoord(x_k)}\le s_*\}$, as in Figure~\ref{fig:stages}. A longitudinal stage leaves it through its top side. Its vertices lie in the block, whose levels run from $pn_*$ to $pn_*+3n_*$ and whose width is $8s_*$. The heavy-path event $\cA^{\pawgt}_{p,h}$ concerns all paths in the block. The crossing event $\cA^{\Elev}_{p,h}$ concerns the levels from $(p+1)n_*$ to $(p+2)n_*$ in the wider dashed region: a longitudinal stage crosses these two levels, and in the proof of Lemma~\ref{lem:chunk-witness} the segment $x_{(p+1)n_*:(p+2)n_*}$ is the candidate witness of $\cA^{\Elev}_{p,h}$.}
\label{fig:chunks}
\end{figure}

Let $x_{0:n}$ be an admissible path with $x_0=0$. Cut it into \emph{stages} as follows, illustrated in Figures~\ref{fig:stages} and~\ref{fig:chunks}.
 The stage starting at level $k$ is the segment $x_{k:l}$ where $l$ is the first level in $(k,n]$ with
\begin{equation}\label{eq:cut}
l-k\ge2n_*\qquad\text{or}\qquad\abs{\rcoord(x_l)-\rcoord(x_k)}\ge s_* ,
\end{equation}
and the next stage starts at level $l$. A stage that ends at a cut specified in \eqref{eq:cut} is a \emph{full} stage.  The stage is \emph{transverse} if $\abs{\rcoord(x_l)-\rcoord(x_k)}\ge s_*$ and \emph{longitudinal} otherwise, in which case $l=k+2n_*$ by minimality of $l$. If no level in $(k,n]$ satisfies \eqref{eq:cut}, the stage is $x_{k:n}$ and it is called a \emph{remainder}. A remainder may have length $0$. Thus for some $N\ge0$,  $x_{0:n}$ is the concatenation of $N$ full stages and a remainder.

Index the stage $x_{k:l}$ by
\begin{equation}\label{eq:chunkindex}
(p,h)=\bigl(\lfloor k/n_*\rfloor,\lfloor \rcoord(x_k)/s_*\rfloor\bigr).
\end{equation}
Since \eqref{eq:cut} fails at every level in $(k,l)$ and $\abs{\Delta\rcoord}\le1$ per step, every vertex $x_i$ of the stage has $0\le i-k\le2n_*$  and $\abs{\rcoord(x_i)-\rcoord(x_k)}\le s_*+1$ (the strict bound $<s_*$ holds for $k<i<l$ because \eqref{eq:cut} fails there, and the last vertex $x_l$ is one step further). 
Since    $pn_*\le k<(p+1)n_*$ and $hs_*\le \rcoord(x_k)<(h+1)s_*$ by \eqref{eq:chunkindex}, every vertex of the stage has $pn_*\le\level\le pn_*+3n_*$ and $\abs{\rcoord-hs_*}\le2s_*+1\le4s_*$. So the stage has all vertices in the block $\cB_{p,h}$ at its index. Consecutive stages $x_{k:l}$ and $x_{l:l'}$ share only the vertex $x_l$, which is collected by the later stage only. So distinct stages have disjoint sets of collected vertices, and hence so do segments of distinct stages.

A full stage $x_{k:l}$ with index $(p,h)$ is \emph{witness-free} if none of its segments is a witness of $\cA_{p,h}$.

\begin{lemma}\label{lem:chunk-witness}
Let $x_{k:l}$ be a witness-free full stage with index $(p,h)$. Then $\pawgt(x_{k:l})\le f_\sig(x_l-x_k)+\eta s_*$ for both $\sigg\in\{-,+\}$. If $x_{k:l}$ is a longitudinal stage, then moreover $\pawgt(x_{k:l})\le f(x_l-x_k)-\eta s_*$.
\end{lemma}

%This last statement is a crucial step where weight goes below $f$.

\begin{proof}
See Figure \ref{fig:chunks}. Since $x_{k:l}$ is witness-free, every segment $x_{i:j}$ of $x_{k:l}$  satisfies $\pawgt(x_{i:j})\le f_\pm(x_j-x_i)+\eta s_*$ and hence also $\pawgt(x_{i:j})\le\tfrac12\bigl(f_-(x_j-x_i)+f_+(x_j-x_i)\bigr)+\eta s_*=f(x_j-x_i)+\eta s_*$. %For $(i,j)=(k,l)$ this is the first assertion.

Consider a longitudinal stage. Then $l=k+2n_*$, and as $pn_*\le k<(p+1)n_*$, the levels $(p+1)n_*$ and $(p+2)n_*$ lie in $[k,l]$. Split $x_{k:l}$ into the three segments $x_{k:(p+1)n_*}$, $x_{(p+1)n_*:(p+2)n_*}$, and $x_{(p+2)n_*:l}$, the last of which may have length $0$. By the first paragraph, the first and the last satisfy $\pawgt\le f+\eta s_*$. The middle segment runs from level $(p+1)n_*$ to level $(p+2)n_*$ and its vertices satisfy $\abs{\rcoord-hs_*}\le2s_*+1$, hence $\abs{\rcoord-\rcoord(v_{p+1,h})}\le2s_*+2\le8s_*$. So it is one of the paths of the family of $\cA^{\Elev}_{p,h}$.  Since the stage is assumed  witness-free, $\pawgt(x_{(p+1)n_*:(p+2)n_*})\le\ghatz n_*-Ks_*$. Its two endpoints have $\abs{\rcoord-\rcoord(x_k)}\le s_*+1$, so its displacement has level $n_*$ and $\abs\rcoord\le2s_*+2$.  By \eqref{t.f(z)}, $f(x_{(p+2)n_*}-x_{(p+1)n_*})\ge\ghatz n_*-\abs a(2s_*+2)$. The two bounds together give 
\[
\pawgt(x_{(p+1)n_*:(p+2)n_*})\le\ghatz n_*-Ks_*\le f(x_{(p+2)n_*}-x_{(p+1)n_*})+\abs a(2s_*+2)-Ks_* .
\]
Add up the bounds over  the three segments and use \eqref{eq:rare-K}: 
\[
\pawgt(x_{k:l})\le f(x_l-x_k)+2\eta s_*+\abs a(2s_*+2)-Ks_*<f(x_l-x_k)-\eta s_* . \qedhere
\]
\end{proof}

\begin{remark}
The transverse case gives no weight shortfall below $f$ by itself. Indeed, if $g$ were differentiable at $\Dir$, then $\halfgap =0$ and $f_-=f_+=f$ would be the tangent plane of $g$ at $\Dir$, so that $f-g=o(\abs\rcoord)$ for displacements of small slope $\rcoord/\level$, by \eqref{eq:contact}. A stage of length $n_*\gg s_*$ that follows a geodesic with transverse displacement $s_*$ would then have $\pawgt\approx g=f-o(s_*)$, that is, the transverse move would cost nothing at the order $\eta s_*$.  When $\Dir$ is a point of non-differentiability, the transverse move costs at least $\halfgap s_*$ relative to $f$, through \eqref{eq:defect-min}, so $\pawgt(x_{k:l})\le f(x_l-x_k)-(\halfgap -\eta)s_*$, which is at most $f(x_l-x_k)-\eta s_*$ once we fix  $\eta\le \halfgap /2$. This is how  the corner hypothesis  enters the proof  in Section~\ref{sub:proofmain}.
\end{remark}

For a stage $x_{k:l}$ with index $(p,h)$ let $\mult(x_{k:l})$ be the maximal number of segments of $x_{k:l}$ that are witnesses of $\cA_{p,h}$ with pairwise disjoint sets of collected vertices. Every witness collects at least one vertex, so $\mult(x_{k:l})\le l-k$ is finite. For a path $x_{0:n}$ with $x_0=0$ let $\mult_{\rm tot}(x_{0:n})$ be the sum of $\mult(x_{k:l})$ over all its stages, the remainder included. Since segments of distinct stages have disjoint sets of collected vertices, the union over all stages of maximal families that attain  $\mult$ consists of $\mult_{\rm tot}(x_{0:n})$ segments of $x_{0:n}$, each a witness of the event $\cA_{p,h}$ of its stage, with pairwise disjoint sets of collected vertices.

\begin{lemma}\label{lem:multiplicity}
Every stage $x_{k:l}$ satisfies
\begin{equation}\label{eq:chunk-defect}
\pawgt(x_{k:l})\le f(x_l-x_k)+\eta s_*+2\eta s_*\,\mult(x_{k:l}).
\end{equation}
\end{lemma}
\begin{proof}
For a single step $x_{i:i+1}$, $\pawgt(x_{i:i+1})=\w_{x_i}\le M$ and $f(x_{i+1}-x_i)=\ghatz+a\tspa\rcoord(x_{i+1}-x_i)$ with $\abs{\rcoord(x_{i+1}-x_i)}\le1$, so by \eqref{eq:rare-step}
\begin{equation}\label{eq:one-step}
\pawgt(x_{i:i+1})\le f(x_{i+1}-x_i)+\abs{\ghatz}+\abs a+M\le f(x_{i+1}-x_i)+\frac{\eta s_*}2 .
\end{equation}
This is the only place where the one-step upper bound $M$ is used.  %Upper bound on weight used here. 

\medskip
Let $x_{k:l}$ be a stage with index $(p,h)$. Define levels $k=d_0<d_1<\dots<d_q\le l$ inductively: given $d_\nu$, let $d_{\nu+1}$ be the first level in $(d_\nu,l]$ with $\pawgt(x_{d_\nu:d_{\nu+1}})>f(x_{d_{\nu+1}}-x_{d_\nu})+\eta s_*$ if there is one,  otherwise stop and set $q=\nu$. Each \emph{split segment} $x_{d_\nu:d_{\nu+1}}$ has, by minimality of $d_{\nu+1}$ and \eqref{eq:one-step},
\begin{align*}
\pawgt(x_{d_\nu:d_{\nu+1}})&=\pawgt(x_{d_\nu:d_{\nu+1}-1})+\pawgt(x_{d_{\nu+1}-1:d_{\nu+1}})\\
&\le f(x_{d_{\nu+1}-1}-x_{d_\nu})+\eta s_*+f(x_{d_{\nu+1}}-x_{d_{\nu+1}-1})+\frac{\eta s_*}2\le f(x_{d_{\nu+1}}-x_{d_\nu})+2\eta s_*
\end{align*}
(for $d_{\nu+1}=d_\nu+1$ the path $x_{d_\nu:d_{\nu+1}-1}$ has length $0$, and $\pawgt=f=0$ for it). Each split segment is a segment of the stage, so all its vertices lie in $\cB_{p,h}$, and by \eqref{eq:defect-min} it has $\pawgt>f+\eta s_*\ge\min_\sig f_\sig+\eta s_*$. So it is a witness of $\cA^{\pawgt}_{p,h}$, hence of $\cA_{p,h}$.
%Here witnesses of $\cA^{\Elev}_{p,h}$ ignored.
The split segments have pairwise disjoint sets of collected vertices, so their number $q$ is at most $\mult(x_{k:l})$. The final piece $x_{d_q:l}$ has $\pawgt\le f+\eta s_*$ by the stopping rule, and if $d_q=l$ it has length $0$ and $\pawgt=f=0$.
Adding up,
\[
\pawgt(x_{k:l})\le f(x_l-x_k)+2\eta s_*\tsp q+\eta s_*\le f(x_l-x_k)+\eta s_*+2\eta s_*\,\mult(x_{k:l}) ,
\]
which is \eqref{eq:chunk-defect}.
\end{proof}

\begin{remark}[Bounded vs.\ unbounded weights]
Under \ref{A3} a collected vertex $z$ can have weight above $M_0$. Let $\exc_z=(\w_z-M_0)^+$ be its \emph{excess}.  Then \eqref{eq:one-step} holds with the additional term $+\exc_{x_i}$ on the right-hand side, and \eqref{eq:chunk-defect} with the additional term $+Z(x_{k:l})$, the total excess of the collected vertices of the stage. See \eqref{eq:one-step-general} and Lemma~\ref{lem:multiplicity-general} in Appendix~\ref{sub:onestep}.
\end{remark}

\subsection{Long paths contain few disjoint witnesses}\label{sub:longpaths}
\begin{proposition}\label{prop:longpaths}
At the fixed scale $(s_*,\clen_*,n_*)$, almost surely there exists a random $N_0<\infty$ such that every admissible path $x_{0:n}$ with $x_0=0$ and $N\ge N_0$ full stages satisfies $\mult_{\rm tot}(x_{0:n})<\lfloor N/32\rfloor$.
\end{proposition}

\gist{For a path with $N$ stages to carry $N/32$ disjoint witnesses, $N/32$ rare events, each of probability at most $2q_0$, must occur disjointly, which by BK has probability at most $(2q_0)^{N/32}$. There are only $60^N$ ways to choose the blocks and to distribute the witnesses among them, and $q_0$ was chosen to beat that.}

\begin{proof}
Let $B_N$ be the event that some admissible path $x_{0:n}$ from $0$ with exactly $N$ full stages has $\mult_{\rm tot}(x_{0:n})\ge\lfloor N/32\rfloor$. It suffices to show $\sum_N\P(B_N)<\infty$.   We describe first how $B_N$ is decomposed according to coarse-grained paths and the distribution of disjoint witnesses among the stages of $x_{0:n}$.

Let $\mathbf b=(p_\tau,h_\tau)_{0\le\tau\le N}$ be the sequence of stage indices of the $N+1$ stages (full stages and remainder) of $x_{0:n}$, in order. We call $\mathbf b$ the \emph{path of blocks} of $x_{0:n}$. Here $(p_0,h_0)=(0,0)$ since $x_0=0$. Consecutive stages $x_{k:l}$, $x_{l:l'}$ have $0<l-k\le2n_*$ and $\abs{\rcoord(x_l)-\rcoord(x_k)}\le s_*+1<2s_*$, so $p_{\tau+1}-p_\tau\in\{0,1,2\}$ and $\abs{h_{\tau+1}-h_\tau}\le2$. Thus there are at most $15^N$ possible paths of blocks. Call a vector $\mathbf m=(m_0,\dots,m_N)\in\Z_{\ge0}^{N+1}$ with $\sum_{\tau=0}^Nm_\tau=\lfloor N/32\rfloor$ an \emph{allocation} that prescribes $m_\tau$  disjoint witness segments to stage $(p_\tau,h_\tau)$, for $\tau=0,\dotsc,N$. There are $\binom{N+\lfloor N/32\rfloor}{N}\le2^{N+\lfloor N/32\rfloor}\le4^N$ allocations. For a path of blocks $\mathbf b$ and an allocation $\mathbf m$ let $D_{\mathbf b,\mathbf m}$ be the disjoint occurrence event $D(\cA_1,\dots,\cA_{\lfloor N/32\rfloor})$ of Section~\ref{sub:bk} for the list in which $\cA_{p_\tau,h_\tau}$ appears $m_\tau$ times  for each $\tau\in\{0,\dotsc,N\}$. By the BK inequality \eqref{eq:bk}  and the upper bound \eqref{eq:rare-prob},
\begin{equation}\label{eq:Dcm}
\P(D_{\mathbf b,\mathbf m})\le(2q_0)^{\lfloor N/32\rfloor}\qquad\text{for every path of blocks $\mathbf b$ and allocation $\mathbf m$.}
\end{equation}

We claim that
\begin{equation}\label{eq:BNincl}
B_N\subset\bigcup_{\mathbf b}\ \bigcup_{\mathbf m}D_{\mathbf b,\mathbf m},
\end{equation}
the union being over the at most $15^N$ paths $\mathbf b$ of blocks and the at most $4^N$ allocations $\mathbf m$. The union on the right-hand side is not necessarily disjoint. 

To argue \eqref{eq:BNincl}, let $\w\in B_N$.   In the environment $\w$, take  any path $x_{0:n}$ from $x_0=0$  that has exactly $N$ full stages and  $\mult_{\rm tot}(x_{0:n})\ge\lfloor N/32\rfloor$. Let $\mathbf b(x_{0:n})=(p_\tau(x_{0:n}),h_\tau(x_{0:n}))_{0\le\tau\le N}$ be its path of blocks.  
Let $\mult_\tau(x_{0:n})$ be the maximal  number of segments of stage $\tau$ of $x_{0:n}$ that  witness  $\cA_{p_\tau(x_{0:n}),\tspb h_\tau(x_{0:n})}$ and have  pairwise disjoint sets of collected vertices.  Since $\sum_{\tau=0}^N\mult_\tau(x_{0:n})=\mult_{\rm tot}(x_{0:n})\ge\lfloor N/32\rfloor$, there exists an allocation $\mathbf m=(m_0,\dots,m_N)$ such that  $m_\tau\le\mult_\tau(x_{0:n})$ for each $\tau\in\{0,\dotsc,N\}$.
In each stage $\tau$ choose $m_\tau$ of the $\mult_\tau(x_{0:n})$ disjoint witnesses. This  gives $\lfloor N/32\rfloor$ segments that are witnesses of the events in the list of $D_{\mathbf b(x_{0:n}),\mathbf m}$. 
  Thus 
$\w\in D_{\mathbf b(x_{0:n}), \mathbf m}$ and   \eqref{eq:BNincl} has been verified.   

By \eqref{eq:BNincl}, then by \eqref{eq:Dcm}, and finally by \eqref{eq:q0}, we get this summable bound:  
\[
\P(B_N)\le15^N\,4^N\,(2q_0)^{\lfloor N/32\rfloor}\le(2q_0)^{-1}\bigl(60\,(2q_0)^{1/32}\bigr)^N\le(2q_0)^{-1}2^{-N}. 
\qedhere
\]
 \end{proof}

\subsection{Proof of Theorem~\ref{thm:main}}\label{sub:proofmain}
Let $\Dir\in\ri\Uset$ be deterministic and suppose $\halfgap =\halfgap (\Dir)>0$. Only now do we choose the parameter of Section~\ref{sub:scale}: set
\begin{equation}\label{eq:etac}
\eta=\min\{1,\halfgap /2\},
\end{equation}
and fix the scale $(s_*,\clen_*,n_*)$ as in Section~\ref{sub:scale} with this $\eta$, so that \eqref{eq:rare-prob}--\eqref{eq:rare-K} hold.  

%The passage from \eqref{eq:chunk-defect} to \eqref{eq:full-defect} is merely arithmetic: the first RHS term switches to negative because $\mult(x_{k:l})\ge1$ allows us to turn  $\eta s_*+2\eta s_*\,\mult(x_{k:l})$ into $-\eta s_*+4\eta s_*\,\mult(x_{k:l})$.

\begin{lemma}\label{lem:fulldefect}
Assume \eqref{eq:etac}. Every full stage $x_{k:l}$ satisfies
\begin{equation}\label{eq:full-defect}
\pawgt(x_{k:l})\le f(x_l-x_k)-\eta s_*+4\eta s_*\,\mult(x_{k:l}).
\end{equation}
Every admissible path $x_{0:n}$ with $x_0=0$ and $N$ full stages satisfies
\begin{equation}\label{eq:total-defect}
\pawgt(x_{0:n})\le f(x_n-x_0)-(N-1)\eta s_*+4\eta s_*\,\mult_{\rm tot}(x_{0:n}).
\end{equation}
\end{lemma}
\begin{proof}
Consider first a witness-free full stage $x_{k:l}$, that is, $\mult(x_{k:l})=0$. If the stage  is longitudinal, Lemma~\ref{lem:chunk-witness} gives $\pawgt(x_{k:l})\le f(x_l-x_k)-\eta s_*$. If the stage is transverse, then $\abs{\rcoord(x_l-x_k)}\ge s_*$, and by Lemma~\ref{lem:chunk-witness}, \eqref{eq:defect-min}, and $\halfgap \ge2\eta$,
\begin{align*}
\pawgt(x_{k:l}) &\le\min_\sig f_\sig(x_l-x_k)+\eta s_*=f(x_l-x_k)-\halfgap \abs{\rcoord(x_l-x_k)}+\eta s_* \\
&\le f(x_l-x_k)-\halfgap s_*+\eta s_*\le f(x_l-x_k)-\eta s_* .
\end{align*}
This is the only place where $\halfgap >0$ is used.

\medskip
%In this part witnesses of $\cA^{\Elev}_{p,h}$ must be included.

Now consider other full stages.  If $\mult(x_{k:l})\ge1$, then $\eta s_*+2\eta s_*\mult(x_{k:l})\le-\eta s_*+4\eta s_*\mult(x_{k:l})$ %because $2\eta s_*\mult\ge2\eta s_*$, 
and \eqref{eq:chunk-defect} gives \eqref{eq:full-defect}.

Finally, by additivity, sum \eqref{eq:full-defect} over the $N$ full stages and \eqref{eq:chunk-defect} (weakened to $\eta s_*+4\eta s_*\mult$) over the remainder to get  \eqref{eq:total-defect}. 
\end{proof}

\begin{remark}[Bounded vs.\ unbounded weights]
Under \ref{A3}, \eqref{eq:total-defect} holds with the additional term $+\sum_{b\in A(x_{0:n})}Z_b$, where $A(x_{0:n})$ is the set of indices of the blocks visited by the path and $Z_b$ is the total excess of the vertices in a fixed neighborhood of the block $\cB_b$, see \eqref{eq:total-defect-general}. The sets $A(x_{0:n})$ are connected sets of about $N$ blocks containing the block of the origin, so the maximum of this sum over all paths with $N$ full stages is a greedy lattice animal. By Lemma~\ref{lem:excess}, which bounds it, and the requirement \eqref{eq:rare-excess}, almost surely this sum is at most $N\eta s_*/32$ for all large $N$, uniformly over the paths with $N$ full stages. Hence \eqref{eq:linear-defect} below holds with the same right-hand side, and the proof of Theorem~\ref{thm:main} goes through unchanged. This and Lemma~\ref{lem:heavy} of Appendix~\ref{sub:staircases} are the two places where $p''>4$ is used. 
\end{remark}
%With the extra term $-N\eta s_*/32$, the bound in \eqref{eq:linear-defect} is $(N-1)\eta s_*-N\eta s_*/8-N\eta s_*/32\ge N\eta s_*/2$ for $N\ge3$, which is assumed there.

\begin{proof}[Proof of Theorem~\ref{thm:main}]
Let $N_0$ be as in Proposition~\ref{prop:longpaths}. For every admissible path $x_{0:n}$ from  $x_0=0$ that has $N\ge\max\{N_0,3\}$ full stages, 
Proposition~\ref{prop:longpaths} guarantees $\mult_{\rm tot}(x_{0:n})<N/32$. 
Substituting this into 
\eqref{eq:total-defect} and using $N\ge8/3$ gives, for any such path, 
\begin{equation}\label{eq:linear-defect}
\pawgt(x_{0:n})<f(x_n-x_0)-(N-1)\eta s_*+\frac{N\eta s_*}8\le f(x_n-x_0)-\frac{N\eta s_*}2 . 
\end{equation}
 
For $n\ge1$ let $z_n$ be the lattice point on level $n$ with $\rcoord(z_n)\in(-1,0]$. 
Let $x^{(n)}_{0:n}$  be a  geodesic 
    from $0$ to $z_n$. Let  $N^{(n)}$ be the number of  full stages in $x^{(n)}_{0:n}$.  Since every full stage has length at most $2n_*$ and the remainder has length less than $2n_*$, $n<2n_*(N^{(n)}+1)$. 
    From this
    \[ N^{(n)}>\frac{n}{2n_*}-1\ge \frac{n}{2\clen_*}-1\underset{n\to\infty}\longrightarrow\infty. \]
 Thus $N^{(n)}\ge\max\{N_0,3\}$   for large enough $n$, and  we can apply   \eqref{eq:linear-defect} to $x^{(n)}_{0:n}$: 
  \[  \pawgt(x^{(n)}_{0:n})< f(z_n)-\frac{N^{(n)}\eta s_*}2  
  \le f(z_n)-n\frac{\eta s_*}2\Bigl(\frac1{2\clen_*}-\frac1n\Bigr) \le f(z_n)-Cn,    \]
  where we can take  $C=\eta s_*/(8\clen_*)$ once  $n\ge 4\clen_*$.  
  
But by 
 \eqref{eq:geodesic-defect} the geodesic satisfies almost surely $\pawgt(x^{(n)}_{0:n})=f(z_n)+o(n)$.
 This contradiction forces $\halfgap =0$, equivalently, that  $g$ is differentiable at $\Dir$. 
  Together with Section~\ref{sec:flat}, which covers the points of $\Uperc\cap\ri\Uset$ by Theorem~\ref{thm:flat-top}, this proves Theorem \ref{thm:main}.   
\end{proof}

\section{Remarks on the proof}
\label{sec:remarks}

\subsection{Relaxing bounded weights}
\label{sub:bdd-chat} 
The lower bound $-M$ of assumption \eqref{ass:bounded} was used only in Hoeffding's inequality \eqref{eq:hoeffding}, for the deterministic attachments of Lemma~\ref{lem:intervals} and the staircases of Lemma~\ref{lem:grid}. The upper bound $M$ was used in \eqref{eq:hoeffding}, in the one-step bound \eqref{eq:one-step}, which says that a single collected vertex cannot lift the weight of a path above $f$ by more than $\eta s_*/2$, and, once the staircases are replaced by random paths, in the trimming of their endpoints (Appendix~\ref{sub:staircases}). 

Without boundedness, two things change. 

(i)  On the lower-tail side, very negative weights can sit on any prescribed vertex in a region of $s^4$ vertices, so the attachments and staircases must be built to avoid such vertices. This requires  a finite moment of order $p'>1$ on $\w_0^-$ (Appendices~\ref{sub:attachments}--\ref{sub:staircases}). 

(ii) On the upper-tail side,  one has to control the excess weight that a path with $N$ stages can collect. This is a greedy lattice animal problem which requires a finite moment of order $p''>4$ on $\w_0^+$  (Appendix~\ref{sub:onestep}).   The exponent $4$ comes from the bound $\clen_k\le s_k^3$ on the characteristic length, since the blocks of Section~\ref{sub:scale} have of order $s_*\clen_*\le s_*^4$ vertices. 

\subsection{Lower tail  moment assumption}\label{sub:lb_mom_chat}

 Assumption \ref{A2} is close to optimal.  
The following example shows that,  with upper-bounded weights,  $\E\w_0^-<\infty$ is necessary for a  finite shape function.

\begin{example}
Fix $p\in(0,1)$ and  $\alpha\in(p\vee\frac12,1)$.
Let the  i.i.d.~weights be nonpositive with
$\P\{\w_0<-t\}=(1+t)^{-\alpha}$ for $t\ge0$.
Then $\E[(\w_0^-)^p]<\infty$ but $\E\w_0^-=\infty$.

Let $\cS=\{(3i,3j):i,j\ge1,\ 1/2\le i/j\le2\}$.
Every up-right path from $0$ to $x\in\N^2$ visits $x-e_1$ or $x-e_2$.
Since all weights are nonpositive,
$\Lpp_{0,x}\le\max\{\w_{x-e_1},\w_{x-e_2}\}$.
For $A>0$, set
\[\cE_x^A=\{\w_{x-e_1}<-A|x|_1,\ \w_{x-e_2}<-A|x|_1\}.\]
These events are independent over $x\in\cS$, and
$\P(\cE_x^A)=(1+A|x|_1)^{-2\alpha}$.
There exist order $k$ points $x\in\cS$ with $\abs{x}_1=3k$,
so $\sum_{x\tsp\in\tsp\cS}\P(\cE_x^A)\asymp\sum_{k\ge1}k^{1-2\alpha}=\infty$.
The Borel--Cantelli lemma, simultaneously for all positive integers $A$, gives
$
\varliminf_{x\tsp\in\tsp\cS :\,  \abs{x}_1\to\infty}
{|x|_1}^{-1}{\Lpp_{0,x}}=-\infty$ 
 almost surely.
\end{example}

\subsection{Upper tail moment assumption}\label{sub:ub_mom-chat}

   Assumption \ref{A3} is not optimal. 
The upper tail of $\mu$ enters the proof at the scale $s_k$ through one requirement: a block of $s_k\times\clen_k$ vertices  must not contain a single vertex with weight of order $s_k$, that is,
\[
s_k\clen_k\,\P\{\w_0>\delta s_k\}\longrightarrow0 .
\]
This is necessary for Proposition~\ref{prop:endpoints}, since one such vertex is by itself a path with $\pawgt>f_\sig+\e s_k$. It is also, up to constants, what Lemmas~\ref{lem:heavy} and~\ref{lem:excess} need: if $\clen_k\le s_k^{\gamma}$, the bound in Lemma~\ref{lem:excess} becomes $CN\bigl(s_*^{(1+\gamma)/p''}+s_*^{2+\gamma-p''}\bigr)$, which is $o(Ns_*)$ if and only if $p''>1+\gamma$. The change of measure of Lemma~\ref{lem:tilt} gives $\clen_k\le s_k^3$,   and therefore the proof requires $p''>4$, as stated in assumption \ref{A3}.

The true order of $\clen_k$ is not known, but scaling heuristics suggest it. Suppose that a path of length $\ell$ has transverse fluctuations of order $\ell^\zeta$ and weight fluctuations of order $\ell^\theta$, with $\theta=2\zeta-1$. A path confined to width $s$ over length $n$ then behaves like a concatenation of about $n/s^{1/\zeta}$ free pieces of length $s^{1/\zeta}$, each falling short of the shape function by about $s^{\theta/\zeta}$. The total shortfall is about $ns^{(\theta-1)/\zeta}=ns^{2-2/\zeta}$, and it reaches order $s$ when $n\asymp s^{2/\zeta-1}$. Under KPZ scaling, $\zeta=\tfrac23$, this gives $\clen_k\asymp s_k^2$, and the requirement above becomes $p''>3$. 

However, KPZ scaling is expected only when $p''>5$. For $2<p''<5$ the Flory argument predicts $\zeta=(1+p'')/(2p''-1)$ and $\theta=3/(2p''-1)$.  See \cite[(6), (7), and Table~I]{Gue-etal-15}, where these exponents are confirmed numerically. This gives $\clen_k\asymp s_k^{3(p''-1)/(p''+1)}$, and the requirement $p''>1+3(p''-1)/(p''+1)$ is equivalent to $(p''-1)(p''-2)>0$, which holds for every $p''>2$. So, heuristically, the scheme of this paper ought to work for every $p''>2$, the threshold of the shape theorem, if the true order of the characteristic length were known, and all of the loss is in the bound $\clen_k\le s_k^3$.  This threshold is close to sharp: by \cite[Theorem~5]{Gan-Gin-Nam-26-}, there are nonnegative i.i.d.\ weights with $\E[\w_0^2]<\infty$ for which $\Lpp_{0,(n,n)}/n\to\infty$ almost surely, so that the shape function is not even finite. Neither the KPZ exponents for confined paths nor the Flory exponents are proved, so this is only conjectural.

Truncation could be a different route to weaker moment conditions. The theorem applies to each truncated law $\mu\wedge M$, and $0\le g_\mu-g_{\mu\wedge M}\le C\int_M^\infty\P(\w_0>t)^{1/2}\,dt\to0$ by the greedy lattice animal bound \cite[Theorem~2.3]{Mar-02-spa}. Differentiability of $g_\mu$ would then follow from a modulus of differentiability of $g_{\mu\wedge M}$ that is uniform in $M$, which the present proof by contradiction does not give.

\subsection{Exclusion of the flat segment}\label{sub:flatexcluded}

As mentioned, the assumption 
  $\Dir\notin\Uperc$ is used in the modification argument to deduce $\clen_k\le s_k^3$ (Proposition~\ref{prop:charlength}). A deeper fact is that, regardless of this particular argument, we expect  the characteristic length  approach to fail on the flat segment.  The reason is that 
known rigorous  results suggest that 
inside the flat segment the characteristic length is much larger than $s_k^3$. 

To describe the picture, 
suppose that $p_1=\mu\{\Mtop\}>\opc$ and that $\Dir$ lies in the relative interior of $\Uperc$ so that $g(\Dir)=\Mtop$. If $x_{0:n}$ satisfies 
 $\pawgt(x_{0:n})=\Mtop n=\ghatz n$,  $x_{0:n}$  witnesses $\Elev(s,n)$ if it stays within transverse distance $8s$ of the line $\rcoord=0$. In a strip of fixed width the property $\pawgt(x_{0:n})=\Mtop n$  fails eventually.  But we still expect a $\lambda_1=\lambda_1(p_1,\Dir)>0$ such that $\P(\Elev(s,n))>q_0$ for all large $s$ and all $n\le e^{\lambda_1s}$, so that the characteristic length of Proposition \ref{prop:charlength} would be  at least exponential in $s$. Lemma~\ref{lem:records} needs $\clen_k2^{-4k}\to0$, which then fails.

 No result in the literature  directly   justifies the sketch above.  Indirect evidence comes from  
 the contact process on an interval of $N$  sites, started from all sites occupied.  The extinction time $\sigma_N$ satisfies $\log\sigma_N/N\to\gamma_2(\lambda)>0$ in probability 
 and the normalized extinction time converges to an exponential random variable \cite{Dur-Liu-88,Dur-Sch-88,Mou-93}.  These results and their proofs have counterparts in oriented percolation when the strip is in the symmetric direction $(\tfrac12,\tfrac12)$  \cite[page~1572]{Dur-Sch-88}. Our direction $\Dir$  is  not necessarily symmetric, but the same contour argument should apply.

\appendix
\section{Oriented percolation and the BK inequality}
\label{app:oriented}

\subsection{Oriented percolation}\label{sub:oriented}

Lemma~\ref{lem:oriented} below is a standard Peierls contour estimate. For $\varkappa=0$ it is the statement that oriented site percolation on $\Z^2$ with site density close to one percolates \cite{Dur-84}. For general $\varkappa$ it also follows from the result of Liggett, Schonmann, and Stacey \cite{Lig-Sch-Sta-97} that a $\varkappa$-dependent field with marginals close to one dominates a product measure of density close to one. We include the short proof. Its geometric part is stated separately, as Lemma~\ref{lem:contour}, because it is used again in the proof of Lemma~\ref{lem:stair-geometry}.

For $y\in\Z^2$ let $Q_y=y+[-\tfrac12,\tfrac12]^2$ be the closed unit square centered at $y$. The sides of these squares are called \emph{edges}. They are the nearest-neighbor edges of the dual lattice $\Z^2+(\tfrac12,\tfrac12)$, whose vertices are the corners of the squares, and each of them crosses exactly one nearest-neighbor edge of $\Z^2$. A \emph{circuit} of \emph{length} $m$ is a closed path of $m$ directed edges of the dual lattice $(\tfrac12,\tfrac12)+\Z^2$, each of which starts at the endpoint of the previous edge, and that does not traverse any edge twice in either direction. A directed edge has a square $Q_x$ on its left  side and a square $Q_y$ on its right side, determined by facing in the   direction of the edge. For example, an upward-pointing  edge has the square to its west on its left.

\begin{lemma}\label{lem:contour}
Let $S\subset\Z^2$ be a finite set that is connected in the nearest-neighbor graph of $\Z^2$, and let $z\in S$. Then there exists a circuit $\Gamma$, of some length $m\ge4$, with the following properties.
\begin{enumerate}[label={\rm(\roman*)}] \itemsep=1pt 
\item\label{con:sides} Every edge of $\Gamma$ has a square $Q_w$ with $w\in S$ on its left and a square $Q_{w'}$ with $w'\notin S$ on its right.
\item\label{con:forward} The squares on the left of the upward and leftward edges of $\Gamma$ are centered at at least $m/4$ distinct points $w$. Each of these points satisfies $w\in S$, and $w+e_1\notin S$ or $w+e_2\notin S$.
\item\label{con:backward} The squares on the right of the downward and rightward edges of $\Gamma$ are centered at at least $m/4$ distinct points $w'$. Each of these points satisfies $w'\notin S$, and $w'+e_1\in S$ or $w'+e_2\in S$.
\item\label{con:size} $\abs S\le m^2/16$.
\item\label{con:count} $\Gamma$ contains the common side of $Q_{z+(j-1)e_1}$ and $Q_{z+je_1}$ for some integer $1\le j\le m/2$. For given $z$ and $m$, there are at most $m\,3^{m-1}$ circuits of length $m$ that contain such an edge.
\end{enumerate}
\end{lemma}

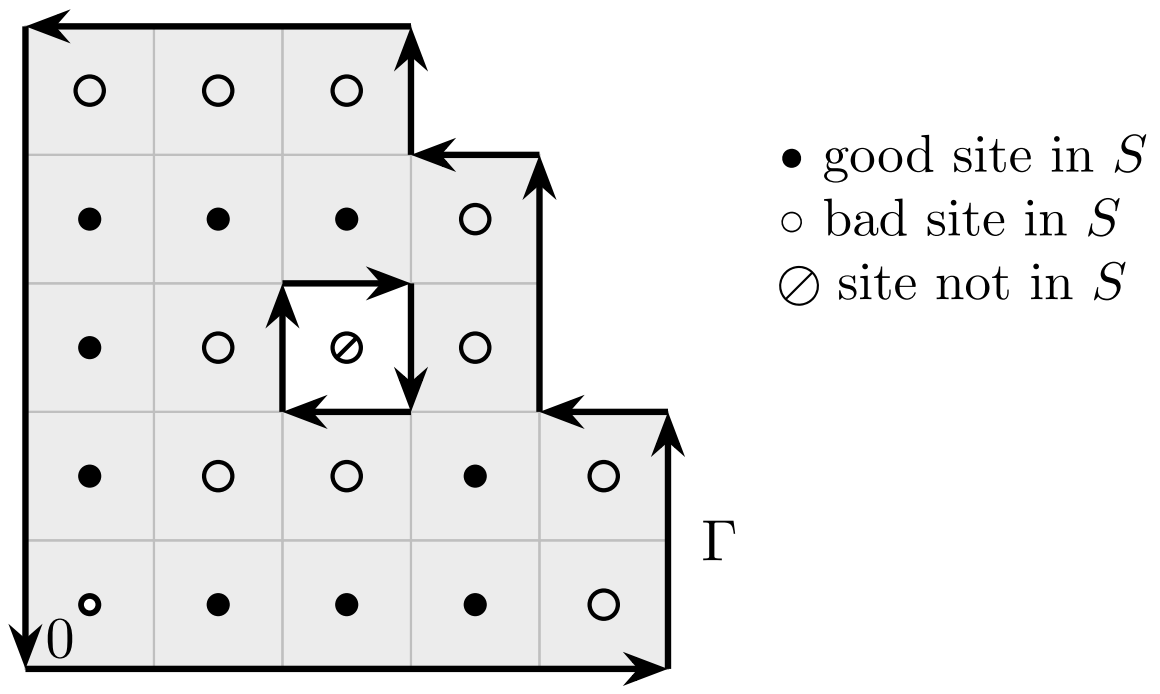
\begin{figure}[htb]
\centering
\begin{tikzpicture}[scale=0.85,>=Stealth]
  %% S = good sites {(0,0),(1,0),(2,0),(3,0),(0,1),(0,2),(0,3),(3,1),(1,3),(2,3)} and their forward neighbors
  \foreach \x/\y in {0/0,1/0,2/0,3/0,4/0,0/1,1/1,2/1,3/1,4/1,0/2,1/2,3/2,0/3,1/3,2/3,3/3,0/4,1/4,2/4} \fill[gray!15] (\x-0.5,\y-0.5) rectangle (\x+0.5,\y+0.5);
  \foreach \x/\y in {0/0,1/0,2/0,3/0,4/0,0/1,1/1,2/1,3/1,4/1,0/2,1/2,3/2,0/3,1/3,2/3,3/3,0/4,1/4,2/4} \draw[gray!50,thin] (\x-0.5,\y-0.5) rectangle (\x+0.5,\y+0.5);
  %% outer contour, counterclockwise
  \draw[very thick,->] (-0.5,-0.5) -- (4.5,-0.5);
  \draw[very thick,->] (4.5,-0.5) -- (4.5,1.5);
  \draw[very thick,->] (4.5,1.5) -- (3.5,1.5);
  \draw[very thick,->] (3.5,1.5) -- (3.5,3.5);
  \draw[very thick,->] (3.5,3.5) -- (2.5,3.5);
  \draw[very thick,->] (2.5,3.5) -- (2.5,4.5);
  \draw[very thick,->] (2.5,4.5) -- (-0.5,4.5);
  \draw[very thick,->] (-0.5,4.5) -- (-0.5,-0.5);
  %% inner contour around the hole at (2,2), clockwise
  \draw[very thick,->] (1.5,1.5) -- (1.5,2.5);
  \draw[very thick,->] (1.5,2.5) -- (2.5,2.5);
  \draw[very thick,->] (2.5,2.5) -- (2.5,1.5);
  \draw[very thick,->] (2.5,1.5) -- (1.5,1.5);
  %% sites
  \foreach \x/\y in {0/0,1/0,2/0,3/0,0/1,0/2,0/3,3/1,1/3,2/3} \fill (\x,\y) circle (0.09);
  \foreach \x/\y in {1/1,2/1,1/2,4/0,4/1,3/2,0/4,1/4,2/4,3/3} \draw[thick] (\x,\y) circle (0.11);
  \draw[thick] (2,2) circle (0.11); \draw[thick] (1.92,1.92) -- (2.08,2.08);
  \fill[white] (0,0) circle (0.045);
  \node[below left=-1pt] at (0,0) {$0$};
  \node[right] at (4.6,0.5) {$\Gamma$};
  \node[right,align=left] at (5.2,3.0) {\small $\bullet$ good site in $S$\\[-1pt] \small $\circ$ bad site in $S$\\[-1pt] \small $\oslash$ site not in $S$};
\end{tikzpicture}
\caption{\small An illustration of Lemma~\ref{lem:contour} with $z=0$, for the set $S$ of the proof of Lemma~\ref{lem:oriented}. The squares $Q_y$, $y\in S$, are shaded. The $1$-chain decomposes into the outer boundary $\Gamma$, oriented counterclockwise, and the boundary of the hole, oriented clockwise, which has winding number $0$ about $0$. Each upward or leftward edge of $\Gamma$ has a square centered at a site of $S$ on its left and a square centered at a site outside $S$ on its right. In the proof of Lemma~\ref{lem:oriented} the site on the left is bad.}
\label{fig:contour}
\end{figure}

\begin{proof}
The proof first constructs $\Gamma$ and then verifies the five items in order. See Figure~\ref{fig:contour}.

Orient the boundary of each square $Q_y$, $y\in S$, counterclockwise, and form the $1$-chain which is the sum of these boundaries: an edge shared by two of the squares appears with both orientations and cancels. Each remaining directed edge is a side of exactly one square $Q_w$ with $w\in S$, which lies on its left by the orientation convention, and the square on its right is centered at a point outside $S$. These remaining directed edges form a finite directed graph in which every vertex has in-degree equal to out-degree (this holds for each square boundary and is preserved under cancellation), so, by Euler's theorem for directed graphs \cite[Theorem~1.6.3]{Ban-Gut-09} applied to each connected component,  the edge set decomposes into circuits.
%The theorem says that a connected directed multigraph in which every vertex has in-degree equal to out-degree has a closed trail that uses every edge exactly once. Applied to each connected component of our graph, it gives one circuit per component. Here is the construction:
%Start with any edge and keep following unused outgoing edges, never repeating an edge, until no continuation is possible. Finiteness guarantees termination. The circuit must terminate at its starting vertex. Indeed, if it ended at another vertex, the circuit would have used one more incoming edge than outgoing edge at that vertex. Since the total incoming and outgoing degrees are equal, an unused outgoing edge would remain, causing a contradiction. Thus we have a closed directed loop. Removing its edges preserves equality of the incoming and outgoing degrees at every vertex. Repeat until all edges have been removed.
The winding number about the point $z$ of the sum of the square boundaries is $1$, because only the boundary of $Q_z$ contributes, and the winding number is additive, so some circuit $\Gamma$ in the decomposition has nonzero winding number about $z$. Let $m$ be its length. Item~\ref{con:sides} holds because it holds for every edge of the chain.

\smallskip
\emph{Items~\ref{con:forward} and~\ref{con:backward}.} Since $\Gamma$ is closed, it has as many upward as downward edges and as many leftward as rightward edges. So $m$ is even, exactly $m/2$ edges of $\Gamma$ are upward or leftward, and exactly $m/2$ are downward or rightward. An upward edge is the right side of the square $Q_w$ on its left, and the square on its right is $Q_{w+e_1}$. A leftward edge is the top side of the square $Q_w$ on its left, and the square on its right is $Q_{w+e_2}$. In both cases $w\in S$ and $w+e_i\notin S$, by item~\ref{con:sides}. A square has one right side and one top side, so the $m/2$ upward or leftward edges give at least $m/4$ distinct points $w$. This is item~\ref{con:forward}. Similarly, a downward edge is the right side of the square $Q_{w'}$ on its right, and the square on its left is $Q_{w'+e_1}$. A rightward edge is the top side of the square $Q_{w'}$ on its right, and the square on its left is $Q_{w'+e_2}$. This gives item~\ref{con:backward}.

\smallskip
\emph{Item~\ref{con:size}.} We first show that $\Gamma$ has nonzero winding number about every point of $S$. Let $y,y'\in S$ be nearest neighbors. The segment $[y,y']$ meets the edges only at its midpoint, which is an interior point of the common side of $Q_y$ and $Q_{y'}$. This side cancels in the chain, so it is not an edge of $\Gamma$, and $[y,y']$ does not meet $\Gamma$. Hence the winding numbers of $\Gamma$ about $y$ and about $y'$ agree. Since $S$ is connected, the winding number of $\Gamma$ about every point of $S$ equals the one about $z$.

Now consider a vertical line that contains a point of $S$, and let $y$ be the highest point of $S$ on it. The upward vertical half-line from $y$ meets the edges at the midpoints of the top sides of the squares $Q_{y+je_2}$, $j\ge0$. For $j\ge1$ this top side is shared by two squares centered outside $S$ and is not in the chain. Since $\Gamma$ has nonzero winding number about $y$, it meets the half-line, so it contains the top side of $Q_y$. In the same way, $\Gamma$ contains the bottom side of the square centered at the lowest point of $S$ on the line. Hence $\Gamma$ has at least $2a$ horizontal edges, where $a$ is the number of vertical lines that contain a point of $S$, and likewise at least $2b$ vertical edges, where $b$ is the number of horizontal lines that contain a point of $S$. Therefore $m\ge2(a+b)$ and $\abs S\le ab\le(a+b)^2/4\le m^2/16$. In particular $m\ge4$.

\smallskip
\emph{Item~\ref{con:count}.} Since $\Gamma$ has nonzero winding number about $z$, it meets the horizontal half-line from $z$ to the right at some point $x$ and the horizontal half-line from $z$ to the left at some point $x'$. These two points divide $\Gamma$ into two arcs, each of length at least $\abs{x-x'}_1>\abs{x-z}_1$, so $\abs{x-z}_1<m/2$. The half-line to the right meets the edges at the points $z+(j-\tfrac12)e_1$, $j\ge1$, which are the midpoints of the common sides of $Q_{z+(j-1)e_1}$ and $Q_{z+je_1}$. So $\Gamma$ contains such a side with $j-\tfrac12<m/2$, that is, with $j\le m/2$, because $m$ is even. A circuit of length $m$ that contains a given edge is determined by the direction in which it traverses this edge and by the choices of the $m-1$ edges that follow, and each of these is one of at most three edges, because it starts at the endpoint of the previous edge and is not the previous edge reversed. So there are at most $\frac m2\cdot2\cdot3^{m-1}$ circuits of length $m$ that contain one of the $m/2$ edges in question.
\end{proof}

\begin{lemma}[high-density oriented paths]\label{lem:oriented}
For every integer $\varkappa\ge0$ there exists $\zeta=\zeta(\varkappa)>0$ with the following property. Let $(G_v)_{v\tsp\in\tsp\Z^2_{\ge0}}$ be events such that every finite family indexed by sites at pairwise sup-distance greater than $\varkappa$ is independent, and $\P(G_v^c)\le\zeta$ for all $v\in\Z^2_{\ge0}$. Then with positive probability there exists an infinite path $v_{0:\infty}$ with $v_0=0$ and steps $v_{i+1}-v_i\in\{e_1,e_2\}$ such that $G_{v_i}$ occurs for every $i\ge0$.
\end{lemma}
\begin{proof}
Call $v$ good if $G_v$ occurs and bad otherwise, and call a path with steps $e_1,e_2$ oriented. Suppose there is no infinite good oriented path from $0$. Let $S$ be the set of endpoints of finite oriented paths from $0$ all of whose sites, except possibly the last, are good. Thus $0\in S$, and if $0$ is bad then $S=\{0\}$. If $S$ were infinite, there would be good oriented paths from $0$ of arbitrary length, since a site of $S$ on level $k$ is reached by a good oriented path of length $k-1$. Identify an oriented path from $0$ with its sequence of steps in $\{e_1,e_2\}$. By compactness of $\{e_1,e_2\}^{\N}$ in the product topology, a sequence of good oriented paths with lengths tending to infinity has a subsequence that converges coordinatewise, and the limit is an infinite oriented path from $0$. Every site of the limit path is a site of all the paths far enough along the subsequence, hence good. This contradicts the assumption, so $S$ is finite.

The set $S$ contains $0$, and it is connected, because every site of an oriented path from $0$ all of whose sites, except possibly the last, are good belongs to $S$. Let $\Gamma$ be the circuit of Lemma~\ref{lem:contour} for this $S$ and $z=0$, and let $m$ be its length. If $w\in S$ and $w+e_i\notin S$ for some $i\in\{1,2\}$, then $w$ is bad: if $w$ were good, then (as $w\in S$ is the endpoint of an oriented path from $0$ whose other sites are good) $w+e_i$ would be the endpoint of an oriented path from $0$ with all sites except possibly the last good, that is, $w+e_i\in S$. Hence, by Lemma~\ref{lem:contour}\ref{con:forward}, the centers of the squares on the left of the upward and leftward edges of $\Gamma$ are at least $m/4$ distinct bad sites. Among any set of sites one can select greedily a subset at pairwise sup-distance greater than $\varkappa$ of size at least the original size divided by $(2\varkappa+1)^2$, so $\Gamma$ determines at least $m/(4(2\varkappa+1)^2)$ bad sites at pairwise sup-distance greater than $\varkappa$. For a fixed circuit these sites are deterministic, so by the independence hypothesis
\[
\P(\text{all sites determined by }\Gamma\text{ are bad})\le\zeta^{m/(4(2\varkappa+1)^2)} .
\]
By Lemma~\ref{lem:contour}\ref{con:count}, there are at most $m\,3^{m-1}$ circuits of length $m$ to consider. Summing over $m\ge4$,
\[
\P(\text{no infinite good oriented path from }0)\le\sum_{m\ge4}m\bigl(3\zeta^{1/(4(2\varkappa+1)^2)}\bigr)^{m-1}\zeta^{1/(4(2\varkappa+1)^2)} ,
\]
which is strictly less than $1$ once $\zeta=\zeta(\varkappa)$ is small enough. \end{proof}

\subsection{Van den Berg--Kesten inequality for path events}\label{sub:bk-appendix}
Recall from Section~\ref{sub:bk} that a path event $\cA$ is the event that $\pawgt(x)>\lambda(x)$ for some path $x$ in a finite family $\Pi$ of admissible paths with all collected vertices in a finite set $\vbad$, where $\lambda:\Pi\to\R$, that a witness of $\cA$ is such a path, and that for path events $\cA_1,\dots,\cA_q$, $D(\cA_1,\dots,\cA_q)$ is the event that there exist witnesses $x^{(1)},\dots,x^{(q)}$ of $\cA_1,\dots,\cA_q$ with pairwise disjoint sets of collected vertices. 

The following lemma is the van den Berg--Kesten inequality \cite{Ber-Kes-85} for path events. Since the weights are real-valued, it does not follow directly from the original statement, which concerns increasing events in a finite product of finite spaces. It is a special case of the inequality in \cite[Theorem~1.3]{Gol-Rin-07}, which allows arbitrary finite families of witnesses. %Indeed, a path event is increasing in $\w$, and a witness of it remains a witness under every change of the weights outside its collected vertices that does not decrease them.
For completeness we include a short self-contained proof, which reduces to the finite case by discretizing the weights.

\begin{lemma}[BK]\label{lem:bk}
$\P(D(\cA_1,\dots,\cA_q))\le\prod_{j=1}^q\P(\cA_j)$.
\end{lemma}
\begin{proof}
Let $\vbad$ be the union of the vertex sets of the $\cA_j$. Realize the weights as $\w_v=F_\mu^{-1}(U_v)$, $v\in \vbad$, with $(U_v)$ i.i.d.\ uniform on $(0,1)$. This has the correct law because $F_\mu^{-1}(U)\le t$ if and only if $U\le F_\mu(t)$. Write $U_v=\sum_{h\ge1}2^{-h}B_{v,h}$ with i.i.d.\ fair bits $B_{v,h}$ (the binary digits of $U_v$, where almost surely $U_v\in(0,1)$ is not a dyadic rational). For a precision $t\in\N$ set $U_v^{(t)}=\sum_{h\le t}2^{-h}B_{v,h}$ and $\w_v^{(t)}=F_\mu^{-1}(U_v^{(t)})$, with the convention $F_\mu^{-1}(0)=-\infty$. Since $U^{(t)}_v\uparrow U_v$ as $t\to\infty$ and $F_\mu^{-1}$ is nondecreasing and left-continuous, $\w_v^{(t)}\uparrow\w_v$. In particular, $\w^{(t)}\le\w$ and $\w^{(t)}_v>-\infty$ for all large $t$. Let $\cA_j^{(t)}$ and $D^{(t)}$ be the events defined as $\cA_j$ and $D$ but with the weights $\w^{(t)}$ (sums involving $-\infty$ are $-\infty$ and never exceed a real bound).

At fixed precision $t$, the events $\cA_j^{(t)}$ are increasing events in the finitely many independent bits $(B_{v,h})_{v\tsp\in\tsp \vbad,h\le t}$, and a witness $x$ of $\cA^{(t)}_j$ occurs ``on'' the set of bits $(v,h)$ with $v$ a collected vertex of $x$ and $h\le t$, in the sense that every configuration of bits agreeing with the present one on this set also has $x$ as a witness. Disjoint sets of collected vertices give disjoint sets of bits, so $D^{(t)}$ is contained in the disjoint occurrence $\cA_1^{(t)}\square\cdots\square\cA_q^{(t)}$. The van den Berg--Kesten inequality \cite{Ber-Kes-85} for increasing events in a finite product of two-point spaces, iterated (disjoint occurrence is associative), gives
\[
\P(D^{(t)})\le\prod_j\P(\cA_j^{(t)})\le\prod_j\P(\cA_j),
\]
the last step because $\w^{(t)}\le\w$ and the $\cA_j$ are increasing in $\w$.

Finally, $D^{(t)}\subset D^{(t+1)}\subset D$: every witness at precision $t$ is a witness at precision $t+1$ and in the true environment, since weights only increase with $t$. Conversely, on $D$ fix witnesses $x^{(1)},\dots,x^{(q)}$ with pairwise disjoint sets of collected vertices. Each satisfies a strict inequality $\pawgt(x^{(j)})>\lambda(x^{(j)})$ with finitely many terms, so it holds at all sufficiently fine precisions, and $\w\in D^{(t)}$ for large $t$. Hence $D^{(t)}\uparrow D$ and $\P(D)=\lim_t\P(D^{(t)})\le\prod_j\P(\cA_j)$.
\end{proof}

\section{Extension to unbounded weights}
\label{app:general}

The main text proves Lemma~\ref{lem:intervals}, Proposition~\ref{prop:endpoints}, and Theorem~\ref{thm:main} under the boundedness assumption \eqref{ass:bounded}, and the remarks titled ``Bounded vs.\ unbounded weights'' record, at each step that uses the bound, what changes under \ref{A2}--\ref{A3}. This appendix fills in those details. The lower bound $-M$ is replaced everywhere as follows: the vertices with very negative weight are observed first,  the paths are then chosen so as to avoid them, and Lemma~\ref{lem:weights} bounds the weights of such paths from below. Section~\ref{sub:attachments} proves Lemma~\ref{lem:intervals} in this way. Section~\ref{sub:staircases} does the construction of Section~\ref{sec:endpoints}  with  deterministic staircases replaced by paths that avoid the very negative vertices. There the upper bound $M$ is removed as well, at the price of trimming a few vertices from the ends of some paths. Section~\ref{sub:onestep} removes the upper bound from the one-step bound of Section~\ref{sec:witnesses} with a  greedy lattice animal  bound.

\subsection{Paths that avoid very bad vertices}\label{sub:verybad}
In this section and the next we assume 
\begin{equation}\label{ass:general}
\begin{split}
&\mu\ \text{satisfies the hypotheses of the shape theorem \cite[Theorem~5.1 and Remark~5.1(ii)]{Mar-04},}\\
&\text{and }\E\bigl[(\w_0^-)^{p'}\bigr]<\infty\ \text{ for some }p'>1 ,
\end{split}
\end{equation}
which holds under \ref{A2}--\ref{A3}. Nothing is assumed about the upper tail of $\mu$ beyond what the shape theorem needs. Sections~\ref{sub:staircases} and~\ref{sub:onestep} assume \ref{A3} as well. Set
\begin{equation}\label{eq:C1}
\Cm=\E\w_0^-\quad\text{and}\quad C_{p'}=\E\bigl[(\w_0^-)^{p'}\bigr],\qquad\text{so that}\quad \P(\w_0<-t)\le C_{p'}t^{-p'}\quad\text{for all }t>0,
\end{equation}
and
\begin{equation}\label{eq:epsprime}
\e'=\frac{p'-1}{2p'},\qquad\text{so that}\quad p'(1-\e')=\frac{p'+1}2>1 .
\end{equation}

The boundedness assumption \eqref{ass:bounded} is used in the main text in three places: the bound \eqref{attachments} for the attachments in the proof of Lemma~\ref{lem:intervals}, the bound \eqref{eq:badstair} for the staircases of Section~\ref{sec:endpoints}, and the one-step bound \eqref{eq:one-step} of Section~\ref{sec:witnesses}. The first two are Hoeffding bounds on the weight of a deterministic path of length $n$, with an error probability small enough to survive a union bound over polynomially many paths, and since a lower bound on the weight is all that is needed there, only the lower bound $-M$ matters. This fails under \eqref{ass:general}. For a fixed $\delta>0$, a deterministic path of length $n$ collects a vertex of weight below $-\delta n$, which by itself lowers the weight of the path by $\delta n$, with a probability that can be of order $n(\delta n)^{-p'}$. This is only polynomially small in $n$, and in general it does not survive the union bound. The replacement follows the same pattern in both places. For a threshold $T>0$, a vertex $z$ is \emph{very bad} if $\w_z<-T$, and $\vbad =\vbad _T$ is the set of very bad vertices. Let $\cF_\vbad$ be the $\sigma$-algebra generated by $\vbad$, that is, by the indicators $\one\{\w_z<-T\}$, $z\in\Z^2$. Conditionally on $\cF_\vbad$, the weights $(\w_z)_{z\notin \vbad}$ are i.i.d.\ with the conditional law
\[
\mu_T=\mu(\,\abullet\,|\,[-T,\infty)) .
\]
The paths of interest are chosen after $\vbad$ is observed, as functions of $\vbad$, so as to collect no very bad vertex. The geometric constructions that produce such paths with the required endpoints are the subject of Sections~\ref{sub:attachments} and~\ref{sub:staircases}. The next lemma bounds the weights of such paths from below. The weights they collect are, conditionally on $\cF_\vbad$, i.i.d.\ with law $\mu_T$ and hence bounded below by $-T$, so a Bernstein bound applies once the few terms with a large negative part are set aside.
Below $\abs\alpha$ denotes the length of a path $\alpha$.

\begin{lemma}\label{lem:weights}
Assume \eqref{ass:general}. Let $T>0$ be such that $\mu([-T,\infty))\ge\tfrac12$, and let $\cP$ be an $\cF_\vbad$-measurable family of at most $N$ admissible paths, each of length at most $n$ and collecting no vertex of $\vbad$. Then for every integer $K_0\ge1$ and all $b,t>0$,
\begin{equation}\label{eq:weights}\begin{aligned}
&\P\Bigl\{\exists\,\alpha\in\cP:\ \pawgt(\alpha)<-2\Cm\abs\alpha-(K_0-1)T-t\Bigr\}\\[3pt] 
&\qquad \le\ N\Bigl[\bigl(2nC_{p'}b^{-p'}\bigr)^{K_0}+\exp\Bigl\{-\frac{t^2}{2b\,(2\Cm n+t/3)}\Bigr\}\Bigr]. 
\end{aligned}\end{equation}
\end{lemma}

\begin{remark}
In both applications below, $n$ is of order $s$ up to powers of $\log s$, $t$ is of order $s$, and $b=s^{1-\e'}$. Then $b^{-p'}=s^{-(p'+1)/2}$ by \eqref{eq:epsprime}, so the first term in \eqref{eq:weights} is polynomially small in $s$, with an exponent that grows with $K_0$, and the second term is $\exp\{-cs^{\e'}\}$. The price is the term $(K_0-1)T$ in the lower bound, which is why $T$ is taken to be a small multiple of the allowed error.
\end{remark}

\begin{proof}[Proof of Lemma \ref{lem:weights}]
Fix $\alpha\in\cP$. Since $\alpha$ is $\cF_\vbad$-measurable and collects $\abs\alpha\le n$ distinct vertices outside $\vbad$, conditionally on $\cF_\vbad$ its weight is a sum $\sum_{i=1}^{\abs\alpha}\w_i$ of $\abs\alpha$ i.i.d.\ $\mu_T$-distributed variables. Since $\mu([-T,\infty))\ge\tfrac12$, $\mu_T(A)\le2\mu(A)$ for every Borel $A\subset[-T,\infty)$. Write $\w_i=\w_i^+-\w_i^-$ and split $\w_i^-=Y_i+Z_i$ with
\[
Y_i=\w_i^-\one\{\w_i^->b\}\quad\text{and}\quad Z_i=\w_i^-\one\{\w_i^-\le b\}.
\]
Then $0\le Y_i\le T$ and $0\le Z_i\le b$, $\mu_T$-almost surely.

By \eqref{eq:C1}, $\mu_T(Y_i>0)\le2\P(\w_0<-b)\le2C_{p'}b^{-p'}$. The number of indices with $Y_i>0$ is binomial with $\abs\alpha\le n$ trials, so 
\[
\P\bigl(\abs{\{i:Y_i>0\}}\ge K_0\,\big|\,\cF_\vbad \bigr)\le\binom{n}{K_0}\bigl(2C_{p'}b^{-p'}\bigr)^{K_0}\le\bigl(2nC_{p'}b^{-p'}\bigr)^{K_0} .
\]
On the complement, $\sum_{i=1}^{\abs\alpha}Y_i\le(K_0-1)T$. 

The $Z_i$ are i.i.d.\ under $\mu_T$ with $0\le Z_i\le b$, $\E_{\mu_T}Z_i\le2\E\w_0^-=2\Cm$, and $\E_{\mu_T}Z_i^2\le b\,\E_{\mu_T}Z_i\le2\Cm b$. Bernstein's inequality \cite[Theorem~2.10]{Bou-Lug-Mas-13} gives
\[
\P\Bigl(\sum_{i=1}^{\abs\alpha}Z_i>2\Cm\abs\alpha+t\,\Big|\,\cF_\vbad \Bigr)\le\exp\Bigl\{-\frac{t^2}{2(2\Cm\abs\alpha\,b+bt/3)}\Bigr\}\le\exp\Bigl\{-\frac{t^2}{2b\,(2\Cm n+t/3)}\Bigr\}.
\]

Thus, outside the above two exceptional events, $\pawgt(\alpha)\ge-\sum_{i=1}^{\abs\alpha}Y_i-\sum_{i=1}^{\abs\alpha}Z_i\ge-2\Cm\abs\alpha-(K_0-1)T-t$. Taking expectations and summing over the at most $N$ paths of $\cP$ gives \eqref{eq:weights}.
\end{proof}

\subsection{Attachments}\label{sub:attachments}
In this section we prove Lemma~\ref{lem:intervals} under \eqref{ass:general}.

\begin{lemma}\label{lem:intervals-general}
Under \eqref{ass:general}, Lemma~\ref{lem:intervals} holds as stated.
\end{lemma}

As explained in Remark~\ref{rk:bounded1}, the boundedness assumption enters the proof of Lemma~\ref{lem:intervals} only through the Hoeffding bound \eqref{attachments}, which controls the weights of the deterministic attachments $\alpha_{u:\hat u}$ and $\alpha_{\hat v:v}$ simultaneously for all $u\in\cJ_1$ and $v\in\cJ_2$. For unbounded weights, a deterministic path of length $m$ contains a vertex of weight $-s$ or worse with a probability that can be of order $ms^{-p'}$, and there are $O(s)$ attachments, so the union bound fails when $p'\le2$. The attachments are therefore chosen after the weights are observed, so as to step around the very bad vertices. Lemma~\ref{lem:attach-general} below provides them. The rest of the proof of Lemma~\ref{lem:intervals} needs only small adjustments, listed after the lemma.

The lemma is formulated for the objects of the proof of Lemma~\ref{lem:intervals}, namely $\cJ_1$ and $\cN_1$ there, see Figure~\ref{fig:intervals}. A level interval $\cJ=\{\level=I,\ \rcoord\in[\rho_1,\rho_2]\}$ with $\rho_2-\rho_1\le Ds$ is given, together with an integer $m\ge1$ and a deterministic set $\cN$ of lattice points on level $I+m$, the \emph{net}, whose transverse coordinates increase with consecutive gaps in $[\theta_0m/4,\theta_0m/2]$, exactly one of them at most $\rho_1-\Dir_1m+\theta_0m/2$ and exactly one at least $\rho_2+(1-\Dir_1)m-\theta_0m/2$.

\begin{lemma}\label{lem:attach-general}
Assume \eqref{ass:general}. Fix $D\ge1$, $\eta>0$, and $\kappa\in(0,\tfrac1{16}]$ with $\kappa\le\eta/(16\Cm+16)$, and let $m=\lceil\kappa s\rceil$. There exists a function $p_{\rm att}(s)\to1$ as $s\to\infty$, depending only on $D,\eta,\kappa,\theta_0,\mu$, such that the following holds for every $s\ge4$ and every $\cJ$ and $\cN$ as above. With probability at least $p_{\rm att}(s)$, for every $u\in\cJ$ there exist a point $\hat u\in\cN$ and an admissible path $\alpha_{u:\hat u}$ from $u$ to $\hat u$ with
\begin{equation}\label{eq:attach-general}
\pawgt(\alpha_{u:\hat u})\ \ge\ -\frac{\eta s}{4},\quad \abs{\hat u-u}_1=m,\quad\text{and}\quad \abs{\rcoord(\hat u-u)}\le m .
\end{equation}
The symmetric statement holds for paths ending at the points of a level interval and starting at a net on the level $m$ below it.
\end{lemma}

\begin{figure}[htb]
\centering
\begin{tikzpicture}[x=0.5cm,y=0.33cm,>=Stealth]
  %% coordinates (r, level) with Dir_1=1/2, so that e_1 and e_2 steps go up-right and up-left; not to scale
  \foreach \y/\lab in {0/$I$, 16/$I+m$} {\draw[gray!50] (-9,\y) -- (9,\y); \node[right,gray] at (9,\y) {\footnotesize \lab};}
  %% the tube Sigma_u around the segment gamma_u, and the box at the net point
  \fill[gray!14] (-1.5,0) -- (1.5,0) -- (3.5,16) -- (0.5,16) -- cycle;
  \fill[gray!32] (2,8) -- (4,12) -- (2,16) -- (0,12) -- cycle;
  \draw[gray] (2,8) -- (4,12) -- (2,16) -- (0,12) -- cycle;
  \draw[dashed] (0,0) -- (2,16);
  %% the level interval J and the net N
  \draw[very thick] (-6,0) -- (6,0);
  \node[below] at (-5,0) {\footnotesize $\cJ$};
  \foreach \x in {-8,-6,...,8} \fill (\x,16) circle (0.08);
  \node[above] at (-7,16) {\footnotesize $\cN$};
  %% the greedy path and the box phase
  \draw[very thick] (0,0) -- (0.5,1) -- (0,2) -- (-0.5,3) -- (0,4) -- (0.5,5) -- (1,6) -- (1.5,7) -- (1,8) -- (1.5,9) -- (1,10) -- (1.5,11) -- (2,12) -- (2.5,13) -- (3,14) -- (2.5,15) -- (2,16);
  %% very bad vertices: inside the tube they force a step, outside it they play no role
  \foreach \p in {(0.5,3),(0.5,7),(-2,6),(3.5,5),(-3,2),(5.5,1),(-4.5,7),(2.5,1),(6.5,5),(-6,4)} \draw[thick] \p ++(-0.17,-0.25) -- ++(0.34,0.5) \p ++(-0.17,0.25) -- ++(0.34,-0.5);
  \fill (0,0) circle (0.1) node[below=2pt] {\footnotesize $u$};
  \fill (2,16) circle (0.1) node[above=2pt] {\footnotesize $\hat u$};
  %% the window of points reachable from u in m steps, and its middle half
  \draw[gray,decorate,decoration={brace,amplitude=4pt}] (-8,20.2) -- (8,20.2) node[midway,above=4pt] {\footnotesize window of $u$};
  \draw[gray] (-4,18.4) -- (-4,18.9) -- (4,18.9) -- (4,18.4);
  \node[gray,right] at (4.1,18.9) {\scriptsize middle half};
  %% labels
  \node[left] (glab) at (-2.6,13.5) {\footnotesize $\gamma_u$};
  \draw[->,gray,shorten >=1pt] (glab) -- (1.65,13.2);
  \node[left] (slab) at (-3.4,10.6) {\footnotesize $\Sigma_u$};
  \draw[->,gray,shorten >=1pt] (slab) -- (-0.1,10.4);
  \node[right] (blab) at (5.3,9.6) {\footnotesize $\hat u-[0,n_{\rm box}]^2$};
  \draw[->,gray,shorten >=1pt] (blab) -- (2.9,10.4);
  \node[right] (flab) at (2.6,2.4) {\footnotesize forced steps};
  \draw[->,gray,shorten >=3pt] (flab) -- (0.5,3);
  \draw[->,gray,shorten >=3pt] (flab) -- (0.5,7);
\end{tikzpicture}
\caption{\small The construction in the proof of Lemma~\ref{lem:attach-general}, drawn in the coordinates $(\rcoord,\level)$ and not to scale. The net point $\hat u$ is chosen in the middle half of the window of lattice points reachable from $u$ in $m$ steps. The very bad vertices (crosses) are observed first. The attachment $\alpha_{u:\hat u}$ (thick) follows the segment $\gamma_u$ (dashed) greedily, and a very bad vertex in its way forces a step to the other side. On $E_1\cap E_2\cap E_3$ there are fewer than $K_0$ forced steps, so the path stays in the tube $\Sigma_u$ of transverse half-width $K_0+2$ (light gray), and its last $n_{\rm box}$ steps lie in the box $\hat u-[0,n_{\rm box}]^2$ (darker gray), which contains no very bad vertex by $E_2$. Very bad vertices outside $\Sigma_u$ play no role.}
\label{fig:attach-general}
\end{figure}
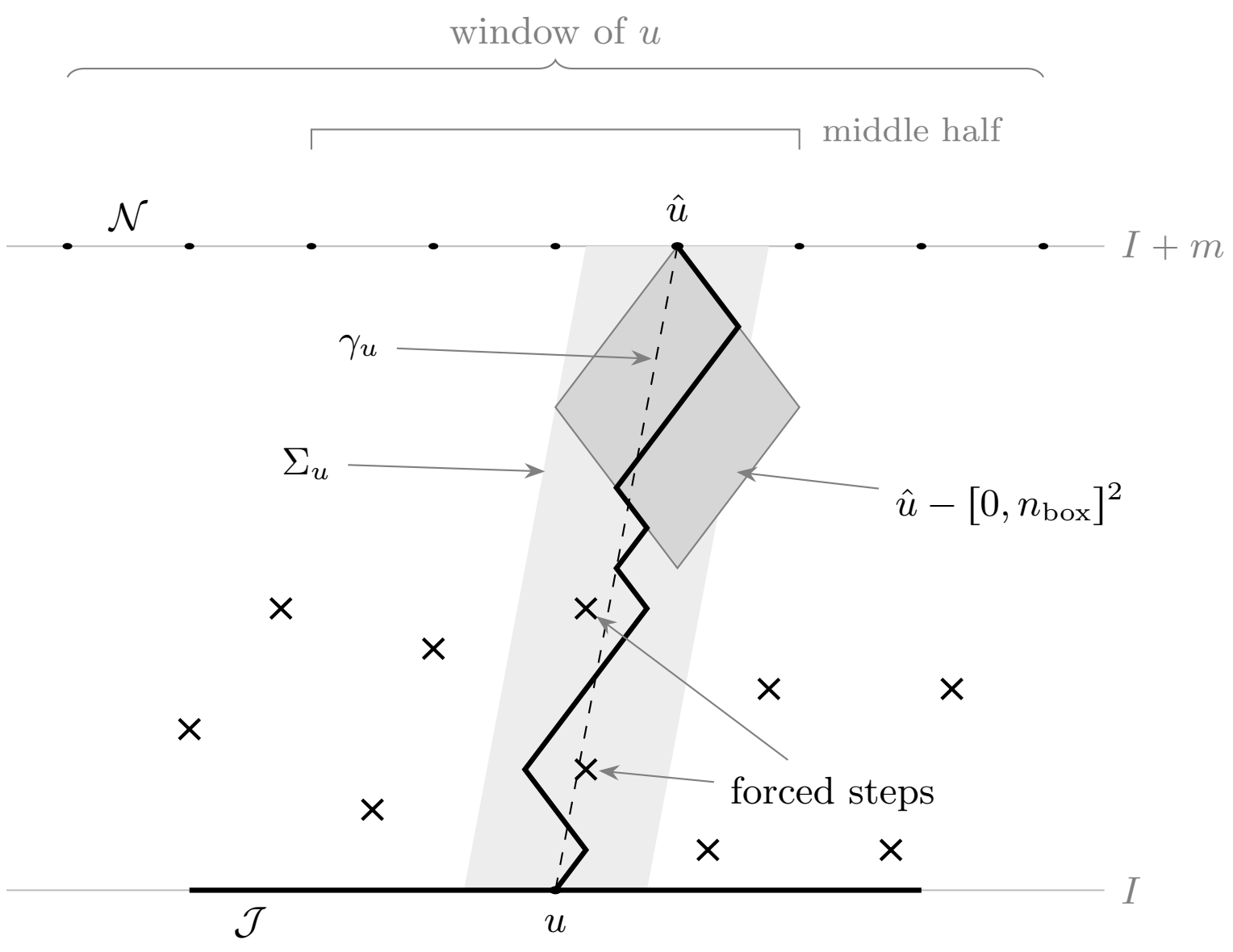

\begin{proof}
Let $K_0\ge2$ be an integer with $K_0(p'-1)\ge4$, let $n_{\rm box}=4K_0+4$, and set
\[
T=\frac{\eta s}{64K_0}\qquad\text{and}\qquad p_T=\P(\w_0<-T)\le C_{p'}\Bigl(\frac{64K_0}{\eta}\Bigr)^{p'}s^{-p'} .
\]
A vertex $z$ is very bad if $\w_z<-T$, and $\vbad$ is the set of very bad vertices. All paths constructed below lie in the region
\[
\cQ=\{z\in\Z^2:\ I\le\level(z)\le I+m,\ \rho_1-m-K_0-2\le\rcoord(z)\le\rho_2+m+K_0+2\},
\]
which has at most $(m+1)(Ds+2m+2K_0+5)\le4Ds^2$ vertices once $s$ is so large that $2K_0+5\le m$.

For $u\in\cJ$, the lattice points reachable from $u$ in $m$ steps are $u+(j,m-j)$, $0\le j\le m$, with transverse coordinates $\rcoord(u)-\Dir_1m+j$. Those with $m/4\le j\le3m/4$ have transverse coordinates in an interval of length $m/2$, the \emph{middle half} of the window $[\rcoord(u)-\Dir_1m,\rcoord(u)+(1-\Dir_1)m]$, at distance $m/4\ge\theta_0m/2$ from both ends of the window. Since the window is contained in $[\rho_1-\Dir_1m,\rho_2+(1-\Dir_1)m]$, the middle half lies in
\[
\cR_1=[\rho_1-\Dir_1m+\theta_0m/2,\ \rho_2+(1-\Dir_1)m-\theta_0m/2],
\]
which is contained in the span of $\cN$, and since consecutive net coordinates are at most $\theta_0m/2\le m/2$ apart, the middle half contains the transverse coordinate of a net point. Fix such a point $\hat u\in\cN$, so that
\begin{equation}\label{eq:middle}
\hat u-u=(j_u,m-j_u)\qquad\text{with}\quad m/4\le j_u\le3m/4 .
\end{equation}
The last two assertions of \eqref{eq:attach-general} are immediate. Let $\gamma_u$ be the straight segment from $u$ to $\hat u$, and for $I\le i\le I+m$ let $r_u(i)=\rcoord(u)+(i-I)\,\rcoord(\hat u-u)/m$ be the transverse coordinate of the point of $\gamma_u$ on level $i$. By \eqref{eq:middle}, the slope $\rcoord(\hat u-u)/m=j_u/m-\Dir_1$ lies in $[\tfrac14-\Dir_1,\tfrac34-\Dir_1]$, so that
\begin{equation}\label{eq:slope}
\rcoord(e_1)-\frac{\rcoord(\hat u-u)}m\in\Bigl[\frac14,\frac34\Bigr]\qquad\text{and}\qquad \frac{\rcoord(\hat u-u)}m-\rcoord(e_2)\in\Bigl[\frac14,\frac34\Bigr].
\end{equation}
Let $\Sigma_u=\{z\in\cQ:\ \abs{\rcoord(z)-r_u(\level(z))}\le K_0+2\}$ be the tube of transverse half-width $K_0+2$ around $\gamma_u$. It has at most $2K_0+5$ vertices on each level, hence at most $(2K_0+5)(m+1)$ vertices.

Define the three events
\begin{align*}
E_1&=\{\text{no two very bad vertices of $\cQ$ are at $\abs\cdot_1$-distance at most $2$}\},\\
E_2&=\{\text{no vertex of $\cJ$ and no vertex of the boxes }\hat u-[0,n_{\rm box}]^2,\ \hat u\in\cN,\text{ is very bad}\},\\
E_3&=\{\text{for every $u\in\cJ$, the tube $\Sigma_u$ contains fewer than $K_0$ very bad vertices}\}.
\end{align*}
A vertex has $12$ other vertices within $\abs\cdot_1$-distance $2$, so by a union bound and independence, \[\P(E_1^c)\le4Ds^2\cdot12\,p_T^2=q_1(s).\] 
The net $\cN$ has at most $4(Ds+m)/(\theta_0m)+2\le8D/(\theta_0\kappa)+2$ points, because all but its two extreme points have coordinates in $\cR_1$, an interval of length at most $Ds+m$, and consecutive gaps are at least $\theta_0m/4$. So 
\[\P(E_2^c)\le\bigl(Ds+1+(8D/(\theta_0\kappa)+2)(n_{\rm box}+1)^2\bigr)p_T=q_2(s).\] 
Finally, for a fixed $u$ the number of very bad vertices in $\Sigma_u$ is binomial with at most $(2K_0+5)(m+1)$ trials and success probability $p_T$, so
\[
\P(E_3^c)\le(Ds+1)\binom{(2K_0+5)(m+1)}{K_0}p_T^{K_0}\le(Ds+1)\bigl((2K_0+5)(m+1)p_T\bigr)^{K_0}=q_3(s) .
\]
Since $p_T\le Cs^{-p'}$ with $p'>1$ and $K_0(p'-1)\ge4$, the three bounds $q_1(s)$, $q_2(s)$, and $q_3(s)$ are $O(s^{2-2p'})$, $O(s^{1-p'})$, and $O(s^{1-K_0(p'-1)})$, all $o(1)$ as $s\to\infty$, with constants depending only on $D,\eta,\kappa,\theta_0,\mu$.\smallskip

We now construct the attachments, see Figure~\ref{fig:attach-general}. On $E_1\cap E_2\cap E_3$, fix $u\in\cJ$ and build $\alpha_{u:\hat u}$ step by step, as a function of $\vbad$ only. Write $x_i$ for its vertex on level $i$ and $(a_i,b_i)=\hat u-x_i\in\Z^2$ for the remaining displacement, so that $a_i+b_i=I+m-i$. Start with $x_I=u$ and proceed in two phases.\smallskip

\emph{Greedy phase}, while $a_i+b_i>n_{\rm box}$: of the two candidates $x_i+e_1$ and $x_i+e_2$ for $x_{i+1}$, the \emph{preferred} one is the one whose transverse coordinate is closer to $r_u(i+1)$ (either one in case of a tie). If the preferred candidate is not very bad, move there. Otherwise the step is \emph{forced}: move to the other candidate, which is not very bad by $E_1$, since the two candidates are at $\abs\cdot_1$-distance $2$.\smallskip

\emph{Box phase}, once $a_i+b_i=n_{\rm box}$: provided $a_i,b_i\ge0$, finish with $a_i$ steps $e_1$ followed by $b_i$ steps $e_2$. These vertices lie in the box $\hat u-[0,n_{\rm box}]^2$, which does not contain any very bad vertices by $E_2$.\smallskip

We claim that, on $E_1\cap E_2\cap E_3$, the greedy phase makes fewer than $K_0$ forced steps, the path stays in $\Sigma_u$, both $a_i$ and $b_i$ are positive throughout the greedy phase (so that the forced alternative is always an admissible step) and nonnegative when it ends, and no vertex of $\alpha_{u:\hat u}$, including $u$ itself, is very bad. The last assertion follows from the others and from $E_2$.

Let $d_i=\rcoord(x_i)-r_u(i)$ be the transverse deviation from $\gamma_u$ and $\phi_i$ the number of forced steps made before level $i$. We show by induction that $\abs{d_i}\le1+\phi_i$ during the greedy phase. This holds at $i=I$. Suppose it holds at $i$. The two candidates have deviations $d_i+\lambda$ and $d_i+\lambda-1$ at level $i+1$, where $\lambda=\rcoord(e_1)-\rcoord(\hat u-u)/m\in[\tfrac14,\tfrac34]$ by \eqref{eq:slope}. If $\abs{d_i}\le1$, then $d_i+\lambda\in[-\tfrac34,\tfrac74]$, so the candidate closer to $0$ has absolute deviation at most $\tfrac34$. If $d_i>1$, then $d_i+\lambda-1\in[d_i-\tfrac34,d_i-\tfrac14]$ is the closer one and has deviation in $(0,d_i-\tfrac14]$. The case $d_i<-1$ is symmetric. Hence the preferred step gives $\abs{d_{i+1}}\le\max\{1,\abs{d_i}-\tfrac14\}\le1+\phi_i=1+\phi_{i+1}$, while a forced step changes the deviation by $1$ relative to the preferred step and gives $\abs{d_{i+1}}\le1+\phi_i+1=1+\phi_{i+1}$. This proves the induction step. Consequently, as long as $\phi_i<K_0$, the vertex $x_i$ and both candidates for $x_{i+1}$ have deviation at most $K_0+1$ and lie in $\Sigma_u$. Each forced step is caused by a very bad vertex of $\Sigma_u$ (the preferred candidate), and distinct forced steps are caused by distinct vertices, so by $E_3$ the greedy phase makes fewer than $K_0$ forced steps and the path stays in $\Sigma_u$.

For the coordinates, note that the point of $\gamma_u$ on level $i$ has remaining displacement $\frac{a_i+b_i}{m}(j_u,m-j_u)$ to $\hat u$, both coordinates of which are at least $(a_i+b_i)/4$ by \eqref{eq:middle}, and that $x_i$ differs from that point by $(d_i,-d_i)$ in the $(e_1,e_2)$ coordinates. Hence $a_i,b_i\ge(a_i+b_i)/4-\abs{d_i}\ge(a_i+b_i)/4-K_0$, which is positive while $a_i+b_i>n_{\rm box}=4K_0+4$ and at least $1$ when $a_i+b_i=n_{\rm box}$. This proves the claim.

\medskip
The events $E_1,E_2,E_3$ and the paths $\alpha_{u:\hat u}$, $u\in\cJ$, are $\cF_\vbad$-measurable, and on $E_1\cap E_2\cap E_3$ each $\alpha_{u:\hat u}$ collects $m$ distinct vertices, none of which are in $\vbad$. Apply Lemma~\ref{lem:weights} to this family of at most $N=Ds+1$ paths (on $(E_1\cap E_2\cap E_3)^c$, to the empty family), with $n=m\le2\kappa s$, the integer $K_0$ fixed above, $b=s^{1-\e'}$, and $t=\eta s/16$, for $s$ so large that $\mu([-T,\infty))\ge\tfrac12$. Since $b^{-p'}=s^{-(p'+1)/2}$, the two terms in \eqref{eq:weights} are at most
\[
\bigl(4\kappa C_{p'}s^{(1-p')/2}\bigr)^{K_0}\le(4\kappa C_{p'})^{K_0}s^{-2}
\qquad\text{and}\qquad
\exp\Bigl\{-\frac{\eta^2}{512(4\Cm\kappa+\eta/48)}\,s^{\e'}\Bigr\}=e^{-c_2s^{\e'}},
\]
the last equality defining $c_2$. So, outside an event of probability at most 
\[q_4(s)=(Ds+1)\bigl[(4\kappa C_{p'})^{K_0}s^{-2}+e^{-c_2s^{\e'}}\bigr],\] 
every $u\in\cJ$ has $\pawgt(\alpha_{u:\hat u})\ge-2\Cm m-(K_0-1)T-\eta s/16$. Using $2\Cm m\le2\Cm(\kappa s+1)\le\eta s/8+2\Cm$ and $(K_0-1)T\le\eta s/64$, we have
\[
\pawgt(\alpha_{u:\hat u})\ \ge\ -\frac{\eta s}8-2\Cm-\frac{\eta s}{64}-\frac{\eta s}{16}\ \ge\ -\frac{\eta s}4
\]
once $2\Cm\le\eta s/64$.

We conclude that the lemma holds with
\[
1-p_{\rm att}(s)\ =\ q_1(s)+q_2(s)+q_3(s)+q_4(s)\ \longrightarrow\ 0
\]
for $s$ beyond the finitely many thresholds used in the proof, and with $p_{\rm att}(s)=0$ for smaller $s$. These four functions depend only on $D,\eta,\kappa,\theta_0,\mu$, and not on $\cJ$ and $\cN$. 
The symmetric statement, for paths ending in a level interval, is proved in the same way with the roles of $e_1,e_2$ and $-e_1,-e_2$ exchanged: the greedy phase runs backwards from $v$ towards a net point $\hat v$ on the level $m$ below, chosen so that $v-\hat v=(j_v,m-j_v)$ with $m/4\le j_v\le3m/4$, and the box phase takes place in $\hat v+[0,n_{\rm box}]^2$.
\end{proof}

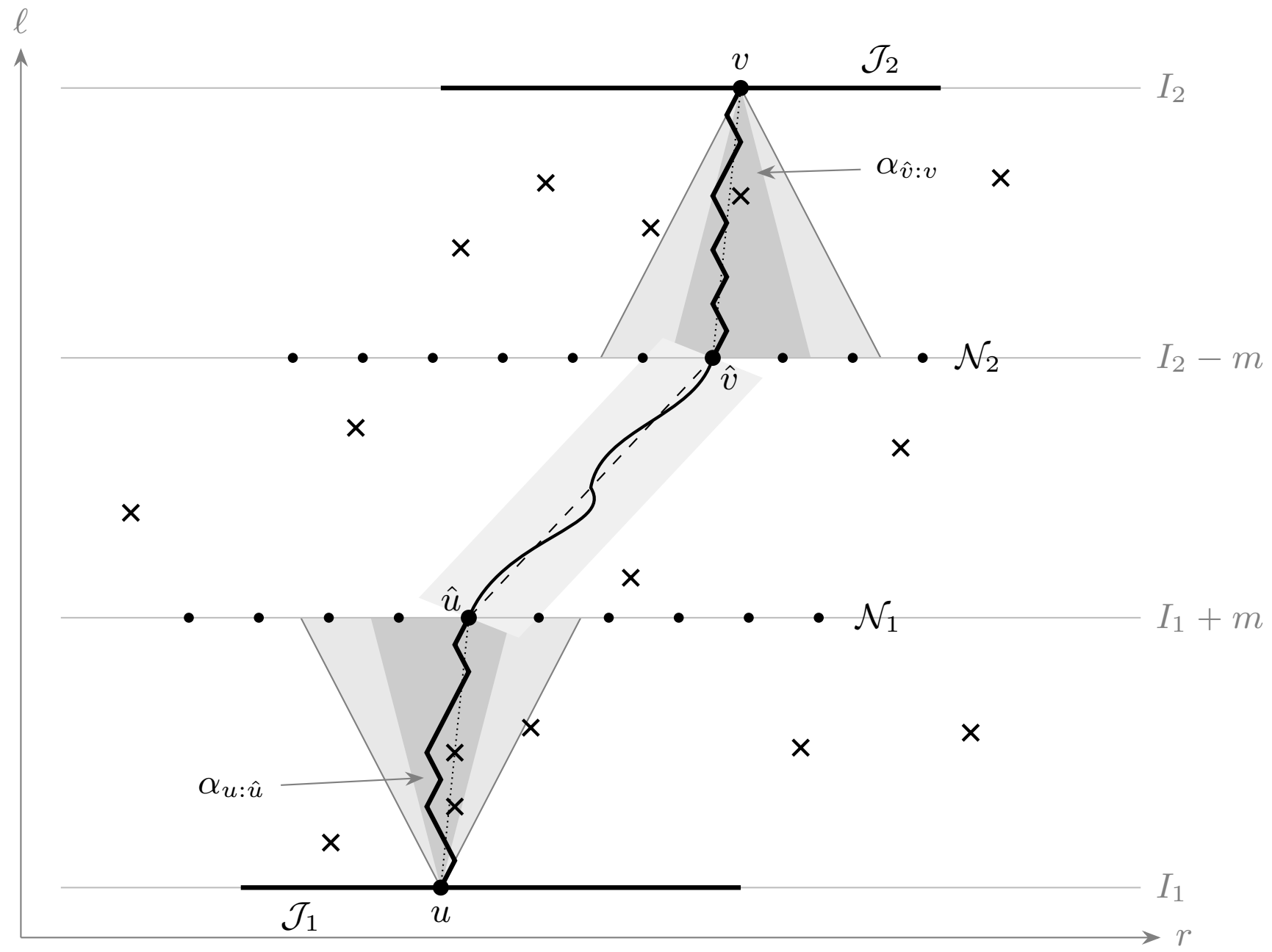
\begin{figure}[htb]
\centering
\begin{tikzpicture}[scale=0.9,>=Stealth]
  %% axes: horizontal = r, vertical = level (not to scale)
  \draw[->,gray] (-2.2,-0.5) -- (-2.2,8.4) node[above] {\footnotesize $\level$};
  \draw[->,gray] (-2.2,-0.5) -- (9.2,-0.5) node[right] {\footnotesize $\rcoord$};
  %% windows of u and v, with their middle halves
  \fill[gray!18] (2,0) -- (0.6,2.7) -- (3.4,2.7) -- cycle;  \fill[gray!40] (2,0) -- (1.3,2.7) -- (2.7,2.7) -- cycle;
  \draw[gray] (0.6,2.7) -- (2,0) -- (3.4,2.7);
  \fill[gray!18] (5,8) -- (3.6,5.3) -- (6.4,5.3) -- cycle;  \fill[gray!40] (5,8) -- (4.3,5.3) -- (5.7,5.3) -- cycle;
  \draw[gray] (3.6,5.3) -- (5,8) -- (6.4,5.3);
  %% connection: tube and path, as in Figure~\ref{fig:intervals}
  \fill[gray!12] (1.78,2.9) -- (4.22,5.5) -- (5.22,5.1) -- (2.78,2.5) -- cycle;
  \draw[dashed] (2.28,2.7) -- (4.72,5.3);
  \draw[thick] (2.28,2.7) .. controls (2.6,3.5) and (3.75,3.6) .. (3.5,4.0) .. controls (3.6,4.6) and (4.55,4.7) .. (4.72,5.3);
  %% level lines
  \foreach \y/\lab in {0/$I_1$, 2.7/$I_1+m$, 5.3/$I_2-m$, 8/$I_2$} {\draw[gray!50] (-1.8,\y) -- (9,\y); \node[right,gray] at (9,\y) {\footnotesize \lab};}
  %% level intervals J_1, J_2
  \draw[very thick] (0,0) -- (5,0);  \node[below] at (0.6,0) {\footnotesize $\cJ_1$};
  \draw[very thick] (2,8) -- (7,8);  \node[above] at (6.4,8) {\footnotesize $\cJ_2$};
  %% nets
  \foreach \x in {-0.52,0.18,...,5.9} \fill (\x,2.7) circle (0.05);
  \node[right] at (6.0,2.7) {\footnotesize $\cN_1$};
  \foreach \x in {0.52,1.22,...,6.9} \fill (\x,5.3) circle (0.05);
  \node[right] at (7.0,5.3) {\footnotesize $\cN_2$};
  %% the segments [u,u-hat], [v-hat,v] and the attachments, which step around the very bad vertices
  \draw[densely dotted] (2,0) -- (2.28,2.7);  \draw[densely dotted] (4.72,5.3) -- (5,8);
  \draw[very thick] (2,0) -- (2.14,0.27) -- (2,0.54) -- (1.86,0.81) -- (2,1.08) -- (1.86,1.35) -- (2,1.62) -- (2.14,1.89) -- (2.28,2.16) -- (2.14,2.43) -- (2.28,2.7);
  \draw[very thick] (5,8) -- (4.86,7.73) -- (5,7.46) -- (4.86,7.19) -- (4.72,6.92) -- (4.86,6.65) -- (4.72,6.38) -- (4.86,6.11) -- (4.72,5.84) -- (4.86,5.57) -- (4.72,5.3);
  \draw[thick] (2.14,0.81) ++(-0.08,-0.08) -- ++(0.16,0.16) (2.14,0.81) ++(-0.08,0.08) -- ++(0.16,-0.16);
  \draw[thick] (2.14,1.35) ++(-0.08,-0.08) -- ++(0.16,0.16) (2.14,1.35) ++(-0.08,0.08) -- ++(0.16,-0.16);
  \draw[thick] (5,6.92) ++(-0.08,-0.08) -- ++(0.16,0.16) (5,6.92) ++(-0.08,0.08) -- ++(0.16,-0.16);
  \draw[thick] (0.9,0.45) ++(-0.08,-0.08) -- ++(0.16,0.16) (0.9,0.45) ++(-0.08,0.08) -- ++(0.16,-0.16);
  \draw[thick] (3.9,3.1) ++(-0.08,-0.08) -- ++(0.16,0.16) (3.9,3.1) ++(-0.08,0.08) -- ++(0.16,-0.16);
  \draw[thick] (1.15,4.6) ++(-0.08,-0.08) -- ++(0.16,0.16) (1.15,4.6) ++(-0.08,0.08) -- ++(0.16,-0.16);
  \draw[thick] (6.6,4.4) ++(-0.08,-0.08) -- ++(0.16,0.16) (6.6,4.4) ++(-0.08,0.08) -- ++(0.16,-0.16);
  \draw[thick] (3.05,7.05) ++(-0.08,-0.08) -- ++(0.16,0.16) (3.05,7.05) ++(-0.08,0.08) -- ++(0.16,-0.16);
  \draw[thick] (7.3,1.55) ++(-0.08,-0.08) -- ++(0.16,0.16) (7.3,1.55) ++(-0.08,0.08) -- ++(0.16,-0.16);
  \draw[thick] (-1.1,3.75) ++(-0.08,-0.08) -- ++(0.16,0.16) (-1.1,3.75) ++(-0.08,0.08) -- ++(0.16,-0.16);
  \draw[thick] (5.6,1.4) ++(-0.08,-0.08) -- ++(0.16,0.16) (5.6,1.4) ++(-0.08,0.08) -- ++(0.16,-0.16);
  \draw[thick] (2.2,6.4) ++(-0.08,-0.08) -- ++(0.16,0.16) (2.2,6.4) ++(-0.08,0.08) -- ++(0.16,-0.16);
  \draw[thick] (7.6,7.1) ++(-0.08,-0.08) -- ++(0.16,0.16) (7.6,7.1) ++(-0.08,0.08) -- ++(0.16,-0.16);
  \draw[thick] (2.9,1.6) ++(-0.08,-0.08) -- ++(0.16,0.16) (2.9,1.6) ++(-0.08,0.08) -- ++(0.16,-0.16);
  \draw[thick] (4.1,6.6) ++(-0.08,-0.08) -- ++(0.16,0.16) (4.1,6.6) ++(-0.08,0.08) -- ++(0.16,-0.16);
  \fill (2,0) circle (0.08) node[below=1pt] {\footnotesize $u$};
  \fill (5,8) circle (0.08) node[above=1pt] {\footnotesize $v$};
  \fill (2.28,2.7) circle (0.08) node[above left=-2pt] {\footnotesize $\hat u$};
  \fill (4.72,5.3) circle (0.08) node[below right=-2pt] {\footnotesize $\hat v$};
  \node[left] (alab) at (0.4,1.0) {\footnotesize $\alpha_{u:\hat u}$};
  \draw[->,gray,shorten >=1pt] (alab) -- (1.87,1.1);
  \node[right] (blab) at (6.2,7.2) {\footnotesize $\alpha_{\hat v:v}$};
  \draw[->,gray,shorten >=1pt] (blab) -- (5.12,7.15);
\end{tikzpicture}
\caption{\small How Figure~\ref{fig:intervals} changes under \eqref{ass:general}, in the proof of Lemma~\ref{lem:intervals-general}. The very bad vertices (crosses) are observed first. The net point $\hat u$ is chosen in the middle half of the window of points reachable from $u$ in $m$ steps (darker shading), and $\hat v$ symmetrically. The attachments $\alpha_{u:\hat u}$ and $\alpha_{\hat v:v}$ (thick) are no longer deterministic: they follow the segments $[u,\hat u]$ and $[\hat v,v]$ (dotted) and step around the very bad vertices, as in Figure~\ref{fig:attach-general}, and their weights are bounded by Lemma~\ref{lem:weights} instead of Hoeffding's inequality. The connections between the nets are as before.}
\label{fig:intervals-general}
\end{figure}

\begin{proof}[Proof of Lemma~\ref{lem:intervals-general}]
Follow the proof of Lemma~\ref{lem:intervals} with
\[
\kappa=\min\Bigl\{\frac{\theta}{16},\ \frac{\delta}{4},\ \frac{\eta}{16(\Lip_{\theta/2}+\Cm+1)}\Bigr\}\quad\text{and}\quad m=\lceil\kappa s\rceil ,
\]
which differs from the $\kappa$ used there only in that $\abs{\E\w_0}$ is replaced by $\Cm$, and which satisfies the hypothesis $\kappa\le\eta/(16\Cm+16)$ of Lemma~\ref{lem:attach-general}. The five requirements on $s_1$ are kept, with the fourth replaced by $s_1\ge16\Lip_{\theta/2}/\eta$. The nets $\cN_1$ and $\cN_2$ are as there. The point $\hat u\in\cN_1$ is now chosen in the middle half of the window of points reachable from $u$ in $m$ steps, as in \eqref{eq:middle}, and $\hat v\in\cN_2$ symmetrically, see Figure~\ref{fig:intervals-general}. The claim on the nets proved there permits this, since the middle half is an interval of length $m/2\ge\theta_0m/2$ contained in $\cR_1$, as shown in the proof of Lemma~\ref{lem:attach-general}, just before \eqref{eq:middle}.  The bounds \eqref{aux397} and \eqref{aux403} hold as before, so the verification that $\tilde z=\hat v-\hat u$ is $\theta/2$-interior with level in $[s/2,Ds]$ is unchanged, and so is the application of Lemma~\ref{lem:travel} to the boundedly many pairs $(\hat u,\hat v)$. The deterministic attachments and the Hoeffding bound \eqref{attachments} are replaced by Lemma~\ref{lem:attach-general}, applied to $(\cJ_1,\cN_1)$ and, in its symmetric form, to $(\cJ_2,\cN_2)$, with the $D$, $\eta$, and $\kappa$ above. With probability at least $2p_{\rm att}(s)-1\to1$, every $u\in\cJ_1$ has an attachment $\alpha_{u:\hat u}$ and every $v\in\cJ_2$ an attachment $\alpha_{\hat v:v}$ of weight at least $-\eta s/4$. On the intersection of this event with the event of Lemma~\ref{lem:travel}, let $x_{u:v}$ be the concatenation of $\alpha_{u:\hat u}$, the connection from $\hat u$ to $\hat v$, and $\alpha_{\hat v:v}$. The final estimate becomes
\[
\pawgt(x_{u:v})\ \ge\ -\frac{\eta s}4+\Bigl(g(\tilde z)-\frac{\eta s}4\Bigr)-\frac{\eta s}4\ \ge\ g(z)-2\tsp\Lip_{\theta/2}m-\frac{3\eta s}4\ \ge\ g(z)-\eta s ,
\]
because $2\tsp\Lip_{\theta/2}m\le2\tsp\Lip_{\theta/2}(\kappa s+1)\le\eta s/8+2\tsp\Lip_{\theta/2}\le\eta s/4$ by the new fourth requirement. The location statement is unchanged.
\end{proof}

\begin{remark}
The above argument needs the exponent $p'>1$ because of the union bound over the $O(s)$ starting points $u$: the failure probability for one $u$ is dominated by the event that the tube $\Sigma_u$ contains $K_0$ very bad vertices, which has probability of order $(mp_T)^{K_0}$, and this is $o(1/s)$ for a fixed $K_0$ only if $mp_T$ decays polynomially in $s$. Under $\E\w_0^-<\infty$ alone one has only $mp_T=o(1)$, possibly very slowly. Letting $K_0$ grow with $s$ (which forces $T=\eta s/(64K_0)$ to decrease) one can presumably reach a condition of the type $\E[\w_0^-\log(1+\w_0^-)]<\infty$, but we have not pursued this.
\end{remark}

\subsection{Staircases}\label{sub:staircases}
In this section we assume \ref{A2} and \ref{A3} and use the notation of Section~\ref{sec:endpoints}: $k\in\cK$ with \eqref{eq:kreq}, $s_k=2^k$, $\clen_k=2^{m_k}\le s_k^3$, $t_k=\e_0s_k/k$, the scale indices $j$ with their meshes $l_{j,k}$ and widths $w_{j,k}=s_{k-j}$, the sets $\cG_{j,k}$ of grid intervals and $\cH_{j,k}$ of bins, and $A,R>0$.

The construction of Section~\ref{sec:endpoints} uses \eqref{ass:bounded} only through the Hoeffding bound \eqref{eq:badstair} on the weights of the deterministic staircases. Under \ref{A2}--\ref{A3} the staircases are replaced by paths that avoid the very bad vertices, in the sense of Section~\ref{sub:verybad}, and the lower bound on their weights comes from Lemma~\ref{lem:weights}. This requires three things. Paths that avoid the very bad vertices must exist with the geometry that the construction needs: through a given vertex, along a given line, and from a given vertex to a given target. This is the content of Lemma~\ref{lem:stair-geometry}. A path that is to be joined to such a path must start and end at suitable vertices, which is arranged by \emph{trimming} a bounded number of vertices from its ends. And the weight of the trimmed vertices must be bounded above, which is where the upper tail \ref{A3} enters, through Lemma~\ref{lem:heavy}. The section ends with the list \ref{ch:gamma}--\ref{ch:endpoints} of the changes to Section~\ref{sec:endpoints}, announced by the remarks there.

All paths of Section~\ref{sec:endpoints} lie in the region
\[
\cQ_k=\{z\in\Z^2:\ \abs{\level(z)}\le(A+5)\clen_k,\ \abs{\rcoord(z)}\le(R+2\cnst_1+10)s_k\},\qquad\abs{\cQ_k}\le\cnst_{\!A,R}\,s_k^4 ,
\]
where $\cnst_1$ is the constant of Lemma~\ref{lem:grid} in its form \ref{ch:cost} below, with a margin that is explained before item~\ref{ch:gamma}. Throughout, $\varrho \ge1$ is a fixed exponent, the families of paths below have at most $s_k^\varrho$ members, and all constants may depend on $A,R,\e_0,\varrho,\theta_0,\mu$. We will see in \ref{ch:family} that $\varrho =11$ suffices.

\medskip
\emph{Parameters.} Fix integers $K_1,K_2,K_3\ge1$ with
\begin{equation}\label{eq:K1K2}
K_3p'>4,\qquad(K_1+1)p'>4,\qquad\text{and}\qquad K_2(p'-1)>2\varrho +2 ,
\end{equation}
and set
\begin{equation}\label{eq:Tk}
T=\frac{t_k}{4K_2}\qquad\text{and}\qquad p_T=\P(\w_0<-T)\le C_{p'}\Bigl(\frac{4K_2k}{\e_0}\Bigr)^{p'}s_k^{-p'} .
\end{equation}
As in Section~\ref{sub:verybad}, a vertex $z$ is \emph{very bad} if $\w_z<-T$, now with the threshold $T$ of \eqref{eq:Tk}, and $\vbad$ is the set of very bad vertices. Let  
\begin{align*}
&\theta_1=\theta_0/5,\quad N_{\rm tr}=K_1(2K_3^2+1)^2,\quad \Delta_0=N_{\rm tr}+1,\quad n_{\rm box}=\bigl\lceil4(\Delta_0+2)/\theta_1\bigr\rceil,\\
&r_0=K_1(2n_{\rm box}+1)+N_{\rm tr},\quad D_1=2\Delta_0/\theta_1+4\Delta_0+2r_0+2n_{\rm box}+8,\quad\text{and}\quad d_0=D_1+2K_3^2 .
\end{align*}
These are constants that depend on $p'$, $\varrho$, and $\theta_0$ only. Their roles are as follows. $K_1$ bounds the size of a cluster of very bad vertices and $K_3^2$ the size of a finite set of vertices reachable without collecting a very bad vertex. $N_{\rm tr}$ bounds the size of a component of trapped vertices, $\Delta_0$ the transverse deviation of a greedy path from the line it follows, $n_{\rm box}$ the side of the box in which an arrival path is completed, and $r_0$ the number of vertices trimmed from an end of a path. $D_1$ and $d_0$ are the distances that define the components and the clusters. $K_2$ is the parameter $K_0$ of Lemma~\ref{lem:weights}. We take $k$ so large that
\begin{equation}\label{eq:wlarge}
w_{j,k}=2^{k-j}\ge2^{(k-3)/4}\ge4(\Delta_0+r_0+n_{\rm box})\qquad\text{for all scale indices }j,
\end{equation}
which is possible by \eqref{t.aux34}. We next introduce some terminology that will facilitate the exposition.
\medskip

\emph{Clusters.} Let $c_0=\Delta_0+2n_{\rm box}+1$ and, for $c>0$, $\cQ_k^{(c)}=\{z\in\Z^2:\ \abs{z-y}_1\le c\text{ for some }y\in\cQ_k\}$. Write $\cQ_k^+=\cQ_k^{(c_0)}$ and $\cQ_k^{++}=\cQ_k^{(c_0+K_3^2)}$. Both have at most $\cnst_{\!A,R}\,s_k^4$ vertices. Two very bad vertices are said to be \emph{neighbors} if their $\abs\cdot_1$-distance is at most $d_0$, and a \emph{cluster} is a connected component of $\vbad \cap\cQ_k^{++}$ for this neighbor relation.\medskip

\emph{Trapped vertices.} For $z\in\Z^2$ let $\Afw(z)$ be the set of vertices reachable from $z$ by admissible paths that collect no very bad vertex ($z$ itself included), and $\Abw(z)$ the set of vertices from which $z$ can be reached by such paths. A vertex $z$ is \emph{trapped} if $\Afw(z)$ is finite and \emph{backward trapped} if $\Abw(z)$ is finite. Write $\Vfw$ and $\Vbw$ for the sets of trapped and backward trapped vertices. By K\H{o}nig's lemma \cite{Kon-27} (see also \cite[Lemma~8.1.2]{Die-18}), 
%https://daiwz.net/course/disc_math/2023/Diestel_Graph_Theory.pdf 
applied to the finitely branching tree of finite admissible paths from $z$ that collect no very bad vertex, $z$ is trapped if and only if every infinite admissible path starting at $z$ collects a very bad vertex, and likewise backward. The two notions are not mirror images of each other, because a path starting at $z$ collects $z$ while a path ending at $z$ does not: $\Vfw$ is the smallest set $S$ such that $z\in S$ whenever $z\in \vbad$ or $z+e_1,z+e_2\in S$, while $\Vbw$ is the smallest set $S$ such that $z\in S$ whenever $z-e_1,z-e_2\in \vbad \cup S$. In particular $\vbad \subset\Vfw$, but a very bad vertex need not be backward trapped, and the set that the backward constructions avoid is $\vbad \cup\Vbw$.\medskip

\emph{Components.} Two vertices \emph{interfere} if their $\abs\cdot_1$-distance is at most $D_1$. A \emph{trapped component} is a connected component of $\Vfw \cap\cQ_k^+$ for the interference relation, and a \emph{backward trapped component} is a connected component of $(\vbad \cup\Vbw)\cap\cQ_k^+$ for the same relation.\medskip

\emph{The event $E_{\rm geo}$} is the event that every cluster has at most $K_1$ vertices and that, for every $z\in\cQ_k^+$, each of the sets $\Afw(z)$ and $\Abw(z)$ is either infinite or has fewer than $K_3^2$ elements. The sets $\vbad,\Vfw ,\Vbw$ and the event $E_{\rm geo}$ are measurable with respect to $\cF_\vbad$, the $\sigma$-algebra generated by the indicators $\one\{\w_z<-T\}$.\medskip

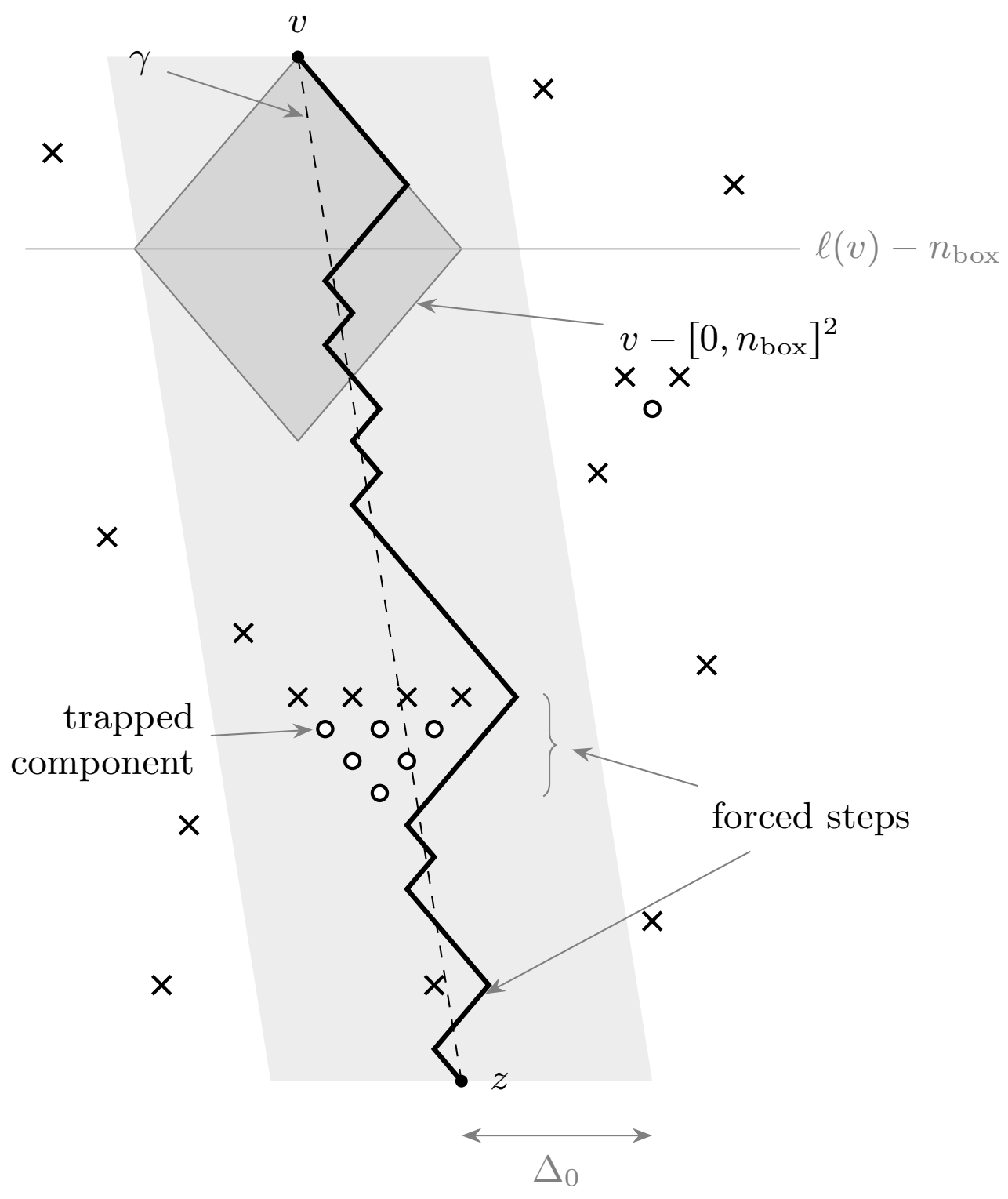
\begin{figure}[htb]
\centering
\begin{tikzpicture}[x=0.46cm,y=0.27cm,>=Stealth]
  %% coordinates (r, level) with Dir_1=1/2, so that e_1 and e_2 steps go up-right and up-left; not to scale
  %% the band of transverse half-width Delta_0 around the segment [z,v], and the box at v
  \fill[gray!14] (-3.5,0) -- (3.5,0) -- (0.5,32) -- (-6.5,32) -- cycle;
  \fill[gray!32] (-3,32) -- (-6,26) -- (-3,20) -- (0,26) -- cycle;
  \draw[gray] (-3,32) -- (-6,26) -- (-3,20) -- (0,26) -- cycle;
  \draw[gray!60] (-8,26) -- (6.2,26); \node[right,gray] at (6.2,26) {\footnotesize $\level(v)-n_{\rm box}$};
  \draw[dashed] (0,0) -- (-3,32);
  %% the arrival path: greedy phase up to level(v)-n_box, then e_1 steps and e_2 steps inside the box
  \draw[very thick] (0,0) -- (-0.5,1) -- (0,2) -- (0.5,3) -- (0,4) -- (-0.5,5) -- (-1,6) -- (-0.5,7) -- (-1,8) -- (-0.5,9) -- (0,10) -- (0.5,11) -- (1,12) -- (0.5,13) -- (0,14) -- (-0.5,15) -- (-1,16) -- (-1.5,17) -- (-2,18) -- (-1.5,19) -- (-2,20) -- (-1.5,21) -- (-2,22) -- (-2.5,23) -- (-2,24) -- (-2.5,25) -- (-2,26) -- (-1.5,27) -- (-1,28) -- (-1.5,29) -- (-2,30) -- (-2.5,31) -- (-3,32);
  %% trapped vertices that are not very bad
  \foreach \p in {(-1.5,9),(-2,10),(-1,10),(-2.5,11),(-1.5,11),(-0.5,11),(3.5,21)} \draw[thick,fill=white] \p circle[radius=1.8pt];
  %% very bad vertices
  \foreach \p in {(-5.5,3),(-0.5,3),(3.5,5),(-5,8),(-3,12),(-2,12),(-1,12),(0,12),(4.5,13),(-4,14),(-6.5,17),(2.5,19),(3,22),(4,22),(5,28),(-7.5,29),(1.5,31)} \draw[thick] \p ++(-0.17,-0.29) -- ++(0.34,0.58) \p ++(-0.17,0.29) -- ++(0.34,-0.58);
  \fill (0,0) circle[radius=1.5pt] node[right=3pt] {\footnotesize $z$};
  \fill (-3,32) circle[radius=1.5pt] node[above=2pt] {\footnotesize $v$};
  %% labels
  \draw[<->,gray] (0,-1.7) -- (3.5,-1.7) node[midway,below=1pt] {\footnotesize $\Delta_0$};
  \node[left] (glab) at (-5.4,31.8) {\footnotesize $\gamma$};
  \draw[->,gray,shorten >=1pt] (glab) -- (-2.82,30.1);
  \node[left,align=right] (tlab) at (-4.6,10.6) {\footnotesize trapped\\[-2pt]\footnotesize component};
  \draw[->,gray,shorten >=3pt] (tlab) -- (-2.5,11);
  \node[right] (flab) at (4.3,8.2) {\footnotesize forced steps};
  \draw[gray,decorate,decoration={brace,amplitude=3pt,mirror}] (1.5,8.9) -- (1.5,12.1);
  \draw[->,gray,shorten >=1pt] (flab) -- (1.95,10.4);
  \draw[->,gray,shorten >=2pt] (flab) -- (0.3,2.6);
  \node[right] (blab) at (2.6,23.2) {\footnotesize $v-[0,n_{\rm box}]^2$};
  \draw[->,gray,shorten >=1pt] (blab) -- (-0.9,24.3);
\end{tikzpicture}
\caption{\small Trapped vertices, a greedy path, and an arrival path, drawn in the coordinates $(\rcoord,\level)$ and not to scale. Crosses are very bad vertices and circles are the trapped vertices that are not very bad: below a row of adjacent very bad vertices on one level there is a triangle of vertices from which every admissible path runs into the row. The thick path is the arrival path from $z$ to $v$. Up to level $\level(v)-n_{\rm box}$ it is the greedy path from $z$ along the line $\gamma$ through $z$ and $v$ (dashed), which moves to the candidate closer to $\gamma$ unless that candidate is trapped, in which case the step is forced. On $E_{\rm geo}$ a trapped component forces at most $N_{\rm tr}$ steps, and the path then returns to $\gamma$ before it meets another component, so it stays within transverse distance $\Delta_0$ of $\gamma$ (light gray). The last $n_{\rm box}$ steps, the $e_1$ steps followed by the $e_2$ steps, lie in the box $v-[0,n_{\rm box}]^2$ (darker gray), which contains no very bad vertex.}
\label{fig:stair-geometry}
\end{figure}

\emph{Greedy paths.} The \emph{$e_1$-fraction} of a nonzero vector $x=(x_1,x_2)\in\R^2_{\ge0}$ is $x_1/(x_1+x_2)$. Let $\gamma$ be a straight line or a segment whose direction is a nonzero vector of $\R^2_{\ge0}$, so that $\gamma$ has at most one point on each level. The \emph{transverse distance} from a point $x$ to $\gamma$ is $\abs{\rcoord(x)-\rcoord(y)}$, where $y$ is the point of $\gamma$ on the level of $x$. Now let $\gamma$ be a line and $z\in\Z^2\setminus\Vfw$ a point of $\gamma$. The \emph{greedy path from $z$ along $\gamma$} is the infinite admissible path $x_{\level(z):\infty}$, indexed by level, that is defined inductively as follows, starting from $x_{\level(z)}=z$. Of the two candidates $x_i+e_1$ and $x_i+e_2$ for $x_{i+1}$, the \emph{preferred candidate} is the one with the smaller transverse distance to $\gamma$ (either one in case of a tie). If the preferred candidate is not trapped, it is $x_{i+1}$. Otherwise the step is \emph{forced}, and $x_{i+1}$ is the other candidate. The greedy path is well defined and none of its vertices is trapped, because the two candidates of a vertex $x_i\notin\Vfw$ are not both trapped, by the characterization of $\Vfw$ above. In particular, it does not collect any very bad vertices. For a point $z\in\Z^2\setminus(\vbad\cup\Vbw)$ of $\gamma$, the \emph{backward greedy path from $z$ along $\gamma$} is defined in the same way, with the candidates $x_i-e_1$ and $x_i-e_2$ for $x_{i-1}$ and with the set $\vbad\cup\Vbw$ in place of $\Vfw$. Thus a backward step is forced if the preferred candidate lies in $\vbad\cup\Vbw$. By the characterization of $\Vbw$, the two candidates of a vertex $x_i\notin\vbad\cup\Vbw$ are not both in $\vbad\cup\Vbw$, so no vertex of the backward greedy path lies in $\vbad\cup\Vbw$, and it does not collect any very bad vertices either. See Figure~\ref{fig:stair-geometry} for an illustration.\medskip

\emph{Arrival paths.} Let $z\in\Z^2\setminus\Vfw$ and $v\in\Z^2$ with $v\ge z$ and $\level(v-z)\ge n_{\rm box}$. Let $x$ be the vertex on level $\level(v)-n_{\rm box}$ of the greedy path from $z$ along the line through $z$ and $v$. If $x\le v$, the \emph{arrival path from $z$ to $v$} is the admissible path that follows this greedy path from $z$ to $x$ and then goes from $x$ to $v$ by $e_1$ steps followed by $e_2$ steps. Its last $n_{\rm box}$ steps stay in the box $v-[0,n_{\rm box}]^2$. If $x\not\le v$, the arrival path from $z$ to $v$ is not defined.\medskip

The next lemma collects the properties of these objects. The greedy paths and the arrival paths are the paths that replace the staircases. Item~\ref{geo:greedy} shows that a greedy path stays close to the line that it follows, item~\ref{geo:arrive} that the arrival path to a suitable target exists, item~\ref{geo:merge} is used to join two such paths, and item~\ref{geo:trim} shows that the vertices at which these constructions can start and end are never far away.

\begin{lemma}\label{lem:stair-geometry}
$\P(E_{\rm geo}^c)\to0$ as $k\to\infty$, and on $E_{\rm geo}$ the following hold.
\begin{enumerate}[label={\rm(\roman*)}]
\item\label{geo:trap} Every trapped component and every backward trapped component has at most $N_{\rm tr}$ vertices.
\item\label{geo:greedy} Let $z\in\cQ_k\setminus\Vfw$ and let $\gamma$ be a straight line through $z$ whose direction has $e_1$-fraction in $[\theta_1,1-\theta_1]$. Then every vertex of the greedy path from $z$ along $\gamma$ is within transverse distance $\Delta_0$ of $\gamma$, for as long as the path stays in $\cQ_k$. The same holds for the backward greedy path from $z$ along $\gamma$ when $z\in\cQ_k\setminus(\vbad \cup\Vbw)$.
\item\label{geo:arrive} Let $z\in\cQ_k\setminus\Vfw$ and $v\ge z$ with $\level(v-z)=n\ge n_{\rm box}$ and $e_1$-fraction of $v-z$ in $[\theta_1,1-\theta_1]$, and suppose that the box $v-[0,n_{\rm box}]^2$ contains no very bad vertex and that every lattice point within transverse distance $\Delta_0$ of the segment $[z,v]$ lies in $\cQ_k$. Then the arrival path from $z$ to $v$ is defined, collects no very bad vertex, and stays within transverse distance $\Delta_0+n_{\rm box}$ of $[z,v]$.
\item\label{geo:merge} Let $x_{m:n}$ and $y_{m:n}$ be admissible paths with $\rcoord(x_m)>\rcoord(y_m)$ and $\rcoord(x_n)<\rcoord(y_n)$. Then $x_i=y_i$ for some $m<i<n$.
\item\label{geo:trim} Among any $r_0+1$ consecutive vertices of an admissible path in $\cQ_k$ there is one whose box $v-[0,n_{\rm box}]^2$ contains no very bad vertex, there is one that is not trapped, and there is one that lies outside $\vbad \cup\Vbw$.
\end{enumerate}
\end{lemma}

\begin{proof}
The proof has two parts. We first show that $\P(E_{\rm geo}^c)\to0$, by bounding the probability that some cluster has more than $K_1$ vertices and then the probability that some set $\Afw(z)$ or $\Abw(z)$ is finite with at least $K_3^2$ elements. Then we prove items~\ref{geo:trap}--\ref{geo:trim} on the event $E_{\rm geo}$.\smallskip

A cluster with more than $K_1$ vertices contains a connected set of exactly $K_1+1$ very bad vertices. The number of connected sets of $K_1+1$ vertices containing a given vertex is at most $K_1!\,(2d_0+1)^{2K_1}$, since each is obtained by $K_1$ additions, each of which chooses an already chosen vertex and one of its at most $(2d_0+1)^2$ neighbors. Hence, by independence and \eqref{eq:Tk}, the probability that some cluster has more than $K_1$ vertices is at most
\[
\abs{\cQ_k^{++}}\,K_1!\,(2d_0+1)^{2K_1}p_T^{K_1+1}\le\cnst\,s_k^4\bigl(k^{p'}s_k^{-p'}\bigr)^{K_1+1}\longrightarrow0
\]
by the second condition in \eqref{eq:K1K2}.

For the sets $\Afw(z)$ and $\Abw(z)$ we use Lemma~\ref{lem:contour}. Suppose that $\Afw(z)$ is finite. It contains $z$ and it is connected, since every vertex of an admissible path from $z$ that collects no very bad vertex lies in $\Afw(z)$. Let $\Gamma$ be the circuit of Lemma~\ref{lem:contour} for $S=\Afw(z)$, and let $m$ be its length. A vertex $w\in\Afw(z)$ with $w+e_1\notin\Afw(z)$ or $w+e_2\notin\Afw(z)$ is very bad, since otherwise both $w+e_1$ and $w+e_2$ could be reached from $z$ through $w$. So, by Lemma~\ref{lem:contour}\ref{con:forward}, $\Gamma$ determines at least $m/4$ vertices, all of which are very bad. By Lemma~\ref{lem:contour}\ref{con:size}, $\abs{\Afw(z)}\ge K_3^2$ forces $m\ge4K_3$, and by Lemma~\ref{lem:contour}\ref{con:count} there are at most $m\,3^{m-1}$ circuits of length $m$ to consider. Hence, for every $z$,
\[
\P\bigl\{K_3^2\le\abs{\Afw(z)}<\infty\bigr\}\le\sum_{m\ge4K_3}m\,3^{m-1}p_T^{m/4}\le\cnst\,p_T^{K_3}\le\cnst\bigl(k^{p'}s_k^{-p'}\bigr)^{K_3}
\]
once $k$ is large, and the union over $z\in\cQ_k^+$ tends to $0$ by the first condition in \eqref{eq:K1K2}. For $\Abw(z)$ the argument is the same, with item~\ref{con:backward} of Lemma~\ref{lem:contour} in place of item~\ref{con:forward}. The set $\Abw(z)$ consists of $z$ and of the vertices $y\notin \vbad$ with $y+e_1\in \Abw(z)$ or $y+e_2\in \Abw(z)$. It is connected, and a vertex $w'\notin\Abw(z)$ with $w'+e_1\in\Abw(z)$ or $w'+e_2\in\Abw(z)$ is very bad. So $\Gamma$ again determines at least $m/4$ very bad vertices, now outside $\Abw(z)$, and the rest of the argument is unchanged. 

We now prove items~\ref{geo:trap}--\ref{geo:trim}. For the rest of the proof we work on the event $E_{\rm geo}$.\smallskip

\emph{\ref{geo:trap}.} Let $z\in\cQ_k^+$ be trapped. Take any infinite admissible path from $z$. Since $z$ is trapped, it collects a very bad vertex. Let $z'$ be the first one. Then $z'$ and all vertices of the path before it lie in $\Afw(z)$, which on $E_{\rm geo}$ has fewer than $K_3^2$ elements, so $\abs{z'-z}_1<K_3^2$. Hence every trapped vertex of $\cQ_k^+$ is within distance $K_3^2$ of a very bad vertex, which lies in $\cQ_k^{++}$, and a very bad vertex has at most $(2K_3^2+1)^2$ trapped vertices within distance $K_3^2$. Now let $y_0,y_1,\dots,y_N$ be a chain of interfering trapped vertices of $\cQ_k^+$, consecutive ones at distance at most $D_1$, and let $y'_j\in \vbad \cap\cQ_k^{++}$ be within distance $K_3^2$ of $y_j$. Consecutive $y'_j$ are at distance at most $D_1+2K_3^2=d_0$, so all $y'_j$ lie in one cluster, which has at most $K_1$ vertices. Therefore a trapped component has at most $K_1(2K_3^2+1)^2=N_{\rm tr}$ vertices. For the backward statement let $z\in(\vbad \cup\Vbw)\cap\cQ_k^+$. If $z\in \vbad$ take $z'=z$. Otherwise $z$ is backward trapped, so every infinite admissible path ending at $z$ collects a very bad vertex. Let $z'$ be the last one on such a path. It is collected, so $z'\notin \Abw(z)$, but the $\abs{z-z'}_1$ vertices after it, including $z$, lie in $\Abw(z)$, so $\abs{z-z'}_1\le\abs{\Abw(z)}<K_3^2$. In both cases $z$ is within distance $K_3^2$ of the very bad vertex $z'$, and the chain argument above applies verbatim to the backward trapped components.\smallskip

\emph{\ref{geo:greedy}.} Let $x_i$ be the vertex on level $i$ of the greedy path from $z$ along $\gamma$, and let $d_i$ be its transverse deviation from $\gamma$, that is, $\rcoord(x_i)$ minus the transverse coordinate of the point of $\gamma$ on level $i$. Recall that a step is forced if the preferred candidate is trapped.  The two candidates have deviations $d_i+\lambda$ and $d_i+\lambda-1$ with $\lambda\in[\theta_1,1-\theta_1]$, by the slope condition. If $\abs{d_i}\le1$, the preferred candidate has deviation of absolute value at most $1-\theta_1$. If $d_i>1$, it has deviation in $(0,d_i-\theta_1]$, and symmetrically if $d_i<-1$. So an unforced step gives $\abs{d_{i+1}}\le\max\{1,\abs{d_i}-\theta_1\}$, and a forced step gives $\abs{d_{i+1}}\le\abs{d_i}+1$.

Each forced step is caused by a distinct trapped vertex (the preferred candidate, on a distinct level), at transverse distance at most $\abs{d_i}+1$ from $\gamma$. This vertex is adjacent to a vertex of the path in $\cQ_k$, hence lies in $\cQ_k^+$, so it belongs to a trapped component. Two such vertices that cause forced steps at levels $i<i'$ with $i'-i\le D_1-4\Delta_0-4$, while $\abs{d_i},\abs{d_{i'}}\le \Delta_0$, are at $\abs\cdot_1$-distance at most $(i'-i)+4\Delta_0+4\le D_1$ (a transverse deviation $d$ costs $\abs\cdot_1$-distance $2\abs d$) and interfere, so they belong to the same trapped component.

We now show by induction over the forced steps that $\abs{d_i}\le \Delta_0$ throughout. Consider a maximal stretch of forced steps caused by one trapped component, and suppose $\abs{d_i}\le1$ at its first forced step (this holds at the start, where $d=0$). By \ref{geo:trap} the stretch contains at most $N_{\rm tr}$ forced steps, and between them the unforced steps do not increase $\abs{d_i}$ beyond $\max\{1,\abs{d_i}\}$, so within the stretch $\abs{d_i}\le1+N_{\rm tr}=\Delta_0$. After the last forced step of the stretch, $\abs{d_i}$ decreases by $\theta_1$ per step until it is at most $1$, which takes at most $\Delta_0/\theta_1$ steps, and no forced step occurs during these steps: a forced step within $\Delta_0/\theta_1+1\le D_1-4\Delta_0-4$ levels of the last one would be caused by a vertex of the same component, contradicting the maximality of the stretch. So $\abs{d_i}\le1$ at the first forced step of the next stretch, and the induction closes. The proof for the backward greedy path is the same, with $\vbad\cup\Vbw$ in place of $\Vfw$ and the backward trapped components in place of the trapped components.
\smallskip

\emph{\ref{geo:arrive}.} Let $x_i$ be the vertex on level $i$ of the greedy path from $z$ along the line through $z$ and $v$, and write $(a_i,b_i)=v-x_i$ for the remaining displacement. The point of $[z,v]$ on level $i$ has remaining displacement $\frac{a_i+b_i}n(v-z)$, both coordinates of which are at least $\theta_1(a_i+b_i)$, and $x_i$ differs from it by $(d_i,-d_i)$ in the $(e_1,e_2)$ coordinates. By \ref{geo:greedy}, $\abs{d_i}\le \Delta_0$, so $a_i,b_i\ge\theta_1(a_i+b_i)-\Delta_0$, which is at least $1$ when $a_i+b_i=n_{\rm box}$, by the choice of $n_{\rm box}$. Hence the vertex $x_i$ on level $\level(v)-n_{\rm box}$ satisfies $x_i\le v$, and the arrival path is defined. Its greedy part collects no very bad vertex, and its final $a_i$ steps $e_1$ and $b_i$ steps $e_2$ stay in the box $v-[0,n_{\rm box}]^2$, which contains no very bad vertex. The transverse bound follows from \ref{geo:greedy} and the fact that the box has transverse width at most $n_{\rm box}$.
\smallskip

\emph{\ref{geo:merge}.} For lattice points on the same level the difference of transverse coordinates is an integer, and a step of either path changes its transverse coordinate by $1-\Dir_1$ or $-\Dir_1$, so $D_i=\rcoord(x_i)-\rcoord(y_i)$ is an integer-valued sequence with increments in $\{-1,0,1\}$, positive at $i=m$ and negative at $i=n$. Hence $D_i=0$ for some $m<i<n$.
\smallskip

\emph{\ref{geo:trim}.} Let $v_0,\dots,v_{r_0}$ be consecutive vertices of an admissible path in $\cQ_k$. A very bad vertex $y$ lies in $v_a-[0,n_{\rm box}]^2$ only if $v_a\in y+[0,n_{\rm box}]^2$, and a path meets this box in at most $2n_{\rm box}+1$ vertices. Two very bad vertices within $\abs\cdot_1$-distance $2n_{\rm box}$ of some $v_a$ and some $v_b$ lie in $\cQ_k^+\subset\cQ_k^{++}$, as $2n_{\rm box}\le c_0$, and are at distance at most $r_0+4n_{\rm box}\le d_0$, so they lie in one cluster, which has at most $K_1$ vertices and rules out at most $K_1(2n_{\rm box}+1)<r_0+1$ of the $v_a$. Similarly, two trapped vertices among the $v_a$, which lie in $\cQ_k\subset\cQ_k^+$, are at distance at most $r_0\le D_1$, so they belong to one trapped component, which by \ref{geo:trap} accounts for at most $N_{\rm tr}<r_0+1$ of the $v_a$. The statement for $\vbad \cup\Vbw$ is proved in the same way, with the backward trapped components.
\end{proof}

Under \ref{A2}, we lower the threshold in the definition of a bad path in Section~\ref{sub:onescale}: a path $\alpha$ of length $\abs\alpha$ is \emph{bad}, and $\Bad(\alpha)$ denotes this event, if
\[
\pawgt(\alpha)<-2\Cm\abs\alpha-\frac{t_k}2 ,
\]
with $\Cm=\E\w_0^-$ as in \eqref{eq:C1}. The following lemma is the precise form of \eqref{eq:badstair-general}.

\begin{lemma}\label{lem:stair-general}
Assume \ref{A2}. For every $k$ let $\cP_k$ be an $\cF_\vbad$-measurable family of at most $s_k^\varrho$ admissible paths in $\cQ_k$, each of length at most $8\cnst_*k^2s_k$ and collecting no very bad vertex. Then
\[
\P\bigl\{\Bad(\alpha)\text{ for some }\alpha\in\cP_k\bigr\}\longrightarrow0\qquad\text{as }k\to\infty .
\]
\end{lemma}

\begin{proof}
Apply Lemma~\ref{lem:weights} with $T$ as in \eqref{eq:Tk}, $N=s_k^\varrho$, $n=8\cnst_*k^2s_k$, $K_0=K_2$, $b=s_k^{1-\e'}$ with $\e'$ from \eqref{eq:epsprime}, and $t=t_k/4$, for $k$ so large that $\mu([-T,\infty))\ge\tfrac12$. Since $(K_2-1)T+t_k/4<t_k/2$ by \eqref{eq:Tk}, it suffices to show that the right-hand side of \eqref{eq:weights} tends to $0$. With $b^{-p'}=s_k^{-(p'+1)/2}$, the first term is
\[
s_k^\varrho \bigl(16\cnst_*C_{p'}k^2s_k^{(1-p')/2}\bigr)^{K_2}=\cnst\,k^{2K_2}s_k^{\,\varrho -K_2(p'-1)/2}\longrightarrow0
\]
by the third condition in \eqref{eq:K1K2}, and the second term is at most
\[
s_k^\varrho \exp\Bigl\{-\frac{t_k^2}{32s_k^{1-\e'}(16\cnst_*\Cm k^2s_k+t_k/12)}\Bigr\}\le s_k^\varrho \exp\Bigl\{-\frac{\e_0^2\,s_k^{\e'}}{\cnst\,k^4}\Bigr\}\longrightarrow0 . \qedhere
\]
\end{proof}

The upper bound $M$ of \eqref{ass:bounded} is needed in the construction below only because of the trimming of \ref{ch:trim}: up to $r_0$ collected vertices are removed from each end of a path whose weight is bounded from below, and the removed weight has to be bounded from above. Under \ref{A3} this is done on the following events. Let $C_{p''}=\E[(\w_0^+)^{p''}]$, so that $\P(\w_0>t)\le C_{p''}t^{-p''}$ for all $t>0$, let
\[
E_{\rm hv}=\{\w_z\le t_k\ \text{ for all }z\in\cQ_k\},
\]
and, for a coarse scale index $j$, $I\in\cG_{j,k}$, $h\in\cH_{j,k}$, and $1\le\nu\le j^2$, let $H^{(j)}_{I,h,\nu}$ be the event that $\w_z\le w_{j,k}$ for every vertex $z$ of the block \eqref{eq:crossblock}.

\begin{lemma}\label{lem:heavy}
Assume \ref{A3}. Then $\P(E_{\rm hv}^c)\to0$ and $\max\P\bigl((H^{(j)}_{I,h,\nu})^c\bigr)\to0$ as $k\to\infty$, the maximum being over the coarse scale indices $j$ and all $I\in\cG_{j,k}$, $h\in\cH_{j,k}$, $1\le\nu\le j^2$. Consequently, for all large $k$,
\begin{equation}\label{eq:heavy-cross}
\P\bigl(\Elev_{z^{(j)}_{I,h,\nu}}\!\!(w_{j,k},l'_{j,k})\cap H^{(j)}_{I,h,\nu}\bigr)\ge\frac{q_0}2
\end{equation}
for all such $(j,I,h,\nu)$, and for fixed $(j,I,h)$ the events in \eqref{eq:heavy-cross} are independent over $\nu$.
\end{lemma}

\begin{proof}
By a union bound, $\P(E_{\rm hv}^c)\le\abs{\cQ_k}\,C_{p''}t_k^{-p''}\le\cnst_{\!A,R}C_{p''}\,s_k^4(k/\e_0)^{p''}s_k^{-p''}\to0$, since $p''>4$. At a coarse scale $l_{j,k}\ge8$, so \eqref{eq:record} applies (see the lines after \eqref{t.aux34}) and gives, with \eqref{eq:Lbounds}, $l_{j,k}\le\tfrac12\clen_{k-j}\le\tfrac12s_{k-j}^3=\tfrac12w_{j,k}^3$. The block \eqref{eq:crossblock} has $l'_{j,k}+1\le2l_{j,k}$ levels and at most $16w_{j,k}+1\le17w_{j,k}$ lattice points on each, hence at most $34l_{j,k}w_{j,k}\le17w_{j,k}^4$ vertices. Therefore
\[
\P\bigl((H^{(j)}_{I,h,\nu})^c\bigr)\le17w_{j,k}^4\,C_{p''}w_{j,k}^{-p''}=17C_{p''}\,w_{j,k}^{4-p''}\le17C_{p''}\,2^{-(p''-4)(k-3)/4}\longrightarrow0
\]
uniformly in $(j,I,h,\nu)$, by \eqref{eq:wlarge} and $p''>4$. Since $\P(\Elev_{z^{(j)}_{I,h,\nu}}\!\!(w_{j,k},l'_{j,k}))>q_0$ by \eqref{eq:trialprob}, \eqref{eq:heavy-cross} follows for all $k$ with $17C_{p''}2^{-(p''-4)(k-3)/4}\le q_0/2$. Both events in \eqref{eq:heavy-cross} depend only on the weights in the block \eqref{eq:crossblock}, and for fixed $(j,I,h)$ the $j^2$ blocks are pairwise disjoint (see the lines before \eqref{eq:allfail}), which gives the independence.
\end{proof}

\medskip
\emph{The changes in Section~\ref{sec:endpoints}.} Under \ref{A2}--\ref{A3} in place of \eqref{ass:bounded}, the construction of Section~\ref{sec:endpoints} goes through with the following modifications, each announced by a remark at the corresponding place in Section~\ref{sec:endpoints}. Throughout, the events $E_{\rm geo}$ of Lemma~\ref{lem:stair-geometry}, $E_{\rm hv}$ of Lemma~\ref{lem:heavy}, and the event of Lemma~\ref{lem:stair-general} for the family $\cP_k$ of \ref{ch:family} are added to the good event $G_{k,J}$.

A word on the region $\cQ_k$. All vertices of the paths of Section~\ref{sec:endpoints} satisfy $\abs\rcoord\le(R+2\cnst_1)s_k$ and $\abs\level\le(A+\tfrac92)\clen_k$ (the grid intervals with $I\in\cG_{j,k}$ reach level $(A+4)\clen_k+l_{j,k}\le(A+\tfrac92)\clen_k$), and the paths constructed below stay within transverse distance $2\Delta_0+n_{\rm box}+r_0+4$ of straight segments joining such vertices. So for large $k$ they lie in $\cQ_k$ together with their $\Delta_0$-neighborhoods, and the hypothesis of Lemma~\ref{lem:stair-geometry}\ref{geo:arrive} on the segment $[z,v]$ holds wherever that item is applied. This is not repeated below.  

\begin{enumerate}[label={\rm(U\arabic*)},itemsep=\medskipamount,listparindent=\parindent]
\item\label{ch:gamma} \emph{Staircases through a vertex.} For a vertex $u\in\cQ_k$ let $\gamma^u$ consist of the backward greedy path from $u$ along the line $\rcoord=\rcoord(u)$, if $u\notin \vbad \cup\Vbw$, and of the forward greedy path from $u$ along the same line, if $u\notin\Vfw$, both run for as long as their vertices stay in $\cQ_k$. Thus $\gamma^u$ consists of whichever of these two parts exist, and a vertex $\gamma^u(I)$ is defined for the levels $I$ covered by these parts. By construction, $\gamma^u$ is an admissible path through $u$ that collects no very bad vertex, and by Lemma~\ref{lem:stair-geometry}\ref{geo:greedy} it stays within transverse distance $\Delta_0$ of the line $\rcoord=\rcoord(u)$. It replaces the deterministic staircase $I\mapsto\hat z(I,\rcoord(u))$: the grid points $\hat z(I,\rcoord(u))$ are replaced by the vertices $\gamma^u(I)$, and the fine pieces $y^{\rcoord(u)}_{I:I+l_{j,k}}$ of the entrance and exit paths of $u$ (defined in the proof of Proposition~\ref{prop:endpoints}) are the segments of $\gamma^u$ between consecutive grid levels. 
\item\label{ch:heavy} \emph{Crossing events.} At a coarse scale, the crossing events $\Elev_{z^{(j)}_{I,h,\nu}}\!\!(w_{j,k},l'_{j,k})$ in \eqref{eq:goodcoarse} are replaced by the events $\Elev_{z^{(j)}_{I,h,\nu}}\!\!(w_{j,k},l'_{j,k})\cap H^{(j)}_{I,h,\nu}$ of \eqref{eq:heavy-cross}. By Lemma~\ref{lem:heavy}, \eqref{eq:allfail} holds for them with $e^{-q_0j^2/2}$ in place of $e^{-q_0j^2}$, which is \eqref{eq:allfail-general}.
\item\label{ch:trim} \emph{Trimming.} A path that is to be reached by an arrival path, Lemma~\ref{lem:stair-geometry}\ref{geo:arrive}, must start at a vertex $v$ whose box $v-[0,n_{\rm box}]^2$ is free of very bad vertices, and a path from which a greedy path is to start must end at a vertex that is not trapped. By Lemma~\ref{lem:stair-geometry}\ref{geo:trim}, this can be arranged by moving the start of the path forward along it by at most $r_0$ steps and its end backward by at most $r_0$ steps. We call this \emph{trimming}. It is applied to the crossings of the events $\Elev_{z^{(j)}_{I,h,\nu}}\!\!(w_{j,k},l'_{j,k})$ at the coarse scales, and to the path $x_{u:u'}$ in the proof of Proposition~\ref{prop:endpoints}, whose trimmed start and end are denoted by $\tilde u$ and $\tilde u'$.

Trimming removes collected vertices, so it lowers the weight of the path by the weight of the removed vertices, and it changes $f_\sig$ of the displacement by at most $\abs{\ghatz}+\abs{\ghatpr\sig}\le\abs{\ghatz}+K$ per removed vertex. So if $\pawgt>f_\sig+b$ holds for a path, then $\pawgt>f_\sig+b-(\abs{\ghatz}+K+\w_z^+)$ holds after the removal of a collected vertex $z$ at one of its ends. On $H^{(j)}_{I,h,\nu}$, the trimming of a crossing of $\Elev_{z^{(j)}_{I,h,\nu}}\!\!(w_{j,k},l'_{j,k})$ therefore lowers its weight by at most $r_0w_{j,k}$ at each end, and on $E_{\rm hv}$ the trimming of $x_{u:u'}$ lowers the margin in $\pawgt(x_{u:u'})>f_\sig(u'-u)+\e s_k$ by at most $r_0(t_k+\abs{\ghatz}+K)$ at each end.

Trimming requires at least $r_0+1$ vertices at each end of the path, with no vertex in common, that is, at least $2r_0+2$ collected vertices. The crossings have $l'_{j,k}\ge l_{j,k}/2\ge4\cnst_*j^2w_{j,k}\ge2r_0+2$ collected vertices, by \eqref{eq:wlarge}. For $x_{u:u'}$, note that a path with $n\le2r_0+1$ collected vertices in $\cQ_k$ has, on $E_{\rm hv}$, weight at most $nt_k$, while its displacement has level $n$ and $\abs\rcoord\le n$, so that
\[
\pawgt\le nt_k\le f_\sig+n(\abs{\ghatz}+K+t_k)\le f_\sig+(2r_0+1)(\abs{\ghatz}+K+t_k)<f_\sig+\e s_k
\]
for large $k$, and it cannot be the path $x_{u:u'}$ of Proposition~\ref{prop:endpoints}.

Finally, consider the trimmed endpoints $\tilde u$ and $\tilde u'$ of $x_{u:u'}$. The box $\tilde u-[0,n_{\rm box}]^2$ is free of very bad vertices, so $\tilde u$ is not very bad, and it is not backward trapped either, because the box contains $(n_{\rm box}+1)^2>K_3^2$ vertices of $\Abw(\tilde u)$ and $E_{\rm geo}$ then makes $\Abw(\tilde u)$ infinite. So $\tilde u\notin \vbad \cup\Vbw$ and the backward part of $\gamma^{\tilde u}$ exists. The vertex $\tilde u'$ is not trapped, so the forward part of $\gamma^{\tilde u'}$ exists. Both parts are used in \ref{ch:endpoints}.
\item\label{ch:coarse} \emph{Coarse pieces.} At a coarse scale, for $\rho=\rcoord(u)$, the path $y^\rho_{I:I+l_{j,k}}$ is the concatenation of four pieces. The first is the arrival path of Lemma~\ref{lem:stair-geometry}\ref{geo:arrive} from $\gamma^u(I)$ to the trimmed start $v$ of the crossing. The second is the trimmed crossing from $v$ to its trimmed end $v'$. The third is the greedy path from $v'$ along the straight line from $v'$ to the point on level $I+l_{j,k}$ with transverse coordinate $\rho-2\Delta_0-3$, run until it meets $\gamma^u$, which happens before level $I+l_{j,k}$ by Lemma~\ref{lem:stair-geometry}\ref{geo:merge}. The fourth is $\gamma^u$ from the meeting point to $\gamma^u(I+l_{j,k})$. The arrival path requires that $\gamma^u(I)$ be defined and not trapped. The coarse pieces are used, in \ref{ch:endpoints}, on grid intervals $[I,I+l_{j,k}]$ that lie either above $u$, where $\gamma^u(I)$ is on the forward greedy path and hence not trapped, or below $u$, where the segment of $\gamma^u$ from $\gamma^u(I)$ to $u$ is an admissible path that collects no very bad vertex and has at least $l_{j,k}+1>K_3^2$ vertices, all of which lie in $\Afw(\gamma^u(I))$, so that $\Afw(\gamma^u(I))$ is infinite on $E_{\rm geo}$ and $\gamma^u(I)$ is not trapped. The slope conditions of Lemma~\ref{lem:stair-geometry} hold for the first and the third piece: their transverse displacements have absolute values in $[21w_{j,k}-r_0-\Delta_0-1,\,(30j^2+8)w_{j,k}+r_0+2\Delta_0+4]\subset(0,40j^2w_{j,k}]$, positive for the first and negative for the third, and their lengths are at least $m'_{j,k}\ge\cnst_*j^2w_{j,k}$, so that the $e_1$-fraction lies in $[\Dir_1-\tfrac45\theta_0,\Dir_1+\tfrac45\theta_0]\subset[\theta_1,1-\theta_1]$.
\item\label{ch:family} \emph{The family $\cP_k$.} It consists of all paths whose weights are collected in the construction, other than the crossings. At the fine scales these are, for every $u\in\cQ_k$ and every fine scale index $j$, the segments of $\gamma^u$ between consecutive grid levels of scale $j$ covered by $\gamma^u$. At the coarse scales these are, for every coarse $(j,I,h)$ and every $u\in\cQ_k$ with $\rcoord(u)$ in the bin $h$ such that $\gamma^u(I)$ is defined and not trapped, the arrival paths from $\gamma^u(I)$ to every vertex $v$ on the levels $I+m'_{j,k},\dots,I+m'_{j,k}+r_0$ of one of the $j^2$ blocks \eqref{eq:crossblock} whose box $v-[0,n_{\rm box}]^2$ is free of very bad vertices, the greedy paths from every non-trapped vertex $v'$ on the levels $I+l_{j,k}-m'_{j,k}-r_0,\dots,I+l_{j,k}-m'_{j,k}$ of one of the blocks up to their meeting level with $\gamma^u$, and the segments of $\gamma^u$ from these meeting levels to $I+l_{j,k}$. All these paths are $\cF_\vbad$-measurable and collect no very bad vertex. Their lengths are at most $8\cnst_*j^2w_{j,k}+r_0\le8\cnst_*k^2s_k$ for large $k$, by \eqref{t.aux32} and $w_{j,k}\le s_k/2$, because at a fine scale the segments of $\gamma^u$ have length $l_{j,k}<8\cnst_*j^2w_{j,k}$, and at a coarse scale each of the three kinds of pieces spans at most $m'_{j,k}+r_0\le2\cnst_*j^2w_{j,k}+r_0$ levels. Their number is at most
\[
4\,\abs{\cQ_k}\cdot k\cdot\max_j\abs{\cG_{j,k}}\abs{\cH_{j,k}}\,j^2(16w_{j,k}+2)(r_0+1)\le\cnst\,k^3s_k^{10}\le s_k^{11}
\]
for large $k$, by $\abs{\cQ_k}\le\cnst s_k^4$, at most $k$ scale indices, $\abs{\cG_{j,k}}\abs{\cH_{j,k}}\le\cnst2^{5j}\le\cnst s_k^5$, $j^2\le k^2$, and $16w_{j,k}+2\le17s_k$. So $\varrho =11$ works in \eqref{eq:K1K2}.
\item\label{ch:cost} \emph{Costs.} At a fine scale, the piece $y^\rho_{I:I+l_{j,k}}$ of \ref{ch:gamma} has length $l_{j,k}$, and its endpoints are within transverse distance $\Delta_0$ of $\rho$, so its displacement has $\abs\rcoord\le2\Delta_0$ and $f_\sig(y^\rho_{I+l_{j,k}}-y^\rho_I)\le\ghatz l_{j,k}+2\Delta_0K$. On the event of Lemma~\ref{lem:stair-general} its weight is at least $-2\Cm l_{j,k}-t_k/2$, so in place of \eqref{eq:finecost},
\begin{equation}\label{eq:finecost-general}
\begin{aligned}
\pawgt(y^\rho_{I:I+l_{j,k}})&\ge f_\sig(y^\rho_{I+l_{j,k}}-y^\rho_I)-(\ghatz+2\Cm)l_{j,k}-2\Delta_0K-\frac{t_k}2\\
&\ge f_\sig(y^\rho_{I+l_{j,k}}-y^\rho_I)-8\cnst_*(\ghatz+2\Cm)j^2w_{j,k}-2\Delta_0K-\frac{t_k}2 ,
\end{aligned}
\end{equation}
where the last step uses $\ghatz+2\Cm\ge\E\w_0+2\E\w_0^-\ge0$, by $\ghatz\ge\E\w_0$ from \ref{S3}. At a coarse scale, the three pieces of \ref{ch:coarse} other than the crossing together span the $2m'_{j,k}$ levels outside the untrimmed crossing plus the at most $2r_0$ trimmed levels, so on the event of Lemma~\ref{lem:stair-general} their total weight is at least $-2\Cm(2m'_{j,k}+2r_0)-3t_k/2$, and the trimmed crossing has weight at least $\ghatz l'_{j,k}-(K+2r_0)w_{j,k}$ by \eqref{success55} and \ref{ch:trim}. The displacement of $y^\rho_{I:I+l_{j,k}}$ again has $\abs\rcoord\le2\Delta_0$, and $l_{j,k}-l'_{j,k}=2m'_{j,k}\le4\cnst_*j^2w_{j,k}$, so in place of \eqref{eq:coarsecost},
\begin{equation}\label{eq:coarsecost-general}
\begin{aligned}
\pawgt(y^\rho_{I:I+l_{j,k}})
&\ge f_\sig(y^\rho_{I+l_{j,k}}-y^\rho_I)-4\cnst_*(\ghatz+2\Cm)j^2w_{j,k}\\
&\qquad\qquad
-(K+2r_0)w_{j,k}-2\Delta_0K-4r_0\Cm-\frac{3t_k}2 .
\end{aligned}
\end{equation}
Since $w_{j,k}\to\infty$ by \eqref{eq:wlarge}, the constant terms $2\Delta_0K+4r_0\Cm$ are at most $j^2w_{j,k}$ for large $k$, so both bounds give
\[\pawgt(y^\rho_{I:I+l_{j,k}})\ge f_\sig(y^\rho_{I+l_{j,k}}-y^\rho_I)-\cnst_1j^2w_{j,k}-2t_k\]
with
\[\cnst_1=8\cnst_*(\ghatz+2\Cm)+2K+2r_0+3\cnst_*+1 ,\]
which replaces \eqref{eq:Ccone}. As for the transverse bounds: at a fine scale all vertices are within $\Delta_0\le\cnst_1j^2w_{j,k}$ of $\rho$ for large $k$. At a coarse scale, the vertices of the crossing are within $(30j^2+9)w_{j,k}$ of $\rho$ as before.  The endpoints of the other pieces, and the target point with transverse coordinate $\rho-2\Delta_0-3$ of the third piece, are within $(30j^2+9)w_{j,k}+r_0+2\Delta_0+3$ of $\rho$. These pieces stay within $\Delta_0+n_{\rm box}$ of the straight segments joining their endpoints or aiming at that target. So all vertices are within $(30j^2+9)w_{j,k}+r_0+3\Delta_0+n_{\rm box}+3\le(30j^2+10)w_{j,k}\le\cnst_*j^2w_{j,k}$ of $\rho$, by \eqref{eq:wlarge}. Hence \eqref{eq:coarseloc} and the second part of \eqref{eq:gridcost} hold as stated.
\item\label{ch:prob} \emph{Probabilities.} There are no families $\Stairs^{(j)}_{I,h}$: $\Good^{(j)}_{I,h}$ is the whole space when scale $j$ is fine, and the union over $\nu$ of the events of \ref{ch:heavy} when scale $j$ is coarse. Accordingly, in \eqref{eq:goodprob} the term $e^{-q_0j^2}$ becomes $e^{-q_0j^2/2}$ and the terms containing $\pbad_k$ disappear, and the good event $G_{k,J}$ of \eqref{eq:Gdef} is intersected with $E_{\rm geo}$, $E_{\rm hv}$, and the event of Lemma~\ref{lem:stair-general} for the family $\cP_k$ of \ref{ch:family}. In \eqref{eq:Gsplit} the crossing term has $e^{-q_0j^2/2}$, so that $\beta_J=\cnst_{\!A,R}\sum_{j\ge J}2^{5j}e^{-q_0j^2/2}$ in \eqref{eq:betaJ}, still finite and tending to $0$ as $J\to\infty$, and the staircase terms are replaced by $\P(E_{\rm geo}^c)+\P(E_{\rm hv}^c)$ plus the probability in Lemma~\ref{lem:stair-general}, all $o(1)$ as $k\to\infty$ and independent of $J$. The conclusion $\P(G_{k,J}^c)\le2\beta_J$ for $k\ge k_1(J)$ is unchanged.
\item\label{ch:endpoints} \emph{Proposition~\ref{prop:endpoints}.} Its statement is unchanged. In its proof, the path $x_{u:u'}$ is first trimmed to $x_{\tilde u:\tilde u'}$ as in \ref{ch:trim}, which lowers the margin in $\pawgt>f_\sig+\e s_k$ by at most $2r_0(t_k+\abs{\ghatz}+K)\le\e s_k/64$ for large $k$, since $t_k=\e s_k/(512k)$. The entrance path of $\tilde u$ is built along $\gamma^{\tilde u}$, from $\gamma^{\tilde u}(I^-_{k,\tilde u})$ to $\tilde u$, and the exit path of $\tilde u'$ along $\gamma^{\tilde u'}$, from $\tilde u'$ to $\gamma^{\tilde u'}(I^+_{k,\tilde u'})$, out of the pieces of Lemma~\ref{lem:grid}\ref{lem:grid:path} in the form of \ref{ch:gamma} and \ref{ch:coarse}, on the grid intervals of the decompositions of $[I^-_{k,\tilde u},\level(\tilde u)]$ and $[\level(\tilde u'),I^+_{k,\tilde u'}]$. These pieces exist: the backward part of $\gamma^{\tilde u}$ and the forward part of $\gamma^{\tilde u'}$ exist by \ref{ch:trim}, and the condition on $\gamma^u(I)$ for the coarse pieces was verified in \ref{ch:coarse}. The bounds \eqref{eq:extcost} and \eqref{eq:extcost-exit} then hold with the new $\cnst_1$, and the loss $\e s_k/64$ from the trimming is absorbed into the margin $\e s_k/2$ in the reduction to crossings between coarse levels. The rest of the proof is unchanged.
\end{enumerate}

\subsection{One-step bound}\label{sub:onestep}
In this section we assume \ref{A3} and use the notation of Section~\ref{sec:witnesses}: the parameter $\eta\in(0,1]$, the scale $(s_*,\clen_*,n_*)$, the blocks $\cB_{p,h}$ of \eqref{eq:Bblock} with the lattice points $v_{p,h}$ (Figure~\ref{fig:blocks}), the stages of a path with their indices \eqref{eq:chunkindex} and multiplicities (Figures~\ref{fig:stages} and~\ref{fig:chunks}), and the paths of blocks of Section~\ref{sub:longpaths}. 

Under \eqref{ass:bounded}, the one-step bound \eqref{eq:one-step} says that a single step cannot lift the weight of a path above $f$ by more than $\eta s_*/2$, because a collected vertex has weight at most $M$ and the scale is fixed so that \eqref{eq:rare-step} holds. This is the only place in Section~\ref{sec:witnesses} where the upper bound $M$ is used. Under \ref{A3} a collected vertex can have an arbitrarily large weight. We keep the one-step bound for the vertices with weight at most
\begin{equation}\label{eq:M0}
M_0=\frac{\eta s_*}4
\end{equation}
and keep track of the \emph{excess}
\begin{equation}\label{eq:excess}
\exc_z=(\w_z-M_0)^+,\qquad z\in\Z^2 ,
\end{equation}
of the other vertices. We will show that, almost surely, for all large $N$ the total excess that a path with $N$ full stages can collect is at most $N\eta s_*/32$, uniformly over these paths. This is a greedy lattice animal problem on the lattice of blocks (Lemma~\ref{lem:excess}), and the moment condition $p''>4$ enters through the size $s_*\clen_*\le s_*^4$ of a block. Let 
\begin{equation}\label{eq:Bprime}\begin{aligned}
    &\cB'=\{z\in\Z^2:\ 0\le\level(z)\le3n_*,\ \abs{\rcoord(z)}\le4s_*+1\} \quad\text{and}\\[3pt] 
    &Z_b=\sum_{z\tsp\in\tsp v_b+\cB'}\exc_z\quad\text{for }b=(p,h)\in\Z^2 .
\end{aligned}
\end{equation}
Since $\level(v_{p,h})=pn_*$ and $\abs{\rcoord(v_{p,h})-hs_*}<1$, the block $\cB_{p,h}$ is contained in $v_{p,h}+\cB'$, and the $Z_b$ are identically distributed. By \eqref{eq:Lbounds}, $n_*\le\clen_*\le s_*^3$, so $\abs{\cB'}\le(3n_*+1)(8s_*+3)\le44s_*^4$.

\medskip
\emph{Fixing the scale.} The scale is fixed as in Section~\ref{sub:scale}, with two changes. In \eqref{eq:rare-step} the constant $M$ is replaced by $M_0$, that is, $k$ is taken so large that $\abs{\ghatz}+\abs a\le\eta s_*/4$, and the requirement
\begin{equation}\label{eq:rare-excess}
C_0\,s_*^{\,2-p''/4}\le\frac{\eta s_*}{64}
\end{equation}
is added, where $C_0$ is the constant of Lemma~\ref{lem:excess} below. Since $p''>4$, \eqref{eq:rare-excess} holds for all large $k$. The condition \eqref{eq:rare-prob} is available under \ref{A2}--\ref{A3}, because Proposition~\ref{prop:endpoints} is, by Section~\ref{sub:staircases}.

\medskip
\emph{The one-step bound with excess.} For a single step $x_{i:i+1}$ we have $\w_{x_i}\le M_0+\exc_{x_i}$, so \eqref{eq:rare-step} with $M_0$ in place of $M$ gives
\begin{equation}\label{eq:one-step-general}
\pawgt(x_{i:i+1})=\w_{x_i}\ \le\ f(x_{i+1}-x_i)+\abs{\ghatz}+\abs a+M_0+\exc_{x_i}\ \le\ f(x_{i+1}-x_i)+\frac{\eta s_*}2+\exc_{x_i} .
\end{equation}
For a stage $x_{k:l}$ let $Z(x_{k:l})=\sum_{i=k}^{l-1}\exc_{x_i}$ be the total excess of its collected vertices, and for an admissible path $x_{0:n}$ with $x_0=0$ let $A(x_{0:n})\subset\Z^2$ be the set of indices of its stages, that is, the set of entries of its path of blocks (defined in the proof of Proposition~\ref{prop:longpaths}). 

\begin{lemma}\label{lem:multiplicity-general}
Assume \ref{A3} and let the scale be fixed as above. Then every stage $x_{k:l}$ satisfies
\begin{equation}\label{eq:chunk-defect-general}
\pawgt(x_{k:l})\le f(x_l-x_k)+\eta s_*+2\eta s_*\,\mult(x_{k:l})+Z(x_{k:l}) .
\end{equation}
Under \eqref{eq:etac}, every full stage satisfies $\pawgt(x_{k:l})\le f(x_l-x_k)-\eta s_*+4\eta s_*\,\mult(x_{k:l})+Z(x_{k:l})$, and every admissible path $x_{0:n}$ with $x_0=0$ and $N$ full stages satisfies
\begin{equation}\label{eq:total-defect-general}
\pawgt(x_{0:n})\le f(x_n-x_0)-(N-1)\eta s_*+4\eta s_*\,\mult_{\rm tot}(x_{0:n})+\sum_{b\tsp\in\tsp A(x_{0:n})}Z_b .
\end{equation}
\end{lemma}

\begin{proof}
Follow the proof of Lemma~\ref{lem:multiplicity} with \eqref{eq:one-step-general} in place of \eqref{eq:one-step}. Each split segment $x_{d_\nu:d_{\nu+1}}$ now satisfies
\[
\pawgt(x_{d_\nu:d_{\nu+1}})\le f(x_{d_{\nu+1}}-x_{d_\nu})+\eta s_*+\frac{\eta s_*}2+\exc_{x_{d_{\nu+1}-1}}\le f(x_{d_{\nu+1}}-x_{d_\nu})+2\eta s_*+\exc_{x_{d_{\nu+1}-1}} ,
\]
it is a witness of $\cA_{p,h}$ as before, and the vertices $x_{d_{\nu+1}-1}$, $0\le\nu<q$, are distinct collected vertices of the stage. Adding up as before gives \eqref{eq:chunk-defect-general}. The bound for full stages follows from \eqref{eq:chunk-defect-general} exactly as \eqref{eq:full-defect} follows from \eqref{eq:chunk-defect} in the proof of Lemma~\ref{lem:fulldefect}: the witness-free case uses Lemma~\ref{lem:chunk-witness} and \eqref{eq:defect-min}, which do not involve the one-step bound, and for $\mult\ge1$ it is arithmetic. Summing over the stages gives \eqref{eq:total-defect-general} with $\sum_{\text{stages}}Z(x_{k:l})$ in place of the last term. This sum is the total excess of the collected vertices of $x_{0:n}$, since distinct stages have disjoint sets of collected vertices. Every collected vertex of a stage with index $b$ lies in $\cB_b\subset v_b+\cB'$, so the collected vertices of all the stages with index $b$ contribute at most $Z_b$, and the sum is at most $\sum_{b\in A(x_{0:n})}Z_b$.
\end{proof}

\begin{lemma}\label{lem:excess}
Assume \ref{A3}. There exists a constant $C_0$, depending only on $\eta$, $p''$, and $\mu$, such that for every scale, almost surely,
\begin{equation}\label{eq:excess-animal}
\varlimsup_{N\to\infty}\ \frac1N\ \max_{A}\ \sum_{b\tsp\in\tsp A}Z_b\ \le\ C_0\,s_*^{\,2-p''/4} ,
\end{equation}
where the maximum is over the sets $A\subset\Z^2$ with at most $N+1$ points that contain $(0,0)$ and are connected for the relation $\abs{b-b'}_\infty\le2$. The set $A(x_{0:n})$ of every admissible path $x_{0:n}$ from $0$ with $N$ full stages is such a set.
\end{lemma}

\begin{proof}
The last assertion holds because the path of blocks of $x_{0:n}$ has $N+1$ entries, starts at $(0,0)$, and consecutive entries satisfy $p_{\tau+1}-p_\tau\in\{0,1,2\}$ and $\abs{h_{\tau+1}-h_\tau}\le2$, see Section~\ref{sub:longpaths}.

Let $C_{p''}=\E[(\w_0^+)^{p''}]$. For $1\le r\le p''$, since $\exc_0\le\w_0^+\one\{\w_0>M_0\}$,
\[
\E\exc_0^{\,r}\le\E\bigl[(\w_0^+)^r\one\{\w_0>M_0\}\bigr]\le M_0^{\,r-p''}\,C_{p''} .
\]
By \cite[Lemma~1]{Ros-70} 
%on page 277: for $1<p<\infty$ and independent $f_1,\dots,f_n\in L^p$, $\lVert f_1+\dots+f_n\rVert_p\le2^p\max\bigl\{(\sum_i\lVert f_i\rVert_p^p)^{1/p},\sum_i\lVert f_i\rVert_1\bigr\}$. No centering is assumed. With $p=4$, raised to the fourth power. 
with exponent $4$, independent nonnegative random variables $Y_1,\dotsc,Y_n$ with finite fourth moments satisfy
\[
\E\Bigl[\Bigl(\sum_{i=1}^nY_i\Bigr)^4\Bigr]\le2^{16}\Bigl(\sum_{i=1}^n\E Y_i^4+\Bigl(\sum_{i=1}^n\E Y_i\Bigr)^4\Bigr) .
\]
We apply this to the i.i.d.\ variables $\exc_x$, $x\in v_b+\cB'$, whose sum is $Z_b$. Their number is $\abs{\cB'}\le44s_*^4$, and $M_0=\eta s_*/4$, so
\[
\E Z_b^4\le2^{16}\Bigl(44s_*^4\,C_{p''}\Bigl(\frac{\eta s_*}4\Bigr)^{4-p''}+\Bigl(44s_*^4\,C_{p''}\Bigl(\frac{\eta s_*}4\Bigr)^{1-p''}\Bigr)^4\Bigr)\le C_4\bigl(s_*^{\,8-p''}+s_*^{\,20-4p''}\bigr)\le2C_4\,s_*^{\,8-p''} ,
\]
with $C_4$ depending on $\eta$, $p''$, and $C_{p''}$, where the last step uses $20-4p''\le8-p''$. If $\E Z_b^4=0$ then $Z_b=0$ almost surely and the left side of \eqref{eq:Ztail} below vanishes. Otherwise let $\varpi=(\E Z_b^4)^{1/4}\in(0,(2C_4)^{1/4}s_*^{\,2-p''/4}]$. 
By Markov's inequality, $\P(Z_b>t)\le\min\{1,\varpi^4t^{-4}\}$, hence
\begin{equation}\label{eq:Ztail}
\int_0^\infty\P(Z_b>t)^{1/2}\,dt\le\int_0^\varpi dt+\int_\varpi^\infty\frac{\varpi^2}{t^2}\,dt=2\varpi\le2(2C_4)^{1/4}\,s_*^{\,2-p''/4} .
\end{equation}

Two boxes $v_b+\cB'$ and $v_{\tilde b}+\cB'$ with $b=(p,h)$, $\tilde b=(\tilde p,\tilde h)$ are disjoint when $\abs{p-\tilde p}\ge4$, since then their level ranges $[pn_*,pn_*+3n_*]$ and $[\tilde pn_*,\tilde pn_*+3n_*]$ are disjoint, and also when $\abs{h-\tilde h}\ge10$, since then $\abs{\rcoord(v_b)-\rcoord(v_{\tilde b})}>10s_*-2>2(4s_*+1)$ for $s_*\ge3$. Partition $\Z^2$ into the $40$ classes $\cC_{\alpha,\beta}=\{(p,h):\ p\equiv\alpha\ (\mathrm{mod}\ 4),\ h\equiv\beta\ (\mathrm{mod}\ 10)\}$. Within a class the boxes are pairwise disjoint, so the variables $Z_b$, $b\in\cC_{\alpha,\beta}$, are i.i.d.

Now we use the following result on greedy lattice animals, which is \cite[Theorem~1.1]{Mar-02-spa} in dimension two.  Let $(X_v)_{v\in\Z^2}$ be i.i.d.\ nonnegative random variables with distribution function $F$ such that $\int_0^\infty(1-F(t))^{1/2}\,dt<\infty$, and let $X^*_n$ be the maximum of $\sum_{v\in\Xi}X_v$ over the lattice animals $\Xi$ of size $n$ containing the origin, that is, the connected subsets of $\Z^2$ (for nearest-neighbor adjacency) with $n$ points that contain $0$. Then $X^*_n/n$ converges almost surely to a constant that is at most $c\int_0^\infty(1-F(t))^{1/2}\,dt$, where $c$ is a universal constant.  Since the $X_v$ are nonnegative and every connected set of size at most $n$ containing $0$ is contained in a connected set of size $n$ containing $0$, $X^*_n$ is also the maximum over the connected sets of size at most $n$ containing $0$. By translation invariance the same conclusion holds for $X^{*,o}_n$, the maximum over the connected sets of size at most $n$ that contain a fixed point $o\in\Z^2$.

Fix a class $\cC=\cC_{\alpha,\beta}$. The map $\phi(p,h)=\bigl((p-\alpha)/4,(h-\beta)/10\bigr)$ is a bijection from $\cC$ onto $\Z^2$, and $\psi(p,h)=\bigl(\lfloor(p-\alpha)/4\rfloor,\lfloor(h-\beta)/10\rfloor\bigr)$ is a map from $\Z^2$ onto $\Z^2$ that agrees with $\phi$ on $\cC$. The variables $X_v=Z_{\phi^{-1}(v)}$, $v\in\Z^2$, are i.i.d.\ with the distribution function $F$ of $Z_{(0,0)}$, and by \eqref{eq:Ztail} the integral condition holds. Let $A$ be as in the lemma. Then
\[
\sum_{b\tsp\in\tsp A\cap\cC}Z_b=\sum_{v\tsp\in\tsp\phi(A\cap\cC)}X_v\le\sum_{v\tsp\in\tsp\psi(A)}X_v .
\]
If $\abs{b-b'}_\infty\le2$ then $\abs{\psi(b)-\psi(b')}_\infty\le1$, so $\psi(A)$ is connected for the relation $\abs{v-v'}_\infty\le1$, has at most $N+1$ points, and contains $o_\cC=\psi(0,0)$. Take a spanning tree of $\psi(A)$ for this relation and insert, on each of its at most $N$ edges that is not a nearest-neighbor edge, one lattice point adjacent to both endpoints. This gives a lattice animal $\Xi\supset\psi(A)$ with at most $2(N+1)$ points containing $o_\cC$, so that $\sum_{v\in\psi(A)}X_v\le X^{*,o_\cC}_{2(N+1)}$. By the quoted result, almost surely,
\begin{align*}
\varlimsup_{N\to\infty}\frac1N\max_A\sum_{b\tsp\in\tsp A\cap\cC}Z_b &\le\varlimsup_{N\to\infty}\frac{X^{*,o_\cC}_{2(N+1)}}N=2\lim_{n\to\infty}\frac{X^{*,o_\cC}_n}n\le2c\int_0^\infty(1-F(t))^{1/2}\,dt\\
&\le4c(2C_4)^{1/4}\,s_*^{\,2-p''/4} .
\end{align*}
Summing over the $40$ classes gives \eqref{eq:excess-animal} with $C_0=160c(2C_4)^{1/4}$.
\end{proof}

\medskip
\emph{Proof of Theorem~\ref{thm:main} under \ref{A2}--\ref{A3}.} Let $N_0$ be as in Proposition~\ref{prop:longpaths}, whose proof does not involve the one-step bound, and let $N_1$ be such that, by Lemma~\ref{lem:excess} and \eqref{eq:rare-excess}, $\sum_{b\in A}Z_b\le N\eta s_*/32$ for all $N\ge N_1$ and all sets $A$ as in that lemma. For every admissible path $x_{0:n}$ with $x_0=0$ and $N\ge\max\{N_0,N_1,3\}$ full stages, \eqref{eq:total-defect-general} and $\mult_{\rm tot}(x_{0:n})<N/32$ give
\[
\pawgt(x_{0:n})<f(x_n-x_0)-(N-1)\eta s_*+\frac{N\eta s_*}8+\frac{N\eta s_*}{32}\le f(x_n-x_0)-\frac{N\eta s_*}2 ,
\]
the last step because $\tfrac{27}{32}N-1\ge\tfrac12N$ for $N\ge3$. This is \eqref{eq:linear-defect}, and the rest of the proof of Theorem~\ref{thm:main} is unchanged.

\section{Glossary of symbols}\label{app:dictionary}
The notation used in the proofs is listed here in order of appearance, grouped by section, with pointers to the definitions.

\newcommand{\dictstart}{\begingroup\renewcommand{\descriptionlabel}[1]{\normalfont##1}\begin{description}[font=\normalfont,style=multiline,leftmargin=9.5em,itemsep=1pt,topsep=3pt]}
\newcommand{\dictend}{\end{description}\endgroup}

\medskip\noindent\emph{Setting and coordinates {\rm(}Section~\ref{sec:setting}{\rm)}.}
\dictstart
\item[$e_1,e_2,\ed,\ead$] the vectors $(1,0)$, $(0,1)$, $(1,1)$, and $(1,-1)$, Section~\ref{sub:cgm}.
\item[$\level(z)$] the antidiagonal level $z_1+z_2$, \eqref{eq:level}.
\item[$x_{m:n}$, $x_{u:v}$] an admissible path, indexed by the levels of its vertices, and its segment from $u$ to $v$. The length of $x_{m:n}$ is $n-m$, Section~\ref{sub:cgm}. In Appendix~\ref{app:general} the length of a path $\alpha$ is also written $\abs\alpha$.
\item[$\pawgt(x_{m:n})$] the weight $\sum_{i=m}^{n-1}\w_{x_i}$ of a path, the terminal vertex excluded. The vertices $x_m,\dotsc,x_{n-1}$ are the collected vertices, \eqref{eq:weight}.
\item[$\Lpp_{u,v}$] the last-passage time from $u$ to $v$. A maximizing path is a geodesic, \eqref{def:LPP}.
\item[$g$] the shape function, \eqref{eq:shape}.
\item[$\Uset$, $\ri\Uset$] the segment $\{(t,1-t):0\le t\le1\}$ of directions and its relative interior, Section~\ref{sub:shape}.
\item[$p'$, $p''$] the moment exponents of \ref{A2} and \ref{A3}.
\item[$\Mtop$, $p_1$] the top of the support of $\mu$ and its mass $\mu\{\Mtop\}$, Section~\ref{sec:flat}.
\item[$\Ltop_{u,v}$, $\gtop$] the maximal number of collected vertices of weight $\Mtop$ along a path from $u$ to $v$, and its shape function, Section~\ref{sec:flat}.
\item[$\Uperc$] the flat segment $\{z\in\Uset:\gtop(z)=1\}=\{z\in\Uset:g(z)=\Mtop\}$, Section~\ref{sec:flat} and \eqref{eq:flat-top}.
\item[$\opc$, $\beta_{p_1}$] the critical density of oriented site percolation, and the edge speed of the cluster of the origin when $p_1>\opc$, Section~\ref{sec:flat}.
\item[$\Dir$] the direction in $(\ri\Uset)\setminus\Uperc$ fixed for the rest of the paper, Section~\ref{sub:coords}.
\item[$\rcoord(z)$] the transverse coordinate $(1-\Dir_1)z_1-\Dir_1z_2$, \eqref{eq:coords}.  $z=\level(z)\Dir+\rcoord(z)\ead$, \eqref{eq:decomp}.
\item[$\hat z(I,\rho)$] the lattice point $z$ with $\level(z)=I$ and $\rcoord(z)\in(\rho-1,\rho]$. The path $I\mapsto\hat z(I,\rho)$ is the staircase at $\rho$, \eqref{stairdef}.
\item[$\theta_0$, $\theta$-interior] $\theta_0=\min\{\Dir_1,1-\Dir_1\}$, and $z$ is $\theta$-interior if $\level(z)>0$ and $\abs{\rcoord(z)}\le(\theta_0-\theta)\level(z)$, Section~\ref{sub:coords}.
\item[$\Lip_\theta$] a Lipschitz constant of $g$ on the cone of $\theta$-interior vectors, Section~\ref{sub:coords} and \ref{S4}.
\item[$\hat g$, $\ghatprm$, $\ghatprp$] $\hat g(h)=g(\Dir+h\ead)$, and its one-sided derivatives $\ghatprm=\hat g'(0-)\ge\ghatprp=\hat g'(0+)$ at $0$, Section~\ref{sub:extreme}.
\item[$f_\sig$, $f$, $a$, $\halfgap $] the extreme supporting functionals $f_\sig(z)=g(\Dir)\level(z)+\ghatpr{\sig}\rcoord(z)$, their average $f$, the average slope $a$, and the half-gap $\halfgap =\tfrac12(\ghatprm-\ghatprp)\ge0$, \eqref{eq:supports}. The function $g$ is differentiable at $\Dir$ if and only if $\halfgap =0$. 
\item[$z_n$, $x^{(n)}_{0:n}$] the lattice point on level $n$ with $\rcoord(z_n)\in(-1,0]$, and a geodesic from $0$ to $z_n$, \eqref{eq:geodesic-defect} and Section~\ref{sub:proofmain}. The point $z_n$ is used with the same meaning in the proof of Proposition~\ref{prop:charlength}. 
\dictend

\medskip\noindent\emph{Localized crossings {\rm(}Section~\ref{sec:crossings}{\rm)}.}
\dictstart
\item[$M$] the bound on the weights in the bounded case, \eqref{ass:bounded}.
\item[block, crossing] a block is a set $\{z:I_1\le\level(z)\le I_2,\ \abs{\rcoord(z)-\rho}\le Rs\}$, of length $I_2-I_1$ and width $2Rs$, and a crossing of it is an admissible path $x_{I_1:I_2}$ with all vertices in it, Section~\ref{sec:crossings}. The blocks $\cB_{p,h}$ of Section~\ref{sec:witnesses} are blocks of a fixed size.
\item[level interval $\cJ$] a set $\{z:\level(z)=I,\ \rcoord(z)\in[\rho_1,\rho_2]\}$, with $\rcoord$-width $\rho_2-\rho_1$, Section~\ref{sec:crossings}.
\item[$\cN_1$, $\cN_2$, $\alpha_{u:\hat u}$] the nets near two level intervals, and the attachment paths of length $m=\lceil\kappa s\rceil$ from the points of a level interval to the net, proof of Lemma~\ref{lem:intervals}.
\dictend

\medskip\noindent\emph{The proof of Proposition~\ref{prop:rarefast}.}
\dictstart
\item[layer] the interval of levels $[I_t,I_{t+1})$, $I_t=t(n+d)$, of $n$ levels for the heavy crossings followed by $d=\lceil D_0s\rceil$ levels for the connectors, step 1.
\item[site, oriented path] a point $v$ of the coarse-grained lattice $\Z^2$, and a path in it with steps $e_1$ and $e_2$, step 1.
\item[$t(v)$, $Y_v$] the layer index $v_1+v_2$ and the transverse center $\sigg(v_1+2v_2)S_0$ of a site $v$, step 1.
\item[$\hcross_{v,j}$] the $j$th heavy-crossing event at $v$: a crossing of the block of width $2Rs$ centered at $z_{v,j}$, in the first $n$ levels of the layer of $v$, with $\pawgt>f_\sig+\e s$, \eqref{vj-block} and \eqref{t.aux89}.
\item[$\cJ^{\rm out}_v$, $\cJ^{\rm in}_{v'}$] the level intervals at which the heavy crossings at $v$ end and those at a forward neighbor $v'$ of $v$ start, step 3.
\item[$\conn_{v,v'}$] the connection event: every point of $\cJ^{\rm out}_v$ is joined to every point of $\cJ^{\rm in}_{v'}$ by a path with $\pawgt\ge f_\sig-\e s/2$ in the region $\{\abs{\rcoord-Y_v}\le4S_0\}$, step 3.
\item[good site, $\cR_v$] $v$ is good if some $\hcross_{v,j}$ and both connection events at $v$ occur. This is determined by the weights in the region $\cR_v$, step 5.
\dictend

\medskip\noindent\emph{Characteristic scale {\rm(}Section~\ref{sec:scales}{\rm)}.}
\dictstart
\item[$F_\mu$, $F_\mu^{-1}$] the distribution function of $\mu$ and its quantile function, Section~\ref{sub:localmod}.
\item[$\w'$, $\tilde\w$, $\chi$, $\delta$] the modified environment, in which, independently at each site of a finite region $\cR$, with probability $\delta$ the weight is replaced by its maximum with an independent copy $\tilde\w$. The variables $\chi_v$ are the Bernoulli$(\delta)$ indicators of replacement, Section~\ref{sub:localmod}.
\item[$\jP$, $\P'$] the joint law of $(\w,\tilde\w,\chi)$, and the law of $\w'$, Section~\ref{sub:localmod}.
\item[$\cI$, $\alpha$, $c'$] the set of weights $t<\Mtop$ that the modification raises by at least $c'$ with conditional probability at least $\delta\alpha/2$, Lemma~\ref{lem:kernel}.
\item[$K$] the constant $100(1+\abs{\ghatprm}+\abs{\ghatprp})$, \eqref{t.Koo}.
\item[$\Elev_z(s,n)$, $\Elev(s,n)$] the crossing event: some path from level $\level(z)$ to level $\level(z)+n$ within transverse distance $8s$ of $\rcoord(z)$ has $\pawgt>g(\Dir)n-Ks$, and $\Elev(s,n)=\Elev_0(s,n)$, \eqref{def:Ens}.
\item[$\varphi$] a lower bound $\tfrac14(1-\gtop(\Dir))$ on the fraction of the vertices of a path in the tube whose weight the modification can raise, proof of Lemma~\ref{lem:tilt}. It is positive because $\Dir\notin\Uperc$.
\item[$q_0$] the threshold probability that defines the characteristic length, fixed in \eqref{eq:q0}.
\item[$s_k$] the transverse width $2^k$ at scale index $k$, Section~\ref{sub:charlength}.
\item[$\clen_k$, $m_k$] the characteristic length $\clen_k=2^{m_k}$, the first dyadic integer at which crossings of width of order $s_k$ become rare, \eqref{eq:clen-def} and \eqref{eq:Lbounds}.
\item[$k_0$] the index from which on $\clen_k$ exists, Proposition~\ref{prop:charlength}.
\item[$\cK$] the infinite set of record scale indices, \eqref{eq:record}.
\item[$n_k$] a length in $(\clen_k/2,\clen_k]$ with $\P(\Elev(s_k,n_k))\le q_0$, \eqref{t.n_k}.
\dictend

\medskip\noindent\emph{Control over a parallelogram {\rm(}Section~\ref{sec:endpoints}{\rm)}.}
\dictstart
\item[$A$, $R$, $J$] the longitudinal and transverse sizes of the region $\cR_k$, and the index of the coarsest scale, Section~\ref{sub:gridscales}.
\item[$j$, $j_*(k)$] a scale index, $J\le j\le j_*(k)=\lceil(m_k-1)/4\rceil$, \eqref{eq:jstar}.
\item[$l_{j,k}$, $w_{j,k}$] the longitudinal mesh $2^{\max\{0,m_k-4j-1\}}$ and the transverse width $2^{k-j}$ of scale $j$, \eqref{eq:meshes}.
\item[$\cG_{j,k}$] the left endpoints $I$ of the grid intervals $[I,I+l_{j,k}]$ of scale $j$, \eqref{eq:gridset}.
\item[$\cnst_*$, coarse, fine] $\cnst_*=50/\theta_0$. Scale $j$ is coarse if $l_{j,k}\ge8\cnst_*j^2w_{j,k}$ and fine otherwise, Section~\ref{sub:gridscales}.
\item[bin, $\cH_{j,k}$] the transverse interval $[hw_{j,k},(h+1)w_{j,k})$ with index $h$, and the set of indices of the bins that meet $[-Rs_k,Rs_k]$, \eqref{eq:bincount}.
\item[$\e_0$, $t_k$] the fluctuation allowance $t_k=\e_0s_k/k$ for staircases, \eqref{eq:tk}.
\item[staircase, $\Bad(\alpha)$, $\pbad_k$] a path fixed deterministically by its endpoints, the event that a staircase $\alpha$ of length $n$ has $\pawgt(\alpha)<n\E\w_0-t_k$, and the bound on its probability, \eqref{eq:badstair}.
\item[$y^\rho_{I:I+l_{j,k}}$] the path from $\hat z(I,\rho)$ to $\hat z(I+l_{j,k},\rho)$ built at scale $j$: a staircase at a fine scale, and staircase, crossing, staircase at a coarse scale, Section~\ref{sub:onescale} and Lemma~\ref{lem:grid}.
\item[$\Stairs^{(j)}_{I,h}$] the family of staircases that may be needed at $(j,I,h)$, \eqref{eq:finecount} and \eqref{eq:staircount}.
\item[$m'_{j,k}$, $l'_{j,k}$] the length $\lceil\cnst_*j^2w_{j,k}\rceil$ of the attaching staircases at a coarse scale, and the length $l_{j,k}-2m'_{j,k}$ of the crossing between them, \eqref{eq:mprime}.
\item[$z^{(j)}_{I,h,\nu}$] the centers of the $j^2$ disjoint blocks that host the crossing events at a coarse $(j,I,h)$, \eqref{eq:rhonu}.
\item[$\Good^{(j)}_{I,h}$] the event that the paths $y^\rho$ can be built at $(j,I,h)$, \eqref{eq:goodfine} and \eqref{eq:goodcoarse}.
\item[$G_{k,J}$] the intersection of the events $\Good^{(j)}_{I,h}$ over all scales, grid intervals, and bins, \eqref{eq:Gdef}.
\item[$\cnst_1$] the constant in the cost and location bounds \eqref{eq:gridcost} of the paths $y^\rho$, \eqref{eq:Ccone}.
\item[$\beta_J$, $k_1(J)$] the bound in $\P(G_{k,J}^c)\le2\beta_J$ and the scale index from which on it holds, \eqref{eq:betaJ} and Lemma~\ref{lem:grid}.
\item[$\cR_k$] the region $\{\abs\level\le A\clen_k,\ \abs\rcoord\le Rs_k\}$ of Proposition~\ref{prop:endpoints}.
\item[$I^\pm_{k,u}$, coarse levels] the levels of $l_{J,k}\Z$ at distance between $l_{J,k}$ and $2l_{J,k}$ below and above $u$, \eqref{eq:Ipm}. The levels of $l_{J,k}\Z$ are the coarse levels.
\item[entrance, exit path] the paths $y^u_{I^-_{k,u}:u}$ and $y^u_{u:I^+_{k,u}}$ that connect $u$ to the coarse levels below and above it, proof of Proposition~\ref{prop:endpoints}.
\dictend

\medskip\noindent\emph{Disjoint witnesses {\rm(}Section~\ref{sec:witnesses}{\rm)}.}
\dictstart
\item[$\eta$] the parameter of the section, fixed as $\min\{1,\halfgap /2\}$ in \eqref{eq:etac}.
\item[$s_*$, $\clen_*$, $n_*$] the fixed scale $(s_k,\clen_k,n_k)$, Section~\ref{sub:scale}.
\item[path event, $D(\cA_1,\dotsc,\cA_q)$] an event that some path of a finite family has weight above a prescribed value, and the event that $\cA_1,\dotsc,\cA_q$ have witnesses with pairwise disjoint sets of collected vertices, Section~\ref{sub:bk}.
\item[$v_{p,h}$, $\cB_{p,h}$] the lattice point with $\level=pn_*$ and $\rcoord\in(hs_*-1,hs_*]$, and the block at $(p,h)$, of length $3n_*$ and width $8s_*$, \eqref{eq:Bblock}.
\item[$\cA^{\pawgt}_{p,h}$] the heavy-path event: some path in $\cB_{p,h}$ has $\pawgt>f_\sig+\eta s_*$ for some $\sigg$, Section~\ref{sub:scale}.
\item[$\cA^{\Elev}_{p,h}$] the crossing event $\Elev_{v_{p+1,h}}(s_*,n_*)$, Section~\ref{sub:scale}.
\item[$\cA_{p,h}$] the union $\cA^{\pawgt}_{p,h}\cup\cA^{\Elev}_{p,h}$, \eqref{eq:Adef}.
\item[stage] a piece of a path from the origin that ends after $2n_*$ levels (longitudinal) or after a transverse move by $s_*$ (transverse). Both kinds are full stages. The last piece, which may end earlier, is the remainder, \eqref{eq:cut}.
\item[index of a stage] the pair $(p,h)$ of \eqref{eq:chunkindex}. A stage with index $(p,h)$ lies in $\cB_{p,h}$.
\item[witness-free] a full stage none of whose segments is a witness of the event $\cA_{p,h}$ of its index, Section~\ref{sub:chunks}.
\item[$\mult(x_{k:l})$, $\mult_{\rm tot}$] the maximal number of witnesses of $\cA_{p,h}$ with pairwise disjoint sets of collected vertices among the segments of a stage, and its sum over the stages of a path, Section~\ref{sub:chunks}.
\item[$\mathbf b$, $\mathbf m$, $D_{\mathbf b,\mathbf m}$] the path of blocks of a path, an allocation of witnesses to its stages, and the corresponding disjoint occurrence event, proof of Proposition~\ref{prop:longpaths}.
\item[$N_0$] the number of full stages beyond which paths carry few disjoint witnesses, Proposition~\ref{prop:longpaths}.
\dictend

\medskip\noindent\emph{Unbounded weights {\rm(}Appendix~\ref{app:general}{\rm)}.}
\dictstart
\item[$\Cm$, $C_{p'}$] $\E\w_0^-$ and $\E[(\w_0^-)^{p'}]$, \eqref{eq:C1}.
\item[very bad, $\vbad$, $\mu_T$] a vertex $z$ with $\w_z<-T$, for a threshold $T$ chosen in each section, the set of very bad vertices, and the law $\mu(\,\cdot\,|\,[-T,\infty))$, Section~\ref{sub:verybad}.
\item[$\gamma^u$] the path through $u$ that replaces the staircase at $\rcoord(u)$ and avoids very bad vertices, item~\ref{ch:gamma} of Appendix~\ref{sub:staircases} and Remark~\ref{rk:stair-general}.
\item[$H^{(j)}_{I,h,\nu}$] the event that no vertex of the block centered at $z^{(j)}_{I,h,\nu}$ has weight above $w_{j,k}$, Remark~\ref{rk:heavy-general}.
\item[$M_0$, $\exc_z$] the level $\eta s_*/4$ and the excess $(\w_z-M_0)^+$ of a vertex, \eqref{eq:M0} and \eqref{eq:excess}.
\item[$Z(x_{k:l})$, $Z_b$] the total excess of the collected vertices of a stage, and of a neighborhood $v_b+\cB'$ of the block with index $b$, Section~\ref{sub:onestep} and \eqref{eq:Bprime}.
\item[$A(x_{0:n})$] the set of indices of the stages of a path, Section~\ref{sub:onestep}.
\dictend

\medskip
\small
\bibliographystyle{aop-no-url}
\bibliography{masterbib,differentiability_extra}

\begin{thebibliography}{48}
\expandafter\ifx\csname natexlab\endcsname\relax\def\natexlab#1{#1}\fi
\expandafter\ifx\csname url\endcsname\relax
  \def\url#1{\texttt{#1}}\fi
\expandafter\ifx\csname urlprefix\endcsname\relax\def\urlprefix{URL }\fi

\bibitem{Auf-Dam-13}
\textsc{Auffinger, A.} and \textsc{Damron, M.} (2013).
\newblock Differentiability at the edge of the percolation cone and related results in first-passage percolation.
\newblock \textit{Probab. Theory Related Fields} \textbf{156} 193--227.

\bibitem{Auf-Dam-14}
\textsc{Auffinger, A.} and \textsc{Damron, M.} (2014).
\newblock A simplified proof of the relation between scaling exponents in first-passage percolation.
\newblock \textit{Ann. Probab.} \textbf{42} 1197--1211.

\bibitem{Auf-Dam-Han-17}
\textsc{Auffinger, A.}, \textsc{Damron, M.} and \textsc{Hanson, J.} (2017).
\newblock \textit{50 years of first passage percolation}, vol.~68 of \textit{University Lecture Series}.
\newblock American Mathematical Society, Providence, RI.

\bibitem{Ban-Gut-09}
\textsc{Bang-Jensen, J.} and \textsc{Gutin, G.} (2009).
\newblock \textit{Digraphs: Theory, Algorithms and Applications}.
\newblock 2nd edn. Springer Monographs in Mathematics, Springer-Verlag London, Ltd., London.

\bibitem{Bou-Lug-Mas-13}
\textsc{Boucheron, S.}, \textsc{Lugosi, G.} and \textsc{Massart, P.} (2013).
\newblock \textit{Concentration inequalities}.
\newblock Oxford University Press, Oxford.
\newblock A nonasymptotic theory of independence, With a foreword by Michel Ledoux.

\bibitem{Cha-13}
\textsc{Chatterjee, S.} (2013).
\newblock The universal relation between scaling exponents in first-passage percolation.
\newblock \textit{Ann. of Math. (2)} \textbf{177} 663--697.

\bibitem{Cha-etal-89}
\textsc{Chayes, J.~T.}, \textsc{Chayes, L.}, \textsc{Grimmett, G.~R.}, \textsc{Kesten, H.} and \textsc{Schonmann, R.~H.} (1989).
\newblock The correlation length for the high-density phase of {B}ernoulli percolation.
\newblock \textit{Ann. Probab.} \textbf{17} 1277--1302.

\bibitem{Cor-12}
\textsc{Corwin, I.} (2012).
\newblock The {K}ardar-{P}arisi-{Z}hang equation and universality class.
\newblock \textit{Random Matrices Theory Appl.} \textbf{1} 1130001, 76.

\bibitem{Cox-Dur-81}
\textsc{Cox, J.~T.} and \textsc{Durrett, R.} (1981).
\newblock Some limit theorems for percolation processes with necessary and sufficient conditions.
\newblock \textit{Ann. Probab.} \textbf{9} 583--603.

\bibitem{Dam-Han-14}
\textsc{Damron, M.} and \textsc{Hanson, J.} (2014).
\newblock Busemann functions and infinite geodesics in two-dimensional first-passage percolation.
\newblock \textit{Comm. Math. Phys.} \textbf{325} 917--963.

\bibitem{Die-18}
\textsc{Diestel, R.} (2018).
\newblock \textit{Graph theory}, vol. 173 of \textit{Graduate Texts in Mathematics}.
\newblock 5th edn. Springer, Berlin.

\bibitem{Dur-84}
\textsc{Durrett, R.} (1984).
\newblock Oriented percolation in two dimensions.
\newblock \textit{Ann. Probab.} \textbf{12} 999--1040.

\bibitem{Dur-Lig-81}
\textsc{Durrett, R.} and \textsc{Liggett, T.~M.} (1981).
\newblock The shape of the limit set in {R}ichardson's growth model.
\newblock \textit{Ann. Probab.} \textbf{9} 186--193.

\bibitem{Dur-Liu-88}
\textsc{Durrett, R.} and \textsc{Liu, X.~F.} (1988).
\newblock The contact process on a finite set.
\newblock \textit{Ann. Probab.} \textbf{16} 1158--1173.

\bibitem{Dur-Sch-88}
\textsc{Durrett, R.} and \textsc{Schonmann, R.~H.} (1988).
\newblock The contact process on a finite set. {II}.
\newblock \textit{Ann. Probab.} \textbf{16} 1570--1583.

\bibitem{Gan-Gin-Nam-26-}
\textsc{Ganguly, S.}, \textsc{Ginsburg, V.} and \textsc{Nam, K.} (2026).
\newblock Last passage percolation in hierarchical environments. Preprint (\href{https://arxiv.org/abs/2411.08018}{\tt arXiv:2411.08018}).

\bibitem{Geo-Ras-Sep-17-ptrf-2}
\textsc{Georgiou, N.}, \textsc{Rassoul-Agha, F.} and \textsc{Sepp{\"a}l{\"a}inen, T.} (2017{\natexlab{a}}).
\newblock Geodesics and the competition interface for the corner growth model.
\newblock \textit{Probab. Theory Related Fields} \textbf{169} 223--255.

\bibitem{Geo-Ras-Sep-17-ptrf-1}
\textsc{Georgiou, N.}, \textsc{Rassoul-Agha, F.} and \textsc{Sepp{\"a}l{\"a}inen, T.} (2017{\natexlab{b}}).
\newblock Stationary cocycles and {B}usemann functions for the corner growth model.
\newblock \textit{Probab. Theory Related Fields} \textbf{169} 177--222.

\bibitem{Geo-Ras-Sep-15}
\textsc{Georgiou, N.}, \textsc{Rassoul-Agha, F.} and \textsc{Sepp\"al\"ainen, T.} (2026).
\newblock Stationary cocycles for the corner growth model. Extended version (\href{https://arxiv.org/abs/1404.7786v4}{\tt arXiv:1404.7786v4}).

\bibitem{Gly-Whi-91}
\textsc{Glynn, P.~W.} and \textsc{Whitt, W.} (1991).
\newblock Departures from many queues in series.
\newblock \textit{Ann. Appl. Probab.} \textbf{1} 546--572.

\bibitem{Gol-Rin-07}
\textsc{Goldstein, L.} and \textsc{Rinott, Y.} (2007).
\newblock Functional {BKR} inequalities, and their duals, with applications.
\newblock \textit{J. Theoret. Probab.} \textbf{20} 275--293.

\bibitem{Gue-etal-15}
\textsc{Gueudr\'e, T.}, \textsc{Le~Doussal, P.}, \textsc{Bouchaud, J.-P.} and \textsc{Rosso, A.} (2015).
\newblock Ground-state statistics of directed polymers with heavy-tailed disorder.
\newblock \textit{Phys. Rev. E} \textbf{91} 062110.

\bibitem{Ham-Wel-65}
\textsc{Hammersley, J.~M.} and \textsc{Welsh, J. A.~D.} (1965).
\newblock First-passage percolation, subadditive processes, stochastic networks, and generalized renewal theory.
\newblock In \textit{Proc. {I}nternat. {R}es. {S}emin., {S}tatist. {L}ab., {U}niv. {C}alifornia, {B}erkeley, {C}alif}. Springer-Verlag, New York, 61--110.

\bibitem{Hoe-63}
\textsc{Hoeffding, W.} (1963).
\newblock Probability inequalities for sums of bounded random variables.
\newblock \textit{J. Amer. Statist. Assoc.} \textbf{58} 13--30.

\bibitem{Hof-08}
\textsc{Hoffman, C.} (2008).
\newblock Geodesics in first passage percolation.
\newblock \textit{Ann. Appl. Probab.} \textbf{18} 1944--1969.

\bibitem{Jan-Ras-Sep-23}
\textsc{Janjigian, C.}, \textsc{Rassoul-Agha, F.} and \textsc{Sepp\"{a}l\"{a}inen, T.} (2023).
\newblock Geometry of geodesics through {B}usemann measures in directed last-passage percolation.
\newblock \textit{J. Eur. Math. Soc. (JEMS)} \textbf{25} 2573--2639.

\bibitem{Jan-Ras-Sep-26}
\textsc{Janjigian, C.}, \textsc{Rassoul-Agha, F.} and \textsc{Sepp\"al\"ainen, T.} (2026).
\newblock Strong existence and uniqueness of the tilt-indexed {B}usemann process in the planar corner growth model.
\newblock \textit{Electron. J. Probab.} \textbf{31} 1--35.

\bibitem{Joc-Pro-Sho-26-}
\textsc{Jockusch, W.}, \textsc{Propp, J.} and \textsc{Shor, P.} (2026).
\newblock Random domino tilings and the arctic circle theorem.
\newblock \textit{Ann. Inst. Henri Poincar\'e Comb. Phys. Interact.} Online first (\href{http://arxiv.org/abs/math/9801068}{\tt arXiv:math/9801068}).

\bibitem{Joh-00}
\textsc{Johansson, K.} (2000).
\newblock Shape fluctuations and random matrices.
\newblock \textit{Comm. Math. Phys.} \textbf{209} 437--476.

\bibitem{Kar-Par-Zha-86}
\textsc{Kardar, M.}, \textsc{Parisi, G.} and \textsc{Zhang, Y.-C.} (1986).
\newblock Dynamic scaling of growing interfaces.
\newblock \textit{Phys. Rev. Lett.} \textbf{56} 889--892.

\bibitem{Kes-86-stflour}
\textsc{Kesten, H.} (1986).
\newblock Aspects of first passage percolation.
\newblock In \textit{\'{E}cole d'{\'e}t{\'e} de probabilit{\'e}s de {S}aint-{F}lour, {XIV}---1984}, vol. 1180 of \textit{Lecture Notes in Math.} Springer, Berlin, 125--264.

\bibitem{Kes-87}
\textsc{Kesten, H.} (1987).
\newblock Scaling relations for {$2$D}-percolation.
\newblock \textit{Comm. Math. Phys.} \textbf{109} 109--156.

\bibitem{Kon-27}
\textsc{K{\H{o}}nig, D.} (1927).
\newblock {\"U}ber eine {S}chlu{\ss}weise aus dem {E}ndlichen ins {U}nendliche.
\newblock \textit{Acta Sci. Math. (Szeged)} \textbf{3} 121--130.

\bibitem{Lig-Sch-Sta-97}
\textsc{Liggett, T.~M.}, \textsc{Schonmann, R.~H.} and \textsc{Stacey, A.~M.} (1997).
\newblock Domination by product measures.
\newblock \textit{Ann. Probab.} \textbf{25} 71--95.

\bibitem{Mar-02}
\textsc{Marchand, R.} (2002).
\newblock Strict inequalities for the time constant in first passage percolation.
\newblock \textit{Ann. Appl. Probab.} \textbf{12} 1001--1038.

\bibitem{Mar-02-spa}
\textsc{Martin, J.~B.} (2002).
\newblock Linear growth for greedy lattice animals.
\newblock \textit{Stochastic Process. Appl.} \textbf{98} 43--66.

\bibitem{Mar-04}
\textsc{Martin, J.~B.} (2004).
\newblock Limiting shape for directed percolation models.
\newblock \textit{Ann. Probab.} \textbf{32} 2908--2937.

\bibitem{Mou-93}
\textsc{Mountford, T.~S.} (1993).
\newblock A metastable result for the finite multidimensional contact process.
\newblock \textit{Canad. Math. Bull.} \textbf{36} 216--226.

\bibitem{New-95}
\textsc{Newman, C.~M.} (1995).
\newblock A surface view of first-passage percolation.
\newblock In \textit{Proceedings of the {I}nternational {C}ongress of {M}athematicians, {V}ol.\ 1, 2 ({Z}{\"u}rich, 1994)}. Birkh{\"a}user, Basel.

\bibitem{OCo-99}
\textsc{O'Connell, N.} (1999).
\newblock Directed percolation and tandem queues.
\newblock \textit{DIAS Technical Report DIAS-APG-9912} .

\bibitem{Ope-26-convdiff}
\textsc{{OpenAI}} (2026).
\newblock Strict convexity and differentiability of the planar exponential first-passage limit shape.
\newblock Preprint, available at \href{https://github.com/openai/math/blob/main/preprints/Strict-convexity-and-differentiability-of-the-planar-exponential-first-passage-limit-shape-September-24-2026/main.pdf}{GitHub}.

\bibitem{Roc-70}
\textsc{Rockafellar, R.~T.} (1970).
\newblock \textit{Convex analysis}.
\newblock Princeton Mathematical Series, No. 28, Princeton University Press, Princeton, N.J.

\bibitem{Ros-70}
\textsc{Rosenthal, H.~P.} (1970).
\newblock On the subspaces of {$L\sp{p}$} {$(p>2)$} spanned by sequences of independent random variables.
\newblock \textit{Israel J. Math.} \textbf{8} 273--303.

\bibitem{Ros-81}
\textsc{Rost, H.} (1981).
\newblock Nonequilibrium behaviour of a many particle process: density profile and local equilibria.
\newblock \textit{Z. Wahrsch. Verw. Gebiete} \textbf{58} 41--53.

\bibitem{Sep-97-aap-1}
\textsc{Sepp{{\"a}}l{{\"a}}inen, T.} (1997).
\newblock A scaling limit for queues in series.
\newblock \textit{Ann. Appl. Probab.} \textbf{7} 855--872.

\bibitem{Sep-98-mprf-2}
\textsc{Sepp{{\"a}}l{{\"a}}inen, T.} (1998{\natexlab{a}}).
\newblock Coupling the totally asymmetric simple exclusion process with a moving interface.
\newblock \textit{Markov Process. Related Fields} \textbf{4} 593--628.
\newblock I Brazilian School in Probability (Rio de Janeiro, 1997).

\bibitem{Sep-98-mprf-1}
\textsc{Sepp{{\"a}}l{{\"a}}inen, T.} (1998{\natexlab{b}}).
\newblock Hydrodynamic scaling, convex duality and asymptotic shapes of growth models.
\newblock \textit{Markov Process. Related Fields} \textbf{4} 1--26.

\bibitem{Ber-Kes-85}
\textsc{van~den Berg, J.} and \textsc{Kesten, H.} (1985).
\newblock Inequalities with applications to percolation and reliability.
\newblock \textit{J. Appl. Probab.} \textbf{22} 556--569.

\end{thebibliography}
\end{document}